\documentclass{amsart}%
\usepackage{amsmath}
\usepackage[dvips]{graphicx}
\usepackage{amsxtra}%
\usepackage{amsfonts}%
\usepackage{amssymb}
\newtheorem{theorem}{Theorem}[section]
\newtheorem{corollary}[theorem]{Corollary}
\newtheorem{lemma}[theorem]{Lemma}

\newtheorem{proposition}[theorem]{Proposition}
\theoremstyle{definition}
\newtheorem{notation}[theorem]{Notation}
\newtheorem{definition}[theorem]{Definition}
\newtheorem{remark}[theorem]{Remark}
\newtheorem{rt}[theorem]{Remark and Terminology}

\newtheorem{observation}[theorem]{Observation}
\newtheorem{example}[theorem]{Example}

\theoremstyle{remark}

\newtheorem*{acknowledgements}{Acknowledgements}
\numberwithin{equation}{section}
\newlength{\displayboxwidth}

\newlength{\bracelength}

\newbox\frogdown
\newlength\frogdrop
\def\hookdownarrow{\setlength{\unitlength}{0.4pt}\setbox\frogdown=\hbox to 0pt{\hss $\displaystyle \downarrow $\hss }\setlength{\frogdrop}{2.5\ht\frogdown}\raisebox{0pt}[5\unitlength][\frogdrop]
{\begin{picture}(0,5)(10,0)
\put(5,0){\oval(10,10)[t]}
\end{picture}\lower\ht\frogdown\box\frogdown}}
\def\openone
{\mathchoice
{\hbox{\upshape \small1\kern-3.3pt\normalsize1}}
{\hbox{\upshape \small1\kern-3.3pt\normalsize1}}
{\hbox{\upshape \tiny1\kern-2.3pt\SMALL1}}
{\hbox{\upshape \Tiny1\kern-2pt\tiny1}}}
\makeatletter
\newbox\ipbox
\newcommand{\ip}[2]{\left\langle #1\mathrel{\mathchoice
{\setbox\ipbox=\hbox{$\displaystyle \left\langle\mathstrut #1#2\right\rangle$}
\vrule height\ht\ipbox width0.25pt depth\dp\ipbox}
{\setbox\ipbox=\hbox{$\textstyle \left\langle\mathstrut #1#2\right\rangle$}
\vrule height\ht\ipbox width0.25pt depth\dp\ipbox}
{\setbox\ipbox=\hbox{$\scriptstyle \left\langle\mathstrut #1#2\right\rangle$}
\vrule height\ht\ipbox width0.25pt depth\dp\ipbox}
{\setbox\ipbox=\hbox{$\scriptscriptstyle \left\langle\mathstrut #1#2\right\rangle$}
\vrule height\ht\ipbox width0.25pt depth\dp\ipbox}
} #2\right\rangle}
\newcommand{\diracb}[1]{\left\langle #1\mathrel{\mathchoice
{\setbox\ipbox=\hbox{$\displaystyle \left\langle\mathstrut #1\right.$}
\vrule height\ht\ipbox width0.25pt depth\dp\ipbox}
{\setbox\ipbox=\hbox{$\textstyle \left\langle\mathstrut #1\right.$}
\vrule height\ht\ipbox width0.25pt depth\dp\ipbox}
{\setbox\ipbox=\hbox{$\scriptstyle \left\langle\mathstrut #1\right.$}
\vrule height\ht\ipbox width0.25pt depth\dp\ipbox}
{\setbox\ipbox=\hbox{$\scriptscriptstyle \left\langle\mathstrut #1\right.$}
\vrule height\ht\ipbox width0.25pt depth\dp\ipbox}
}\right. }
\newcommand{\dirack}[1]{\left. \mathrel{\mathchoice
{\setbox\ipbox=\hbox{$\displaystyle \left.\mathstrut #1\right\rangle$}
\vrule height\ht\ipbox width0.25pt depth\dp\ipbox}
{\setbox\ipbox=\hbox{$\textstyle \left.\mathstrut #1\right\rangle$}
\vrule height\ht\ipbox width0.25pt depth\dp\ipbox}
{\setbox\ipbox=\hbox{$\scriptstyle \left.\mathstrut #1\right\rangle$}
\vrule height\ht\ipbox width0.25pt depth\dp\ipbox}
{\setbox\ipbox=\hbox{$\scriptscriptstyle \left.\mathstrut #1\right\rangle$}
\vrule height\ht\ipbox width0.25pt depth\dp\ipbox}
} #1\right\rangle}
\newcommand{\rip}[2]{\left( #1\mathrel{\mathchoice
{\setbox\ipbox=\hbox{$\displaystyle \left(\mathstrut #1#2\right)$}
\vrule height\ht\ipbox width0.25pt depth\dp\ipbox}
{\setbox\ipbox=\hbox{$\textstyle \left(\mathstrut #1#2\right)$}
\vrule height\ht\ipbox width0.25pt depth\dp\ipbox}
{\setbox\ipbox=\hbox{$\scriptstyle \left(\mathstrut #1#2\right)$}
\vrule height\ht\ipbox width0.25pt depth\dp\ipbox}
{\setbox\ipbox=\hbox{$\scriptscriptstyle \left(\mathstrut #1#2\right)$}
\vrule height\ht\ipbox width0.25pt depth\dp\ipbox}
} #2\right)}
\let\subsubsubsectionname\@empty
\newcounter{subsubsubsection}[subsubsection]

\def\l@subsubsubsection{\@tocline{4}{0pt}{1pc}{9pc}{}}

\def\subsubsubsection{\@startsection{subsubsubsection}{4}%
  \z@{.5\linespacing\@plus.7\linespacing}{-.5em}%
  {\normalfont\itshape}}
\AtBeginDocument{%
  \@for\@tempa:=-1,0,1,2,3,4\do{%
    \@ifundefined{r@tocindent\@tempa}{%
      \@xp\gdef\csname r@tocindent\@tempa\endcsname{0pt}}{}%
  }%
}
\def\@writetocindents{%
  \begingroup
  \@for\@tempa:=-1,0,1,2,3,4\do{%
    \immediate\write\@auxout{%
      \string\newlabel{tocindent\@tempa}{%
        \csname r@tocindent\@tempa\endcsname}}%
  }%
  \endgroup}
\renewcommand*\subjclass[2][2000]{%
  \def\@subjclass{#2}%
  \@ifundefined{subjclassname@#1}{%
    \ClassWarning{\@classname}{Unknown edition (#1) of Mathematics
      Subject Classification; using '2000'.}%
  }{%
    \@xp\let\@xp\subjclassname\csname subjclassname@#1\endcsname
  }%
}
\makeatother
\def\LaTeXparent#1{}%
\def\ChildStyles#1{}%

\newcommand{\matpic}{\begin{picture}
(12,3.25)(0,-0.625) \put(0,0){\makebox(1,1){$0$}} \put(0,1){\makebox(1,1){$0$}}
\put(1,0){\makebox(1,1){$B_{1,0}^{\left(  0\right)  }$}} \put(1,1){\makebox
(1,1){$B_{0,0}^{\left(  0\right)  }$}} \put(2,0){\makebox(1,1){$B_{1,1}
^{\left(  0\right)  }$}} \put(2,1){\makebox(1,1){$B_{0,1}^{\left(  0\right)
}$}} \put(3,0){\makebox(1,1){$B_{1,0}^{\left(  1\right)  }$}} \put
(3,1){\makebox(1,1){$B_{0,0}^{\left(  1\right)  }$}} \put(4,0){\makebox
(1,1){$B_{1,1}^{\left(  1\right)  }$}} \put(4,1){\makebox(1,1){$B_{0,1}
^{\left(  1\right)  }$}} \put(5,0){\makebox(2,2){$\cdots$}}
\put(7,0){\makebox(1,1){$B_{1,0}^{\scriptscriptstyle\left(  g-2\right)  }$}}
\put(7,1){\makebox(1,1){$B_{0,0}^{\scriptscriptstyle\left(  g-2\right)  }$}}
\put(8,0){\makebox(1,1){$B_{1,1}^{\scriptscriptstyle\left(  g-2\right)  }$}}
\put(8,1){\makebox(1,1){$B_{0,1}^{\scriptscriptstyle\left(  g-2\right)  }$}}
\put(9,0){\makebox(1,1){$B_{1,0}^{\scriptscriptstyle\left(  g-1\right)  }$}}
\put(9,1){\makebox(1,1){$B_{0,0}^{\scriptscriptstyle\left(  g-1\right)  }$}}
\put(10,0){\makebox(1,1){$B_{1,1}^{\scriptscriptstyle\left(  g-1\right)  }$}}
\put(10,1){\makebox(1,1){$B_{0,1}^{\scriptscriptstyle\left(  g-1\right)  }$}}
\put(11,0){\makebox(1,1){$0$}} \put(11,1){\makebox(1,1){$0$}} \multiput
(0,0)(0,2){2}{\line(1,0){12}} \multiput(0,0)(2,0){3}{\line(0,1){2}}
\multiput(8,0)(2,0){3}{\line(0,1){2}} \multiput(1.05,0.05)(2,0){2}
{\dashbox{0.11176}(1.9,1.9){}} \multiput(7.05,0.05)(2,0){2}{\dashbox
{0.11176}(1.9,1.9){}} \put(0,2.05){\makebox(2,1)[b]{$\overbrace{\hbox
to\bracelength{}}^{\displaystyle A^{\left(  0\right)  }}$}} \put
(2,2.05){\makebox(2,1)[b]{$\overbrace{\hbox to\bracelength{}}^{\displaystyle
A^{\left(  1\right)  }}$}} \put(8,2.05){\makebox(2,1)[b]{$\overbrace{\hbox
to\bracelength{}}^{\displaystyle A^{\left(  g-1\right)  }}$}} \put
(10,2.05){\makebox(2,1)[b]{$\overbrace{\hbox to\bracelength{}}^{\displaystyle
A^{\left(  g\right)  }}$}} \put(1,-1.05){\makebox(2,1)[t]{$\underbrace{\hbox
to\bracelength{}}_{\displaystyle B^{\left(  0\right)  }}$}} \put
(3,-1.05){\makebox(2,1)[t]{$\underbrace{\hbox to\bracelength{}}_{\displaystyle
B^{\left(  1\right)  }}$}} \put(7,-1.05){\makebox(2,1)[t]{$\underbrace{\hbox
to\bracelength{}}_{\displaystyle B^{\left(  g-2\right)  }}$}} \put
(9,-1.05){\makebox(2,1)[t]{$\underbrace{\hbox to\bracelength{}}_{\displaystyle
B^{\left(  g-1\right)  }}$}}
\end{picture}}

\mathcode`\;="303B
\begin{document}
\title{Minimality of the data in wavelet filters}
\author{Palle E. T. Jorgensen}
\address{Department of Mathematics\\
The University of Iowa\\
14 MacLean Hall\\
Iowa City, IA 52242-1419\\
U.S.A.}
\email{jorgen@math.uiowa.edu}
\thanks{Partially supported by the National Science Foundation}
\dedicatory{\upshape with an Appendix by Brian Treadway}

\begin{abstract}
Orthogonal wavelets, or wavelet frames, for $L^{2}\left(  \mathbb{R}\right)  $
are associated with quadrature mirror filters (QMF), a
set of complex numbers which relate the dyadic scaling of functions on
$\mathbb{R}$ to the $\mathbb{Z}$-translates.
In this paper, we show that generically, the data in the QMF-systems of
wavelets is minimal, in the sense that it cannot be nontrivially reduced. The
minimality property is given a geometric formulation in the Hilbert space
$\ell^{2}\left(  \mathbb{Z}\right)  $, and it is then shown that minimality
corresponds to irreducibility of a wavelet representation of the algebra
$\mathcal{O}_{2}$; and so our result is that this family of representations of
$\mathcal{O}_{2}$ on the Hilbert space $\ell^{2}\left(  \mathbb{Z}\right)  $
is irreducible for a generic set of values of the parameters which label the
wavelet representations.
\end{abstract}
\maketitle

\section{\label{Int}Introduction}

\setlength{\displayboxwidth}{\textwidth}\addtolength{\displayboxwidth
}{-2\leftmargini}%
Let $L^{2}\left(  \mathbb{R}\right)  $ be the Hilbert space of all $L^{2}%
$-functions. For $\psi\in L^{2}\left(  \mathbb{R}\right)  $, set
\begin{equation}
\psi_{n,k}\left(  x\right)  :=2^{\frac{n}{2}}\psi\left(  2^{n}x-k\right)
\text{\qquad for }x\in\mathbb{R}\text{, and }n,k\in\mathbb{Z}. \label{eqInt.1}%
\end{equation}
We say that $\psi$ is a wavelet (in \emph{the strict sense\/}) if
$\left\{  \psi_{n,k};n,k\in\mathbb{Z}\right\}  $
constitutes an orthonormal basis in $L^{2}\left(  \mathbb{R}\right)  $; and we
say that $\psi$ is a wavelet in \emph{the frame sense} (tight frame) if%
\begin{equation}
\left\|  f\right\|  _{L^{2}\left(  \mathbb{R}\right)  }^{2}=\sum
_{n,k\in\mathbb{Z}}\left|  \ip{\psi_{n,k}}{f}\right|  ^{2} \label{eqInt.3}%
\end{equation}
holds for all $f\in L^{2}\left(  \mathbb{R}\right)  $, where $\ip{\,\cdot
\,}{\,\cdot\,}$ is the usual $L^{2}\left(  \mathbb{R}\right)  $-inner product,
i.e.,
$\ip{\psi_{n,k}}{f}=\int_{\mathbb{R}}\overline{\psi_{n,k}\left(  x\right)
}f\left(  x\right)  \,dx=c_{n,k}$.  The numbers $c_{n,k}$
are the wavelet coefficients.
It is known \cite{Dau92,Hor95} that a given wavelet $\psi$ in the
sense of frames is a (strict) wavelet if and only if $\left\|  \psi\right\|
_{L^{2}\left(  \mathbb{R}\right)  }=1$.
We shall have occasion to consider scaling on $\mathbb{R}$ other than the
dyadic one, say $x\mapsto Nx$ where $N\in\mathbb{N}$, $N>2$. Then the analogue
of (\ref{eqInt.1}) is%
\begin{equation}
\psi_{n,k}\left(  x\right)  :=N^{\frac{n}{2}}\psi\left(  N^{n}x-k\right)
,\qquad x\in\mathbb{R},\;n,k\in\mathbb{Z}. \label{eqInt.5}%
\end{equation}
However, in that case, it is generally not enough to consider only one
$\psi$ in $L^{2}\left(  \mathbb{R}\right)  $: If the wavelet is derived
from an $N$-subband wavelet filter as in \cite{BrJo00}, then we
construct $\psi^{\left(  1\right)  },\psi^{\left(  2\right)  },\dots
,\psi^{\left(  N-1\right)  }$ in $L^{2}\left(  \mathbb{R}\right)  $ such that
the functions in (\ref{eqInt.5}) have the basis property, either in the strict
sense, or in the sense of frames. Then the system%
\begin{equation}
\left\{  \psi_{n,k}^{\left(  i\right)  };1\leq i<N,\;n,k\in\mathbb{Z}\right\}
\label{eqInt.6}%
\end{equation}
constitutes an orthonormal basis of $L^{2}\left(  \mathbb{R}\right)  $,
or, alternatively, a tight frame, as in (\ref{eqInt.3}) but with the
$\psi_{n,k}^{\left(  i\right)  }$ functions in place of $\psi_{n,k}$.

Our main point is to show how the notion of
\emph{irreducibility} for representations
of the Cuntz algebra $\mathcal{O}_{N}$ corresponds to
\emph{optimality} of the
corresponding wavelet filters. Since we are addressing two different audiences
(wavelets vs.\ representation theory), a few more details are included in this
paper than might otherwise be customary. Our main result is that the
irreducibility of the representation (equivalently, minimality of the filter)
is \emph{generic} for the wavelet representations; see Theorems \ref{ThmMin.7}
and \ref{theorem5.6}. In addition, we show that generically, two different
filters yield inequivalent representations, i.e., the corresponding two
representations are not unitarily equivalent. This was known earlier only in
very restrictive special cases \cite{BrJo00}, and the general case treated here
has not previously been discussed in the literature. Moreover, the methods used
for the special cases in fact do \emph{not at all} carry over to the general
case.
We are concerned with the wavelet filters which enter into the
construction of $\psi^{(1)},\psi^{(2)},\dots,\psi^{(N-1)}$
in (\ref{eqInt.6}).
These filters (see (\ref{eqInt.7})--(\ref{eqInt.9}) and (\ref{eq2.15}) below)
are just a finite set of numbers which relate the $\mathbb{Z}%
$-translates of these functions to the corresponding scalings by $x\mapsto
Nx$. Hence the analysis may be discretized via the filters,
but the question arises whether or not the data which go into the wavelet
filters are \emph{minimal}. Representation theory is ideally
suited to make the minimality question mathematically precise. (This is a
QMF-multiresolution construction, and it is its minimality and efficiency
which concern us here. While it is true, see, e.g.,
\cite{Gab98}, \cite{FPT99}, \cite{Bag00}, \cite{BaMe99}, and \cite{DaLa98},
that there are other and different possible wavelet constructions, it is not
yet clear how our present techniques might adapt to the alternative
constructions, although the approach in \cite{DaLa98} is also based on
operator-theoretic considerations.)

To explain the \emph{minimality issue}
for multiresolution quadrature mirror
(QMF) wavelet filters, we recall the
\emph{scaling function} $\varphi$ of a
resolution in $L^{2}\left(  \mathbb{R}\right)  $.
Let $g\in\mathbb{N}$, and let $a_{0},a_{1},\dots,a_{2g-1}$ be given complex
numbers such that
\begin{equation}
\sum_{k=0}^{2g-1}a_{k}=2, \label{eqInt.7}%
\end{equation}
and%
\begin{equation}
\sum_{k}a_{k+2l}\bar{a}_{k}=%
\begin{cases}
2 & \text{\quad if }l=0,\\
0 & \text{\quad if }l\neq0.
\end{cases}
\label{eqInt.8}%
\end{equation}
In the summation (\ref{eqInt.8}), and elsewhere, we adopt the convention
that terms are defined to be zero when the index is not in the
specified range. Hence, in (\ref{eqInt.8}), it is understood that $a_{k+2l}=0$
whenever $k$ and $l$ are such that $k+2l$ is not in $\left\{  0,1,\dots
,2g-1\right\}  $. It is known
\cite{BrJo00,BEJ00,Mal99} that there is a $\varphi\in L^{2}\left(
\mathbb{R}\right)  \setminus\left\{  0\right\}  $
of compact support,
unique up to a constant
multiple, such that
\begin{equation}
\varphi\left(  x\right)  =\sum_{k=0}^{2g-1}a_{k}\varphi\left(  2x-k\right)
,\qquad x\in\mathbb{R}\mathord{;} \label{eqInt.9}%
\end{equation}
in fact,
$\operatorname*{supp}\left(  \varphi\right)  \subset\left[  0,2g-1\right]  $.
(If $H$ denotes the Hilbert transform of $L^{2}\left(  \mathbb{R}\right)  $,
and $\varphi$ solves (\ref{eqInt.9}), then $H\varphi$ does as well; but
$H\varphi$ will not be of compact support if $\varphi$ is.) In
finding $\varphi$ in (\ref{eqInt.9}), there are methods based on iteration
(see Appendix), on random matrix products, and on Fourier transform, see
\cite{BrJo00}, \cite{BEJ00}, \cite{BrJo99b}, \cite{Coh92b}, \cite{CoRy95}, and
\cite{Dau92}; and the various methods intertwine in the analysis of $\varphi$,
i.e., in deciding when $\varphi\left(  x\right)  $ is continuous, or not, or
if it is differentiable. This issue will be resumed in the Appendix below,
which is based on \cite{BrJo99b}. But the next two sections will deal with the
minimality question alluded to above.

Let $\varphi$ be as in (\ref{eqInt.9}), and let $\mathcal{V}_{0}$ be the
closed subspace in $\mathcal{H}$ ($:=L^{2}\left(  \mathbb{R}\right)  $)
spanned by $\left\{  \varphi\left(  x-k\right)  ;k\in\mathbb{Z}\right\}  $,
i.e., by the integral translates of the scaling function $\varphi$. Let $U$
($:=U_{N}$) be
\begin{equation}
Uf\left(  x\right)  :=N^{-\frac{1}{2}}f\left(  \frac{x}{N}\right)  ,\qquad
f\in L^{2}\left(  \mathbb{R}\right)  , \label{eqInt.11}%
\end{equation}
the unitary scaling operator in $\mathcal{H}=L^{2}\left(  \mathbb{R}\right)
$. Then if $N=2$,
\begin{equation}
U\mathcal{V}_{0}\subset\mathcal{V}_{0} \label{eqInt.12}%
\end{equation}
is a proper subspace, and
\begin{equation}
\bigwedge_{n}U^{n}\mathcal{V}_{0}=\left\{  0\right\}  ; \label{eqInt.13}%
\end{equation}
see \cite{BEJ00} and \cite[Ch.~5]{Dau92}. Setting
$\mathcal{V}_{n}:=U^{n}\mathcal{V}_{0}$
and%
\begin{equation}
\mathcal{W}_{n}:=\mathcal{V}_{n-1}\ominus\mathcal{V}_{n}, \label{eqInt.15}%
\end{equation}
we arrive at the resolution%
\begin{equation}
\mathcal{V}_{0}=\sideset{}{^{\smash{\oplus}}}{\sum}_{n\geq1}\mathcal{W}_{n},
\label{eqInt.16}%
\end{equation}
and the wavelet function $\psi$ is picked in $\mathcal{W}_{0}$;
see Table \ref{tablegeomappr}. We will set up
an isomorphism between the resolution subspace $\mathcal{V}_{0}$ and $\ell
^{2}\left(  \mathbb{Z}\right)  $, and associate operators in $\ell^{2}\left(
\mathbb{Z}\right)  $ with the wavelet operations in $\mathcal{V}_{0}\subset
L^{2}\left(  \mathbb{R}\right)  $. This is of practical significance given
that the operators in $\ell^{2}\left(  \mathbb{Z}\right)  $ are those which
are defined directly from the wavelet filters, and it is the digital filter
operations which lend themselves to algorithms.
Generalizing (\ref{eqInt.15}), in the case of scale $N$ ($>2$) the space $\mathcal{V}_{0}\ominus
U_{N}\mathcal{V}_{0}$ splits up as a sum of orthogonal spaces $\mathcal{W}_{1}
^{\left(  i\right)  }$, $i=1,2,\dots,N-1$;
see (\ref{eq2.25})--(\ref{eq2.25bis}).

\begin{table}[tbp]
\caption{Discrete vs.\ continuous wavelets, i.e., $\ell\sp2$ vs.\ $L\sp
2\left(  \mathbb{R}\right)  $}%
\label{tablegeomappr}
$\renewcommand{\arraystretch}{0.25}%
\begin{tabular}
[c]{ccccrrrcl}%
$\left\{  0\right\}  $ & $\longleftarrow$ & $\cdots$ & $\longleftarrow$ &
$\mathcal{V}_{2}\raisebox{-12pt}{$\searrow\vphantom{\raisebox{-2pt}{$\searrow
$}}$\hskip-8pt}$ & $\mathcal{V}_{1}\raisebox{-12pt}{$\searrow\vphantom
{\raisebox{-2pt}{$\searrow$}}$\hskip-8pt}$ & $\mathcal{V}_{0}\raisebox
{-12pt}{$\searrow\vphantom{\raisebox{-2pt}{$\searrow$}}$\hskip-8pt}$ &  &
finer scales\\\hline
&  &  &  & \multicolumn{1}{c}{} & \multicolumn{1}{c}{} & \multicolumn{1}{c}{}
&  & \\\cline{1-7}
&  &  &  & \multicolumn{1}{c}{} & \multicolumn{1}{c}{} & \multicolumn{1}{c}{}
& \multicolumn{1}{|c}{} & \\\cline{1-6}
&  &  &  & \multicolumn{1}{c}{} & \multicolumn{1}{c}{} & \multicolumn{1}{|c}{}
& \multicolumn{1}{|c}{} & \\\cline{1-5}
&  &  &  & \multicolumn{1}{c}{} & \multicolumn{1}{|c}{} &
\multicolumn{1}{|c}{} & \multicolumn{1}{|c}{} & \\\cline{1-4}
&  & $\cdots$ &  & \multicolumn{1}{|c}{$\mathcal{W}_{3}$} &
\multicolumn{1}{|c}{$\mathcal{W}_{2}$} & \multicolumn{1}{|c}{$\mathcal{W}_{1}%
$} & \multicolumn{1}{|c}{$\cdots$} & rest of $L^{2}\left(  \mathbb{R}\right)
$\\\cline{1-4}
&  &  &  & \multicolumn{1}{c}{} & \multicolumn{1}{|c}{} &
\multicolumn{1}{|c}{} & \multicolumn{1}{|c}{} & \\\cline{1-5}
&  &  &  & \multicolumn{1}{c}{} & \multicolumn{1}{c}{} & \multicolumn{1}{|c}{}
& \multicolumn{1}{|c}{} & \\\cline{1-6}
&  &  &  & \multicolumn{1}{c}{} & \multicolumn{1}{c}{} & \multicolumn{1}{c}{}
& \multicolumn{1}{|c}{} & \\\cline{1-7}
&  &  &  & \multicolumn{1}{c}{} & \multicolumn{1}{c}{} & \multicolumn{1}{c}{}
&  & \\\hline
\rule[6pt]{0pt}{6pt} &  & $\cdots$ & \multicolumn{1}{r}{$\rlap{$%
\underset{\textstyle U}{\longleftarrow}$}\hskip3pt$} & $\rlap{$\underset
{\textstyle U}{\longleftarrow}$}\hskip3pt$ & $\rlap{$\underset{\textstyle
U}{\longleftarrow}$}\hskip3pt$ & \multicolumn{1}{c}{} &  & \\
&  & $\llap{$W$}\makebox[6pt]{\raisebox{4pt}{\makebox[0pt]{\hss$\uparrow$\hss
}}\raisebox{-4pt}{\makebox[0pt]{\hss$|$\hss}}}$ &  &
\multicolumn{1}{c}{$\makebox[6pt]{\raisebox{4pt}{\makebox[0pt]{\hss$%
\uparrow$\hss}}\raisebox{-4pt}{\makebox[0pt]{\hss$|$\hss}}}$} &
\multicolumn{1}{c}{$\makebox[6pt]{\raisebox{4pt}{\makebox[0pt]{\hss$%
\uparrow$\hss}}\raisebox{-4pt}{\makebox[0pt]{\hss$|$\hss}}}$} &
\multicolumn{1}{c}{$\makebox[6pt]{\raisebox{4pt}{\makebox[0pt]{\hss$%
\uparrow$\hss}}\raisebox{-4pt}{\makebox[0pt]{\hss$|$\hss}}}\rlap{$W$}$} &  &
\\
\rule[-6pt]{0pt}{6pt}$\left\{  0\right\}  $ & $\longleftarrow$ & $\cdots$ &
\multicolumn{1}{r}{$\rlap{$\overset{\textstyle S_0}{\longleftarrow}$}%
\hskip3pt$} & $\rlap{$\overset{\textstyle S_0}{\longleftarrow}$}\hskip3pt$ &
$\rlap{$\overset{\textstyle S_0}{\longleftarrow}$}\hskip3pt$ &
\multicolumn{1}{c}{} &  & \\\cline{1-7}
&  &  &  & \multicolumn{1}{c}{} & \multicolumn{1}{c}{} & \multicolumn{1}{c}{}
& \multicolumn{1}{|c}{} & \\\cline{1-6}
&  &  &  & \multicolumn{1}{c}{} & \multicolumn{1}{c}{} & \multicolumn{1}{|c}{}
& \multicolumn{1}{|c}{} & \\\cline{1-5}
&  &  &  & \multicolumn{1}{c}{} & \multicolumn{1}{|c}{} &
\multicolumn{1}{|c}{} & \multicolumn{1}{|c}{} & \\\cline{1-4}
&  & $\cdots$ &  & \multicolumn{1}{|c}{$S_{0}^{2}\mathcal{L}$} &
\multicolumn{1}{|c}{$S_{0}\mathcal{L}$} & \multicolumn{1}{|c}{$\mathcal{L}%
=S_{1}\ell^{2}$} & \multicolumn{1}{|c}{} & \\\cline{1-4}
&  &  &  & \multicolumn{1}{c}{} & \multicolumn{1}{|c}{} &
\multicolumn{1}{|c}{} & \multicolumn{1}{|c}{} & \\\cline{1-5}
&  &  &  & \multicolumn{1}{c}{} & \multicolumn{1}{c}{} & \multicolumn{1}{|c}{}
& \multicolumn{1}{|c}{} & \\\cline{1-6}
&  &  &  & \multicolumn{1}{c}{} & \multicolumn{1}{c}{} & \multicolumn{1}{c}{}
& \multicolumn{1}{|c}{} & \\\cline{1-7}
&  &  &  & \makebox[36pt]{\hfill}\llap{$S_{0}^{2}\ell^{2}$}\raisebox
{12pt}{$\nearrow\vphantom{\raisebox{2pt}{$\nearrow$}}$\hskip-8pt} &
\makebox[36pt]{\hfill}\llap{$S_{0}\ell^{2}$}\raisebox{12pt}{$\nearrow
\vphantom{\raisebox{2pt}{$\nearrow$}}$\hskip-8pt} & $\ell^{2}\raisebox{12pt}%
{$\nearrow\vphantom{\raisebox{2pt}{$\nearrow$}}$\hskip-8pt}$ &  &
\end{tabular}
\ $\end{table}

\section{\label{Rep}Representations of $\mathcal{O}_{N}$ and Table
\ref{tablegeomappr} (Discrete vs.\ continuous wavelets)}

The computational significance of the operator system in
Table \ref{tablegeomappr}
(scale $N=2$) is that the operators which generate wavelets in $L^{2}\left(
\mathbb{R}\right)  $ become modeled by an associated system of operators in
the \emph{sequence space} $\ell^{2}$ ($:=\ell^{2}\left(  \mathbb{Z}\right)
\cong L^{2}\left(  \mathbb{T}\right)  $). (We will do the discussion here in
Section \ref{Rep} just for $N=2$, but this is merely for simplicity; it easily
generalizes to arbitrary $N$.) Then the algorithms are implemented in
$\ell^{2}$ by basic discrete operations, and only in the end are the results
then ``translated'' back to the space $L^{2}(\mathbb{R})$. The space
$L^{2}(\mathbb{R})$ is not amenable (in its own right) to \emph{discrete}
computations. This is made precise by the frame operator
$W\colon\ell^{2}$ ($\cong L^{2}\left(  \mathbb{T}\right)  $) $\rightarrow
\mathcal{V}_{0}$ ($\subset L^{2}\left(  \mathbb{R}\right)  $)
defined as%
\begin{equation}
W\colon\ell^{2}\ni\left(  \xi_{k}\right)  \longmapsto\sum_{k\in\mathbb{Z}}%
\xi_{k}\varphi\left(  x-k\right)  \in L^{2}\left(  \mathbb{R}\right)  .
\label{eqRep.1}%
\end{equation}
If $\varphi$ has orthogonal
translates, then $W$ will be an isometry of $\ell^{2}$ onto $\mathcal{V}_{0}$
($\subset L^{2}\left(  \mathbb{R}\right)  $). Even if the functions $\left\{
\varphi\left(  x-k\right)  \right\}  _{k\in\mathbb{Z}}$ formed from $\varphi$
by $\mathbb{Z}$-translates only constitute a frame in $\mathcal{V}_{0}$, then
we will have the following estimates:%
\begin{equation}
c_{1}^{1/2}\cdot\left\|  \xi\right\|  _{\ell^{2}}\leq\left\|  W\xi\right\|
_{L^{2}\left(  \mathbb{R}\right)  }\leq c_{2}^{1/2}\cdot\left\|  \xi\right\|
_{\ell^{2}}\,, \label{eqRep.3}%
\end{equation}
where $c_{1}$ and $c_{2}$ are positive constants depending only on $\varphi$.

\begin{lemma}
\label{LemRep.1}If the coefficients $\left\{  a_{k};k=0,1,\dots,2g-1\right\}
$ from \textup{(\ref{eqInt.9})} satisfy the conditions in
\textup{(\ref{eqInt.8}),} then the corresponding operator $S_{0}\colon\ell
^{2}\rightarrow\ell^{2}$, given by%
\begin{equation}
\left(  S_{0}\xi\right)  _{k}=\frac{1}{\sqrt{2}}\sum_{l\in\mathbb{Z}}%
a_{k-2l}\xi_{l}=\frac{1}{\sqrt{2}}\sum_{\substack{p\in\mathbb{Z}%
\colon\\p\equiv k\operatorname{mod}2}}a_{p}\xi_{\frac{k-p}{2}},\qquad
k\in\mathbb{Z}, \label{eqRep.4}%
\end{equation}
is \emph{isometric} and satisfies the following \emph{intertwining identity:}%
\begin{equation}
WS_{0}=UW, \label{eqRep.5}%
\end{equation}
where $U$ is the dyadic scaling operator in $L^{2}\left(  \mathbb{R}\right)  $
introduced in \textup{(\ref{eqInt.11}).} \textup{(}Here we restrict
attention to $N=2$, but just for notational simplicity!\/\textup{)} Setting
$b_{k}:=\left(  -1\right)  ^{k}\bar{a}_{2g-1-k}$,
and defining a second isometric operator $S_{1}\colon\ell^{2}\rightarrow
\ell^{2}$ by formula \textup{(\ref{eqRep.4})} with the only modification that
$\left(  b_{k}\right)  $ is used in place of $\left(  a_{k}\right)  $, we get%
\begin{equation}
S_{j}^{\ast}S_{k}^{{}}=\delta_{j,k}\openone_{\ell^{2}} \label{eqRep.7}%
\end{equation}
and%
\begin{equation}
\sum_{j}S_{j}^{{}}S_{j}^{\ast}=\openone_{\ell^{2}}\,, \label{eqRep.8}%
\end{equation}
which are the Cuntz identities from operator theory \cite{Cun77}, and the
operators $S_{0}$ and $S_{1}$ satisfy the identities indicated in Table
\textup{\ref{tablegeomappr}.}
\end{lemma}

\begin{remark}
\label{RemRepNew.pound}For understanding the second line in Table
\textup{\ref{tablegeomappr},} note that $S_{0}$ is a \emph{shift} as an
isometry, in the sense of \cite{SzFo70}, and $\mathcal{L}:=S_{1}\ell^{2}$ is a
wandering subspace for $S_{0}$, in the sense that the spaces $\mathcal{L}$,
$S_{0}^{{}}\mathcal{L}$, $S_{0}^{2}\mathcal{L},\;\dots$ are mutually
orthogonal in $\ell^{2}$. To see this, note that \textup{(\ref{eqRep.8})}
implies that
($\mathcal{L}:=$) $S_{1}\ell^{2}=\ell^{2}\ominus S_{0}
\ell^{2}=\ker\left(  S_{0}^{\ast}\right) $.
As a result, we get the following:
\end{remark}

\begin{corollary}
\label{corollary2new}The projections onto the orthogonal subspaces in the
second line of Table $1$ corresponding to the $\mathcal{W}_{1}%
,\mathcal{W}_{2},\dots$ subspaces of the first line \textup{(}see
\textup{(\ref{eqInt.15}))} are
\begin{align*}
\operatorname*{proj}\mathcal{L}  &  =S_{1}S_{1}^{\ast}=I-S_{0}S_{0}^{\ast},\\
\vdots &  \qquad\vdots\\
\operatorname*{proj}S_{0}^{n-1}\mathcal{L}  &  =S_{0}^{n-1}S_{0}^{\ast
n-1}-S_{0}^{n}S_{0}^{\ast n}.
\end{align*}
\end{corollary}

\begin{proof}
Immediate from Lemma \ref{LemRep.1}, Remark \ref{RemRepNew.pound},
and (\ref{eqInt.16}).
\end{proof}

\begin{remark}
\label{RemRep.2}Any system of operators $\left\{  S_{j}\right\}  $ satisfying
\textup{(\ref{eqRep.7})--(\ref{eqRep.8})} is said to be a
\emph{representation} of
the $C^{\ast}$-algebra $\mathcal{O}_{2}$, and there is a similar notion for
$\mathcal{O}_{N}$ when $N>2$, with $\mathcal{O}_{N}$ having generators
$S_{0},S_{1},\dots,S_{N-1}$, but otherwise also satisfying the operator
identities \textup{(\ref{eqRep.7})--(\ref{eqRep.8}).}
The power and the usefulness of the multiresolution subband filters for the
analysis of wavelets and their algorithms was first demonstrated forcefully in
\cite{CoWi93} and \cite{Wic93}; see especially \cite[p.~140]{CoWi93} and
\cite[p.~157]{Wic93}, where the $\mathcal{O}_{N}$-relations
(\ref{eqRep.7})--(\ref{eqRep.8})
are identified, and analyzed in the case $N=2$. Around the
same time, A. Cohen \cite{Coh92b} identified and utilized the interplay
between $\ell^{2}$ and $L^{2}\left(  \mathbb{R}\right)  $ which, as noted in
Section \ref{Rep} above, is implied by the $\mathcal{O}_{N}$-relations and
their representations. But neither of those prior references takes up the
construction of $\mathcal{O}_{N}$-representations in a systematic fashion.
Of course the quadrature mirror filters (QMF's) have a long history in
electrical engineering (speech coding problems), going back to long before
they were used in wavelets, but the form in which we shall use them here is
well articulated, for example, in \cite{CEG77}. Some more of the history of
and literature on wavelet filters is covered well in \cite{Mey93} and
\cite{Ben00}.
\end{remark}

\begin{definition}
\label{DefRep.3}A representation of $\mathcal{O}_{N}$ on the Hilbert space
$\ell^{2}$ is said to be \emph{irreducible} if there are no closed subspaces
$\left\{  0\right\}  \subsetneqq\mathcal{H}_{0}\subsetneqq\ell^{2}$ which
reduce the representation, i.e., which yield a representation of
\textup{(\ref{eqRep.7})--(\ref{eqRep.8})} on each of the two subspaces in the
decomposition%
\begin{equation}
\ell^{2}=\mathcal{H}_{0}\oplus\left(  \ell^{2}\ominus\mathcal{H}_{0}\right)  ,
\label{eqRep.9}%
\end{equation}
where
$\ell^{2}\ominus\mathcal{H}_{0}=\left(  \mathcal{H}_{0}\right)  ^{\perp
}=\left\{  \xi\in\ell^{2};\ip{\xi}{\eta}=0,\;\forall\,\eta\in\mathcal{H}
_{0}\right\}  $.
\end{definition}

\begin{proof}
[Proof of Lemma \textup{\ref{LemRep.1}}]Most of the details of the proof are
contained in \cite{BrJo97b} and \cite{BrJo00}, so we only sketch points not
already covered there. The essential step (for the present applications) is
the formula (\ref{eqRep.5}), which shows that $W$ intertwines the isometry
$S_{0}\,$with the restriction of the unitary operator $U\colon f\mapsto
\frac{1}{\sqrt{2}}f\left(  x/2\right)  $ to the resolution subspace
$\mathcal{V}_{0}\subset L^{2}\left(  \mathbb{R}\right)  $. We have:%
\begin{align*}
\left(  UW\xi\right)  \left(  x\right)   &  =\frac{1}{\sqrt{2}}\left(
W\xi\right)  \left(  \frac{x}{2}\right)  &  & \\
&  =\frac{1}{\sqrt{2}}\sum_{k\in\mathbb{Z}}\xi_{k}\varphi\left(  \frac{x}%
{2}-k\right)  &  &  \text{\qquad(by (\ref{eqRep.1}))}\\
&  =\frac{1}{\sqrt{2}}\sum_{k\in\mathbb{Z}}\sum_{l\in\mathbb{Z}}\xi_{k}%
a_{l}\varphi\left(  x-2k-l\right)  &  &  \text{\qquad(by (\ref{eqInt.9}))}\\
&  =\frac{1}{\sqrt{2}}\sum_{p\in\mathbb{Z}}\left(  \sum_{k\in\mathbb{Z}}%
\xi_{k}a_{p-2k}\right)  \varphi\left(  x-p\right)  &  & \\
&  =\sum_{p\in\mathbb{Z}}\left(  S_{0}\xi\right)  _{p}\varphi\left(
x-p\right)  &  &  \text{\qquad(by (\ref{eqRep.4}))}\\
&  =\left(  WS_{0}\xi\right)  \left(  x\right)  &  &  \text{\qquad(by
(\ref{eqRep.1}))}%
\end{align*}
for all $\xi\in\ell^{2}$, and all $x\in\mathbb{R}$. This proves (\ref{eqRep.5}).
\end{proof}

For later use, we record the operators on the
respective Hilbert spaces $L^{2}\left(  \mathbb{T}\right)  \cong\ell^{2}$ and
$L^{2}\left(  \mathbb{R}\right)  $, and the corresponding transformation rules
with respect to the operator $W$. Let $N$ be the scale number, and let
$\left(  a_{k}\right)  _{k=0}^{Ng-1}$ be given satisfying%
\begin{equation}
\sum_{k\in\mathbb{Z}}a_{k+Nl}\bar{a}_{k}=\delta_{0,l}N \label{eqRep.18}%
\end{equation}
and set $m_{0}\left(  z\right)  :=\frac{1}{\sqrt{N}}\sum_{k=0}^{Ng-1}%
a_{k}z^{k}$, $z\in\mathbb{T}$. 
The following summary table of transformation
rules may clarify the proof.
\begin{equation}%
\begin{array}
[c]{llll}
& \quad\;\text{\textbf{SCALING}} & \quad\;\text{\textbf{TRANSLATION}} & \\
\displaystyle\vphantom{\frac{1}{N}}L^{2}\left(  \mathbb{R}\right)  \colon &
\displaystyle\quad F\longmapsto\smash{\frac{1}{\sqrt{N}}}F\left(  \frac{x}%
{N}\right)  & \displaystyle\quad F\left(  x\right)  \longmapsto F\left(
x-1\right)  & \text{\quad real wavelets}\\
\displaystyle\uparrow\rlap{$\scriptstyle W$} &  &  & \\
\displaystyle\vphantom{\frac{1}{N}}\ell^{2}\colon & \displaystyle\quad
\xi\longmapsto\smash{\sum_{l}}a_{k-Nl}\xi_{l} & \displaystyle\quad\left(
\xi_{k}\right)  \longmapsto\left(  \xi_{k-1}\right)  & \text{\quad discrete
model}\\
\displaystyle\uparrow\rlap{\Small Fourier transform} &  &  & \\
\displaystyle\vphantom{\frac{1}{N}}L^{2}\left(  \mathbb{T}\right)  \colon &
\displaystyle\quad f\longmapsto m_{0}\left(  z\right)  f\left(  z^{N}\right)
& \displaystyle\quad f\left(  z\right)  \longmapsto zf\left(  z\right)  &
{\normalsize \quad\setlength{\arraycolsep}{0pt}%
\begin{array}
[t]{l}%
\text{periodic model,}\\
\quad\mathbb{T}=\mathbb{R}/2\pi\mathbb{Z}%
\end{array}
}%
\end{array}
\label{eqRep.19}%
\end{equation}

\begin{remark}
\label{RemRep.4}The significance of irreducibility \textup{(}when
satisfied\/\textup{)} is that the \emph{wavelet subbands} which are indicated
in Table \textup{\ref{tablegeomappr}} are then the \emph{only subbands} of the
corresponding multiresolution. We will show that in fact irreducibility holds
\emph{generically,} but it does not hold, for example, for the Haar wavelets.
In the simplest case, the Haar wavelet has $g=2=N$, and the numbers
from Lemma \textup{\ref{LemRep.1}} are
\begin{equation}%
\begin{pmatrix}
a_{0} & a_{1}\\
b_{0} & b_{1}%
\end{pmatrix}
=%
\begin{pmatrix}
1 & 1\\
1 & -1
\end{pmatrix}
. \label{eqRep.11}%
\end{equation}
Hence, for this representation of $\mathcal{O}_{2}$ on $\ell^{2}$, we may take
$\mathcal{H}_{0}=\ell^{2}\left(  0,1,2,\dots\right)  $, and therefore
$\mathcal{H}_{0}^{\perp}=\ell^{2}\left(  \dots,-3,-2,-1\right)  $. Returning
to the multiresolution diagram in Table \textup{\ref{tablegeomappr},} this
means that we get additional subspaces of $L^{2}\left(  \mathbb{R}\right)  $,
on top of the standard ones which are listed in Table
\textup{\ref{tablegeomappr}.} Specifically, in addition to
\[
\mathcal{V}_{n}=U^{n}\mathcal{V}_{0}=WS_{0}^{n}\ell^{2}
\text{\quad and\quad }
\mathcal{W}_{n}=\mathcal{V}_{n-1}\ominus\mathcal{V}_{n}
=WS_{0}^{n-1}S_{1}^{{}}\ell^{2},
\]
we get a new system with ``twice as many'', as follows: $\mathcal{V}%
_{n}^{\left(  \pm\right)  }$ and $\mathcal{W}_{n}^{\left(  \pm\right)  }$,
where
\[
\mathcal{V}_{n}^{\left(  +\right)  }    =WS_{0}^{n}\left(  \mathcal{H}
_{0}\right)  ,\qquad
\mathcal{W}_{n}^{\left(  +\right)  }    =WS_{0}^{n-1}S_{1}^{{}}\left(
\mathcal{H}_{0}\right)  ;
\]
and
\[
\mathcal{V}_{n}^{\left(  -\right)  }    =WS_{0}^{n}\left(  \mathcal{H}%
_{0}^{\perp}\right)  ,\qquad
\mathcal{W}_{n}^{\left(  -\right)  }    =WS_{0}^{n-1}S_{1}^{{}}\left(
\mathcal{H}_{0}^{\perp}\right)  .
\]
For the case of the Haar wavelet, see \textup{(\ref{eqRep.11}),}%
\[
\mathcal{V}_{0}^{\left(  +\right)  }\subset L^{2}\left(  0,\infty\right)
,\qquad\mathcal{V}_{0}^{\left(  -\right)  }\subset L^{2}\left(  -\infty
,0\right)  ,
\]
or rather, $\mathcal{V}_{0}$ consists of finite linear combinations of
$\mathbb{Z}$-translates of
\begin{equation}
\varphi\left(  x\right)  =%
\begin{cases}
1 & \text{ if }0\leq x<1,\\
0 & \text{ if }x\in\mathbb{R}\setminus\left[  0,1\right)  ,
\end{cases}
\label{eqRep.16}%
\end{equation}
i.e., functions in $L^{2}\left(
\mathbb{R}\right)  $ which are constant between $n$ and $n+1$ for all
$n\in\mathbb{Z}$; and%
\begin{equation}
\mathcal{V}_{0}^{\left(  +\right)  }=\mathcal{V}_{0}^{{}}\cap L^{2}\left(
0,\infty\right)  ,\qquad\mathcal{V}_{0}^{\left(  -\right)  }=\mathcal{V}%
_{0}^{{}}\cap L^{2}\left(  -\infty,0\right)  . \label{eqRep.17}%
\end{equation}
Hence we get two separate wavelets, but with translations built on $\left\{
0,1,2,\dots\right\}  $ and $\left\{  \dots,-3,-2,-1\right\}  $. In view of the
graphics in the Appendix below, it is perhaps surprising that other wavelets
\textup{(}different from the Haar wavelets\/\textup{)} do not have the
corresponding additional ``positive vs.\ negative'' splitting into subbands
within the Hilbert space $L^{2}\left(  \mathbb{R}\right)  $.
\end{remark}

\begin{remark}
\label{RemRep.5}There are other dyadic Haar wavelets
(mock Haar wavelets), in addition to
\textup{(\ref{eqRep.16}).} For example, let%
\begin{equation}
\varphi_{k}\left(  x\right)  =%
\begin{cases}
\frac{1}{\sqrt{2k+1}} & \text{ if }0\leq x<2k+1,\\
0 & \text{ if }x\in\mathbb{R}\setminus\left[  0,2k+1\right)  .
\end{cases}
\label{eqRep.17bis}%
\end{equation}
Then it follows that there is a splitting of $\mathcal{V}_{0}$ into orthogonal
subspaces which is analogous to \textup{(\ref{eqRep.17}),} but it has many
more subbands than the two, ``positive vs.\ negative'',
which are special to the standard Haar wavelet
\textup{(\ref{eqRep.16}).} For details on these other Haar wavelets, and their
decompositions, see \cite[Proposition 8.2]{BrJo99a}.
They are only tight frames, and the
$m$-functions of \textup{(\ref{eqRep.17bis})} are%
\begin{equation}
m_{0}\left(  z\right)  =\frac{1}{\sqrt{2}}\left(  1+z^{2k+1}\right)  ,\qquad
m_{1}\left(  z\right)  =\frac{1}{\sqrt{2}}\left(  1-z^{2k+1}\right)  ,\qquad
z\in\mathbb{T}. \label{eqRep.17ter}%
\end{equation}
Hence, after adjusting the $\mathcal{O}_{2}$-representation $T$ with a
rotation $V\in\mathrm{U}_{2}\left(  \mathbb{C}\right)
$, we have%
\begin{equation}
T_{0}f\left(  z\right)  =f\left(  z^{2}\right)  ,\qquad T_{1}f\left(
z\right)  =z^{2k+1}f\left(  z^{2}\right)  ,\qquad f\in L^{2}\left(
\mathbb{T}\right)  \cong\ell^{2}, \label{eqRep.17tetra}%
\end{equation}
and the two new operators $T_{0},T_{1}$ will satisfy the
$\mathcal{O}_{2}$-identities (\ref{eqRep.7})--(\ref{eqRep.8}); the
representation will have the same reducing subspaces as the one
defined directly from $m_{0}$ and $m_{1}$. The explicit decomposition of the
multiresolution subspaces corresponding to \textup{(\ref{eqRep.17})} may be
derived, via $W$ in Table \textup{\ref{tablegeomappr},} from the
decomposition into sums of irreducibles for the $\mathcal{O}_{2}%
$-representation on $\ell^{2}$ which corresponds to \textup{(\ref{eqRep.17}).}
This means that the decomposition \textup{(\ref{eqRep.9})} associated
with \textup{(\ref{eqRep.17bis})} and \textup{(\ref{eqRep.17tetra})} has
\emph{more than two} terms in its subspace configuration.
\end{remark}

\section{\label{Wav}Wavelet filters and subbands}

The operators of wavelet filters may be realized on either one
of the two Hilbert spaces $\ell^{2}(\mathbb{Z})$ or $L^{2}(\mathbb{T})$,
$\mathbb{T}=\mathbb{R}/2\pi\mathbb{Z}$, and $L^{2}(\mathbb{T})$ defined from
the normalized Haar measure $\mu$ on $\mathbb{T}$. But, of course, $\ell
^{2}(\mathbb{Z})$ $\cong L^{2}(\mathbb{T})$ via the Fourier series.
For a given sequence $a_{0},a_{1},\dots,a_{Ng-1}$, consider the operator
$S_{0}$ in $\ell^{2}(\mathbb{Z})$ given by%

\begin{equation}
\xi\longmapsto S_{0}\xi \text{\quad and\quad }
(S_{0}\xi)_{k}=\frac{1}{\sqrt{N}}\sum_{l}a_{k-lN}\xi_{l}. \label{eq2.4}%
\end{equation}
Setting
$m_{0}(z)=\frac{1}{\sqrt{N}}\sum_{k=0}^{Ng-1}a_{k}z^{k}$
and%
\begin{equation}
(\hat{S}_{0}f)(z)=m_{0}(z)f(z^{N}),\qquad f\in L^{2}(\mathbb{T}),\ 
z\in\mathbb{T}, \label{eq2.6}%
\end{equation}
we note that $S_{0}$ and $\hat{S}_{0}$ are really two versions of the same
operator, i.e.,
that $(\hat{S}_{0}f)\sphat=S_{0}(\hat{f})$
when $\hat{f}=(\xi_{k})$ from
the Fourier series.
(The first one is the discrete
model, and the second, the periodic model, referring to the diagram
(\ref{eqRep.19}).) Hence, we shall simply use the same notation $S_{0}$ in
referring to this operator in either one of its incarnations. It is the
(\ref{eq2.4}\textbf{) }version which is used in algorithms, of course.

Let $\varphi\in L^{2}(\mathbb{R})$ be the compactly supported scaling function
solving
\begin{equation}
\varphi(x)=\sum_{k=0}^{Ng-1}a_{k}\varphi(Nx-k). \label{eq2.8}%
\end{equation}
Then define the operator
$W\colon\ell^{2}(\mathbb{Z})\rightarrow L^{2}(\mathbb{R})$
by (\ref{eqRep.1}).
The conditions on the wavelet filter $\{a_{k}\}$ in
(\ref{eqInt.7})--(\ref{eqInt.8}) and
(\ref{eqRep.18}) may now be restated in terms of $m_{0}(z)$ in
(\ref{eq2.6}) as follows:
\begin{equation}
\sum_{k=0}^{N-1}\left|  m_{0}(ze^{i\frac{k2\pi}{N}})\right|  ^{2}=N,\text{ }
\label{eq2.11}%
\end{equation}
and%
\begin{equation}
m_{0}(1)=\sqrt{N}. \label{eq2.12}%
\end{equation}
It then follows from Lemma \ref{LemRep.1} that $W$ in (\ref{eqRep.1}) maps
$\ell^{2}(\mathbb{Z})$ onto the resolution subspace $\mathcal{V}_{0}$
($\subset
L^{2}(\mathbb{R})$), and that
\begin{equation}
U_{N}W=WS_{0} \label{eq2.13}%
\end{equation}
where
$U_{N}f(x)=N^{-1/2}f\left(  x/N\right)  $, $f\in
L^{2}(\mathbb{R})$, $x\in\mathbb{R}$.
We showed in \cite{BrJo00} that there are
$L^{\infty}$-functions $m_{1},\dots,m_{N-1}$ such
that the $N$-by-$N$ complex matrix
\begin{equation}
\frac{1}{\sqrt{N}}\left(  m_{j}(e^{i\frac{k2\pi}{N}}z)\right)  _{j,k=0}^{N-1}
\label{eq2.15}%
\end{equation}
is unitary for all $z\in\mathbb{T}$. If we define
\begin{equation}
S_{j}f(z)=m_{j}(z)f(z^{N}),\qquad f\in L^{2}(\mathbb{T}),\qquad z\in
\mathbb{T}\text{, } \label{eq2.16}%
\end{equation}
then
\begin{equation}
S_{j}^{\ast}S_{k}=\delta_{j,k}I_{L^{2}(\mathbb{T})}, \label{eq2.17}%
\end{equation}
and%
\begin{equation}
\sum_{j=0}^{N-1}S_{j}S_{j}^{\ast}=I_{L^{2}(\mathbb{T})}. \label{eq2.18}%
\end{equation}
((\ref{eq2.12}) is not needed for this, only for the algorithmic
operations of the Appendix.)

\begin{lemma}
\label{LemWavNew.1}
The solutions $(m_{j}%
)_{j=0}^{N-1}$ to
\textup{(\ref{eq2.15})}
are in $1$--$1$ correspondence with the semigroup
of all polynomial functions
\begin{equation}
A\colon\mathbb{T}\longrightarrow\mathrm{U}_{N}(\mathbb{C}), \label{eq2.19}%
\end{equation}
where $\mathrm{U}_{N}(\mathbb{C})$ denotes the unitary $N\times N$ matrices.
\end{lemma}

\begin{proof}
The correspondence is $m\leftrightarrow A$ with
\begin{equation}
m_{j}(z)=\sum_{k=0}^{N-1}A_{j,k}(z^{N})z^{k}, \label{eq2.20}%
\end{equation}
and in the reverse direction,
\begin{equation}
A_{j,k}(z)=\frac{1}{N}\sum_{w^{N}=z}w^{-k}m_{j}(w) \label{eq2.21}%
\end{equation}
does the job, as can be checked by direct substitution.
\end{proof}

We also showed
in \cite{BrJo00} that if $m_{0}$ is given, and if it satisfies (\ref{eq2.11}),
then it is possible to construct $m_{1},\dots,m_{N-1}$ such that the extended
system $m_{0},m_{1},\dots,m_{N-1}$ will satisfy (\ref{eq2.15}). As a
consequence, $A$ in (\ref{eq2.21}) will be a $\mathrm{U}_{N}(\mathbb{C}%
)$-loop, and
the original $m_{0}$ is then recovered from (\ref{eq2.20}) for $j=0$.
To stress the dependence of the operators in (\ref{eq2.16}) on the loop group
element $A$ we will denote the corresponding operators $T_{i}^{\left(
A\right)  }$, and it follows that, if $A=\openone_{N}$, then the operators
$S_{i}$ of (\ref{eq2.16}) are
\begin{equation}
f\left(  z\right)  \longmapsto z^{i}f\left(  z^{N}\right)  \text{,\qquad where
}i=0,1,\dots,N-1, \label{eqMin.6}%
\end{equation}
and we will reserve the notation $S_{i}$ for those special ones, i.e.,
$S_{i}:=T_{i}^{\left(  \openone_{N}\right)  }$.

Let $s_{j}\mapsto T_{j}^{\left(A\right)}$ be an arbitrary
wavelet representation.
By virtue of (\ref{eq2.17})--(\ref{eq2.18}), $L^{2}(\mathbb{T})$, or
equivalently $\ell^{2}(\mathbb{Z})$, splits up as an orthogonal sum%
\begin{equation}
T_{j}^{\left(A\right)}
(\ell^{2}(\mathbb{Z}))\text{,\qquad}j=0,1,\dots,N-1. \label{eq2.22}%
\end{equation}
We saw that the wavelet transform $W$ of (\ref{eqRep.1}) maps
$\ell^{2}(\mathbb{Z})$ onto $\mathcal{V}_{0}$, and from (\ref{eq2.13}) we
conclude that $W$ maps
$T_{0}^{\left(A\right)}(\ell^{2}(\mathbb{Z}))$ onto $U_{N}%
(\mathcal{V}_{0})$ ($=:\mathcal{V}_{1}$). Hence, in the $N$-scale wavelet
case, $W$ transforms the spaces
$T_{j}^{\left(A\right)}(\ell^{2}(\mathbb{Z}))$ ($\subset
\ell^{2}(\mathbb{Z})$) onto orthogonal subspaces $\mathcal{W}_{1}^{(j)}$,
$j=1,\dots,N-1$ in $L^{2}(\mathbb{R})$, and
\begin{equation}
\mathcal{W}_{1}=\mathcal{V}_{0}\ominus\mathcal{V}_{1}=\sideset{}%
{^{\smash{\oplus}}}{\sum}_{j=1}^{N-1\,}\mathcal{W}_{1}^{\left(  j\right)  },
\label{eq2.25}%
\end{equation}
where
\begin{equation}
\mathcal{W}_{1}^{\left(  j\right)  }
=T_{j}^{\left(A\right)}\ell^{2},\qquad j=1,\dots,N-1.
\label{eq2.25bis}%
\end{equation}
Each of the spaces $\mathcal{V}_{1}$ and $\mathcal{W}_{1}^{\left(  j\right)
}$ is split further into orthogonal subspaces corresponding to iteration of
the operators
$T_{0}^{\left(A\right)},T_{1}^{\left(A\right)},\dots,T_{N-1}^{\left(A\right)}$
of (\ref{eq2.17})--(\ref{eq2.18})$.$
It is the system
$\{T_{j}^{\left(A\right)}\}_{j=0}^{N-1}$ which is called a wavelet
representation, and it follows that the wavelet decomposition may be recovered
from the representation. Moreover, the variety of all wavelet representations
is in $1$--$1$ correspondence with the semigroup of polynomial functions $A$
in (\ref{eq2.19}). Operators
$\{T_{j}^{\left(A\right)}\}$ satisfying (\ref{eq2.17})--(\ref{eq2.18})
are said to constitute a representation of the $C^{\ast}%
$-algebra $\mathcal{O}_{N}$, the Cuntz algebra \cite{Cun77}, and it is the
irreducibility of the representations from (\ref{eq2.16}) which will concern
us. If a representation (\ref{eq2.16}) is reducible (Definition \ref{DefRep.3}%
), then there is a subspace
\begin{equation}
0\subsetneqq\mathcal{H}_{0}\subsetneqq L^{2}(\mathbb{T}) \label{eq2.24}%
\end{equation}
which is invariant under all the operators
$T_{j}^{\left(A\right)}$ and $T_{j}^{\left(A\right)\,\ast}$, and so
the data going into the wavelet filter system $\{m_{j}\}$ are then not minimal.

\section{\label{Lem}A lemma about projections}

Our main result is that for a generic set within the class of all wavelet
representations, we do have irreducibility, i.e., there is no reduction as
indicated in (\ref{eq2.24}) in Section \ref{Wav}. In proving this, we will
first reduce the question to a
\emph{finite-dimensional matrix problem}. We will also,
using \cite{BJKW00},
show that every wavelet representation, if it is reducible, decomposes into a
\emph{finite} orthogonal sum of irreducible representations, i.e., if the
$S_{j}$ operators from (\ref{eq2.16}) are given, then there is a finite
orthogonal splitting
\begin{equation}
\ell^{2}(\mathbb{Z})=\sideset{}{^{\smash{\oplus}}}{\sum}_{p}\mathcal{H}_{p}
\label{eq3.1}%
\end{equation}
such that each of the subspaces $\mathcal{H}_{p}$ reduces the representation,
each of the restricted representations of $\mathcal{O}_{N}$ is irreducible,
and moreover that the irreducible subrepresentations which do occur are
\emph{mutually inequivalent} (and therefore disjoint). It is this last
property of inequivalence of the irreducible subrepresentation which amounts
to the fact that the commutant of the original representation from
(\ref{eq2.16}) is \emph{abelian}. Let $\mathcal{H}:=\ell^{2}(\mathbb{Z})$,
let $\mathcal{B}(\mathcal{H})$ denote the algebra of all bounded operators in
$\mathcal{H}$, and
let $s_{i}\mapsto T_{i}=T_{i}^{\left(A\right)}$ be an arbitrary
wavelet representation.
Then the commutant is
\begin{equation}
\mathcal{O}_{N}^{\prime}  =\left\{  X\in\mathcal{B}(\mathcal{H});
T_{i}X=XT_{i}\;\forall\,i\right\}
=\left\{  X\in\mathcal{B}(\mathcal{H}); \sum_{i=0}^{N-1}T_{i}%
XT_{i}^{\ast}=X\right\}  . \label{eq3.2}
\end{equation}
Of course, there are many representations of $\mathcal{O}_{N}$ such that the
corresponding commutant $\mathcal{O}_{N}^{\prime}$ is not abelian, see for
example \cite{DKS99}, but the
\emph{abelian property}
(i.e., that the decomposition
into irreducibles is
\emph{multiplicity-free\/}), is specific to the wavelet
representations;
see Sections \ref{Irr}
and \ref{Exp} below.
The proof of the abelian property is based on a lemma regarding a certain
matrix which turns out to be diagonal with respect to a basis which is a
finite subset of the Fourier basis
\begin{equation}
\{z^{n};n\in\mathbb{Z}\}\text{\qquad(also denoted }e_{n}\left(  z\right)
:=z^{n}\text{),} \label{eq3.3}%
\end{equation}
or equivalently the canonical basis vectors in $\ell^{2}(\mathbb{Z})$. This
lemma in turn depends on a sublemma about a finite set of projections
$P_{1},\dots,P_{g}$ in Hilbert space $\mathcal{H}$. Recall $P\in
\mathcal{B}(\mathcal{H})$ is a projection iff $P=P^{\ast}=P^{2}$. However,
there are more details to the full argument, and they will be taken up in
Sections \ref{Irr} and \ref{Exp} below.

\begin{lemma}
\label{lemma3.1}Let $P_{1},\dots,P_{g}$ be projections. Suppose the operator
\begin{equation}
R=P_{g}P_{g-1}\cdots P_{2}P_{1}P_{2}\cdots P_{g-1}P_{g} \label{eq3.4}%
\end{equation}
is nonzero. Then $R$ is a projection if and only if the $P_{i}$'s are mutually commuting.
\end{lemma}

\begin{proof}
It is clear that the operator $R$ in (\ref{eq3.4}) is a projection if the
family $P_{1},\dots,P_{g}$ consists of mutually commuting projections. We now
prove the converse by induction starting with two given projections
$P_{1},P_{2}$ such that $R:=P_{2}P_{1}P_{2}$ is given to be a projection. Then
the commutator $S:=P_{1}P_{2}-P_{2}P_{1}$ satisfies $S^{\ast}=-S$. Using that
$R^{2}=R$ we conclude that $S^{3}=0$, and therefore $S=0$; in other words, the
two projections $P_{1},P_{2}$ commute.

Suppose the lemma holds for fewer than $g$ projections. If $R$ is given as in
(\ref{eq3.4}), then
\begin{equation}
R=P_{g}TP_{g} \label{eq3.5}%
\end{equation}
where
\begin{equation}
T=P_{g-1}\cdots P_{2}P_{1}P_{2}\cdots P_{g-1}. \label{eq3.6}%
\end{equation}
Writing the operator $T$ in matrix form relative to the two projections
$P_{g}$ and $P_{g}^{\perp}=I-P_{g}$, we get
\begin{equation}
T    =\left(
\begin{array}
[c]{cc}%
R & P_{g}TP_{g}^{\perp}\\
P_{g}^{\perp}TP_{g} & P_{g}^{\perp}TP_{g}^{\perp}%
\end{array}
\right) 
  =(T_{ij})_{i,j=0}^{1}\label{eq3.7}
\end{equation}
with $T_{0,0}=R$, etc. But then (\ref{eq3.8}%
)--(\ref{eq3.9}) yield the conclusion:
\begin{equation}
(T^{2})_{0,0}=(T_{0,0})^{2}+T_{0,1}T_{1,0}=R+T_{0,1}(T_{0,1})^{\ast},
\label{eq3.8}%
\end{equation}
and%
\begin{equation}
(T^{2})_{0,0}\leq T_{0,0}=R \label{eq3.9}%
\end{equation}
imply $T_{0,1}(T_{0,1})^{\ast}=0$, and therefore $T_{0,1}=0$.
As a result, the block matrix in (\ref{eq3.7}) reduces to
\[
T=\left(
\begin{array}
[c]{cc}%
R & 0\\
0 & P_{g}^{\perp}TP_{g}^{\perp}%
\end{array}
\right)  .
\]
A further calculation shows that $T$ must then itself be a projection. From
the definition of $T$ in (\ref{eq3.6}), and the induction hypothesis, we then
conclude that the family $\{P_{i}\}_{i=1}^{g}$ is indeed commutative.
\end{proof}

\begin{remark}
\label{RemLem.2}In the special case when all the projections $\{P_{i}%
\}_{i=1}^{g}$ are one-dimen\-sion\-al, i.e.,
$P_{i}=\left|  v_{i}\right\rangle \left\langle v_{i}\right| $
in the Dirac notation, and $\left\|  v_{i}\right\|  =1$, there is a simpler
proof based on the Schwarz inequality, as follows: Let
$R$ in (\ref{eq3.4}) be given to be a projection,
i.e., $R^{2}=R\neq0$. We also have
$R=\left|  \lambda_{1,2}\lambda_{2,3}\cdots\lambda_{g-1,g}\right|  ^{2}P_{g}$
with $\lambda_{i,j}:=\ip{v_{i}}{v_{j}}$. We then conclude that
$\left|  \lambda_{1,2}\lambda_{2,3}\cdots\lambda_{g-1,g}\right|  =1$,
and therefore by Schwarz, there are constants $\zeta_{i}\in\mathbb{C}$,
$\left|  \zeta_{i}\right|  =1$, such that $v_{2}=\zeta_{1}v_{1}$, $v_{3}%
=\zeta_{2}v_{2},\dots,$ and the commutativity of the family $\{P_{i}%
\}_{i=1}^{g}$ is immediate.
But in this case we find, in addition, that the projections all coincide.
\end{remark}

\section{\label{Min}Minimality and representations}

The representations of the $C^{\ast}$-algebra $\mathcal{O}_{N}$ are used in
other parts of mathematics, in addition to wavelet analysis. While it is known
that in general the irreducible representations of $\mathcal{O}_{N}$ cannot be
given a measurable labeling, see, e.g., \cite{BrJo97a}, \cite{BrJo97b},
\cite{Cun77}, and \cite{BJKW00}, there are various families of $\mathcal{O}%
_{N}$-representations which do admit labeling of their irreducibles, and their
decomposition into sums of irreducibles. We show that
the decomposition into sums of irreducibles occurs only for the special
(permutative) representations
\cite{BrJo99a}
which generalize those derived from the Haar
wavelets. When decompositions do occur, the irreducibles have multiplicity at
most one; see Section \ref{Exp} below.
The basis for our
analysis is the presence of certain finite-dimensional subspaces $\mathcal{K}$
which are invariant under the operators $S_{i}^{\ast}$ when the representation
is defined from the $S_{i}$'s with relations (\ref{eq2.17})--(\ref{eq2.18}).
For related $\mathcal{O}_{N}$-representations which arise in statistical
mechanics, see \cite{FNW92}, \cite{FNW94},
and Section \ref{Irr} below. These finite-dimensional subspaces
have the the significance of labeling the correlations of the sites in the
quantum spin chain model. If it is an infinite spin model on a one-dimensional
lattice, then $\mathcal{K}$ describes the correlations of spin observables
$\sigma_{0},\sigma_{1},\dots$ with those on the other side, $\dots,\sigma
_{-2},\sigma_{-1}$.

We say that a representation of $\mathcal{O}_{N}$ in a Hilbert space
$\mathcal{H}$ is a \emph{wavelet representation} if $\mathcal{H}=L^{2}\left(
\mathbb{T}\right)  $ ($\cong\ell^{2}\left(  \mathbb{Z}\right)  $) and if the
corresponding operators $S_{i}$ are given by (\ref{eq2.16}) for some QMF
functions $\left\{  m_{i}\right\}  _{i=0}^{N-1}$. By (\ref{eq2.20}%
)--(\ref{eq2.21}) that is equivalent to using polynomial functions
$A\colon\mathbb{T}\rightarrow\mathrm{U}_{N}\left(  \mathbb{C}\right)  $ for
labeling the representations. We will let $\mathcal{P}\left(  \mathbb{T}%
,\mathrm{U}_{N}\left(  \mathbb{C}\right)  \right)  $ be the semigroup of such
\emph{polynomial loops}, loops because they may be viewed as
loops in the unitary group $\mathrm{U}_{N}\left(  \mathbb{C}\right)  $, see
\cite{PrSe86}.
We will use the notation $A\left(  z\right)  =\left(  A_{i,j}\left(  z\right)
\right)  _{i,j=0}^{N-1}$ for the loop-group element $A\colon\mathbb{T}%
\rightarrow\mathrm{U}\left(  N\right)  $. Since the Fourier expansion is
finite, there is a $g$ such that $A\left(  z\right)  $ has the form%
\begin{equation}
A\left(  z\right)  =\sum_{k=0}^{g-1}z^{k}A^{\left(  k\right)  }\qquad
(A^{(g-1)}\neq0) \label{eqBrJo00Rep.18}%
\end{equation}
where $A^{\left(  k\right)  }\in\mathcal{B}\left(  \mathbb{C}^{N}\right)  $
for $k=0,\dots,g-1$. The factorization in \cite[Lemma 3.3]{BrJo00} motivates
the name \emph{genus} for $g$.

\begin{lemma}
\label{CorBrJo00Rep.2}If $A\left(  z\right)  $ is a general polynomial of $z$
with values in $\mathcal{B}\left(  \mathbb{C}^{N}\right)  $ of the form
\textup{(\ref{eqBrJo00Rep.18}),} the following four conditions
\textup{(\ref{eqBrJo00Rep.19})--(\ref{eqBrJo00Rep.22})} are equivalent:%
\begin{align}
&
\begin{minipage} [t]{\displayboxwidth}\raggedright$A\left( z\right) ^{\ast}
A\left( z\right) =\openone_{N}$, $z\in\mathbb{T}%
$, i.e., $A$ takes values in $\mathrm{U}\left( N\right) $; \end{minipage}%
\label{eqBrJo00Rep.19}\\
&
\begin{minipage} [t]{\displayboxwidth}\raggedright$\sum_{k}A^{\left
( k\right) \,\ast}A^{\left( k+n\right) }= \begin{cases} \openone_{N}
& \text{ if }n=0, \\ 0 & \text{ if }n\in\mathbb{Z}\setminus\left
\{ 0\right\} , \end{cases}
$ with the convention that $A^{\left( m\right) }=0$ if $m\notin\left
\{ 0,1,\dots,g-1\right\} $; \end{minipage}%
\label{eqBrJo00Rep.20}\\
&
\begin{minipage} [t]{\displayboxwidth}\raggedright
there are projections $P_{1},\dots,P_{s}$ in $\mathcal{B}\left( \mathbb{C}%
^{N}\right) $, positive integers $r_{1},\dots,r_{s}$, and a unitary $W\in
\mathrm{U}\left( N\right) $ such that $A\left( z\right) =\left( \prod
_{j=1}^{s}\left( \openone_{N}-P_{j}+z^{r_{j}}P_{j}\right) \right
) W$; \end{minipage}
\label{eqBrJo00Rep.21}%
\end{align}
and%
\begin{equation}%
\begin{minipage} [t]{\displayboxwidth}\raggedright
there are projections $Q_{0},Q_{1},\dots,Q_{g-2}$ and a unitary $V\in
\mathrm{U}\left( N\right) $ such that \end{minipage}
\label{eqBrJo00Rep.22}%
\end{equation}%
\begin{align*}
A^{\left(  0\right)  }  &  =V\prod_{j=0}^{g-2}\left(  \openone_{N}%
-Q_{j}\right)  ,\\
A^{\left(  1\right)  }  &  =V\sum_{j=0}^{g-2}\left(  \openone_{N}%
-Q_{0}\right)  \cdots\\
&  \qquad\cdots\left(  \openone_{N}-Q_{j-1}\right)  Q_{j}\left(  \openone
_{N}-Q_{j+1}\right)  \cdots\\
&  \qquad\qquad\cdots\left(  \openone_{N}-Q_{g-2}\right)  ,\\
\vdots &  \mathrel{\phantom{=W}}\vdots\\
\qquad A^{\left(  g-1\right)  }  &  =V\prod_{j=0}^{g-2}Q_{j}.
\end{align*}
\end{lemma}

\begin{proof}
We refer the reader to \cite[Proposition 3.2]{BrJo00}.
\end{proof}

\begin{remark}
\label{RemBrJo00Rep.3}The case $g=2=N$ includes the family of wavelets
introduced by Daubechies \cite{Dau92} and studied further in \cite{BEJ00}.
Note that $g=2$ yields the representation%
\begin{equation}
A^{\left(  0\right)  }=V\left(  \openone_{N}-Q\right)  ,\qquad A^{\left(
1\right)  }=VQ, \label{eqBrJo00Rep.23}%
\end{equation}
by \textup{(\ref{eqBrJo00Rep.22}).} But then \textup{(\ref{eqBrJo00Rep.20})}
takes the form%
\begin{equation}
A^{\left(  0\right)  \,\ast}A^{\left(  0\right)  }=\openone_{N}-Q,\qquad
A^{\left(  1\right)  \,\ast}A^{\left(  1\right)  }=Q, \label{eqBrJo00Rep.24}%
\end{equation}
which will be used in the Sections \ref{Fil} and \ref{Exp} below.
\end{remark}

In the general case, we will need the operators (alias matrices) $R\left(
k,l\right)  :=A^{\left(  l\right)  \,\ast}A^{\left(  k\right)  }$, and the
representation (\ref{eqBrJo00Rep.22}) then yields%
\begin{align}
R\left(  0,0\right)   &  =Q_{g-2}^{\perp}\cdots Q_{1}^{\perp}Q_{0}^{\perp
}Q_{1}^{\perp}\cdots Q_{g-2}^{\perp},\label{eqBrJo00Rep.24bis}\\
\vdots &  \qquad\vdots\nonumber\\
R\left(  g-1,g-1\right)   &  =Q_{g-2}\cdots Q_{1}Q_{0}Q_{1}\cdots
Q_{g-2},\nonumber
\end{align}
which were introduced in Lemma \ref{lemma3.1} above.

A loop $A\in\mathcal{P}(\mathbb{T},\mathrm{U}_{N}(\mathbb{C}))$ is viewed
as an entire analytic matrix function, $\mathbb{C}\rightarrow
M_{N}(\mathbb{C})$, and we consider (\ref{eqBrJo00Rep.18}) also as a
representation for this extended (entire) function. The (unique) entire
extension will be denoted $A(z)$ as well.
The estimates in the next corollary translate into a \emph{stability property}
for the corresponding wavelet filters, the significance of which will be
established in Section \ref{Exp} below.

\begin{corollary}
\label{corollary5new}If $g$ is the genus, then we have the following estimate
relative to the order on the positive operators on $\mathbb{C}^{N}$:%
\[
\left(  \min\left(  1,\left|  z\right|  ^{2}\right)  \right)  ^{g-1}%
\cdot\openone_{N}\leq A(z)^{\ast}A(z)\leq\left(  \max\left(  1,\left|
z\right|  ^{2}\right)  \right)  ^{g-1}\cdot\openone_{N},
\]
valid for all $z\in\mathbb{C}$,
where $\openone_{N}$ is the identity matrix.
\end{corollary}

\begin{proof}
The corollary is applied in Section \ref{Exp} below, so we postpone its proof
to Section \ref{Exp}. The argument is in Observation \ref{ObsExpNew.2}, and it
is based on the
ordered
factorizations (\ref{eqBrJo00Rep.21})--(\ref{eqBrJo00Rep.22}) in Lemma
\ref{CorBrJo00Rep.2}.
\end{proof}

It follows from the corollary and (\ref{eq2.20}) that the
system $m_{0},m_{1},\dots,m_{N-1}$ of polynomials that makes up the
multiresolution filter cannot have any other common zeroes than $z=0$, i.e.,
if some $z_{0}\in\mathbb{C}$ satisfies $m_{i}\left(  z_{0}\right)  =0$ for all
$i$, then $z_{0}=0$.

We now turn to some representation theory for the $C^{\ast}$-algebra
$\mathcal{O}_{N}$ which will be needed in the following sections. Some
background references for this are \cite{BrJo99a}, \cite{Eva80}, \cite{Pop92},
and \cite{ReWe98}. Our references for wavelets and
filters are \cite{Hor95}, \cite{Pol90}, and \cite{Vai93}.

Let $P$ ($\in\mathcal{B}\left(  \mathcal{H}\right)  $) be a projection. We say
that it is \emph{co-invariant} for some (fixed) representation $\left\{
T_{i}\right\}  _{i=0}^{N-1}$ of $\mathcal{O}_{N}$ if%
\begin{equation}
T_{i}^{\ast}P=PT_{i}^{\ast}P\text{\qquad for all }i. \label{eqMin.7}%
\end{equation}
Let $\mathcal{H}_{-}$ be the closed span of $\left\{  z^{-n};n=0,1,\dots
\right\}  $, and let $P_{-}$ be the projection onto $\mathcal{H}_{-}$. Then
(\ref{eqMin.7}) is satisfied for $P_{-}$ and all wavelet representations
$T^{\left(  A\right)  }$, as follows from
(\ref{eq2.16}), (\ref{eq2.20}), and
the formula%
\begin{equation}
T_{i}^{\left(  A\right)  \,\ast}=\sum_{j=0}^{N-1}\overline{A_{i,j}\left(
z\right)  }\,S_{j}^{\ast}, \label{eqMin.8}%
\end{equation}
where $S_{j}^{\ast}$ are the adjoints of the respective operators $S_{j}$ in
(\ref{eqMin.6}). Specifically,%
\begin{equation}
\left(  S_{j}^{\ast}f\right)  \left(  z\right)  =\frac{1}{N}\sum_{w^{N}%
=z}w^{-j}f\left(  w\right)  ,\qquad f\in L^{2}\left(  \mathbb{T}\right)  .
\label{eqMin.9}%
\end{equation}

\begin{lemma}
\label{LemMin.3}Let $E$ and $P$ be co-invariant projections for a fixed
representation $T^{\left(  A\right)  }$. Suppose $E\leq P\leq P_{-}$, and
further that for some $r\in\mathbb{N}$,
\begin{equation}
P\mathcal{H}=\operatorname*{span}\left\{  z^{-k};0\leq k\leq r\right\}  .
\label{eqMin.10}%
\end{equation}
Finally assume that%
\begin{equation}
T_{i}^{\left(  A\right)  \,\ast}E=ET_{i}^{\left(  A\right)  \,\ast
}P\text{\qquad for all }i=0,1,\dots,N-1. \label{eqMin.11}%
\end{equation}
Then we have the following identities:%
\begin{equation}
S_{j}^{\ast}ES_{k}^{{}}P=\sum_{i=0}^{N-1}A_{i,j}E\bar{A}_{i,k}P_{k}
\label{eqMin.12}%
\end{equation}
for all $j,k=0,1,\dots,N-1$, where the functions $A_{i,j}$ are the matrix
entries of the given loop, a function is identified with the corresponding
multiplication operator in $\mathcal{H}=L^{2}\left(  \mathbb{T}\right)  $, and
$P_{k}:=S_{k}^{\ast}PS_{k}^{{}}$ are projections.
\end{lemma}

\begin{proof}
It is given that both $E$ and $P$ satisfy (\ref{eqMin.7}) relative to
$T^{\left(  A\right)  }$, and further that
$T_{i}^{\left(  A\right)  \,\ast}E=ET_{i}^{\left(  A\right)  \,\ast}P$.
Equivalently, by (\ref{eqMin.8}),
$\sum_{l}\bar{A}_{i,l}^{{}}S_{l}^{\ast}E=E\sum_{k}\bar{A}_{i,k}^{{}}S_{k}
^{\ast}P$.
Using $\sum_{i}A_{i,j}\bar{A}_{i,l}=\delta_{j,l}$, we get%
\begin{equation}
S_{j}^{\ast}E=\sum_{i,k}A_{i,j}^{{}}E\bar{A}_{i,k}^{{}}S_{k}^{\ast}P.
\label{eqMin.15}%
\end{equation}
Now multiplying through from the right with $S_{k}P$ on both sides in
(\ref{eqMin.15}), the conclusion of the lemma follows. To see this, notice
first from (\ref{eqMin.9})--(\ref{eqMin.10}) that%
\begin{equation}
S_{k}^{\ast}PS_{l}^{{}}P=0\text{\qquad if }k\neq l. \label{eqMin.16}%
\end{equation}
The proof of (\ref{eqMin.16}) is based on the observation that the
representation $S$ ($=T^{\left(  \openone_{N}\right)  }$)
in (\ref{eqMin.6}) is permutative, see
\cite{BrJo99a}. Specifically, $S_{k}^{\ast}\left(  z^{j-nN}\right)
=\delta_{k,j}z^{-n}$ if $0\leq j<N$, and $n\in\mathbb{Z}$, and $S_{l}\left(
z^{-n}\right)  =z^{l-nN}$.
The desired formula (\ref{eqMin.12}) now follows from this and $S_{k}^{\ast
}PS_{k}^{{}}P=S_{k}^{\ast}PS_{k}^{{}}=P_{k}$, since $P$ is relatively
co-invariant for the representation $S=T^{\left(  \openone_{N}\right)  }$ by assumption.
\end{proof}

\begin{rt}
\label{RemMin.4}The proof shows more generally that the implication
\textup{(\ref{eqMin.11})} $\Rightarrow$ \textup{(\ref{eqMin.12})} holds for
any operator $E\in\mathcal{B}\left(  P\mathcal{H}\right)  $ when $P$ is
specified as in the statement of the lemma. In $\mathcal{B}\left(
P\mathcal{H}\right)  $, we may then introduce the basis $e_{-k,-l}:=\left|
z^{-k}\right\rangle \left\langle z^{-l}\right|  $, and coordinates%
\begin{equation}
E=\sum_{k,l}X_{k,l}e_{-k,-l}. \label{eqRemMin.4pound}%
\end{equation}
\end{rt}

If the loop $A\left(  z\right)  $ is 
given by (\ref{eqBrJo00Rep.18}),
then the operators%
\begin{equation}
R\left(  k,l\right)  :=A^{\left(  l\right)  \,\ast}A^{\left(  k\right)  }
\label{eqMinNew.2}%
\end{equation}
of Lemma \ref{CorBrJo00Rep.2} go into the calculation of the right-hand side
in (\ref{eqMin.12}) as follows: The $\left(  r,s\right)  $-matrix entry of the
matrix
($S_{j}^{\ast}ES_{k}^{{}}=$)  $\sum_{i=0}^{N-1}A_{i,j}
E\bar{A}_{i,k}P$
is given by the following matrix product:%
\begin{equation}
\sum_{\substack{p,q\\p\geq r,\,q\geq s}}X_{p,q}R\left(  p-r,q-s\right)
_{k,j}, \label{eqMinNew.3}%
\end{equation}
again with the convention that the summation indices restrict to the range
where the terms in the sum are defined and nonvanishing.

\begin{definition}
\label{DefMin.4}Let $\left\{  T_{i}\right\}  _{i=0}^{N-1}$ be a representation
of $\mathcal{O}_{N}$ on a Hilbert space $\mathcal{H}$, and let $\mathcal{K}$
be a finite-dimensional subspace which satisfies
\begin{equation}
T_{i}^{\ast}\mathcal{K}\subset\mathcal{K}\text{\qquad for all }i.
\label{eqMin.16bis}%
\end{equation}
Hence the projection $P$ onto $\mathcal{K}$ satisfies%
\begin{equation}
PT_{i}=PT_{i}P\text{\qquad for all }i. \label{eqMin.17}%
\end{equation}
We say that $\mathcal{K}$ is \emph{cyclic} if it is cyclic for the
$\mathcal{O}_{N}$-representation, i.e., if%
\begin{equation}
\bigvee_{i_{1},i_{2},\dots,i_{n}}T_{i_{1}}T_{i_{2}}\cdots T_{i_{n}}%
\mathcal{K}=\mathcal{H}. \label{eqMin.18}%
\end{equation}
\end{definition}

For the wavelet representations, $\mathcal{H}=L^{2}\left(  \mathbb{T}\right)
$, the Fourier basis $\left\{  z^{n};n\in\mathbb{Z}\right\}  $ has the
following property: There is an $r_{0}\in\mathbb{N}$ such that, for all
$n\in\mathbb{Z}$, there is a $p\in\mathbb{N}$ satisfying
\begin{equation}
T_{i_{q}}^{\ast}\cdots T_{i_{2}}^{\ast}T_{i_{1}}^{\ast}\left(  z^{n}\right)
\in\operatorname*{span}\left\{  z^{-k};0\leq k\leq r_{0}\right\}  \label{eqMin.19}%
\end{equation}
for all multi-indices $i_{1},\dots,i_{q}$ and $q\geq p$. We showed
\cite{BrJo00} that $r_{0}$ may be taken
\begin{equation}
r_{0}=\left\lfloor \frac{gN-1}{N-1}\right\rfloor \label{eqMin.20}%
\end{equation}
where $g$ is the genus, $N$ is the scale, and $\left\lfloor x\right\rfloor $
is the largest integer $\leq x$. We also showed that, whenever (\ref{eqMin.19}%
) holds, then%
\begin{equation}
\mathcal{K}:=\operatorname*{span}\left\{  z^{-k};0\leq k\leq r_{0}\right\}
\label{eqMin.21}%
\end{equation}
is cyclic.
It is known in general that, if some $\mathcal{K}$ is \emph{minimal} with
respect to the two properties, (\ref{eqMin.16bis}) and $\mathcal{O}_{N}%
$-cyclicity, then
\begin{equation}
\mathcal{B}\left(  \mathcal{K}\right) ^{\sigma^{\left(  A\right)  }}:=
\left\{  X\in\mathcal{B}\left(  \mathcal{K}\right)  ;\sum\nolimits_{i}%
PT_{i}^{\left(  A\right)  }XT_{i}^{\left(  A\right)  \,\ast}P=X\right\}
\label{eqMin.22}%
\end{equation}
is an \emph{algebra.} The set (\ref{eqMin.22}) is the fixed-point set for
the completely positive map
\begin{equation}
\sigma_{\mathcal{K}}^{\left(  A\right)  }\left(  \,\cdot\,\right)  =\sum
_{i}V_{i}^{{}}\left(  \,\cdot\,\right)  V_{i}^{\ast}\text{,\qquad where }%
V_{i}:=PT_{i}^{\left(  A\right)  }. \label{eqMin.22bis}%
\end{equation}
We further showed in \cite{BJKW00} that $T^{\left(  A\right)  }$ is
irreducible if and only if $\mathcal{B}\left(  \mathcal{K}\right)
^{\sigma^{\left(  A\right)  }}=\mathbb{C}\,\openone_{\mathcal{K}}$.
In general, this set is \emph{not} an
algebra, but the above minimality on $\mathcal{K}$ forces it to be an algebra,
see \cite{DKS99}.

We shall need,
in the later proofs,
the following two results from \cite{BJKW00}.
We include the statements here since they seem not to be well known
in the wavelet community.
Let $\pi$ be
a representation of $\mathcal{O}_{N}$ on a Hilbert space $\mathcal{H}$.

\begin{theorem}
\label{ThmBrJo00Rep.1} \cite[Section~6]{BJKW00} There is a positive
norm-preserving linear isomorphism between the commutant algebra%
\begin{equation}
\pi\left(  \mathcal{O}_{N}\right)  ^{\prime}=\left\{  A\in\mathcal{B}\left(
\mathcal{H}\right)  ;A\pi\left(  x\right)  =\pi\left(  x\right)  A\text{ for
all }x\in\mathcal{O}_{N}\right\}  \label{eqBrJo00Rep.10}%
\end{equation}
and the fixed-point set%
\begin{equation}
\mathcal{B}\left(  \mathcal{K}\right)  ^{\sigma}=\left\{  A\in\mathcal{B}%
\left(  \mathcal{K}\right)  ;\sigma\left(  A\right)  =A\right\}
\label{eqBrJo00Rep.11}%
\end{equation}
given by%
\begin{equation}
\pi\left(  \mathcal{O}_{N}\right)  ^{\prime}\ni A\longrightarrow PAP,
\label{eqBrJo00Rep.12}%
\end{equation}
where $P$ is the projection of $\mathcal{H}$ onto $\mathcal{K}$. In
particular, $\pi$ is irreducible if and only if $\sigma$ is ergodic
\textup{(}where $\sigma$ is the mapping $\mathcal{B}\left(  \mathcal{K}%
\right)  \rightarrow\mathcal{B}\left(  \mathcal{K}\right)  $ defined in
\textup{(\ref{eqMin.22bis})).}

More generally, if $\mathcal{K}_{1}$, $\mathcal{K}_{2}$ \textup{(}with
corresponding projections $P^{\left(  1\right)  }$ and $P^{\left(  2\right)
}$\textup{)} are $T^{\ast}$-invariant cyclic subspaces for two representations
$\pi_{1}$, $\pi_{2}$ of $\mathcal{O}_{N}$ on $\mathcal{H}_{1}$, $\mathcal{H}%
_{2}$, and%
\begin{equation}
V_{i}^{\left(  j\right)  }=P_{{}}^{\left(  j\right)  }\pi_{j}^{{}}\left(
s_{i}\right)  |_{\mathcal{K}_{j}} \label{eqBrJo00Rep.13}%
\end{equation}
for $j=1,2$, $i=0,\dots,N-1$, define $\rho$ on $\mathcal{B}\left(
\mathcal{K}_{1},\mathcal{K}_{2}\right)  $ by%
\begin{equation}
\rho\left(  A\right)  =\sum_{i}V_{i}^{\left(  2\right)  }AV_{i}^{\left(
1\right)  \,\ast}. \label{eqBrJo00Rep.14}%
\end{equation}
Then there is an isometric linear isomorphism between the set of intertwiners%
\begin{equation}
\left\{  A\in\mathcal{B}\left(  \mathcal{H}_{1},\mathcal{H}_{2}\right)
;A\pi_{1}\left(  x\right)  =\pi_{2}\left(  x\right)  A\text{ for all }%
x\in\mathcal{O}_{N}\right\}  \label{eqBrJo00Rep.15}%
\end{equation}
and the fixed-point set%
\begin{equation}
\left\{  B\in\mathcal{B}\left(  \mathcal{K}_{1},\mathcal{K}_{2}\right)
;\rho\left(  B\right)  =B\right\}  \label{eqBrJo00Rep.16}%
\end{equation}
given by%
\begin{equation}
A\longrightarrow B=P^{\left(  2\right)  }AP^{\left(  1\right)  }.
\label{eqBrJo00Rep.17}%
\end{equation}
\end{theorem}

\begin{theorem}
\label{ThmBJKW003.5} \cite[Theorem~3.5]{BJKW00} Let $\varphi
=\sum_{i}V_{i}^{{}}\,\cdot\,V_{i}^{\ast}$ be a normal unital completely
positive map of $\mathcal{B}\left(  \mathcal{K}\right)  $. Then%
\[
\left\{  V_{i}^{{}},V_{i}^{\ast}\right\}  ^{\prime}\subset\mathcal{B}\left(
\mathcal{H}\right)  ^{\varphi}.
\]
Furthermore, the space $\mathcal{B}\left(  \mathcal{H}\right)  ^{\varphi}$
contains a largest $\ast$-subalgebra, and this algebra is $\left\{  V_{i}^{{}%
},V_{i}^{\ast}\right\}  ^{\prime}$.
\end{theorem}

We are now ready to state and prove the main result of the present section.
Its significance becomes more clear when it is seen in the light of the two
previous Theorems \ref{ThmBrJo00Rep.1} and \ref{ThmBJKW003.5}. In particular,
we will show in Section \ref{Exp} below that Theorem \ref{ThmBJKW003.5} is
applicable in verifying irreducibility, as we will show that $\mathcal{B}%
\left(  \mathcal{K}\right)  ^{\sigma}$ is generically an algebra for the
wavelet representations.

\begin{theorem}
\label{ThmMin.7}Let $T^{\left(  A\right)  }$ be a wavelet representation of
$\mathcal{O}_{N}$ on $\mathcal{H}=L^{2}\left(  \mathbb{T}\right)  $, and
assume the genus of $A$ is $g$. Let $r_{0}$
be as in \textup{(\ref{eqMin.20}),}
and let $P$ be the projection onto $\mathcal{K}%
:=\operatorname*{span}\left\{  z^{-k};0\leq k\leq r_{0}\right\}  $. Suppose
there is a second projection $E\in\mathcal{B}\left(  \mathcal{K}\right)  $
such that $0\neq E\neq P$, and $E$ commutes with $T_{i}^{\left(  A\right)
\,\ast}P$ for all $i=0,1,\dots,N-1$. Then it follows that $E$ is diagonal with
respect to the basis $\left\{  z^{-k};k=0,1,\dots,r_{0}\right\}  $ in
$\mathcal{K}$. Moreover, $A\left(  z\right)  $ has a matrix corner of the
form
\begin{equation}
V%
\begin{pmatrix}
z^{n_{0}} & 0 & \cdots & 0\\
0 & z^{n_{1}} & \cdots & 0\\
\vdots & \vdots & \ddots & \vdots\\
0 & 0 & \cdots & z^{n_{M-1}}%
\end{pmatrix}
, \label{eqThmMin.7AzVm}%
\end{equation}
where $V\in\mathrm{U}_{M}\left(  \mathbb{C}\right)  $, and where the exponents
$n_{i}$ of the diagonal corner in $A\left(  z\right)  $ satisfy $0\leq
n_{i}\leq g-1$ for all $i=0,1,\dots,M-1$.
\end{theorem}

\begin{remark}
\label{RemThmMin.7}The loops $A\colon\mathbb{T}\rightarrow\mathrm{U}%
_{N}\left(  \mathbb{C}\right)  $ which do admit nontrivial projections $E$ as
in the statement of the theorem are described in detail in Definition
\textup{\ref{definition5.5}} below, to which we refer. Hence, the existence of
such a projection $E$ means that it is possible to ``split off'' a matrix
block in $A\left(  z\right)  $ which is in diagonal form.
\end{remark}

\begin{proof}
[Proof of Theorem \textup{\ref{ThmMin.7}}]Let $V_{i}:=PT_{i}^{\left(
A\right)  }$. Suppose $E\in\mathcal{B}\left(  P\mathcal{H}\right)  $ satisfies
$EV_{i}^{\ast}=V_{i}^{\ast}E$, or equivalently $ET_{i}^{\left(  A\right)
\,\ast}P=T_{i}^{\left(  A\right)  \,\ast}E$. Then by Lemma \ref{LemMin.3} and
Remark \ref{RemMin.4}, we have%
\begin{equation}
\left(  S_{j}^{\ast}ES_{k}^{{}}\right)  _{r,s}=\sum_{\substack{p,q\\p\geq
r,\,q\geq s}}X_{p,q}R\left(  p-r,q-s\right)  _{k,j}=\sum_{\substack{p,q\\p\geq
0,\,q\geq0}}X_{p+r,q+s}R\left(  p,q\right)  _{k,j}. \label{eqThmMin.7proof.1}%
\end{equation}
The $j,k$-indices are in $\left\{  0,1,\dots,N-1\right\}  $. For the matrices
$R\left(  p,q\right)  $, we have the identities%
\begin{equation}
\sum_{p}R\left(  p,p\right)  =\openone_{N}
\text{\quad and\quad }
\sum_{p}R\left(  p,p+l\right)  =0\text{\qquad if }l\neq0.
\label{eqThmMin.7proof.2}%
\end{equation}
See Lemma \ref{CorBrJo00Rep.2} above. The terms on the
left-hand side in (\ref{eqThmMin.7proof.1}) are%
\begin{equation}
\left(  S_{j}^{\ast}ES_{k}^{{}}\right)  _{r,s}=X_{rN-j,sN-k}\,,
\label{eqThmMin.7proof.4}%
\end{equation}
again with the convention that the terms are defined to be zero when the
subscript indices are not in the prescribed range.

If $E\neq0$, consider the lexicographic order on the subscript indices of the
corresponding matrix entries $X_{p,q}$ (in (\ref{eqRemMin.4pound})). The range
on both indices $p,q$ is $\left\{  0,1,2,\dots,r_{0}\right\}  $ where $r_{0}$
is determined as in Lemma \ref{LemMin.3}, see also (\ref{eqMin.20}).
Then pick the last
(relative to lexicographic order) nonzero $X_{r,s}$, i.e., $r,s$ are
determined such that
\begin{equation}
X_{p+r,q+s}=0\text{\qquad if }p>0\text{ or }q>0. \label{eqThmMin.7proof.5}%
\end{equation}
It follows that there are only the following possibilities for this $\left(
r,s\right)  $:%
\[%
\begin{array}
[c]{cll}%
\left(  0,0\right)  &  & E=X_{0,0}\left|  z^{0}\right\rangle \left\langle
z^{0}\right|  ,\\
&  & \\
\left(  1,1\right)  &  & E=X_{0,0}\left|  z^{0}\right\rangle \left\langle
z^{0}\right|  +X_{1,1}\left|  z^{-1}\right\rangle \left\langle z^{-1}\right|
,\\
&  & \\
\left(  2,2\right)  &  & E=X_{0,0}\left|  z^{0}\right\rangle \left\langle
z^{0}\right|  +X_{1,1}\left|  z^{-1}\right\rangle \left\langle z^{-1}\right|
+X_{2,2}\left|  z^{-2}\right\rangle \left\langle z^{-2}\right|  ,\\
\vdots &  & \;\vdots\;.
\end{array}
\]
If $\left(  r,s\right)  =\left(
0,0\right)  $, then, using (\ref{eqThmMin.7proof.1}) and
(\ref{eqThmMin.7proof.4}), we arrive at the matrix identity%
\begin{equation}
X_{0,0}R\left(  0,0\right)  =\left(
\begin{array}
[c]{c|ccc}%
X_{0,0} & 0 & \cdots & 0\\\hline
0 & 0 & \cdots & 0\\
\vdots & \vdots & \ddots & \vdots\\
0 & 0 & \cdots & 0
\end{array}
\right)  \in M_{N}\left(  \mathbb{C}\right)  \text{,\qquad where }X_{0,0}%
\neq0, \label{eqThmMin.7proof.6}%
\end{equation}
and therefore $R\left(  0,0\right)  =\left|  \varepsilon_{0}\right\rangle
\left\langle \varepsilon_{0}\right|  $ where $\varepsilon_{0}$ is the first
canonical basis vector in $\mathbb{C}^{N}$. By Lemma \ref{lemma3.1},
(\ref{eqBrJo00Rep.24bis}), and Remark \ref{RemLem.2}, we conclude that
$Q_{i}^{\perp}\geq\left|  \varepsilon_{0}\right\rangle \left\langle
\varepsilon_{0}\right|  $ for all $i$, and therefore%
\begin{equation}
A\left(  z\right)  =V\left(
\begin{array}
[c]{c|ccc}%
1 & 0 & \cdots & 0\\\hline
0 &  &  & \\
\vdots &  & B(z) & \\
0 &  &  &
\end{array}
\right)  \label{eqThmMin.7proof.7}%
\end{equation}
for some $V\in\mathrm{U}_{N}\left(  \mathbb{C}\right)  $ and $B\in
\mathcal{P}\left(  \mathbb{T},\mathrm{U}_{N-1}\left(  \mathbb{C}\right)
\right)  $; see Lemma \ref{CorBrJo00Rep.2}.

If $\left(  r,s\right)  =\left(  1,1\right)  $, then, using again
(\ref{eqThmMin.7proof.1}) and (\ref{eqThmMin.7proof.4}), we arrive at the
matrix identity%
\begin{equation}
X_{1,1}R\left(  0,0\right)  =\left(
\begin{array}
[c]{ccc|c}%
0 & \cdots & 0 & 0\\
\vdots & \ddots & \vdots & \vdots\\
0 & \cdots & 0 & 0\\\hline
0 & \cdots & 0 & X_{1,1}%
\end{array}
\right)  \in M_{N}\left(  \mathbb{C}\right)  \text{,\qquad where }X_{1,1}%
\neq0, \label{eqThmMin.7proof.8}%
\end{equation}
and therefore $R\left(  0,0\right)  =\left|  \varepsilon_{N-1}\right\rangle
\left\langle \varepsilon_{N-1}\right|  $. Using again Lemma \ref{lemma3.1},
(\ref{eqBrJo00Rep.24bis}), and Remark \ref{RemLem.2}, we conclude that
$Q_{i}^{\perp}\geq\left|  \varepsilon_{N-1}\right\rangle \left\langle
\varepsilon_{N-1}\right|  $ for all $i$, and therefore%
\begin{equation}
A\left(  z\right)  =V\left(
\begin{array}
[c]{ccc|c}
&  &  & 0\\
& C\left(  z\right)  &  & \vdots\\
&  &  & 0\\\hline
0 & \cdots & 0 & 1
\end{array}
\right)  \label{eqThmMin.7proof.9}%
\end{equation}
for some $V\in\mathrm{U}_{N}\left(  \mathbb{C}\right)  $ and $C\in
\mathcal{P}\left(  \mathbb{T},\mathrm{U}_{N-1}\left(  \mathbb{C}\right)
\right)  $.

The reason for why the matrix of $E$ has diagonal form relative to the natural
Fourier basis is as follows: Let $N=2$, for simplicity. (The argument is the
same, \emph{mutatis mutandis,} in the general case.) Then pick the last term
$\left(  r,s\right)  $, $r\neq s$, with $X_{r,s}\neq0$, where again ``last''
refers to the lexicographic order of the matrix-entry indices, see
(\ref{eqThmMin.7proof.5}). We then get, using (\ref{eqThmMin.7proof.1}) and
(\ref{eqThmMin.7proof.4}), the following matrix-identity (where we specialize
to $\left(  r,s\right)  =\left(  0,1\right)  $):%
\[
X_{0,1}R\left(  0,0\right)  =%
\begin{pmatrix}
0 & 0\\
X_{0,1} & 0
\end{pmatrix}
\in M_{2}\left(  \mathbb{C}\right)
\]
and $X_{0,1}\neq0$, as mentioned. This forces
$R\left(  0,0\right)  =\left( 
\begin{smallmatrix}
0 & 0\\
1 & 0
\end{smallmatrix}
\right) $, which is impossible by Lemma \ref{CorBrJo00Rep.2}, since $R\left(
0,0\right)  $ is positive and $\left( 
\begin{smallmatrix}
0 & 0\\
1 & 0
\end{smallmatrix}
\right) $ is not.

Let $N\geq2$, and suppose $\left(  r,s\right)  =\left(  2,2\right)  $, i.e.,
assume that $X_{2,2}\neq0$, and $X_{2+p,2+q}=0$ if $p>0$ or $q>0$, referring
to the lexicographic order. Then by the same argument which we used in the
earlier cases,%
\begin{equation}
X_{2,2}\left(  R\left(  0,0\right)  \right)  _{k,j}=\left(  X_{2N-j,2N-k}%
\right)  _{j,k=0}^{N-1}\,. \label{eqThmMin.7proof.10}%
\end{equation}
But all the double indices $\left(  2N-j,2N-k\right)  $ of the matrix on the
right are strictly bigger than $\left(  2,2\right)  $ in the lexicographic
order, and we conclude that%
\begin{equation}
R\left(  0,0\right)  =0\text{\qquad in }M_{N}\left(  \mathbb{C}\right)  .
\label{eqThmMin.7proof.11}%
\end{equation}
The formula for $R\left(  0,0\right)  $ then yields $Q_{0}^{\perp}Q_{1}%
^{\perp}\cdots Q_{g-1}^{\perp}=0$. Moreover, the additional restrictions are:%
\[
X_{0,0}R\left(  0,0\right)  +X_{1,1}R\left(  1,1\right)  +X_{2,2}R\left(
2,2\right)  =X_{0,0}\left|
\varepsilon_{0}\right\rangle \left\langle \varepsilon_{0}\right|
\in M_{N}\left(  \mathbb{C}\right)
\]
as in (\ref{eqThmMin.7proof.6}), and%
\[
X_{1,1}R\left(  0,0\right)  +X_{2,2}R\left(  1,1\right)  =\left(
\begin{array}
[c]{ccc|cc}%
0 & \cdots & 0 & 0 & 0\\
\vdots & \ddots & \vdots & \vdots & \vdots\\
0 & \cdots & 0 & 0 & 0\\\hline
0 & \cdots & 0 & X_{2,2} & 0\\
0 & \cdots & 0 & 0 & X_{1,1}%
\end{array}
\right)  .
\]
Now substituting $R\left(  0,0\right)  =0$, we arrive at
\begin{equation}
X_{1,1}R\left(  1,1\right)  +X_{2,2}R\left(  2,2\right)  =X_{0,0}\left|
\varepsilon_{0}\right\rangle \left\langle \varepsilon_{0}\right|
\label{eqThmMin.7proof.12}%
\end{equation}
and%
\begin{equation}
X_{2,2}R\left(  1,1\right)  =X_{2,2}\left|  \varepsilon_{N-2}\right\rangle
\left\langle \varepsilon_{N-2}\right|  +X_{1,1}\left|  \varepsilon
_{N-1}\right\rangle \left\langle \varepsilon_{N-1}\right|  .
\label{eqThmMin.7proof.13}%
\end{equation}
Since $E=X_{0,0}\left|  1\right\rangle \left\langle 1\right|  +X_{1,1}\left|
z^{-1}\right\rangle \left\langle z^{-1}\right|  +X_{2,2}\left|  z^{-2}%
\right\rangle \left\langle z^{-2}\right|  $ is a projection, and $X_{2,2}%
\neq0$, we must have $X_{2,2}=1$ and $X_{0,0}$ and $X_{1,1}\in\left\{
0,1\right\}  $. The conclusion of the theorem can then be checked case by
case, using (\ref{eqThmMin.7proof.12}) and (\ref{eqThmMin.7proof.13}).

In general, let $X_{s,s}\neq0$ be the last (in lexicographic order) nonzero
term, and assume $s\geq2$. Then by (\ref{eqThmMin.7proof.2}%
)--(\ref{eqThmMin.7proof.4}), we get $X_{s,s}R\left(  0,0\right)  =0$, and
therefore $R\left(  0,0\right)  =0$ as before. Using this, the equation for
the $\left(  s-1,s-1\right)  $ term is then
\[
X_{s,s}R\left(  1,1\right)  =\left(  X_{\left(  s-1\right)  N-j,\left(
s-1\right)  N-k}\right)  _{j,k=0}^{N-1}\,.
\]
If $\left(  s-1\right)  \left(  N-1\right)  \leq N$, then all the entries in
the matrix on the right must vanish, and we get $R\left(  1,1\right)  =0$. If
not, we proceed as in (\ref{eqThmMin.7proof.13}). If $R\left(  1,1\right)
=0$, we go to the $\left(  s-2,s-2\right)  $ term, viz.,%
\begin{equation}
X_{s,s}R\left(  2,2\right)  =\left(  X_{\left(  s-2\right)  N-j,\left(
s-2\right)  N-k}\right)  _{j,k=0}^{N-1}\,. \label{eqThmMin.7proof.14}%
\end{equation}
Eventually the matrix on the right will have nonzero terms, starting with
$X_{s,s}$, and terms before that in the lexicographic order. Suppose, for
example, that the matrix on the right in (\ref{eqThmMin.7proof.14}) has
nonzero entries. Then the equation for the $\left(  s-3,s-3\right)  $ term is%
\[
X_{s-1,s-1}R\left(  2,2\right)  +X_{s,s}R\left(  3,3\right)  =\left(
X_{\left(  s-3\right)  N-j,\left(  s-3\right)  N-k}\right)  _{j,k=0}^{N-1}\,,
\]
and the argument is done by a case-by-case check, using that the coordinates
$X_{0,0},X_{1,1},\dots$ are in $\left\{  0,1\right\}  $ while $X_{s,s}=1$.

There is a similar argument, based on the reversed lexicographic order,
starting with $\left(  N-1,N-1\right)  $, which will account for a possible
lower right matrix corner of diagonal form.
This completes the proof of the theorem.
\end{proof}

\begin{remark}
\label{RemMin.8}(\emph{Permutative Representations\/}) The form
\begin{equation}
A\left(  z\right)  =V%
\begin{pmatrix}
z^{n_{0}} & 0 & \cdots & 0\\
0 & z^{n_{1}} & \cdots & 0\\
\vdots & \vdots & \ddots & \vdots\\
0 & 0 & \cdots & z^{n_{N-1}}%
\end{pmatrix}
, \label{eqRemMin.8.15}%
\end{equation}
$V\in\mathrm{U}_{N}\left(  \mathbb{C}\right)  $, in the conclusion of Theorem
\textup{\ref{ThmMin.7}} corresponds to the representations of $\mathcal{O}%
_{N}$ which permute the basis elements $\left\{  z^{n};n\in\mathbb{Z}\right\}
$ for $\mathcal{H}=L^{2}\left(  \mathbb{T}\right)  $; they are studied more
generally in \cite{BrJo99a} under the name \emph{permutative representations.}
We also met them, in a special case, in Remark \ref{RemRep.5} above, in
connection with the ``stretched out'' Haar wavelets. So the conclusion of
Theorem \textup{\ref{ThmMin.7}} is that the wavelet representations which are
not of this form are irreducible.

Now for the details: Let $T^{\left(  A\right)  }$ be the representation of
$\mathcal{O}_{N}$ corresponding to $A\left(  z\right)  $ in
\textup{(\ref{eqRemMin.8.15}).} The element $V\in\mathrm{U}_{N}\left(
\mathbb{C}\right)  $ defines an automorphism of $\mathcal{O}_{N}$, denoted
$\alpha_{V}$ or $\operatorname*{Ad}\left(  V\right)  $. Let $D\left(
z\right)  =V^{-1}A\left(  z\right)  $ be the diagonal factor in
\textup{(\ref{eqRemMin.8.15}).} If $\pi^{\left(  A\right)  }\left(
s_{i}\right)  =T_{i}^{\left(  A\right)  }$ and $\pi^{\left(  D\right)
}\left(  s_{i}\right)  =T_{i}^{\left(  D\right)  }$ are the corresponding
representations, then it follows that%
\begin{equation}
\pi^{\left(  A\right)  }=\pi^{\left(  D\right)  }\circ\alpha_{V}\,,
\label{eqRemMin.8.16}%
\end{equation}
which means that $\pi^{\left(  A\right)  }$ and $\pi^{\left(  D\right)  }$
have the same decomposition into sums of irreducibles, corresponding to
irreducible subspaces of $L^{2}\left(  \mathbb{T}\right)  $. The formulas for
the operators $T_{i}^{\left(  D\right)  }$ are as follows:%
\begin{equation}
T_{i}^{\left(  D\right)  }\left(  z^{k}\right)     =z^{N\left(
n_{i}+k\right)  +i},\qquad k\in\mathbb{Z},\ i=0,\dots ,N-1,
\label{eqRemMin.8.17}
\end{equation}
which justifies the ``permutative'' label; in other
words, both the operators $T_{i}^{\left(  D\right)  }$ and their adjoints
permute the basis elements of the Fourier basis $\left\{  z^{k};k\in
\mathbb{Z}\right\}  $ for $L^{2}\left(  \mathbb{T}\right)  $. The
decomposition structure of these representations is worked out in
\cite{BrJo99a}; see also \cite{DKS99}.
\end{remark}

\begin{remark}
\label{RemMinNew.12}Note that if $N>2$, then some representation $T^{\left(
A\right)  }$ may be reducible even if $A$ is not itself of the form
\textup{(\ref{eqRemMin.8.15});} it may only have a matrix corner of this form.
Take, for example,%
\begin{equation}
A\left(  z\right)  =\left(
\begin{array}
[c]{c|cc}%
1 & 0 & 0\\\hline
0 & 1/\sqrt{2} & z/\sqrt{2}^{\mathstrut}\\
0 & z/\sqrt{2} & -z^{2}/\sqrt{2}%
\end{array}
\right)  =\openone_{0}\oplus\frac{1}{\sqrt{2}}\left(
\begin{array}
[c]{cc}%
1 & z\\
z & -z^{2}%
\end{array}
\right)  , \label{eqRemMinNew.12.1}%
\end{equation}
i.e., $N=3$, $g=3$. Then $T^{\left(  A\right)  }$ is a \emph{reducible}
representation of $\mathcal{O}_{3}$ acting on $L^{2}\left(  \mathbb{T}\right)
$, and, in fact, the Hardy subspace $H^{2}\subset L^{2}\left(  \mathbb{T}%
\right)  $ reduces this representation. For more details, see Section
\textup{\ref{Irr}} below.
\end{remark}

\section{\label{Irr}Irreducibility}

In this section we consider
as in Lemma \ref{LemMin.3} (\ref{eqMin.21})
the finite-dimensional subspace
\begin{equation}
\mathcal{K}:=\operatorname*{span}\{z^{-k};0\leq k\leq r_{0}\}\subset
L^{2}(\mathbb{T}),\qquad 
r_{0}=\left\lfloor \frac{gN-1}{N-1}\right\rfloor ,\label{eqIrrNew.1}
\end{equation}
defined from a polynomial loop $A(z)$ of scale size $N$ and genus $g$, and we
show that the irreducibility property for the corresponding representation
$T^{(A)}$ of $\mathcal{O}_{N}$ is \emph{generic}, i.e., it holds for all $A$
except for a subvariety of smaller dimension, once $N$ and $g$ are fixed. In
order to apply the results in Sections \ref{Min} and \ref{Irr}, some more
details are needed regarding the subspace $\mathcal{K}$, and they are taken up
in Section \ref{Exp}.

We begin with some notation and a lemma:

\begin{notation}
\label{notation5.2}Let $e_{n}(z):=z^{n},n\in\mathbb{Z}$, denote the Fourier
basis in $L^{2}(\mathbb{T})$. For finite subsets $J\subset\mathbb{Z}$, set
$\left\langle J\right\rangle :=\operatorname*{span}\left\{  e_{j};j\in
J\right\}  \subset L^{2}(\mathbb{T}).$ If
\begin{equation}
J_{0}=\left\{  0,-1,-2,\dots,-r_{0}\right\}  , \label{eq5.16}%
\end{equation}
set $\mathcal{K}:=\left\langle J_{0}\right\rangle $. If $T=T^{(A)}$ is a
wavelet representation and $r_{0}$ is as above, we note \cite[Proposition
5.5]{BrJo00} that $\mathcal{K}$ is cyclic.
\end{notation}

\begin{lemma}
\label{lemma5.1}Let $A\in\mathcal{P}(\mathbb{T},\mathrm{U}_{N}(\mathbb{C}))$
be a \textup{(}polynomial\/\textup{)}
loop of genus $g$. Then
\begin{equation}
\mathcal{K}=\left\langle \{0,-1,\dots,-r_{0}\}\right\rangle
=\operatorname*{span}\left\{  z^{-k};0\leq k\leq r_{0}\right\} ,\qquad
r_{0}=\left\lfloor \frac{gN-1}{N-1}\right\rfloor ,
\label{eqlemma5.1.1}%
\end{equation}
contains no one-dimensional subspace which is both $T_{i}^{(A)\,\ast}%
$-invariant, and also cyclic for the representation of $\mathcal{O}_{N}$ on
$L^{2}(\mathbb{T})$.
\end{lemma}

\begin{proof}
To show that a subspace $\mathcal{K}$ is minimal in the sense specified in the
lemma, we must check that whenever
\begin{equation}
(0\neq)\qquad\mathcal{K}_{0}\subsetneqq\mathcal{K} \label{eq5.2}%
\end{equation}
is a subspace satisfying
\begin{equation}
T_{i}^{(A)\,\ast}\mathcal{K}_{0}\subset\mathcal{K}_{0}\text{\qquad for
}i=0,1,\dots,N-1, \label{eq5.3}%
\end{equation}
then $\mathcal{K}_{0}$ cannot be \emph{cyclic} for the representation
$T^{(A)}$ of $\mathcal{O}_{N}$, i.e., it generates a cyclic subspace which is
a \emph{proper} subspace of $L^{2}(\mathbb{T}).$ The cyclic subspace generated
by $\mathcal{K}_{0}$ is the closed subspace spanned by
\begin{equation}
T_{i_{1}}^{(A)}\cdots T_{i_{n}}^{(A)}\mathcal{K}_{0}\text{\qquad for
}n=0,1,\dots, \label{eq5.3a}%
\end{equation}
and all multi-indices $(i_{1},\dots,i_{n})$. This follows from (\ref{eq5.3}),
and we will denote this space $\left[  \mathcal{O}_{N}\mathcal{K}_{0}\right]
.$ We will prove the assertion by checking that if (\ref{eq5.2})--(\ref{eq5.3}%
) hold, then there is an $m_{0}\in\mathbb{Z}$ such that $\left[
\mathcal{O}_{N}\mathcal{K}_{0}\right]  $ is contained in the closed span of
$\{z^{k};k\geq m_{0}\}$. Note that this integer $m_{0}$ might be negative, and
also that $\left[  \mathcal{O}_{N}\mathcal{K}_{0}\right]  $ might well be a
\emph{proper} subspace.

Now let $\mathcal{K}_{0}$ be given subject to conditions (\ref{eq5.2}%
)--(\ref{eq5.3}), and suppose (as in the lemma) that $\dim
\mathcal{K}_{0}=1$. Let $\xi\in\mathcal{K}_{0}$, $\left\|  \xi\right\|  =1$.
Then by (\ref{eq5.3}) there are $\lambda_{i}\in\mathbb{C}$ with
\begin{equation}
T_{i}^{(A)\,\ast}\xi=\overline{\lambda_{i}}\xi, \label{eq5.4}%
\end{equation}
or equivalently,%
\begin{equation}
\xi(z)=\sum_{i}\overline{\lambda_{i}}m_{i}^{(A)}(z)\xi(z^{N})\text{.}
\label{eq5.5}%
\end{equation}
Using \cite{Jor99a}, \cite{BrJo97b} we conclude that
\begin{equation}
\xi(z)=z^{-k}\text{\qquad for some }k\text{,} \label{eq5.6}%
\end{equation}
after adjusting with a constant multiple, and moreover that
\begin{equation}
\sum_{i}\overline{\lambda_{i}}m_{i}^{(A)}(z)=z^{k(N-1)}. \label{eq5.7}%
\end{equation}
Setting $\alpha(z):=\left(
1,
z,
\dots,
z^{N-1}
\right)  ^{\operatorname*{tr}}$, this may be rewritten as
\begin{equation}
\ip{\lambda}{A(z^{N})\alpha\left(  z\right)  }=z^{k(N-1)}\text{.}
\label{eq5.8}%
\end{equation}
Now pick $j\in\{0,1,\dots,N-1\}$ such that $-k\equiv j\operatorname*{mod}N$,
and apply the operators $S_{l}^{\ast}$, $l=0,1,\dots,N-1$, to both sides in
(\ref{eq5.8}). It follows that there is some $m\in\mathbb{Z}$ such that
\begin{equation}
\ip{\lambda}{A(z)\varepsilon_{j}}=z^{m} 
\text{\quad and\quad }
\ip{\lambda}{A(z)\varepsilon_{l}}=0 
\text{\quad if }l\neq j.
\label{eq5.9and10}
\end{equation}
Since
$\left|  \ip{\lambda}{A(z)\varepsilon_{j}}\right|  \leq\left\|  \lambda
\right\|  \left\|  A(z)\varepsilon_{j}\right\|  =1
$, the first part of (\ref{eq5.9and10})
implies equality in a Schwarz inequality. Then
(\ref{eq5.9and10}) yields $A(z)^{\ast}\lambda=z^{-m}\varepsilon_{j}%
$, or equivalently,
\begin{equation}
A(z)\varepsilon_{j}=z^{m}\lambda. \label{eq5.11}%
\end{equation}
Using the formula in Lemma \ref{CorBrJo00Rep.2} for the coefficients in
$A(z)$, and Lemma \ref{lemma3.1}, we note that (\ref{eq5.11}) implies
\begin{equation}
A(z)\varepsilon_{j}=z^{m}V\varepsilon_{j}, \label{eq5.12}%
\end{equation}
and in particular $\lambda=V\varepsilon_{j}$, where $V\in\mathrm{U}%
_{N}(\mathbb{C})$ is as in Lemma \ref{CorBrJo00Rep.2}.

For the convenience of the reader, we sketch the argument, but only in the
simplest case $m=0$. With the notation of Lemma \ref{CorBrJo00Rep.2}, we get%
\begin{equation}
A^{(0)}\varepsilon_{j}=\lambda,\qquad
A^{(k)}\varepsilon_{j}=0,\qquad 1\leq k<g, \label{eq5.13aandb}%
\end{equation}
where $j$ is the (fixed) number determined from (\ref{eq5.8}) as described.
Introducing the projections $Q_{0},Q_{1},\dots\in\mathcal{B}(\mathbb{C}^{N})$
from Lemma \ref{CorBrJo00Rep.2},
the first part of (\ref{eq5.13aandb}) then reads%
\begin{equation}
VQ_{0}^{\perp}Q_{1}^{\perp}\cdots Q_{g-1}^{\perp}\varepsilon_{j}=\lambda,
\label{eq5.13c}%
\end{equation}
and since
\begin{equation}
\left\|  \lambda\right\|  =\left(  \sum_{i}\left|  \lambda_{i}\right|
^{2}\right)  ^{\frac{1}{2}}=1,\qquad\left\|  Q_{0}^{\perp}Q_{1}^{\perp}\cdots
Q_{g-1}^{\perp}\varepsilon_{j}\right\|  =1. \label{eq5.13d}%
\end{equation}
Hence $Q_{1}^{\perp}\cdots Q_{g-1}^{\perp}\varepsilon_{j}$ is in the range of
$Q_{0}^{\perp}$, and $Q_{0}^{\perp}(Q_{1}^{\perp}\cdots Q_{g-1}^{\perp
}\varepsilon_{j})=Q_{1}^{\perp}\cdots Q_{g-1}^{\perp}\varepsilon_{j}$. We get
\begin{equation}
\left\|  Q_{1}^{\perp}\cdots Q_{g-1}^{\perp}\varepsilon_{j}\right\|  =1,
\label{eq5.13e}%
\end{equation}
and by induction,
\begin{equation}
Q_{g-1}^{\perp}\varepsilon_{j}=Q_{g-2}^{\perp}\varepsilon_{j}=\cdots
=Q_{1}^{\perp}\varepsilon_{j}=Q_{0}^{\perp}\varepsilon_{j}=\varepsilon_{j}
\label{eq5.13f}%
\end{equation}
and therefore $Q_{k}\varepsilon_{j}=0$, and $A(z)\varepsilon_{j}%
=V\varepsilon_{j}=\lambda$. This proves the claim (\ref{eq5.12}). Since we
wish to prove that $\left[  \mathcal{O}_{N}\mathcal{K}_{0}\right]  $ is
contained in $z^{m}\mathcal{H}_{+}$ for some $m\in\mathbb{Z}$, where
$\mathcal{H}_{+}:=\overline{\operatorname*{span}}\left\{  z^{k};0\leq
k\right\} $
is the Hardy space in $L^{2}(\mathbb{T})$, we may assume that $m$ in
(\ref{eq5.12}) is taken as $m=0$. Since the invariant subspaces for a
representation $\pi$ of $\mathcal{O}_{N}$ are the same as for $\pi\circ
\alpha_{V}$ where $\alpha_{V}=\operatorname*{Ad}V$, we may replace $A(z)$ with
$V^{-1}A(z)$, or equivalently, reduce to the special case $A(z)\varepsilon
_{j}=\varepsilon_{j}$ of formula (\ref{eq5.12}). Since the matrix of the basis
permutation $\varepsilon_{0}\leftrightarrow\varepsilon_{j}$ is in
$\mathrm{U}_{N}(\mathbb{C})$, the same argument leaves us with the simpler
case $A(z)\varepsilon_{0}=\varepsilon_{0}$, or equivalently,
\begin{equation}
A(z)=\left(
\begin{array}
[c]{c|ccc}%
1 & 0 & \cdots & 0\\\hline
0 & A_{1,1}(z) & \cdots & A_{1,N-1}(z)\\
\vdots & \vdots & \ddots & \vdots\\
0 & A_{N-1,1}(z) & \cdots & A_{N-1,N-1}(z)
\end{array}
\right)  . \label{eq5.15}%
\end{equation}
Then \cite[Theorem 6.2]{BrJo00} implies $\xi(z)=z^{0}=1$, or equivalently
$k=0$ in (\ref{eq5.6}), and so $\left[  \mathcal{O}_{N}\xi\right]
\subset\mathcal{H}_{+}$. Since $\mathcal{K}_{0}=\mathbb{C}\xi$, we have proved
that $\mathcal{K}_{0}$ is not cyclic, i.e., the representation on the single
vector $\xi$ does not generate all of $L^{2}(\mathbb{T})$. This concludes the
proof of the lemma.
\end{proof}

\begin{remark}
\label{remarknew}If $A\in\mathcal{P}(\mathbb{T},\mathrm{U}_{N}(\mathbb{C}))$,
then the wavelet representation $T^{(A)}$ may or may not be irreducible, as a
representation of $\mathcal{O}_{N}$ on $L^{2}(\mathbb{T})$. It is not
irreducible, for example, for the Haar wavelet. Nonetheless, we will show that
when $N$ and $g$ are given, then irreducibility holds
\emph{generically} for
$T^{(A)}$ as $A$ ranges over all $\mathcal{P}_{g}(\mathbb{T},\mathrm{U}%
_{N}(\mathbb{C}))$, i.e., has scale number $N$ and genus $g$.

Now for $T^{(A)}$ irreducible, of course every $\xi\in L^{2}(\mathbb{T})$,
$\xi\neq0$, will be cyclic, so cyclicity will then not be an issue, but rather
the question of when $\xi$ satisfies
\begin{equation}
T_{i}^{(A)\ast}\xi\in\mathbb{C}\xi\qquad\text{ for all }i\text{. }
\label{eq5new1}%
\end{equation}
The lemma states that (\ref{eq5new1}) cannot hold for $\xi\neq0$ if $T^{(A)}$
is irreducible. The issue is then instead to find \emph{minimal} subspaces
$\mathcal{K}_{\operatorname*{red}}$ such that
\begin{equation}
T_{i}^{(A)\ast}\left(  \mathcal{K}_{\operatorname*{red}}\right)
\subset\mathcal{K}_{\operatorname*{red}}\text{. } \label{eq5new2}%
\end{equation}
These subspaces are relevant for the algorithms which are used in the
construction of and the analysis of wavelets, as we sketched in Section
\ref{Int}, as the subspaces $\mathcal{K}_{\operatorname*{red}}$ in
$L^{2}(\mathbb{T})\cong\ell^{2}$ correspond to subspaces in the associated
multiresolution subspaces in $L^{2}(\mathbb{R})$ \textup{(}See Table
\textup{\ref{tablegeomappr}).} We first addressed (\ref{eq5new2}) in
\cite{BrJo00}, but the minimality was not considered there. We also note that
(\ref{eq5new2}) has applications for different representations of
$\mathcal{O}_{N}$, and is there connected with finitely correlated states in
statistical mechanics, see \cite{BrJo97a,FNW94,FNW92}. The
minimality issue was also considered in \cite{DKS99} in a different context.
We noted in (\ref{eqMin.21}) that
$\mathcal{K}=\operatorname*{span}\{z^{0},z^{-1},\dots,z^{-r_{0}}\}$,
with $r_{0}$ as in (\ref{eqMin.20}),
satisfies (\ref{eq5new2}). Let $P$ or $P_{\mathcal{K}}$ denote the projection
onto the subspace $\mathcal{K}$. We will consider subspaces of
$\mathcal{K}$ which still satisfy (\ref{eq5new2}), are cyclic, and
minimal with respect to (\ref{eq5new2}) and cyclicity. If, for example, $g=3$
and $N=2$, we will show that for some $A\in\mathcal{P}_{3}(\mathbb{T}%
,\mathrm{U}_{2}(\mathbb{C}))$ we may have
\begin{equation}
\mathcal{K}_{\operatorname*{red}}=\operatorname*{span}\left\{  z^{-2}%
,z^{-3}\right\}  . \label{eq5new3bis}%
\end{equation}
This is a little surprising since then $r_{0}=5$, and so $\mathcal{K}$ in
(\ref{eqMin.21}) is $6$-dimensional.
\end{remark}

\begin{corollary}
\label{corollary5.3}Let $J_{0}$ be as above in $(\ref{eq5.16})$, and consider
the two finite-dimensional subspaces $\mathcal{K}_{0}=\left\langle
J_{0}\backslash\{0\}\right\rangle $ and $\mathcal{K}_{1}=\left\langle
J_{0}\backslash\{-r_{0}\}\right\rangle $ in $L^{2}(\mathbb{T})$. Then
$\mathcal{K}_{0}$ is non-cyclic if
\begin{equation}
T_{i}^{(A)\ast}e_{0}\in\mathbb{C}e_{0}\qquad\text{for all }i. \label{eq5.17}%
\end{equation}
Suppose $N-1$ divides $gN-1$. Then $\mathcal{K}_{1}$ is non-cyclic if
\begin{equation}
T_{i}^{(A)\ast}\left(  e_{-r_{0}}\right)  \in\mathbb{C}e_{-r_{0}}\text{\qquad
for all }i. \label{eq5.18}%
\end{equation}
Moreover, $\mathcal{K}_{0}$ is cyclic if $\lambda_{0}=\sum_{i}\left|
A_{i,0}^{(0)}\right|  ^{2}=0$.
\end{corollary}

\begin{proof}
The two vectors $e_{0}$ and $e_{-r_{0}}$, corresponding to the endpoints in
$J_{0}$, are special in that
\begin{equation}
PT_{i}^{(A)}e_{0}\in\mathbb{C}e_{0}\text{, } \label{eq5.18bis}%
\end{equation}
and when $N-1$ divides $gN-1$,%
\begin{equation}
PT_{i}^{(A)}e_{-r_{0}}\in\mathbb{C}e_{-r_{0}}\text{\qquad for all }i.
\label{eq5.18ter}%
\end{equation}
It follows that, when (\ref{eq5.17}) holds, then
\begin{equation}
\left[  \mathcal{O}_{N}\left\langle \{0\}\right\rangle \right]  \oplus\left[
\mathcal{O}_{N}\left\langle J_{0}\backslash\{0\}\right\rangle \right]
=L^{2}(\mathbb{T}) \label{eq5.19}%
\end{equation}
where $\left\langle \{0\}\right\rangle =\mathbb{C}z^{0}=$ the one-dimensional
space of the constants, and neither of the two subspaces in this orthogonal
sum is zero. Hence (\ref{eq5.17}) implies that $\mathcal{K}_{0}=\left\langle
J_{0}\backslash\{0\}\right\rangle $ is not cyclic. The same argument proves
that $\mathcal{K}_{1}=\left\langle J_{0}\backslash\{-r_{0}\}\right\rangle $ is
not cyclic if (\ref{eq5.18}) holds.

We also note that (\ref{eq5.17}) holds if and only if%
\begin{equation}
\sum_{i}\left|  A_{i,0}^{(0)}\right|  ^{2}=1\text{,\qquad i.e., }\lambda
_{0}\left(  A\right)  =1, \label{eq5.20}%
\end{equation}
or equivalently,%
\begin{equation}
R(0,0)_{0,0}\left(  =\left(  A^{(0)\ast}A^{(0)}\right)  _{0,0}\right)  =1.
\label{eq5.21}%
\end{equation}
(The issue is resumed in Remark \ref{remark5.4} below.)

We now turn to the converse implications in the corollary, doing the details
only for $\mathcal{K}_{0}=\left\langle J_{0}\backslash\{0\}\right\rangle .$ If
(\ref{eq5.17}) does not hold, then
\begin{equation}
\lambda_{0}:=\sum_{i}\left|  A_{i,0}^{(0)}\right|  ^{2} \label{eq5.22}%
\end{equation}
satisfies $\lambda_{0}<1.$ Using the following formula,
\begin{align}
\ip{e_{0}}{T_{i}^{(A)}(e_{-k})}  &  =\ip{T_{i}^{(A)\,\ast} e_{0}}{z^{-k}}
=\ip{\overline{A_{i,0}(z)\vphantom{z^{k}}}}{\overline{z^{k}}}\label{eq5.23}\\
&  =\ip{z^{k}}{A_{i,0}(z)}=A_{i,0}^{(k)},\nonumber
\end{align}
we therefore have the conclusion: For each $i$, $T_{i}^{(A)\ast}e_{0}$ is in
$\mathcal{K}$, and it splits according to the sum $\mathcal{K}=\mathbb{C}%
e_{0}+\mathcal{K}_{0}$ as follows:%
\begin{equation}
T_{i}^{(A)\ast}e_{0}=\overline{A_{i,0}^{(0)}}e_{0}+\xi_{i},\text{\qquad}0\leq
i<N, \label{eq5.24}%
\end{equation}
where $\xi_{i}\in\mathcal{K}_{0}$ is computed according to formula
(\ref{eq5.23}). Applying $PT_{i}^{(A)}$ to both sides in (\ref{eq5.24}), we
conclude that
\[
e_{0}=\sum_{i}PT_{i}^{(A)}T_{i}^{(A)\ast}e_{0}\in\left(  \sum_{i}\left|
A_{i,0}^{(0)}\right|  ^{2}\right)  e_{0}+P\left[  \mathcal{O}_{N}%
\mathcal{K}_{0}\right]  ,
\]
or equivalently,
\begin{equation}
\left(  1-\lambda_{0}\right)  e_{0}\in P\left[  \mathcal{O}_{N}\mathcal{K}%
_{0}\right]  . \label{eq5.25}%
\end{equation}
Since in the second part of the corollary, $\lambda_{0}=0$ by assumption, we
conclude from (\ref{eq5.24}) that $e_{0}\in\left[  \mathcal{O}_{N}%
\mathcal{K}_{0}\right]  $, and the inclusion
\begin{equation}
\underset{=J_{0}}{\langle\underbrace{\{0,-1,\dots,-r_{0}\}}\rangle}%
\subset\left[  \mathcal{O}_{N}\mathcal{K}_{0}\right]  \label{eq5.25bis}%
\end{equation}
follows. Finally, we get
\begin{equation}
L^{2}(\mathbb{T})=\left[  \mathcal{O}_{N}\left\langle J_{0}\right\rangle
\right]  \subset\left[  \mathcal{O}_{N}\mathcal{K}_{0}\right]  \subset
L^{2}(\mathbb{T}), \label{eq5.26}%
\end{equation}
which proves that the reduced space $\mathcal{K}_{0}$ is then cyclic, and the
corollary is established.
\end{proof}

\begin{remark}
\label{remark5.4}To summarize, $T^{(A)}$ is given by some
$A\in\mathcal{P}(\mathbb{T},\mathrm{U}_{N}(\mathbb{C}))$,
where the coefficients $A^{(0)},A^{(1)},\dots\in M_{N}(\mathbb{C})$ satisfy
$A^{(g-1)}\neq0$ and
the conditions in
(\ref{eqBrJo00Rep.18}) and
Lemma \ref{CorBrJo00Rep.2}. The conclusions about cyclicity
in the previous corollary may be restated as follows in terms of these
matrices:
\begin{equation}
e_{0}\not \in \left[  \mathcal{O}_{N}\left\langle J_{0}\backslash
\{0\}\right\rangle \right]  \Longleftrightarrow A\text{ has the form }
\label{eq5.28}%
\end{equation}%
\[
\setlength{\unitlength}{15pt}\setlength{\bracelength}{3.8\unitlength}%
\begin{picture}(24,6.25)(0,-1.75)
\put(0,0){\makebox(1,1){$\ast$}}
\put(0,1){\makebox(1,1){$\raisebox{4pt}{\vdots}$}}
\put(0,2){\makebox(1,1){$\ast$}}
\put(0,3){\makebox(1,1){$\ast$}}
\multiput(1,0)(1,0){3}{\makebox(1,1){$\cdot$}}
\multiput(1,2)(1,0){3}{\makebox(1,1){$\cdot$}}
\multiput(1,3)(1,0){3}{\makebox(1,1){$\cdot$}}
\put(1,1){\makebox(3,1){$M$}}
\put(4,0){\makebox(1,1){$0$}}
\put(4,1){\makebox(1,1){$\raisebox{4pt}{\vdots}$}}
\put(4,2){\makebox(1,1){$0$}}
\put(4,3){\makebox(1,1){$0$}}
\multiput(5,0)(1,0){3}{\makebox(1,1){$\cdot$}}
\multiput(5,2)(1,0){3}{\makebox(1,1){$\cdot$}}
\multiput(5,3)(1,0){3}{\makebox(1,1){$\cdot$}}
\put(5,1){\makebox(3,1){$I$}}
\put(8,0){\makebox(1,1){$0$}}
\put(8,1){\makebox(1,1){$\raisebox{4pt}{\vdots}$}}
\put(8,2){\makebox(1,1){$0$}}
\put(8,3){\makebox(1,1){$0$}}
\multiput(9,0)(1,0){3}{\makebox(1,1){$\cdot$}}
\multiput(9,2)(1,0){3}{\makebox(1,1){$\cdot$}}
\multiput(9,3)(1,0){3}{\makebox(1,1){$\cdot$}}
\put(9,1){\makebox(3,1){$S$}}
\put(12,0){\makebox(4,4){$\cdots$}}
\put(17,1){\makebox(3,1){$C$}}
\put(20,0){\makebox(1,1){$0$}}
\put(20,1){\makebox(1,1){$\raisebox{4pt}{\vdots}$}}
\put(20,2){\makebox(1,1){$0$}}
\put(20,3){\makebox(1,1){$0$}}
\multiput(21,0)(1,0){3}{\makebox(1,1){$\cdot$}}
\multiput(21,2)(1,0){3}{\makebox(1,1){$\cdot$}}
\multiput(21,3)(1,0){3}{\makebox(1,1){$\cdot$}}
\put(21,1){\makebox(3,1){$L$}}
\multiput(0,0)(0,4){2}{\line(1,0){24}}
\multiput(0,0)(4,0){4}{\line(0,1){4}}
\multiput(20,0)(4,0){2}{\line(0,1){4}}
\multiput(0,0)(4,0){3}{\dashbox{0.235294}(1,4){}}
\put(20,0){\dashbox{0.235294}(1,4){}}
\put(0,-2.1){\makebox(4,2)[t]{$\underbrace{\hbox to\bracelength{}
}_{\displaystyle A^{\left(  0\right)  }}$}}
\put(4,-2.1){\makebox(4,2)[t]{$\underbrace{\hbox to\bracelength{}
}_{\displaystyle A^{\left(  1\right)  }}$}}
\put(20,-2.1){\makebox(4,2)[t]{$\underbrace{\hbox to\bracelength{}
}_{\displaystyle A^{\left(  g-1\right) \rlap{\quad ,}}}$}}
\end{picture}%
\]
i.e., all in the first of the first columns; and
\begin{equation}
e_{-r_{0}}\not \in \left[  \mathcal{O}_{N}\left\langle J_{0}\backslash
\{-r_{0}\}\right\rangle \right]  \Longleftrightarrow A\text{ has the form }
\label{eq5.29}%
\end{equation}%
\[
\setlength{\unitlength}{15pt}\setlength{\bracelength}{3.8\unitlength}%
\begin{picture}(24,6.25)(0,-1.75)
\multiput(0,0)(1,0){3}{\makebox(1,1){$\cdot$}}
\multiput(0,2)(1,0){3}{\makebox(1,1){$\cdot$}}
\multiput(0,3)(1,0){3}{\makebox(1,1){$\cdot$}}
\put(0,1){\makebox(3,1){$M$}}
\put(3,0){\makebox(1,1){$0$}}
\put(3,1){\makebox(1,1){$\raisebox{4pt}{\vdots}$}}
\put(3,2){\makebox(1,1){$0$}}
\put(3,3){\makebox(1,1){$0$}}
\multiput(4,0)(1,0){3}{\makebox(1,1){$\cdot$}}
\multiput(4,2)(1,0){3}{\makebox(1,1){$\cdot$}}
\multiput(4,3)(1,0){3}{\makebox(1,1){$\cdot$}}
\put(4,1){\makebox(3,1){$I$}}
\put(7,0){\makebox(1,1){$0$}}
\put(7,1){\makebox(1,1){$\raisebox{4pt}{\vdots}$}}
\put(7,2){\makebox(1,1){$0$}}
\put(7,3){\makebox(1,1){$0$}}
\multiput(8,0)(1,0){3}{\makebox(1,1){$\cdot$}}
\multiput(8,2)(1,0){3}{\makebox(1,1){$\cdot$}}
\multiput(8,3)(1,0){3}{\makebox(1,1){$\cdot$}}
\put(8,1){\makebox(3,1){$S$}}
\put(11,0){\makebox(1,1){$0$}}
\put(11,1){\makebox(1,1){$\raisebox{4pt}{\vdots}$}}
\put(11,2){\makebox(1,1){$0$}}
\put(11,3){\makebox(1,1){$0$}}
\put(12,0){\makebox(4,4){$\cdots$}}
\multiput(16,0)(1,0){3}{\makebox(1,1){$\cdot$}}
\multiput(16,2)(1,0){3}{\makebox(1,1){$\cdot$}}
\multiput(16,3)(1,0){3}{\makebox(1,1){$\cdot$}}
\put(16,1){\makebox(3,1){$C$}}
\put(19,0){\makebox(1,1){$0$}}
\put(19,1){\makebox(1,1){$\raisebox{4pt}{\vdots}$}}
\put(19,2){\makebox(1,1){$0$}}
\put(19,3){\makebox(1,1){$0$}}
\multiput(20,0)(1,0){3}{\makebox(1,1){$\cdot$}}
\multiput(20,2)(1,0){3}{\makebox(1,1){$\cdot$}}
\multiput(20,3)(1,0){3}{\makebox(1,1){$\cdot$}}
\put(20,1){\makebox(3,1){$L$}}
\put(23,0){\makebox(1,1){$\ast$}}
\put(23,1){\makebox(1,1){$\raisebox{4pt}{\vdots}$}}
\put(23,2){\makebox(1,1){$\ast$}}
\put(23,3){\makebox(1,1){$\ast$}}
\multiput(0,0)(0,4){2}{\line(1,0){24}}
\multiput(0,0)(4,0){7}{\line(0,1){4}}
\multiput(3,0)(4,0){3}{\dashbox{0.235294}(1,4){}}
\multiput(19,0)(4,0){2}{\dashbox{0.235294}(1,4){}}
\put(0,-2.1){\makebox(4,2)[t]{$\underbrace{\hbox to\bracelength{}
}_{\displaystyle A^{\left(  0\right)  }}$}}
\put(4,-2.1){\makebox(4,2)[t]{$\underbrace{\hbox to\bracelength{}
}_{\displaystyle A^{\left(  1\right)  }}$}}
\put(20,-2.1){\makebox(4,2)[t]{$\underbrace{\hbox to\bracelength{}
}_{\displaystyle A^{\left(  g-1\right) \rlap{\quad ,}}}$}}
\end{picture}%
\]
i.e., all in the last of the last columns. It follows from Theorem
\textup{\ref{ThmMin.7}} that the representation $T^{\left(  A\right)  }$ is
\emph{reducible} if either one of the conditions \textup{(\ref{eq5.28})} or
\textup{(\ref{eq5.29})} holds.

We can show, using \cite[Theorem 6.2]{BrJo00}, that a wavelet
representation $T^{(A)}$ satisfies \textup{(\ref{eq5.28})} if and only if
there are $V\in\mathrm{U}_{N}(\mathbb{C})$ and $B\in\mathcal{P}(\mathbb{T}%
,\mathrm{U}_{N-1}(\mathbb{C}))$ such that
$A$ has the form
(\ref{eqThmMin.7proof.7}).
There is a similar conclusion concerning the other condition
\textup{(\ref{eq5.29}).}

Let us say that, for fixed $N$ and $g$, a property is \emph{generic} if it
holds for all loops $A(z)$ of scale $N$ and genus $g$, except for $A$ in a
variety of lower dimension. Then we conclude from \textup{(\ref{eq5.28}%
)--(\ref{eq5.29})} that $\mathcal{K}_{\operatorname*{red}}=\left\langle
\left\{  -1,\dots,-(r_{0}-1)\right\}  \right\rangle $ or $\left\langle
\{-1,\dots,-(r_{0}-1),-r_{0}\}\right\rangle $ is cyclic for a generic set of
loops, when $g$ and $N$ are fixed.
The process described in (\ref{eq5.22}) of elimination starting with the
elimination of, if possible, $0$ and $-r_{0}$ from $\left\langle
\{0,-1,\dots,-r_{0}\}\right\rangle $ to get a smaller space, say
$\mathcal{K}_{0}$ such that%
\begin{equation}
T_{i}^{(A)\ast}(\mathcal{K}_{0})\subset\mathcal{K}_{0}\text{\qquad for all }i,
\label{eq5.32a}%
\end{equation}
and
\begin{equation}
\mathcal{K}_{0}\text{ is }\mathcal{O}_{N}\text{-cyclic in }L^{2}(\mathbb{T}),
\label{eq5.32b}%
\end{equation}
may be continued, subject to certain spectral conditions on the given loop
$A$. These conditions are \emph{generic} in the same sense; for example, as
noted in (\ref{eq5.22}), $e_{0}$ can be eliminated (so that the smaller
$\mathcal{K}_{0}$ will still satisfy (\ref{eq5.32a})--(\ref{eq5.32b})) if and
only if $R(0,0)_{0,0}<1$. There is a similar spectral condition for the
elimination of two vectors $e_{0}$ \emph{and} $e_{-1}$, i.e., for getting
$\mathcal{K}_{0}=\left\langle \left\{  -2,-3,\dots,-r_{0}\right\}
\right\rangle $ to also satisfy (\ref{eq5.32a})--(\ref{eq5.32b}): For the
$N=2$ case, this condition is that the $2$-by-$2$ matrix%
\begin{equation}
\left(
\begin{array}
[c]{ll}%
R(1,1)_{0,0} & R(0,1)_{0,1}\\
R(1,0)_{1,0} & R(0,0)_{1,1}%
\end{array}
\right)  \label{eq5brian-give-number}%
\end{equation}
does not have $1$ in its spectrum (recall $R(k,l):=A^{(l)\ast}A^{(k)}$). The
argument is the same as before, even if $N>2$, \emph{mutatis mutandis}. If $1$
is not in the spectrum, and if $R(0,0)_{0,0}=0$, then we show that both the
vectors $e_{0}$ and $e_{-1}$ are in the cut-down of the cyclic space, i.e., in%
\begin{equation}
P\left[  \mathcal{O}_{2}\left\langle \left\{  -2,-3,\dots,-r_{0}\right\}
\right\rangle \right]  , \label{eq5brianbis}%
\end{equation}
which is then $T_{i}^{(A)\ast}$-invariant. \textup{(}Here we use the symbol
$P$ for the projection onto the subspace $\mathcal{K}$ spanned by $\left\{
z^{-k};k=0,1,\dots,r_{0}\right\}  $ where $r_{0}=\left\lfloor \frac{Ng-1}%
{N-1}\right\rfloor $. If $N=2$, then, of course, $\mathcal{K}$ is of dimension
$2g$.\/\textup{)} Hence, the smaller subspace $\left\langle \{-2,-3,\dots
,-r_{0}\}\right\rangle \subsetneqq\left\langle \left\{  0,-1,-2,\dots
,-r_{0}\right\}  \right\rangle $ will also satisfy conditions (\ref{eq5.32a}%
)--(\ref{eq5.32b}).

To illustrate
the spectral condition more explicitly, we need $g>2$.
In the case $g=3,N=2$, there are $V\in\mathrm{U}_{2}(\mathbb{C})$, and
projections $P,Q$ in $\mathbb{C}^{2}$, such that%
\begin{equation}
A(z)=V(P^{\perp}+zP)(Q^{\perp}+zQ), \label{eqSpe.9}%
\end{equation}
and we get
\begin{equation}
\begin{gathered}
R(0,0)  =Q^{\perp}P^{\perp}Q^{\perp},\qquad
R(0,1)  =QP^{\perp}Q^{\perp},\qquad
R(1,0)  =Q^{\perp}P^{\perp}Q,\\
R(1,1)  =QP^{\perp}Q+Q^{\perp}PQ^{\perp},\qquad
R(2,2)  =QPQ.
\end{gathered}
\label{eqSpe.10to14}
\end{equation}
Hence, $R(0,0)_{0,0}=0$ holds if and only if $P^{\perp}Q^{\perp}%
\varepsilon_{0}=0$, or equivalently,
\begin{equation}
PQ\varepsilon_{0}=P\varepsilon_{0}+Q\varepsilon_{0}-\varepsilon_{0}.
\label{eqSpe.15}%
\end{equation}
The entries of the matrix (\ref{eq5brian-give-number}) are then
\begin{equation}
\left(
\begin{array}
[c]{clc}%
\left\|  P^{\perp}Q\varepsilon_{0}\right\|  ^{2}+\left\|  PQ^{\perp
}\varepsilon_{0}\right\|  ^{2} &  & \ip{\varepsilon_{0}}{P^{\perp}Q^{\perp
}\varepsilon_{1}}\\
\ip{P^{\perp}Q^{\perp}\varepsilon_{1}}{\varepsilon_{0}} &  & \left\|
P^{\perp}Q^{\perp}\varepsilon_{1}\right\|  ^{2}%
\end{array}
\right)  . \label{eqSpe.16}%
\end{equation}
Since
\begin{equation}
R(0,0)+R(1,1)+R(2,2)=\left(
\begin{array}
[c]{ll}%
1 & 0\\
0 & 1
\end{array}
\right)  , \label{eqSpe.17}%
\end{equation}
the restriction $R(0,0)_{0,0}=0$ therefore implies that
\begin{equation}
R(1,1)_{0,0}=1-R(2,2)_{0,0}=1-\left\|  PQ\varepsilon_{0}\right\|  ^{2}\text{.}
\label{eqSpe.18}%
\end{equation}
Using this, we get that $1$ is in the spectrum of
(\ref{eq5brian-give-number}), so reduced, if and only if
\begin{equation}
\left\|  PQ\varepsilon_{0}\right\|  ^{2}\cdot\left(  1-\left\|  P^{\perp
}Q^{\perp}\varepsilon_{1}\right\|  ^{2}\right)  =\left|  \ip{\varepsilon_{0}%
}{P^{\perp}Q^{\perp}\varepsilon_{1}}\right|  ^{2}. \label{eqSpe.19}%
\end{equation}
To solve this, let for example $P$ and $Q$ be the respective projections onto
\begin{equation}
\left(
\begin{array}
[c]{l}%
\cos\theta\\
\sin\theta
\end{array}
\right)  
\text{\quad and\quad }
\left(
\begin{array}
[c]{l}%
\cos\rho\\
\sin\rho
\end{array}
\right)  \label{eqSpe.20and21}%
\end{equation}
and solve for $\theta$ and $\rho$. For these examples (i.e., in the
complementary region of the $(\theta,\rho)$-plane), we will then have
$\mathcal{K}_{\operatorname*{red}}$ :=$\left\langle \left\{  -2,-3\right\}
\right\rangle $ satisfy the covariance condition as well as the cyclicity. All
of the cases $N=2$, $g=3$, will be taken up again in the Appendix, where the
algorithmic properties of \textup{(\ref{eqInt.9})} for the scaling function
$\varphi$ are displayed in detail. This is an iteration based on
\textup{(\ref{eqInt.9}),} and the regularity of the corresponding
$x\mapsto\varphi_{\theta,\rho}\left(  x\right)  $ turns out to depend on the
spectral properties of the operators in
\textup{(\ref{eqSpe.10to14}).}

We now turn to the distinction between the diagonal elements $A(z)$ in
\linebreak $\mathcal{P}(\mathbb{T},\mathrm{U}_{N}(\mathbb{C}))$, and the
non-diagonal ones. We say that $A$ is \emph{diagonal} if it maps into the
diagonal matrices in $\mathrm{U}_{N}(\mathbb{C})$, except for a constant
factor, i.e., if there is some $V\in\mathrm{U}_{N}(\mathbb{C})$, and
$n_{0},\dots,n_{N-1}\geq0$, such that
$A$ has the form (\ref{eqRemMin.8.15}).
The variety of these diagonal loops will be called $\mathcal{P}%
_{\operatorname{diag}}(\mathbb{T},\mathrm{U}_{N}(\mathbb{C}))$. Note that this
definition includes
\begin{equation}
\left(
\begin{array}
[c]{ll}%
0 & z\\
1 & 0
\end{array}
\right)  =\left(
\begin{array}
[c]{ll}%
0 & 1\\
1 & 0
\end{array}
\right)  \left(
\begin{array}
[c]{ll}%
1 & 0\\
0 & z
\end{array}
\right)  \label{eqSpe.23}%
\end{equation}
in $\mathcal{P}_{\operatorname{diag}}(\mathbb{T},\mathrm{U}_{2}(\mathbb{C}))$.
\end{remark}

\begin{definition}
\label{definition5.5}We say that a loop $A\in\mathcal{P}(\mathbb{T}%
,\mathrm{U}_{N}(\mathbb{C}))$ is \emph{purely non-diagonal} if there is not a
decomposition $N=d_{0}+b+d_{1}$ with $d_{0}>0$, or $d_{1}>0$, diagonal
elements $D_{i}(z)\in\mathcal{P}_{\operatorname{diag}}(\mathbb{T}%
,\mathrm{U}_{d_{i}}(\mathbb{C}))$, $i=0,1$, $B(z)\in\mathcal{P}(\mathbb{T}%
,\mathrm{U}_{b}(\mathbb{T}))$, and $V\in\mathrm{U}_{N}(\mathbb{C})$ such that
\begin{equation}
A(z)=V\left(
\begin{array}
[c]{ccc}%
D_{0}(z) & 0 & 0\\
0 & B(z) & 0\\
0 & 0 & D_{1}(z)
\end{array}
\right)  . \label{eqSpe.24}%
\end{equation}
\end{definition}

Easy examples of loop matrices for $N=2$ and $g=3$ which are not diagonal,
i.e., do not have the representation (\ref{eqRemMin.8.15}) or (\ref{eqSpe.24}) for
any $V$, are%
\begin{equation}
\frac{1}{\sqrt{2}}\left(
\begin{array}
[c]{cc}%
1 & z\\
z & -z^{2}%
\end{array}
\right)  \text{\quad and\quad}\frac{1}{\sqrt{2}}\left(
\begin{array}
[c]{cc}%
z & z^{2}\\
1 & -z
\end{array}
\right)  . \label{eqSpe.24bis}%
\end{equation}
Both have $\lambda_{0}\left(  A\right)  =1/2$. Both correspond to Haar
wavelets, and both are exceptional cases in the wider family of the Appendix
below. Both have irreducible wavelet representations, by the next theorem.
For more details about matrix factorizations in the loop groups, we refer the
reader to \cite{PrSe86}, and the paper \cite{AlPe99}.

We now turn to our first explicit result about minimal subspaces
$\mathcal{L}\subset\mathcal{K}$, i.e., subspaces $\mathcal{L}$ which are
$T^{\ast}$-invariant, cyclic, and which do not contain proper $T^{\ast}%
$-invariant subspaces which are also cyclic. In Theorem \ref{ThmExp.2} below,
we shall then further give a formula for the (unique) minimal such space
$\mathcal{L}$. We stress that these results are special for the wavelet
representations, and that they do not hold for other kinds of representations
of $\mathcal{O}_{N}$.

\begin{theorem}
\label{theorem5.6}$\mathrm{(a)}$ Let $A\in\mathcal{P}(\mathbb{T}%
,\mathrm{U}_{N}(\mathbb{C}))$ be given. Suppose it is purely non-diagonal, and
let $T^{(A)}$ be the corresponding wavelet representation of $\mathcal{O}_{N}$
on $L^{2}(\mathbb{T})$. Then it follows that $T^{(A)}$ is irreducible.

$\mathrm{(b)}$ Let
$r_{0}$ be as in \textup{(\ref{eqMin.20}).}
Let optimal numbers $p,q$, $0\leq p\leq q\leq r_{0}$
be determined by the spectral condition in Remark \textup{\ref{remark5.4}}
such that
\begin{equation}
\mathcal{K}_{\operatorname*{red}}=\left\langle \{-p,-(p+1),\dots
,-q\}\right\rangle \label{eqSpe.26}%
\end{equation}
is $T_{i}^{(A)\ast}$-invariant for all $i$, and further satisfies
\begin{equation}
\left\langle \{0,-1,\dots,-p+1,-q-1,\dots,-r_{0}\}\right\rangle \subset\left[
\mathcal{O}_{N}\mathcal{K}_{\operatorname*{red}}\right]  . \label{eqSpe.27}%
\end{equation}
Then the following three properties hold:

\begin{enumerate}
\item[(i)] $T_{i}^{(A)\ast}(\mathcal{K}_{\operatorname*{red}})\subset
\mathcal{K}_{\operatorname*{red}}$ for all $i$,

\item[(ii)] $\mathcal{K}_{\operatorname*{red}}$ is cyclic \textup{(}for
$L^{2}(\mathbb{T})$\textup{),}

\item[(iii)] $\mathcal{K}_{\operatorname*{red}}$ is minimal with respect to
properties $\mathrm{(i)}$--$\mathrm{(ii)}$.
\end{enumerate}

$\mathrm{(c)}$ The minimal space $\mathcal{K}_{\operatorname*{red}}$ from
$\mathrm{(b)}$ is reduced from the right if $N-1$ divides $gN-1$,
where $g$ is the genus,
and if not,
it is $\left\langle \{-p,\dots,-r_{0}\}\right\rangle $; so it is only
``truncated'' at one end when $N-1$ does not divide $gN-1$.
\end{theorem}

\begin{proof}
Once $\mathcal{K}_{\operatorname*{red}}$ has been chosen as in the statement
(b) of the theorem, the three properties (i)--(iii) follow from Theorem
\ref{ThmBJKW003.5} and \ref{ThmMin.7}. The significance of (i)--(iii) is that
they imply that if
\begin{equation}
\sigma(\;\cdot\;):=\sum_{i}P_{\mathcal{K}_{\operatorname*{red}}}T_{i}%
^{(A)}(\;\cdot\;)T_{i}^{(A)\ast}P_{\mathcal{K}_{\operatorname*{red}}},
\label{eqSpe.28}%
\end{equation}
then the fixed-point set $\mathcal{B}(\mathcal{K}_{\operatorname*{red}%
})^{\sigma}$ is in fact an algebra. This is a result of Davidson et al.\
\cite{DKS99}. Using Theorem \ref{ThmBJKW003.5}, we conclude that the
projections in $\mathcal{B}(\mathcal{K}_{\operatorname*{red}})^{\sigma}$ are
characterized by the condition of Lemma \ref{LemMin.3}. Now, by \cite{DKS99},
there are projections $E_{j}\in\mathcal{B}(\mathcal{K}_{\operatorname*{red}})$
such that, for each $i,j$, we have the covariance properties
\begin{equation}
E_{j}V_{i}^{\ast}E_{j}=V_{i}^{\ast}E_{j}, \label{eq5.36}%
\end{equation}
where $V_{i}^{\ast}=T_{i}^{(A)\ast}P_{\mathcal{K}_{\operatorname*{red}}}$, or
equivalently,
\begin{equation}
V_{i}=P_{\mathcal{K}_{\operatorname*{red}}}T_{i}^{(A)}\mathpunct{;}
\label{eq5.37}%
\end{equation}
and in addition, we have
\begin{equation}
\sum_{j}E_{j}=\openone_{\mathcal{K}_{\operatorname*{red}}}, \label{eq5.38}%
\end{equation}
and each subspace $\left[  \mathcal{O}_{N}E_{j}\mathcal{K}%
_{\operatorname*{red}}\right]  $ \emph{irreducible}, in the sense that each
$\left[  \mathcal{O}_{N}E_{j}\mathcal{K}_{\operatorname*{red}}\right]  $
reduces the representation $\mathcal{O}_{N}$ to one which is irreducible on
the subspace.

It follows from (\ref{eq5.36})--(\ref{eq5.38}) that the complementary
projection
\begin{equation}
\openone_{\mathcal{K}_{\operatorname*{red}}}-E_{j}=\sum_{l:l\neq j}E_{l}
\label{eq5.38bis}%
\end{equation}
then also satisfies (\ref{eq5.36}), and so in particular $E_{j}$ must commute
with each $V_{i}$ ($=P_{\mathcal{K}_{\operatorname*{red}}}T_{i}^{(A)}$). Then
by Theorems \ref{ThmBJKW003.5} and \ref{ThmMin.7}, we conclude that each
$E_{j}$ has a matrix which is diagonal with respect to the Fourier basis
$\left\{  z^{-k}\right\}  $. Since the loop $A(z)$ is picked to be purely
non-diagonal, we finally conclude that the decomposition $\{E_{j}\}$ of
(\ref{eq5.38}) can only have one term, and the proof is concluded.
\end{proof}

\begin{remark}
\label{remark5.7}Even if the assumption in Theorem \ref{theorem5.6}, to the
effect that $A$ be purely non-diagonal, is removed, we have the decomposition
into irreducibles, and these irreducibles $\left[  \mathcal{O}_{N}%
E_{j}\mathcal{K}_{\operatorname*{red}}\right]  $ are mutually disjoint, i.e.,
inequivalent representations when $j\neq j^{\prime}$ for two possible terms
$j,j^{\prime}$ in a decomposition. This follows from the Theorems
\ref{ThmBJKW003.5} and \ref{ThmMin.7}, which state that the projections
$E_{j}$ are all diagonal relative to the same basis (see also Theorem
\ref{ThmExp.2} below!). So in particular, $\mathcal{B}(\mathcal{K}%
_{\operatorname*{red}})^{\sigma}$ is abelian when $\mathcal{K}%
_{\operatorname*{red}}$ is chosen subject to conditions (i)--(iii) in the
statement of Theorem \ref{theorem5.6}. This means that the corresponding
decomposition of $T^{(A)}$ into a sum of irreducible representations of
$\mathcal{O}_{N}$ is multiplicity-free.
\end{remark}

\begin{example}
[An Application\/]\label{ExaIrr.application}Even though we list only the
scaling functions $\varphi\left(  x\right)  $ in the examples in the Appendix,
the \emph{wavelet generator} $\psi\left(  x\right)  $ is significant. But it
is not unique: We can have a loop $A$ of genus $2$, and a different one $B$ of
genus $3$, which have the same $\varphi$. Then, of course, there are different
wavelet generators, say $\psi_{A}^{{}}$ and $\psi_{B}^{{}}$. To see this, take
$\varphi$ as follows \textup{(}see also Remark \textup{\ref{RemRep.5}):}%
\begin{equation}%
\begin{array}
[c]{ccc}
& \setlength{\unitlength}{1bp}%
\begin{picture}
(265,129)(-9,0) \put(0,0){\includegraphics[bb=8 0 337 164,
height=123bp,width=247bp] {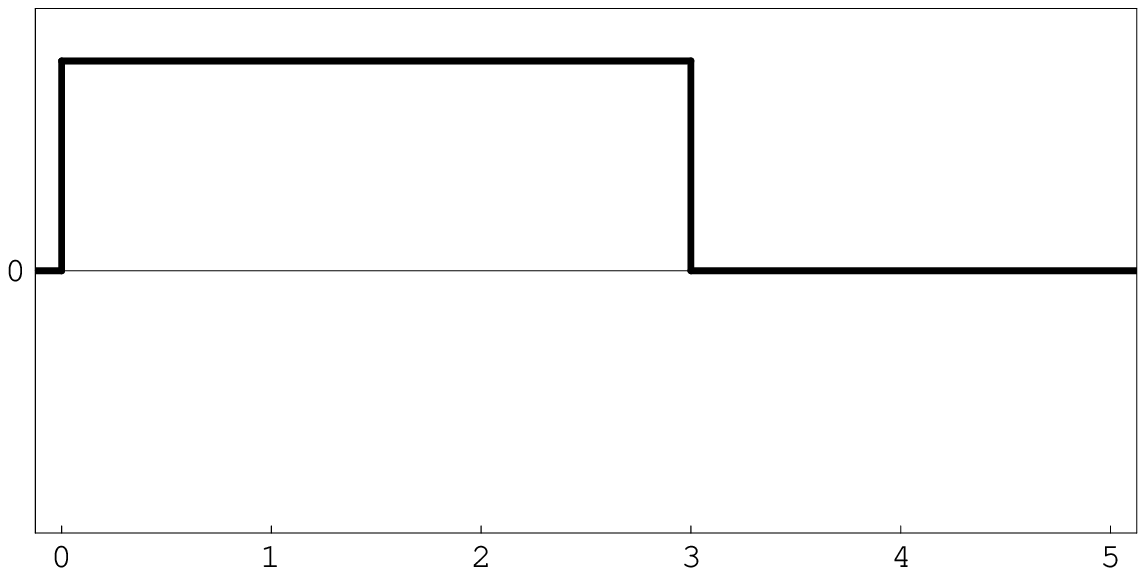}} \put(248,2){\makebox(0,12)[l]{$x$}}
\put(159,98){\makebox(0,12)[l]{$\varphi\left(x\right)$}}
\end{picture}%
&
\end{array}
\label{figExaIrr.application.1}%
\end{equation}
A loop $A$ in diagonal form giving this $\varphi$ is
\[
\frac{1}{\sqrt{2}}\left(
\begin{array}
[c]{cc}%
1 & 1\\
1 & -1
\end{array}
\right)  \left(
\begin{array}
[c]{cc}%
1 & 0\\
0 & z
\end{array}
\right)  \text{,\qquad genus }g=2.
\]
This is of the form \textup{(\ref{eqRemMin.8.15}).} The corresponding wavelet
generator $\psi_{A}^{{}}$ is then%
\begin{equation}%
\begin{array}
[c]{ccc}
& \setlength{\unitlength}{1bp}%
\begin{picture}
(265,129)(-9,0) \put(0,0){\includegraphics[bb=8 0 337 164,
height=123bp,width=247bp] {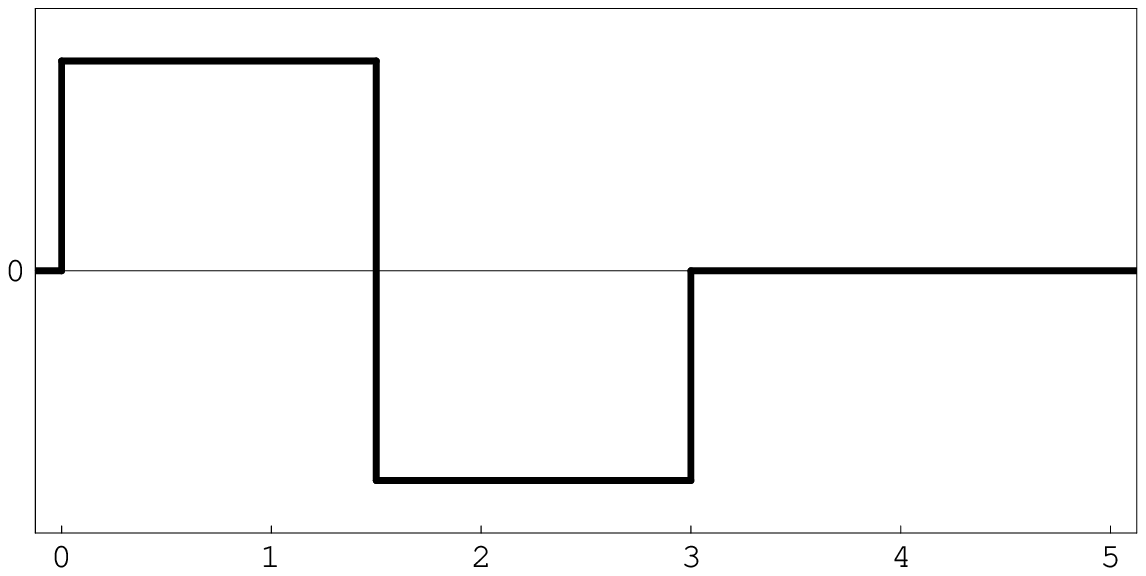}} \put(248,2){\makebox(0,12)[l]{$x$}}
\put(92,98){\makebox(0,12)[l]{$\psi_{A}^{}\left(x\right)$}}
\end{picture}%
&
\end{array}
\label{figExaIrr.application.2}%
\end{equation}
But setting%
\begin{equation}
B\left(  z\right)  =\frac{1}{\sqrt{2}}\left(
\begin{array}
[c]{cc}%
1 & z\\
z & -z^{2}%
\end{array}
\right)  , \label{eqExaIrr.application.1}%
\end{equation}
then this loop has the same $\varphi$. Since
\begin{equation}
m_{0}^{\left(  B\right)  }\left(  z\right)     =\frac{1}{\sqrt{2}}\left(
1+z^{3}\right)  =m_{0}^{\left(  A\right)  }\left(  z\right)
,\qquad
m_{1}^{\left(  B\right)  }\left(  z\right)   
=z^{2}m_{1}^{\left(  A\right)  }\left(  z\right)  ,
\label{eqExaIrr.application.2}%
\end{equation}
the corresponding wavelet generator $\psi_{B}^{{}}$ is now different from
$\psi_{A}^{{}}$ only by a translation.%
\begin{equation}%
\begin{array}
[c]{ccc}
& \setlength{\unitlength}{1bp}%
\begin{picture}
(265,129)(-9,0) \put(0,0){\includegraphics[bb=8 0 337 164,
height=123bp,width=247bp] {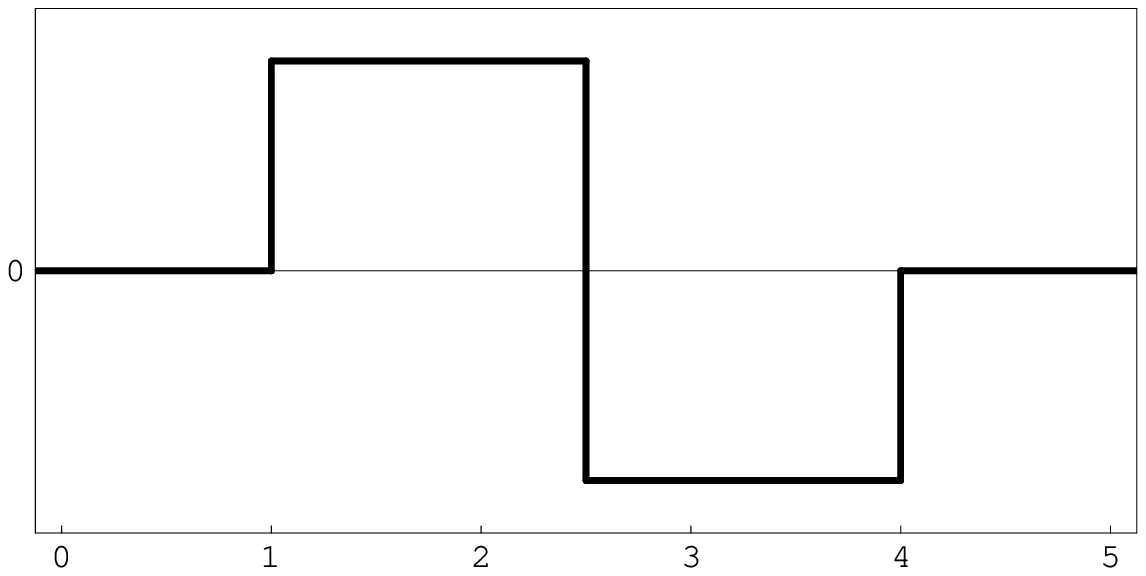}} \put(248,2){\makebox(0,12)[l]{$x$}}
\put(137,98){\makebox(0,12)[l]{$\psi_{B}^{}\left(x\right)$}}
\end{picture}%
&
\end{array}
\label{figExaIrr.application.3}%
\end{equation}
In fact, $\psi_{B}^{{}}\left(  x\right)  =\psi_{A}^{{}}\left(  x-1\right)  $.
However, the most striking contrast between the two loops $A$ and $B$ is that
the minimality question comes out differently from one to the next: The
representation $T^{\left(  A\right)  }$ of $\mathcal{O}_{2}$ on $L^{2}\left(
\mathbb{T}\right)  $ is \emph{reducible,} while $T^{\left(  B\right)  }$ is
\emph{irreducible,} i.e., there are no nonzero closed subspaces of
$L^{2}\left(  \mathbb{T}\right)  $, other than $L^{2}\left(  \mathbb{T}%
\right)  $, which are invariant under all $T_{i}^{\left(  B\right)  }$ and
$T_{i}^{\left(  B\right)  \,\ast}$. \textup{(}Or, stated equivalently, by
\textup{(\ref{eq3.2})} we have the implication $\sum_{i}T_{i}^{\left(
B\right)  }XT_{i}^{\left(  B\right)  \,\ast}=X$, $X\in\mathcal{B}\left(
L^{2}\left(  \mathbb{T}\right)  \right)  $, $\Rightarrow X\in\mathbb{C}%
\,\openone_{L^{2}\left(  \mathbb{T}\right)  }$.\textup{)} The two conclusions
for $T^{\left(  A\right)  }$ and $T^{\left(  B\right)  }$ follow from Lemma
\textup{\ref{lemma5.1},} Corollary \textup{\ref{corollary5.3},} and Theorems
\textup{\ref{ThmMin.7}} and \textup{\ref{theorem5.6},} respectively; but
Theorem \textup{\ref{ThmExp.2}} is also used. What is perhaps more surprising
is that the matrix loop
\begin{equation}
\mathcal{A}\left(  z\right)  :=B\oplus B=\frac{1}{\sqrt{2}}\left(
\begin{array}
[c]{cc|cc}%
1 & z & 0 & 0\\
z & -z^{2} & 0 & 0\\\hline
0 & 0 & 1 & z\\
0 & 0 & z & -z^{2}%
\end{array}
\right)  \label{eqExaIrr.application.4}%
\end{equation}
\textup{(}see \textup{(\ref{eqExaIrr.application.1}))} in $\mathrm{U}%
_{4}\left(  \mathbb{C}\right)  $, i.e., $N=4$ and $g=3$, defines a
representation $T^{\left(  \mathcal{A}\right)  }$ of $\mathcal{O}_{4}$ which
acts \emph{irreducibly} on $L^{2}\left(  \mathbb{T}\right)  $.
\end{example}

Our general result in this paper is that the wavelet representations are
\emph{irreducible,} except for isolated examples of Haar type, such as
$\psi_{A}^{{}}$ in \textup{(\ref{figExaIrr.application.2})}. But
\textup{(\ref{eqExaIrr.application.2})}
above shows that even for the reducible ones, irreducibility can still be
achieved, if only $\mathbb{Z}$-translations are allowed; see
\textup{(\ref{figExaIrr.application.3}).}
The following result is a corollary of Theorem \ref{theorem5.6}, and it helps
to distinguish the wavelet representations $T^{\left(  A\right)  }$ from the
more general representations of \cite{FNW92,FNW94} associated with finitely
correlated states in statistical mechanics. It is
a crucial distinction, and it is
concerned with the
completely positive maps $\sigma$ which are described in Theorems
\ref{ThmBrJo00Rep.1} and \ref{ThmBJKW003.5}. In \cite{FNW94}, the
representations are determined by maps $\sigma$ which possess
\emph{faithful}
invariant states, and these states play a role in the proofs of the results
there. Our next corollary asserts that such faithful invariant states do
\emph{not} exist for the wavelet representations.

\begin{corollary}
\label{CorIrr.next}Let $T^{\left(  A\right)  }$ be a wavelet representation of
$\mathcal{O}_{N}$ on $L^{2}\left(  \mathbb{T}\right)  $ which satisfies the
conditions in Theorem \textup{\ref{theorem5.6},} and let $\sigma_{\mathcal{K}%
}^{\left(  A\right)  }\left(  \,\cdot\,\right)  =\sum_{i}V_{i}\left(
\,\cdot\,\right)  V_{i}^{\ast}$ be the corresponding completely positive
mapping of Theorems \textup{\ref{ThmBrJo00Rep.1}} and
\textup{\ref{ThmBJKW003.5}.} Then there is no \emph{faithful} state $\rho$ on
$\mathcal{B}\left(  \mathcal{K}\right)  $ which leaves $\sigma_{\mathcal{K}%
}^{\left(  A\right)  }$ invariant, i.e., which satisfies%
\begin{equation}
\rho\circ\sigma_{\mathcal{K}}^{\left(  A\right)  }=\rho.
\label{eqCorIrr.next.1}%
\end{equation}
\end{corollary}

\begin{proof}
We will restrict to the case $N=2$, although for $g=2$, we cover arbitrary $N$
in \cite{BrJo00}. (If $g=2$, then $\mathcal{K}=\left\langle e_{0}%
,e_{-1},e_{-2},e_{-3}\right\rangle $. Setting $E_{-k,-l}:=\left|
e_{-k}\right\rangle \left\langle e_{-l}\right|  $, we showed in \cite{BrJo00}
that the density matrix $D$ given by $D=\lambda_{N-2}E_{-1,-1}+\left(
1-\lambda_{N-1}\right)  E_{-2,-2}$ satisfies $\sigma^{\ast}\left(  D\right)
=D$, where $\sigma=\sigma_{\mathcal{K}}^{\left(  A\right)  }$ and $\lambda
_{i}:=R\left(  0,0\right)  _{i,i}$, and where $\sigma^{\ast}$ is the adjoint
of $\sigma\colon\mathcal{B}\left(  \mathcal{K}\right)  \rightarrow
\mathcal{B}\left(  \mathcal{K}\right)  $ with respect to the trace inner
product. Defining the state $\rho$ on $\mathcal{B}\left(  \mathcal{K}\right)
$ by%
\begin{equation}
\rho\left(  X\right)  :=\operatorname*{trace}\left(  DX\right)  ,\qquad
\forall\,X\in\mathcal{B}\left(  \mathcal{K}\right)  ,
\label{eqCorIrr.next.proof2}%
\end{equation}
we check that $\rho$ satisfies (\ref{eqCorIrr.next.1}). We know from
\cite{BJKW00} that $\ker\left(  \sigma-\openone\right)  $ and $\ker\left(
\sigma^{\ast}-\openone\right)  $ have the same dimension. But $T^{\left(
A\right)  }$ is irreducible by Theorem \ref{theorem5.6} when $0<\lambda_{0}%
<1$. Hence $\ker\left(  \sigma-\openone\right)  $ is one-dimensional by
Theorem \ref{ThmBrJo00Rep.1}, and there are therefore no other states $\rho$
satisfying (\ref{eqCorIrr.next.1}). But the state $\rho$ in
(\ref{eqCorIrr.next.proof2}) is clearly not faithful, and the proof is
complete, in the special case $g=2$.)

We now turn to the details for $N=2$, $g=3$, and it will be clear that they
generalize to arbitrary $g$. If $N=2$, $g=3$, we get $\mathcal{K}=\left\langle
e_{0},e_{-1},e_{-2},e_{-3},e_{-4},e_{-5}\right\rangle \cong\mathbb{C}^{6}$,
and $T_{i}^{\left(  A\right)  \,\ast}e_{-k}$ may easily be computed; see,
e.g., the details in Section \ref{Exp} below, especially (\ref{eqExp.23a}%
)--(\ref{eqExp.23f}). As a result, we get $\sigma^{\ast}\left(  E_{-k,-l}%
\right)  =\sum_{i}\left|  T_{i}^{\ast}e_{-k}\right\rangle \left\langle
T_{i}^{\ast}e_{-l}\right|  $, and therefore%
\begin{align}
\sigma^{\ast}\left(  E_{-1,-1}\right)   &  =\sum_{k,l}R\left(  l,k\right)
_{1,1}E_{-1-k,-1-l}\,,\label{eqCorIrr.next.proof2a}\\
\sigma^{\ast}\left(  E_{-2,-2}\right)   &  =\sum_{k,l}R\left(  l,k\right)
_{0,0}E_{-1-k,-1-l}\,,\label{eqCorIrr.next.proof2b}\\
\sigma^{\ast}\left(  E_{-3,-3}\right)   &  =\sum_{k,l}R\left(  l,k\right)
_{1,1}E_{-2-k,-2-l}\,,\label{eqCorIrr.next.proof2c}\\
\sigma^{\ast}\left(  E_{-4,-4}\right)   &  =\sum_{k,l}R\left(  l,k\right)
_{0,0}E_{-2-k,-2-l}\,, \label{eqCorIrr.next.proof2d}%
\end{align}
where the $k,l$ summations are both over $\left\{  0,1,2\right\}  $. In
addition, by (\ref{eq5.18bis}) and (\ref{eq5.18ter}),
\begin{equation}
\sigma\left(  E_{0,0}\right)  =\lambda_{0}E_{0,0}\text{,\quad and\quad}%
\sigma\left(  E_{-5,-5}\right)  =\lambda_{0}E_{-5,-5},
\label{eqCorIrr.next.proof2e}%
\end{equation}
where $\lambda_{0}=\lambda_{0}\left(  A\right)  =R\left(  0,0\right)  _{0,0}$.
So the complement of $\left\langle E_{0,0},E_{-5,-5}\right\rangle $ in
$\mathcal{B}\left(  \mathcal{K}\right)  $ is invariant under $\sigma^{\ast}$,
and the element $D$ which is fixed by $\sigma^{\ast}$ must be diagonal in the
Fourier basis, by Theorem \ref{ThmMin.7}. Using Lemma \ref{CorBrJo00Rep.2} and
the argument from the previous step, we then check that a density matrix $D$
may be found in the form
\begin{equation}
D=\delta_{1}E_{-1,-1}+\delta_{2}E_{-2,-2}+\delta_{3}E_{-3,-3}+\delta
_{4}E_{-4,-4},\qquad\delta_{i}\geq0,\;\sum_{i}\delta_{i}=1,
\label{eqCorIrr.next.proof3}%
\end{equation}
such that the state $\rho\left(  \,\cdot\,\right)  =\operatorname*{trace}%
\left(  D\,\cdot\,\right)  $ on $\mathcal{B}\left(  \mathbb{C}^{6}\right)  $
will satisfy (\ref{eqCorIrr.next.1}). But if $0<\lambda_{0}\left(  A\right)
<1$, the wavelet representation $T^{\left(  A\right)  }$ is irreducible, and
so (\ref{eqCorIrr.next.1}) has no other state solutions. Finally, it is clear
from (\ref{eqCorIrr.next.proof3}) that $\rho$ is not faithful.
\end{proof}

\section{\label{Fil}Filtrations in $\mathcal{P}(\mathbb{T},\mathrm{U}%
_{2}(\mathbb{C}))$ as\\
factorizations of quadrature mirror filters}

Since $\mathcal{P}(\mathbb{T},\mathrm{U}_{N}(\mathbb{C}))$ has multiplicative
structure, it has ideals, and since the unimodular polynomials, i.e.,
$\mathbb{T}\rightarrow\mathbb{T}$, are monomials, we may reduce the
consideration to the ideals $z^{k}\mathcal{P}(\mathbb{T},\mathrm{U}%
_{N}(\mathbb{C}))$, $k=0,1,2,\dots.$

In view of the examples, we specialize the discussion to the case $N=2$, but
the arguments work generally.

In this section, we explain how the subspace
$\mathcal{K}:=\left\langle \{0,-1,\dots,-(2g-1)\}\right\rangle$
in (\ref{eqMin.21}) is reduced first to the smaller one
$\left\langle \{-1,-2,\dots,-(2g-2)\}\right\rangle $,
and then further to $\left\langle \{-2,\dots,-(2g-3)\}\right\rangle $,
in the case $N=2$. Returning to the semigroup $\mathcal{P}(\mathbb{T}%
,\mathrm{U}_{2}(\mathbb{C}))$, we note that it has a natural filtration of
ideals:
\begin{equation}
z\mathcal{P}(\mathbb{T},\mathrm{U}_{2}(\mathbb{C}))\supset z^{2}%
\mathcal{P}(\mathbb{T},\mathrm{U}_{2}(\mathbb{C}))\supset\cdots. \label{eq6.2}%
\end{equation}
A loop $A(z)$ is in $z^{k}\mathcal{P}(\mathbb{T},\mathrm{U}_{2}(\mathbb{C}))$
if and only if there is some $B(z)$ $\in\mathcal{P}(\mathbb{T},\mathrm{U}%
_{2}(\mathbb{C}))$ such that
\begin{equation}
A(z)=z^{k}B(z)\text{,\qquad}z\in\mathbb{T}\text{.} \label{eq6.3}%
\end{equation}
Since
\begin{equation}
m_{i}^{(B)}(z)=\sum_{j}B_{i,j}(z^{2})z^{j}, \label{eq6.4}%
\end{equation}
we get
$m_{i}^{(A)}(z)=z^{2k}m_{i}^{(B)}(z)$,
and for the representations
\begin{equation}
T_{i}^{(A)}=M_{z^{2k}}T_{i}^{(B)} \label{eq6.6}%
\end{equation}
where $M_{z^{2k}}$ denotes multiplication by $z^{2k}$ on the Hilbert space
$L^{2}(\mathbb{T})$. Despite this simple relationship between $T^{(A)}$ and
$T^{(B)}$, the irreducibility question can come out differently from one to
the other. \begin{figure}[tbp]%
\[
\setlength{\unitlength}{30pt}\setlength{\bracelength}{1.9\unitlength}\matpic
\]
\caption{$B\left(  z\right)  \in\mathcal{P}_{g}\left(  \mathbb{T}%
,\mathrm{U}_{2}\left(  \mathbb{C}\right)  \right)  $ vs.\ $m_{i}^{\left(
A\right)  }\left(  z\right)  =zm_{i}^{\left(  B\right)  }\left(  z\right)
\in\mathcal{P}_{g+1}\left(  \mathbb{T},\mathrm{U}_{2}\left(  \mathbb{C}%
\right)  \right)  $}%
\label{fig3}%
\end{figure}

If $A(z)=zB(z)$, and $B$ is of genus $g$, then $A$ is of genus $g+1$, but it
has vanishing first and last columns in its representation, as is clear from
Figure \ref{fig3}. Specifically, suppose $m_{i}^{\left(  A\right)  }\left(
z\right)  =zm_{i}^{\left(  B\right)  }\left(  z\right)  $ for all $i$; then we
have the following system of identities:%
\begin{equation}
A_{i,0}^{(0)}\equiv0\text{, }A_{i,1}^{(0)}=B_{i,0}^{(0)}\text{, }A_{i,0}%
^{(1)}=B_{i,1}^{(0)}\text{, }A_{i,1}^{(1)}=B_{i,0}^{(1)},\dots,A_{i,0}%
^{(g)}=B_{i,1}^{(g-1)}\text{, }A_{i,1}^{(g)}\equiv0 \label{eq6.7}%
\end{equation}
for $i=0,1$, and so the matrix
\[
\left(
\begin{array}
[c]{ll}%
R(1,1)_{0,0} & R(0,1)_{0,1}\\
R(1,0)_{1,0} & R(0,0)_{1,1}%
\end{array}
\right)
\]
of Section \ref{Irr} takes the following form:
\begin{equation}
\left(
\begin{array}
[c]{ll}%
R_{A}(1,1)_{0,0} & R_{A}(0,1)_{0,1}\\
R_{A}(1,0)_{1,0} & R_{A}(0,0)_{1,1}%
\end{array}
\right)  =\left(
\begin{array}
[c]{ll}%
R_{B}(0,0)_{1,1} & R_{B}(0,0)_{1,0}\\
R_{B}(0,0)_{0,1} & R_{B}(0,0)_{0,0}%
\end{array}
\right)  . \label{eq6.8}%
\end{equation}
Moreover, a given $A\in\mathcal{P}_{g+1}(\mathbb{T},\mathrm{U}_{2}%
(\mathbb{C}))$ has the form $m_{i}^{\left(  A\right)  }\left(  z\right)
=zm_{i}^{\left(  B\right)  }\left(  z\right)  $ for some $B\in\mathcal{P}%
(\mathbb{T},\mathrm{U}_{2}(\mathbb{C}))$ if and only if
\begin{equation}
\lambda_{0}(A)(:=R_{A}(0,0)_{0,0})=0. \label{eq6.9}%
\end{equation}

Putting this together, we get the following result:

\begin{proposition}
\label{proposition6.1}

\raggedright
\begin{enumerate}
\item[(a)] Let
\[
A\in\mathcal{P}_{g+1}(\mathbb{T},\mathrm{U}_{2}(\mathbb{C})),
\]
and let $P$ be the projection onto the subspace $\mathcal{K}$. Then the
following three conditions, $\mathrm{(i)}$, $\mathrm{(ii)}$, and
\textup{(iii)}, are equivalent:

\begin{enumerate}
\item[(i)] $\lambda_{0}(A)=0\mathpunct{;}$

\item[(ii)] $m_{i}^{\left(  A\right)  }\left(  z\right)  =zm_{i}^{\left(
B\right)  }\left(  z\right)  \mathpunct{;}$

\item[\textup{(iii)}] $e_{0}\in\left[  \mathcal{O}_{2}\left\langle \left\{
-1,-2,\dots,-2\left(  g-1\right)  \right\}  \right\rangle \right]  $.
\end{enumerate}

\item[(b)] The following two conditions are equivalent:

\begin{enumerate}
\item[(i)] $\lambda_{0}(A)=1\mathpunct{;}$

\item[(ii)] there is some
\[
V\in\mathrm{U}_{2}(\mathbb{C})\text{,\qquad}b\in\mathbb{T},
\]
such that
\begin{equation}
A(z)=V\left(
\begin{array}
[c]{lc}%
1 & 0\\
0 & bz^{g}%
\end{array}
\right)
\end{equation}
for all $z\in\mathbb{T}$.
\end{enumerate}

\item[(c)] The following two conditions are equivalent:

\begin{enumerate}
\item[(i)] $\lambda_{0}(A)<1\mathpunct{;}$

\item[(ii)]
\begin{align}
e_{0}  &  \in P\left[  \mathcal{O}_{2}\left\langle \{-1,-2,\dots
,-2g\}\right\rangle \right]  ,\qquad\text{and}\label{eq6.11}\\
e_{-(2g+1)}  &  \in P\left[  \mathcal{O}_{2}\left\langle \{-1,-2,\dots
,-2g\}\right\rangle \right]  .\nonumber
\end{align}
\end{enumerate}

\item[(d)] Suppose $\lambda_{0}(A)=0$. Then the following three conditions are equivalent:

\begin{enumerate}
\item[(i)] $1$ is \emph{not} in the spectrum of the matrix $(\ref{eq6.8}%
)\mathpunct{;}$

\item[(ii)]
\begin{align}
e_{0},e_{-1}  &  \in P\left[  \mathcal{O}_{2}\left\langle \{-2,\dots
,-2g+1\}\right\rangle \right]  ,\qquad\text{and}\label{eq6.12}\\
e_{-2g-1},e_{-2g}  &  \in P\left[  \mathcal{O}_{2}\left\langle \{-2,\dots
,-2g+1\}\right\rangle \right]  \mathpunct{;}\nonumber
\end{align}

\item[(iii)] the loop $B$ in $m_{i}^{\left(  A\right)  }\left(  z\right)
=zm_{i}^{\left(  B\right)  }\left(  z\right)  $ has the two vectors $\left(
B_{i,0}^{(0)}\right)  _{i}$ and $\left(  B_{i,1}^{(0)}\right)  _{i}$ linearly
independent in $\mathbb{C}^{2}$. \textup{(}Hence, given the factorization
$m_{i}^{\left(  A\right)  }\left(  z\right)  =zm_{i}^{\left(  B\right)
}\left(  z\right)  $, cyclicity of the reduced subspace
\[
\left\langle \left\{  -2,-3,\dots,-2g+1\right\}  \right\rangle ,
\]
i.e.,
\[
\operatorname*{span}\left\{  z^{-k};2\leq k\leq2g-1\right\}  ,
\]
holds for a generic subfamily $\{B\}$ in $\mathcal{P}_{g}(\mathbb{T}%
,\mathrm{U}_{2}(\mathbb{C}))$.\textup{)}
\end{enumerate}
\end{enumerate}
\end{proposition}

\begin{proof}
(a), (i) $\Rightarrow$ (ii): If $\lambda_{0}(A)=0$, then $A_{i,0}^{(0)}%
\equiv0$, and therefore%
\begin{equation}
A_{0,1}^{(g)}=\overline{A_{1,0}^{(0)}},\qquad A_{1,1}^{(g)}=-\overline
{A_{0,0}^{(0)}}, \label{eq6.13}%
\end{equation}
i.e., $A_{i,1}^{(g)}\equiv0$. This means that the coefficient matrices in the
expansion
\[
A(z)=A^{(0)}+A^{(1)}z+\cdots+A^{(g)}z^{g}%
\]
satisfy the conditions in Figure \ref{fig3}; and, if we define matrices
$B^{(0)},B^{(1)},\dots,B^{(g-1)}$ by (\ref{eq6.7}) above, then it follows that
$A(z)=zB(z)$ where
\begin{equation}
B(z)=B^{(0)}+B^{(1)}z+\cdots+B^{(g-1)}z^{g-1}. \label{eq6.14}%
\end{equation}
Hence (ii) holds.

(ii) $\Rightarrow$ (i): This is clear from reading (\ref{eq6.7}) in reverse.

The equivalence (i) $\Leftrightarrow$ (iii) follows from the observation that
the following sum representation
\[
e_{0}=\sum_{i_{1},\dots,i_{n}}T_{i_{1}}T_{i_{2}}\cdots T_{i_{n}}l_{i_{1}%
,\dots,i_{n}}%
\]
holds for some $n$ and $l_{i_{1},\dots,i_{n}}\in\left\langle \left\{
-1,-2,\dots\right\}  \right\rangle $ if and only if%
\[
T_{i_{n}}^{\ast}\cdots T_{i_{2}}^{\ast}T_{i_{1}}^{\ast}e_{0}\in\left\langle
\left\{  -1,-2,\dots\right\}  \right\rangle \text{\qquad for all }i_{1}%
,\dots,i_{n}.
\]
The conclusion can therefore be read off from the following general fact:
\[
T_{i_{n}}^{\ast}\cdots T_{i_{2}}^{\ast}T_{i_{1}}^{\ast}e_{0}\in\overline
{A_{i_{1},0}^{\left(  0\right)  }}\cdots\overline{A_{i_{n},0}^{\left(
0\right)  }}e_{0}+\left\langle \left\{  -1,-2,\dots\right\}  \right\rangle .
\]

(b), (i) $\Leftrightarrow$ (ii): If $\lambda_{0}(A)=1$, then $A_{i,0}%
^{(k)}\equiv0$ for $k>0$, and conversely. This follows from the identity%
\begin{equation}
A^{(0)\ast}A^{(0)}+\cdots+A^{(g)\ast}A^{(g)}=\left(
\begin{array}
[c]{ll}%
1 & 0\\
0 & 1
\end{array}
\right)  \text{,} \label{eq6.15}%
\end{equation}
which is part of the defining axiom system for $A$. Hence, the result follows
from \cite[Theorem 6.2]{BrJo00}, once we note that the only polynomials $b(z)$
such that $\left|  b(z)\right|  =1$ for all $z\in\mathbb{T}$ are the
monomials; see also \cite[Lemma 3.1]{BrJo00}.

(c): We already showed in Section \ref{Irr} that (i) \emph{implies} the first
of the conditions in (ii). The second one then follows from (\ref{eq6.13}),
i.e., the second line in (\ref{eq6.12}) follows from the first one.

(ii) $\Rightarrow$ (i): This follows from (b) above. For if $\lambda_{0}%
(A)=1$, then it follows from (b) that
\begin{equation}
L^{2}(\mathbb{R})=\left[  \mathcal{O}_{2}(\mathbb{C}e_{0})\right]
\oplus\left[  \mathcal{O}_{2}\left\langle \{-1,-2,\dots\}\right\rangle
\right]  , \label{eq6.16}%
\end{equation}
and so $e_{0}$ is not in the subspace
\[
\left[  \mathcal{O}_{2}\left\langle \{-1,-2,\dots\}\right\rangle \right]  ;
\]
and the same argument, based on (\ref{eq6.13}), shows that $e_{-(2g+1)}$ is
not in
\[
\left[  \mathcal{O}_{2}\left\langle \{\dots,-2g+1,-2g\}\right\rangle \right]
,
\]
concluding the proof of (c).

(d): We already saw that if $\lambda_{0}(A)=0$, then the conditions in
(\ref{eq6.12}) hold if and only if $1$ is not in the spectrum of the matrix
from (\ref{eq6.8}). Having now $m_{i}^{\left(  A\right)  }\left(  z\right)
=zm_{i}^{\left(  B\right)  }\left(  z\right)  $ from (a) above, we can use the
identity (\ref{eq6.8}) relating the $R_{A}(\,\cdot\,,\,\cdot\,)$-numbers
to the $R_{B}(\,\cdot\,,\,\cdot\,)$-numbers. But $1$ is in the spectrum of the
matrix
\[
\left(
\begin{array}
[c]{ll}%
d_{1} & c\\
\bar{c} & d_{0}%
\end{array}
\right)
\]
if and only if
\begin{equation}
(1-d_{0})(1-d_{1})=\left|  c\right|  ^{2}. \label{eq6.17}%
\end{equation}
The matrix on the right-hand side in (\ref{eq6.8}) is of this form, and
\begin{equation}
\left|  c\right|  ^{2}\leq d_{0}d_{1} \label{eq6.18}%
\end{equation}
by Schwarz's inequality. Here
\begin{equation}
c    :=\sum_{i}\overline{B_{i,1}^{(0)}}B_{i,0}^{(0)},\qquad
d_{0}    :=\sum_{i}\left|  B_{i,0}^{(0)}\right|  ^{2}\leq1,\qquad
d_{1}    :=\sum_{i}\left|  B_{i,1}^{(0)}\right|  ^{2}\leq1.\label{eq6.19}
\end{equation}
But using (\ref{eq6.13}) and (\ref{eq6.15}), we also get $d_{0}+d_{1}\leq1$.
Now (\ref{eq6.17})--(\ref{eq6.18}) yield $1\leq d_{0}+d_{1}$, and therefore
$d_{0}+d_{1}=1$. Substituting this back into formula (\ref{eq6.17}) then
yields $d_{0}d_{1}=\left|  c\right|  ^{2}$, which amounts to ``equality'' in
Schwarz's inequality (\ref{eq6.18}); and so the corresponding vectors (d)(iii)
are proportional. We already noted the equivalence of (i) and (ii) in (d); and
we just established that the negation of (i) amounts to linear dependence of
the vectors in (d)(iii). So (i) is equivalent to the linear independence, as
claimed in (d)(iii). This completes the proof of the proposition.
\end{proof}

\begin{remark}
\label{RemFilNew.2}Let loops $A$ and $B$ be as in Proposition
\textup{\ref{proposition6.1}(a),} see also Figure \textup{\ref{fig3},} and let
$T^{\left(  A\right)  }$ and $T^{\left(  B\right)  }$ be the corresponding
wavelet representations. Then, as a result of the theorems in Sections
\textup{\ref{Irr}} and \textup{\ref{Exp},} we conclude that $T^{\left(
A\right)  }$ is irreducible if and only if $T^{\left(  B\right)  }$ is. Since
the factorization in Proposition \textup{\ref{proposition6.1}(a)} corresponds
to $\lambda_{0}\left(  A\right)  =0$, we conclude that the general
irreducibility question has therefore been reduced to the case $\lambda
_{0}\left(  A\right)  >0$, which is the subject of the next section.
\end{remark}

\section{\label{Exp}An explicit formula for the minimal subspace}

Given $N=2$, and $A\in\mathcal{P}_{g}\left(  \mathbb{T}\text{,}\mathrm{U}%
_{2}\left(  \mathbb{C}\right)  \right)  $, we considered the wavelet
representation $T^{\left(  A\right)  }$ of $\mathcal{O}_{2}$ on $L^{2}\left(
\mathbb{T}\right)  \cong\ell^{2}$. We showed that
\begin{equation}
\mathcal{K}:=\left\langle z^{0},z^{-1},\dots,z^{-\left(  2g-1\right)
}\right\rangle \label{eqExp.1}%
\end{equation}
is $T^{\left(  A\right)  \,\ast}$-invariant and cyclic (in $L^{2}\left(
\mathbb{T}\right)  $) for the representation $T^{\left(  A\right)  }$. But we
also showed that the first one of the basis vectors, $z^{0}$, \emph{is} in the
cyclic space generated by $z^{-1},z^{-2},\dots,z^{-2\left(  g-1\right)  }$ and
the representation, if and only if $A_{i,0}^{\left(  0\right)  }\equiv0$ for
all $i$. Specifically, setting%
\begin{equation}
\lambda_{0}\left(  A\right)  :=\sum_{i}\left|  A_{i,0}^{\left(  0\right)
}\right|  ^{2}, \label{eqExp.2}%
\end{equation}
we showed in Corollary \ref{corollary5.3} and Proposition \ref{proposition6.1}%
(a) that%
\begin{equation}
e_{0}\in\left[  \mathcal{O}_{2}\left\langle e_{-1},e_{-2},\dots,e_{-2\left(
g-1\right)  }\right\rangle \right]  \label{eqExp.3}%
\end{equation}
if and only if $\lambda_{0}\left(  A\right)  =0$. Hence, it follows that
$\mathcal{K}$ is \emph{not minimal} (in the sense of the following definition)
if $\lambda_{0}\left(  A\right)  =0$. In other words, if $\lambda_{0}\left(
A\right)  =0$, $\mathcal{K}$ then contains a strictly smaller subspace which
is both $T^{\left(  A\right)  \,\ast}$-invariant and cyclic. In this section,
we show the converse implication. But first a definition:

\begin{definition}
\label{DefExp.1}We say that a subspace $\mathcal{L}\subset\mathcal{K}$ is
\emph{minimal} if it is $T^{\left(  A\right)  \,\ast}$-invariant,
cyclic
for the representation $T^{\left(  A\right)  }$, and minimal with respect
to the two properties, i.e., it does not contain a proper subspace which is
also $T^{\left(  A\right)  \,\ast}$-invariant and cyclic.
\end{definition}

We will now prove the converse to the above-mentioned result, showing, in
particular, that if $\lambda_{0}\left(  A\right)  >0$, then $\mathcal{K}$ is
\emph{generically} minimal; see Corollary \ref{CorExp.6}.

\begin{theorem}
\label{ThmExp.2}Let $A\in\mathcal{P}_{g}\left(  \mathbb{T},\mathrm{U}%
_{2}\left(  \mathbb{C}\right)  \right)  $, and let $T^{\left(  A\right)  }$ be
its wavelet representation. Let $\mathcal{K}=\left\langle
z^{0},z^{-1},\dots,z^{-\left(  2g-1\right)  }\right\rangle $ as in
\textup{(\ref{eqExp.1}),} and assume $\lambda_{0}\left(  A\right)  >0$. Then
$\mathcal{K}$ contains a unique minimal subspace $\mathcal{L}$, i.e.,
$\mathcal{L}$ is $T^{\left(  A\right)  \,\ast}$-invariant, cyclic, and
minimal. It is spanned by the complex conjugates of the following family of
$4g$ functions:%
\begin{equation}
A_{i,j}\left(  z\right)  z^{k+j}\text{,\qquad where }i,j\in\left\{
0,1\right\}  \text{ and }k\in\left\{  0,1,\dots,g-1\right\}  .
\label{eqExpNew.pound}%
\end{equation}
\textup{(}Here $k$ varies
independently of \emph{both} $i$ and $j$.\textup{)} Moreover, within the class%
\begin{equation}
\mathcal{P}_{g}\left(  \mathbb{T},\mathrm{U}_{2}\left(  \mathbb{C}\right)
\right)  ,\qquad\lambda_{0}\left(  A\right)  >0, \label{eqExpNew.poundbis}%
\end{equation}
the dimension of $\mathcal{L}$ is $2g$, for a generic subfamily, and so, for
this subfamily, $\mathcal{L}=\mathcal{K}$, and $\mathcal{K}$ itself is minimal.
\end{theorem}

\begin{remark}
\label{remark8new}An immediate consequence of the definition of the
subspace\linebreak $\mathcal{L}\subset\mathcal{K}\subset L^{2}(\mathbb{T})$ is
the following formula for the ``deficiency space'':
\begin{equation}
\mathcal{K}\ominus\mathcal{L}=\bigwedge_{i}\ker(PT_{i}P)
\label{eqExpNew.poundter}%
\end{equation}
where $P$ denotes the projection onto $\mathcal{K}$. We will show below, using
\textup{(\ref{eqExpNew.poundter})} and Corollary \textup{\ref{corollary5new},}
that $\mathcal{L}=\mathcal{K}$ if $\lambda_{0}\left(  A\right)  >0$. This
means that $\mathcal{K}$ itself is then the unique minimal subspace when the
loop $A$ does not have $m_{i}^{\left(  A\right)  }\left(  z\right)  \in
z\mathbb{C}\left[  z\right]  $, as in Section \textup{\ref{Fil}.} But before
arriving at the conclusion, we must first derive several \emph{a priori}
properties of $\mathcal{L}$.
\end{remark}

\begin{proof}
[Proof of Theorem \textup{\ref{ThmExp.2}}]The details are somewhat technical,
and it seems more practical to first do them for the special case when $g=3$,
and then comment at the end on the (relatively minor) modifications needed in
the proof for the case when $g$ is arbitrary $g\geq2$.

Using the terminology of (\ref{eqBrJo00Rep.18}), we then get
\begin{equation}
A\left(  z\right)  =A^{\left(  0\right)  }+A^{\left(  1\right)  }z+A^{\left(
2\right)  }z^{2}, \label{eqExp.4}%
\end{equation}
where $A^{(2)}\neq0$ and $A^{\left(  0\right)  }$, $A^{\left(
1\right)  }$, and $A^{\left(  2\right)  }$ are $2$-by-$2$ complex matrices
satisfying%
\begin{equation}
\sum_{k=0}^{2}A^{\left(  k\right)  \,\ast}A^{\left(  k+l\right)  }%
=\delta_{0,l}\openone_{2}. \label{eqExp.5}%
\end{equation}
When $A$ is given, we denote that corresponding wavelet representation by
$T^{\left(  A\right)  }$, or just $T$ for simplicity. Recall%
\begin{equation}
\left(  T_{i}f\right)  \left(  z\right)  =\sum_{j}A_{i,j}\left(  z^{2}\right)
z^{j}f\left(  z^{2}\right)  , \label{eqExp.6}%
\end{equation}
or simply%
\begin{equation}
T_{i}f\left(  z\right)  =m_{i}^{\left(  A\right)  }\left(  z\right)  f\left(
z^{2}\right)  , \label{eqExp.7}%
\end{equation}
where%
\begin{equation}
m_{i}^{\left(  A\right)  }\left(  z\right)  =\sum_{j}A_{i,j}\left(
z^{2}\right)  z^{j}. \label{eqExp.8}%
\end{equation}
As we saw in (\ref{eqExp.1}), the subspace%
\begin{equation}
\mathcal{K}=\left\langle z^{0},z^{-1},z^{-2},z^{-3},z^{-4},z^{-5}\right\rangle
\label{eqExp.9}%
\end{equation}
is then $T^{\ast}$-invariant, and also cyclic for the representation. But the
issue is when $\mathcal{K}$ is \emph{minimal} with respect to these two
properties. The minimality of some subspace $\mathcal{L}\subset\mathcal{K}$
then means that $\mathcal{L}$ is $T^{\ast}$-invariant, and cyclic, and that no
proper $T^{\ast}$-invariant subspace of $\mathcal{L}$ is cyclic.

In working out details on $\mathcal{L}$, we use (\ref{eqExp.6}%
)--(\ref{eqExp.8}) in conjunction with (\ref{eqMin.8})--(\ref{eqMin.9}), and
it is more helpful to work with the complex conjugates%
\begin{equation}
\mathcal{M}:=\overline{\mathcal{L}}=\left\{  \overline{f\left(  x\right)
};f\in\mathcal{L}\right\}  , \label{eqExp.10}%
\end{equation}
and so $\mathcal{M}$ consists of polynomials of degree at most $5$. It follows
from (\ref{eqExp.6})--(\ref{eqExp.8}) that $\mathcal{M}$ is then spanned by
the functions (polynomials) in the following list:
\begin{equation}
A_{i,j}\left(  z\right)  z^{k+j},\qquad i=0,1,\;j=0,1,\;k=0,1,2.
\label{eqExp.11}%
\end{equation}
Hence, by (\ref{eqExp.10}), $\mathcal{L}$ consists of the space spanned by the
complex conjugates of these functions.\renewcommand{\qed}{} By (\ref{eqExp.9})
it is clear that $\mathcal{L}\subset\mathcal{K}$.
\end{proof}

The proof of Theorem \ref{ThmExp.2} will now be split up into several lemmas:

\begin{lemma}
\label{LemExp.3}The space $\mathcal{L}$ is $T^{\ast}$-invariant.
\end{lemma}

\begin{proof}
Now, the functions $A_{i,j}\left(  z\right)  $ in (\ref{eqExp.11})
are the matrix elements of the loop $A\left(  z\right)
$, and so it follows from (\ref{eqExp.4}) that each of them is a polynomial of
degree at most $2$, say%
\begin{equation}
a\left(  z\right)  =c_{0}+c_{1}z+c_{2}z^{2} \label{eqExp.12}%
\end{equation}
(since $A\left(  z\right)  $ itself has degree $2$ when $g=3$). Hence,%
\begin{equation}
T_{i}^{\ast}\left(  \bar{a}\right)  =\bar{c}_{0}\overline{A_{i,0}\left(
z\right)  }+\bar{c}_{1}\overline{A_{i,1}\left(  z\right)  }z^{-1}+\bar{c}%
_{2}\overline{A_{i,0}\left(  z\right)  }z^{-1}, \label{eqExp.13}%
\end{equation}
or equivalently,%
\begin{equation}
\overline{T_{i}^{\ast}\left(  \bar{a}\right)  }=c_{0}A_{i,0}\left(  z\right)
+c_{1}A_{i,1}\left(  z\right)  z^{1}+c_{2}A_{i,0}\left(  z\right)  z^{1}.
\label{eqExp.14}%
\end{equation}
Now set $b\left(  z\right)  :=z^{2p}a\left(  z\right)  $ where $a$ is as in
(\ref{eqExp.12}). From (\ref{eqExp.8}), we then get%
\[
T_{i}^{\ast}\left(  \bar{b}\right)  =z^{-p}T_{i}^{\ast}\left(  \bar{a}\right)
,
\]
or equivalently,%
\begin{equation}
\overline{T_{i}^{\ast}\left(  \bar{b}\right)  }=z^{p}\left(  c_{0}%
A_{i,0}\left(  z\right)  +c_{1}A_{i,1}\left(  z\right)  z+c_{2}A_{i,0}\left(
z\right)  z\right)  , \label{eqExp.15}%
\end{equation}
using (\ref{eqExp.14}). In view of (\ref{eqExp.11}), we
need then only to compute the following:%
\[
T_{i}^{\ast}\left(  \overline{za\left(  z\right)  }\right)  =\bar{c}%
_{0}\overline{A_{i,1}\left(  z\right)  }z^{-1}+\bar{c}_{1}\overline
{A_{i,0}\left(  z\right)  }z^{-1}+\bar{c}_{2}\overline{A_{i,1}\left(
z\right)  }z^{-2},
\]
or equivalently,%
\[
\overline{T_{i}^{\ast}\left(  \overline{za\left(  z\right)  }\right)  }%
=c_{0}A_{i,1}\left(  z\right)  z+c_{1}A_{i,0}\left(  z\right)  z+c_{2}%
A_{i,1}\left(  z\right)  z^{2}.
\]
Now putting the formulas together, we get the value of $T_{i}^{\ast}$ on each
of the functions (\ref{eqExp.11}) which go into the
definition of $\mathcal{L}$, and the conclusion of the lemma follows.
\end{proof}

\begin{lemma}
\label{LemExp.4}The space $\mathcal{L}$ is cyclic.
\end{lemma}

\begin{proof}
Since both $\mathcal{L}$ and $\mathcal{K}$ are $T^{\ast}$-invariant, the
conclusion will follow if we check the inclusion
\begin{equation}
\mathcal{K}\subset\operatorname*{span}_{i}\left(  T_{i}\mathcal{L}\right)  .
\label{eqExp.16}%
\end{equation}
For the space on the right-hand side in (\ref{eqExp.16}), we shall use the
terminology$\left[  \mathcal{O}_{2}^{1}\mathcal{L}\right]  $, and similarly,
the space spanned by all the spaces%
\begin{equation}
T_{i_{1}}T_{i_{2}}\cdots T_{i_{n}}\mathcal{L} \label{eqExp.17}%
\end{equation}
will be denoted $\left[  \mathcal{O}_{2}^{n}\mathcal{L}\right]  $. In
(\ref{eqExp.17}), we vary the multi-index $\left(  i_{1},i_{2},\dots
,i_{n}\right)  $ over all the $2^{n}$ possibilities. It follows from the
$T^{\ast}$-invariance of $\mathcal{K}$ (in (\ref{eqExp.1})) and $\mathcal{L}$
(in (\ref{eqExpNew.pound})) that we get different families of nested
finite-dimensional subspaces:%
\begin{equation}
\mathcal{K}\subset\left[  \mathcal{O}_{2}^{1}\mathcal{K}\right]
\subset\left[  \mathcal{O}_{2}^{2}\mathcal{K}\right]  \subset\left[
\mathcal{O}_{2}^{3}\mathcal{K}\right]  \subset\dots\subset\left[
\mathcal{O}_{2}^{n}\mathcal{K}\right]  \subset\left[  \mathcal{O}_{2}%
^{n+1}\mathcal{K}\right]  \subset\cdots, \label{eqExp.18}%
\end{equation}
and a similar sequence for $\mathcal{L}$. Since $\mathcal{K}$ is cyclic, we
have
\begin{equation}
\bigvee_{n}\left[  \mathcal{O}_{2}^{n}\mathcal{K}\right]  =L^{2}\left(
\mathbb{T}\right)  \qquad(\cong\ell^{2}). \label{eqExp.19}%
\end{equation}
But $\mathcal{K}\subset\left[  \mathcal{O}_{2}^{n}\mathcal{L}\right]  $, for
some $n$, so $\mathcal{L}$ is also cyclic. The conclusion of the lemma follows
from (\ref{eqExp.18}) and (\ref{eqExp.19}), once we check that
\begin{equation}
\mathcal{K}\subset\left[  \mathcal{O}_{2}^{1}\mathcal{L}\right]  ,
\label{eqExp.20}%
\end{equation}
and so $n=1$ works, and%
$\left[  \mathcal{O}_{2}^{p}\mathcal{K}\right]  \subset\left[  \mathcal{O}
_{2}^{p+1}\mathcal{L}\right]  $ for all $p$.

Turning now to the details: Since $\mathcal{K}$ is spanned by $z^{-k}$,
$k=0,1,\dots,5$, we must check that each of these basis functions has the
representation%
\begin{equation}
z^{-k}=\sum_{i}T_{i}l_{i} \label{eqExp.21}%
\end{equation}
for $l_{0},l_{1}\in\mathcal{L}$, where we refer to (\ref{eqExp.11})
(see also (\ref{eqExpNew.pound})) for the
characterization of the space $\mathcal{L}$, or rather $\mathcal{M}%
:=\overline{\mathcal{L}}$.

But (\ref{eqExp.21}) is equivalent to the assertion that%
\begin{equation}
T_{i}^{\ast}\left(  z^{-k}\right)  \in\mathcal{L} \label{eqExp.22}%
\end{equation}
for all $i=0,1$ and $k=0,1,\dots,5$; and (\ref{eqExp.22}) can be checked by a
direct calculation, which is very similar to the one going into the proof of
Lemma \ref{LemExp.3}. Specifically, using (\ref{eqMin.8})--(\ref{eqMin.9}) we
get the following:%
\begin{align}
T_{i}^{\ast}\left(  z^{0}\right)   &  =\overline{A_{i,0}\left(  z\right)  }%
\in\mathcal{L},\label{eqExp.23a}\\
T_{i}^{\ast}\left(  z^{-1}\right)   &  =\overline{A_{i,1}\left(  z\right)
}z^{-1}\in\mathcal{L},\label{eqExp.23b}\\
T_{i}^{\ast}\left(  z^{-2}\right)   &  =\overline{A_{i,0}\left(  z\right)
}z^{-1}\in\mathcal{L},\label{eqExp.23c}\\
T_{i}^{\ast}\left(  z^{-3}\right)   &  =\overline{A_{i,1}\left(  z\right)
}z^{-2}\in\mathcal{L},\label{eqExp.23d}\\
T_{i}^{\ast}\left(  z^{-4}\right)   &  =\overline{A_{i,0}\left(  z\right)
}z^{-2}\in\mathcal{L},\label{eqExp.23e}\\%
\intertext{and finally}%
T_{i}^{\ast}\left(  z^{-5}\right)   &  =\overline{A_{i,1}\left(  z\right)
}z^{-3}\in\mathcal{L}. \label{eqExp.23f}%
\end{align}
Recall that the complex conjugates of the functions on the right-hand side in
this list are precisely the ones from (\ref{eqExp.11}), or
equivalently, (\ref{eqExpNew.pound}). This proves (\ref{eqExp.20}), and
therefore the cyclicity of $\mathcal{L}$, which was claimed in the lemma. As a
bonus, we get from (\ref{eqExp.23a})--(\ref{eqExp.23f}) that the inclusion
$\mathcal{L}\subset\left\langle \{-1,-2,-3,-4\}\right\rangle $ holds if and
only if $\lambda_{0}(A)=0$. To see this, use the fact (for $g=3$) that
\[
A_{i,1}^{(2)}=(-1)^{i}\overline{A_{1-i,0}^{(0)}}.
\]
\end{proof}

\begin{lemma}
\label{LemExp.5}The space $\mathcal{L}$ is minimal in the sense of Definition
\textup{\ref{DefExp.1}.}
\end{lemma}

\begin{proof}
We will establish the conclusion by proving that if $\mathcal{L}_{1}$ is any
subspace of $\mathcal{K}$ which is both $T^{\ast}$-invariant and cyclic, then
$\mathcal{L}\subset\mathcal{L}_{1}$. So in particular, $\mathcal{L}$ does not
contain a \emph{proper} subspace which is both $T^{\ast}$-invariant and cyclic.

Now suppose that some space $\mathcal{L}_{1}$ has the stated properties. Since
it is cyclic, we must have%
\begin{equation}
\mathcal{K}\subset\left[  \mathcal{O}_{2}^{n}\mathcal{L}_{1}^{{}}\right]
\label{eqExp.24}%
\end{equation}
satisfied for \emph{some} $n\in\mathbb{N}$. As noted in the proof of Lemma
\ref{LemExp.4}, this is equivalent to%
\begin{equation}
T_{i_{n}}^{\ast}\cdots T_{i_{2}}^{\ast}T_{i_{1}}^{\ast}\left(  z^{-k}\right)
\in\mathcal{L}_{1} \label{eqExp.25}%
\end{equation}
for all $i_{1},\dots,i_{n}\in\left\{  0,1\right\}  $, and all $k\in\left\{
0,1,\dots,5\right\}  $. But we also saw in the proof of Lemma \ref{LemExp.3}
that the functions on the left-hand side in (\ref{eqExp.25}) are precisely
those which are listed in (\ref{eqExp.11}). Note that the
functions in (\ref{eqExp.11}), or (\ref{eqExpNew.pound}),
are those given by%
\begin{equation}
T_{i}^{\ast}\left(  z^{-k}\right)  ,\qquad i=0,1,\qquad k=0,1,\dots,5.
\label{eqExp.26}%
\end{equation}
But $\lambda_{0}\left(  A\right)  >0$ by assumption, so for some $i$, we have
$A_{i,0}^{\left(  0\right)  }\neq0$, and the calculation in the proof of Lemma
\ref{LemExp.3}, and in the previous two sections, then shows that the families
of functions in (\ref{eqExp.26}) and (\ref{eqExp.25}) are the same, i.e., we
get the same functions in (\ref{eqExp.25}) for $n>1$ as the ones which are
already obtained for $n=1$ in (\ref{eqExp.26}). This is the step which uses
the assumption $\lambda_{0}\left(  A\right)  >0$. Since $\mathcal{L}$ is
spanned by the vectors in (\ref{eqExp.26}), the desired inclusion
$\mathcal{L}\subset\mathcal{L}_{1}$ follows. More details are worked out in
Remark \ref{RemExpNew.1} below.
\end{proof}

\begin{proof}
[Proof of Theorem \textup{\ref{ThmExp.2}} concluded]The result in the theorem
is now immediate from the three lemmas, and we need only comment on the size
of the genus $g$. We argued the case $g=3$; but, for the general case,
$\mathcal{K}$ is spanned by $z^{-k}$, $k=0,1,\dots,2g-1$, and the functions
from the list (\ref{eqExp.11}), or equivalently
(\ref{eqExpNew.pound}), will then be
\begin{equation}
A_{i,j}\left(  z\right)  z^{k+j},\qquad i=0,1,\;j=0,1,\;k=0,1,\dots,g-1.
\label{eqExp.27}%
\end{equation}
Otherwise, all the arguments from the proofs of the lemmas carry over. See
Remark \ref{RemExpNew.1} for more details.
\end{proof}

\begin{corollary}
\label{CorExp.6}When $g$ is given, and $\lambda_{0}:=\sum_{i}\left|
A_{i,0}^{\left(  0\right)  }\right|  ^{2}=R\left(  0,0\right)  _{0,0}>0$, then
$\mathcal{L}=\mathcal{K}$ for a generic set of loops $A$ in $\mathcal{P}%
_{g}\left(  \mathbb{T},\mathrm{U}_{2}\left(  \mathbb{C}\right)  \right)  $.
\end{corollary}

\begin{proof}
The proof comes down to a dimension count. Since $\mathcal{K}=\left\langle
z^{0},\dots,z^{-\left(  2g-1\right)  }\right\rangle $ is of dimension $2g$, we
just need to check that the space $\mathcal{L}$ ($\subset\mathcal{K}$),
spanned by the $4g$ functions in (\ref{eqExp.27}), is of
dimension $2g$ for a generic set of loops $A$ in $\mathcal{P}_{g}\left(
\mathbb{T},\mathrm{U}_{2}\left(  \mathbb{C}\right)  \right)  $, and that can
be checked by a determinant argument based on the conditions for the matrices
$A^{\left(  0\right)  },A^{\left(  1\right)  },\dots,A^{\left(  g-1\right)  }$
defining $A\left(  z\right)  $; see (\ref{eqExp.4})--(\ref{eqExp.5}) above.

The above-mentioned dimension count is based on the following consideration
(which we only sketch in rough outline). A possible linear relation among the
functions from (\ref{eqExpNew.pound}) takes the form%
\begin{equation}
\sum_{i}\sum_{j}\sum_{k}C_{i,j,k}A_{i,j}\left(  z\right)  z^{j+k}\equiv0,
\label{eqExp.28}%
\end{equation}
where the $i,j$ summation indices are $0,1$, and the $k$ summation is over
$0,1,\dots,g-1$. As a result, we get the following system of relations:%
\begin{equation}
\sum_{i=0}^{1}C_{i,j,k}A_{i,j}\left(  z\right)  \equiv0\qquad\left(
\operatorname{mod}z^{g-j-k}\right)  \label{eqExp.29}%
\end{equation}
for all $j=0,1$, and all $k=0,1,\dots,g-1$. Note that (\ref{eqExp.29}) is a
matrix multiplication. Using finally
\begin{equation}
\sum_{i}\left|  A_{i,0}^{\left(  0\right)  }\right|  ^{2}>0,
\label{eqExp.30bis}%
\end{equation}
we see that the dimension of the space spanned by $\left\{  A_{i,j}\left(
z\right)  z^{j+k}\right\}  $ is $2g$, as claimed. See Remark \ref{RemExpNew.1}
and Observation \ref{ObsExpNew.2} for details.
\end{proof}

\begin{remark}
\label{RemExpNew.1}A more detailed study of the space $\mathcal{L}$ will be
postponed to a later paper, but one point is included here: The function
$z^{0}$ \textup{(}$=e_{0}\equiv1$\textup{)} is in $\mathcal{L}$ if and only if
the polynomials%
\begin{equation}
\left\{  A_{i,0}\left(  z\right)  ,A_{j,1}\left(  z\right)  z\right\}  _{i,j}
\label{eqRemExpNew.1.1}%
\end{equation}
do not have a common divisor. This follows from \textup{(\ref{eqExp.23a}%
)--(\ref{eqExp.23f}).} Indeed, for $e_{0}$ to be in $\mathcal{L}$, we must
have the existence of $h_{i,j}\left(  z\right)  \in\mathbb{C}\left[  z\right]
$ such that
\begin{equation}
1=\sum_{i}h_{i,0}\left(  z\right)  A_{i,0}\left(  z\right)  +\sum_{j}%
h_{j,1}\left(  z\right)  zA_{j,1}\left(  z\right)  . \label{eqRemExpNew.1.2}%
\end{equation}
But by algebra, this amounts to the assertion that the family of polynomials
listed in \textup{(\ref{eqRemExpNew.1.1})} is mutually prime within the ring
$\mathbb{C}\left[  z\right]  $. Also note that, by the result in Section
\textup{\ref{Fil},} monomials such as $d\left(  z\right)  =z$ are not common
divisors in the polynomials of \textup{(\ref{eqRemExpNew.1.1})} if
$d_{0}\left(  A\right)  >0$. In fact, a possible common divisor $d\left(
z\right)  \in\mathbb{C}\left[  z\right]  $ for \textup{(\ref{eqRemExpNew.1.1}%
)} yields the following factorization:%
\begin{equation}
A_{i,0}\left(  z\right)  =d\left(  z\right)  k_{i,0}\left(  z\right)  ,\qquad
A_{j,1}\left(  z\right)  z=d\left(  z\right)  k_{j,1}\left(  z\right)
\label{eqRemExpNew.1.3}%
\end{equation}
($k_{i,0}\left(  z\right)  ,k_{j,1}\left(  z\right)  \in\mathbb{C}\left[
z\right]  $). Hence:

\begin{observation}
\label{ObsExpNew.2}If $d_{0}\left(  A\right)  >0$, then $e_{0}\in\mathcal{L}$.
\end{observation}

\begin{proof}
For if not, the greatest common divisor $d\left(  z\right)  $ of the family
\textup{(\ref{eqRemExpNew.1.1})} would have a root $\gamma\in\mathbb{C}%
\setminus\left\{  0\right\}  $, i.e., $d\left(  \gamma\right)  =0$. By
\textup{(\ref{eqRemExpNew.1.3}),} we would then have%
\begin{equation}
A\left(  \gamma\right)  =0, \label{eqRemExpNew.1.4}%
\end{equation}
where $A\in\mathcal{P}\left(  \mathbb{T},\mathrm{U}_{2}\left(  \mathbb{C}%
\right)  \right)  $ is the originally given loop. Recall that, by
\textup{(\ref{eqBrJo00Rep.18}),} we may view $A\left(  z\right)  $ as an
entire analytic matrix function, i.e., an entire analytic function,
$\mathbb{C}\rightarrow M_{2}\left(  \mathbb{C}\right)  $, whose restriction to
$\mathbb{T}$ takes values in $\mathrm{U}_{2}\left(  \mathbb{C}\right)  $.
\textup{(}These are also called \emph{inner} matrix functions \cite{PrSe86}%
.\textup{)} But \textup{(\ref{eqRemExpNew.1.4})} is impossible \textup{(}for
$\gamma\neq0$\textup{)} in view of Lemma \textup{\ref{CorBrJo00Rep.2}} and its
corollary. We will give the details for $g=3$, but they apply with the obvious
modifications to the general case of $g\geq2$. If
\textup{(\ref{eqRemExpNew.1.4})} holds, then by Lemma
\textup{\ref{CorBrJo00Rep.2},}%
\begin{equation}
V^{-1}A\left(  \gamma\right)  =\left(  Q_{0}^{\perp}+\gamma Q_{0}^{{}}\right)
\left(  Q_{1}^{\perp}+\gamma Q_{1}^{{}}\right)  =0, \label{eqObsExpNew.2.5}%
\end{equation}
where we use the projections $Q_{0},Q_{1}$ in $\mathbb{C}^{2}$ from
\textup{(\ref{eqBrJo00Rep.22}).} Setting $A_{j}\left(  z\right)  =Q_{j}%
^{\perp}+zQ_{j}^{{}}$, $j=0,1$, \textup{(\ref{eqObsExpNew.2.5})} then yields
the following estimate:%
\begin{align}
0  &  =A_{1}\left(  \gamma\right)  ^{\ast}A_{0}\left(  \gamma\right)  ^{\ast
}A_{0}\left(  \gamma\right)  A_{1}\left(  \gamma\right)  =A_{1}\left(
\gamma\right)  ^{\ast}\left(  Q_{0}^{\perp}+\left|  \gamma\right|  ^{2}%
Q_{0}^{{}}\right)  A_{1}\left(  \gamma\right) \label{eqObsExpNew.2.5bis}\\
&  \geq\min\left(  1,\left|  \gamma\right|  ^{2}\right)  \cdot\left(
A_{1}\left(  \gamma\right)  ^{\ast}A_{1}\left(  \gamma\right)  \right)
  =\min\left(  1,\left|  \gamma\right|  ^{2}\right)  \cdot\left(
Q_{1}^{\perp}+\left|  \gamma\right|  ^{2}Q_{1}^{{}}\right) \nonumber\\
&  \geq\left(  \min\left(  1,\left|  \gamma\right|  ^{2}\right)  \right)
^{2}\cdot\openone_{2},\nonumber
\end{align}
where the order $\geq $ is that of positive operators
on $\mathbb{C}^{2}$. But this is impossible, since $\gamma\neq0$. The latter
is from the assumption $d_{0}\left(  A\right)  >0$.
\end{proof}
\end{remark}

This last argument in this proof also serves as a proof of
Corollary \ref{corollary5new}, and this corollary is again the basis for the
following stronger result, which we now sketch:

\begin{theorem}
\label{ThmExpNew.10}If $\lambda_{0}\left(  A\right)  >0$, it follows that
$\mathcal{L}=\mathcal{K}$.
\end{theorem}

\begin{proof}
By Corollary \ref{corollary5new} and (\ref{eq2.20}), we have the estimate%
\begin{equation}
\sum_{i=0}^{1}\left|  m_{i}\left(  z\right)  \right|  ^{2}\geq\left(
\min\left(  1,\left|  z\right|  ^{2}\right)  \right)  ^{g-1}\left(  1+\left|
z\right|  ^{2}\right)  \text{\qquad for all }z\in\mathbb{C}.
\label{eqThmExpNew.10proof.2}%
\end{equation}
If $k\in\mathcal{K}\ominus\mathcal{L}$, then by (\ref{eqExpNew.poundter}), we
get%
\begin{equation}
Pm_{i}\left(  z\right)  k\left(  z^{2}\right)  =0,\qquad i=0,1,\;z\in
\mathbb{T}. \label{eqThmExpNew.10proof.3}%
\end{equation}
If $\lambda_{0}\left(  A\right)  >0$, then the value $z=0$ is not a common
root of the two complex polynomials $m_{0}\left(  z\right)  $, $m_{1}\left(
z\right)  $, and so by (\ref{eqThmExpNew.10proof.2}) the two polynomials
$m_{0}$ and $m_{1}$ have no common roots at all, by Proposition
\ref{proposition6.1}(a). Therefore, when (\ref{eqThmExpNew.10proof.2}) and
(\ref{eqThmExpNew.10proof.3}) are combined, we get $k=0$. Hence $\mathcal{K}%
\ominus\mathcal{L}=0$, and the proof is completed.
\end{proof}

\begin{remark}
\label{RemExpNew.pound}Not everything that works out easily in the case $g=3$
generalizes immediately to $g>3$: The case $g=3$, and arbitrary $N$, amounts
to a choice of two projections, say $P$ and $Q$, in $\mathbb{C}^{N}$, and the
finite-dimensional representations of the algebra generated by two projections
are completely known by folklore \textup{(}see, e.g., \cite{JSW95}\textup{).}
In fact, it can be easily checked that this is the same as displaying the
finite-dimensional representations of the Clifford algebra with two
generators, $A$, $B$, say. The relations between $A$, $B$ are:%
\begin{equation}
A^{\ast}=A,\qquad B^{\ast}=B,\qquad A^{2}+B^{2}=\openone_{N}\text{,\quad
and\quad}AB+BA=0. \label{eqRemExpNew.pound.1}%
\end{equation}
For any such pair, set%
\begin{equation}
P=\frac{1}{2}\left(  \openone+A-B\right)  ,\qquad Q=\frac{1}{2}\left(
\openone-A-B\right)  . \label{eqRemExpNew.pound.2}%
\end{equation}
Then it is immediate that $P$ and $Q$ are projections, i.e., $P=P^{\ast}%
=P^{2}$, etc. Conversely, if $P$, $Q$ are any projections, set%
\begin{equation}
A=P-Q,\qquad B=\openone_{N}-P-Q, \label{eqRemExpNew.pound.3}%
\end{equation}
and an easy calculation shows that \textup{(\ref{eqRemExpNew.pound.1})} is
then satisfied. Since the finite-di\-men\-sion\-al representations of
\textup{(\ref{eqRemExpNew.pound.1}),} the Clifford algebra $\mathcal{C}_{2}$,
are known \cite{JSW95}, we then get a useful classification of $\mathcal{P}%
_{3}\left(  \mathbb{T},\mathrm{U}_{N}\left(  \mathbb{C}\right)  \right)  $.
But these comments do not carry over to the case $g>3$. \textup{(}A good
reference on the Clifford algebra and its representations is \cite[Ch.~1, \S
5]{LaMi89}.\textup{)}
\end{remark}

It is perhaps a little early to identify regions in the parameters $\theta$,
$\rho$ where the scaling function $x\mapsto\varphi_{\theta,\rho}^{{}}\left(
x\right)  $ is regular, and where it is not, but a primitive test would be the
vanishing-moment condition of Daubechies \cite[Ch.~6]{Dau92}. We would look
for values $\theta$, $\rho$ such that $m_{0}^{\left(  \theta,\rho\right)
}\left(  z\right)  $ is divisible by $\left(  1+z\right)  ^{2}$ and by
$\left(  1+z\right)  ^{3}$. These are the conditions which ensure that%
\begin{equation}
\left(  \frac{d\;}{d\xi}\right)  ^{k}\hat{\psi}\left(  \xi\right)  |_{\xi
=0}=0,\qquad k=0,1,\dots, \label{eqMom.null}%
\end{equation}
starting with
\begin{align}
0  &  =\hat{\psi}\left(  0\right)  =\int_{\mathbb{R}}\psi\left(  x\right)
\,dx,\qquad\dots,\label{eqMom.0}\\
0  &  =\int_{\mathbb{R}}x^{k}\psi\left(  x\right)  \,dx. \label{eqMom.1}%
\end{align}
Here (\ref{eqMom.0}) is automatic since $m_{0}\left(  -1\right)  =0$. Recall
the coordinates $z=e^{-i\xi}$, $\xi\in\mathbb{R}$. The second one
(\ref{eqMom.1}) corresponds to
$\smash[b]{\left(  \frac{d\;}{d\xi}\right)  ^{k}}\hat
{\psi}\left(  \xi\right)  |_{\xi=0}=0$, or alternatively, $\left(  \frac
{d\;}{dz}\right)  ^{k}m_{0}\left(  z\right)  |_{z=-1}=0$, or in yet another
form, the condition that $\left(  1+z\right)  ^{k+1}$ is a factor of
$m_{0}\left(  z\right)  $, etc.

\begin{proposition}
\label{ProMom.2}\textup{(a)} The polynomial $m_{0}^{\left(  \theta
,\rho\right)  }\left(  z\right)  $ is divisible by $\left(  1+z\right)  ^{2}$
\textup{(}see \textup{(\ref{eqMom.1}))} if and only if
\begin{equation}
\cos\left(  2\theta\right)  +\cos\left(  2\rho\right)  =\frac{1}{2}
\label{eqMom.3}%
\end{equation}
\textup{(}shown as curves in the four corners of Figure
\textup{\ref{FigSmooth}} in the Appendix\/\textup{).}

\textup{(b)} The polynomial $m_{0}^{\left(  \theta,\rho\right)  }\left(
z\right)  $ is divisible by $\left(  1+z\right)  ^{3}$ if and only if%
\begin{equation}
\cos2\theta+\cos2\rho=\frac{1}{2}\text{\quad and\quad}\sin2\theta+\sin
2\rho=2\sin\left(  2\theta-2\rho\right)  , \label{eqMom.4}%
\end{equation}
i.e., when
\[
\theta=\cos^{-1}\sqrt[4]{\frac{5}{32}}\approx0.89\approx0.28\,\pi ,\qquad
\rho=\cos^{-1}\sqrt{\frac{5}{4}-\sqrt{\frac{5}{32}}}\approx0.39\approx
0.12\,\pi
\]
or when $\left(  \theta,\rho\right)  $ is related to this pair by
$\left(  \theta,\rho\right)  \rightarrow\left(  \theta+m\pi
,\rho+n\pi\right)  $, $m,n\in\mathbb{Z}$, or $\left(
\theta,\rho\right)  \rightarrow\left(  -\theta,-\rho\right)  $, or both.
\textup{(}The points in $\left[  0,\pi\right]  \times\left[  0,\pi\right]  $
are shown in Figure \textup{\ref{FigSmooth}.)}
\end{proposition}

\begin{proof}
(a) In this example, $N=2$ and $g=3$ . Let $m\left(  z\right)  :=\left(
\begin{smallmatrix}
m_{0}\left(  z\right) \\
m_{1}\left(  z\right)
\end{smallmatrix}
\right)  $ be the usual QMF-polynomials in $z$, viewed as a column vector. By
Lemma \ref{CorBrJo00Rep.2}, we have
\begin{equation}
m\left(  z\right)  =V\left(  P^{\perp}+z^{2}P\right)  \left(  Q^{\perp}%
+z^{2}Q\right)  \alpha\left(  z\right)  , \label{eqMom.7}%
\end{equation}
where $V\in\mathrm{U}_{2}\left(  \mathbb{C}\right)  $, $P$ and $Q$ are
projections, and as before, $\alpha\left(  z\right)  :=\left(
\begin{smallmatrix}
1\\
z
\end{smallmatrix}
\right)  $. Here $V=\frac{1}{\sqrt{2}}\left(
\begin{smallmatrix}
1 & 1\\
1 & -1
\end{smallmatrix}
\right)  $ and $P$, $Q$ are the $\theta$, $\rho$ projections
specified in (\ref{eqSpe.20and21}) (see also (\ref{eqA.9})
below). We need only check that $\frac{d\;}{dz}m\left(  z\right)
|_{z=-1}=\left(
\begin{smallmatrix}
0\\
\xi_{0}^{{}}%
\end{smallmatrix}
\right)  $ for some number $\xi_{0}^{{}}$. But%
\begin{multline}
\frac{d\;}{dz}m\left(  z\right)  =V2zP\left(  Q^{\perp}+z^{2}Q\right)
\alpha\left(  z\right)  +V\left(  P^{\perp}+z^{2}P\right)  2zQ\alpha\left(
z\right) \label{eqMom.8}\\
+V\left(  P^{\perp}+z^{2}P\right)  \left(  Q^{\perp}+z^{2}Q\right)  \left(
\begin{array}
[c]{c}%
0\\
1
\end{array}
\right)  .
\end{multline}
Substitution of $z=-1$ yields%
\begin{equation}
m^{\prime}\left(  -1\right)  =-2V\left(  P+Q\right)  \left(
\begin{array}
[c]{c}%
1\\
-1
\end{array}
\right)  +V\left(
\begin{array}
[c]{c}%
0\\
1
\end{array}
\right)  . \label{eqMom.9}%
\end{equation}
Hence%
\begin{equation}
2\left(  P+Q\right)  \left(
\begin{array}
[c]{c}%
1\\
-1
\end{array}
\right)  -\left(
\begin{array}
[c]{c}%
0\\
1
\end{array}
\right)  =\left(
\begin{array}
[c]{c}%
\xi_{1}^{{}}\\
-\xi_{1}^{{}}%
\end{array}
\right)  , \label{eqMom.10}%
\end{equation}
for some number $\xi_{1}^{{}}$, and therefore, with%
\begin{equation}
P+Q=\left(
\begin{array}
[c]{cc}%
1 & 0\\
0 & 1
\end{array}
\right)  +\frac{1}{2}\left(
\begin{array}
[c]{cc}%
\cos2\theta+\cos2\rho & \sin2\theta+\sin2\rho\\
\sin2\theta+\sin2\rho & -\left(  \cos2\theta+\cos2\rho\right)
\end{array}
\right)  \label{eqMom.11}%
\end{equation}
we have
\[
\left(
\begin{array}
[c]{c}%
2\\
-2
\end{array}
\right)  +\left(
\begin{array}
[c]{cc}%
\cos2\theta+\cos2\rho & \sin2\theta+\sin2\rho\\
\sin2\theta+\sin2\rho & -\left(  \cos2\theta+\cos2\rho\right)
\end{array}
\right)  \left(
\begin{array}
[c]{c}%
1\\
-1
\end{array}
\right)  -\left(
\begin{array}
[c]{c}%
0\\
1
\end{array}
\right)  =\left(
\begin{array}
[c]{c}%
\xi_{1}^{{}}\\
-\xi_{1}^{{}}%
\end{array}
\right)  .
\]
The result (a) follows.

Part (b) follows upon solving%
\begin{equation}
\frac{d^{2}\;}{dz^{2}}m\left(  z\right)  |_{z=-1}=\left(
\begin{array}
[c]{c}%
0\\
\xi_{2}^{{}}%
\end{array}
\right)  , \label{eqMom.13}%
\end{equation}
for some number $\xi_{2}^{{}}$, in addition to the conditions in (a). But we
only need to work out the next derivative, using (\ref{eqMom.3}) to eliminate
the cosine terms where possible:%
\begin{align}
&  V^{-1}\frac{d^{2}\;}{dz^{2}}m\left(  z\right)  |_{z=-1}\label{eqMom.14}\\
&  \qquad=\left(  2\left(  P+Q\right)  +8PQ\right)  \left(
\begin{array}
[c]{c}%
1\\
-1
\end{array}
\right)  -4\left(  P+Q\right)  \left(
\begin{array}
[c]{c}%
0\\
1
\end{array}
\right) \nonumber\\
&  \qquad=\left(
\begin{array}
[c]{c}%
\frac{5}{2}-3\left(  \sin2\theta+\sin2\rho\right) \\
-\frac{9}{2}+\left(  \sin2\theta+\sin2\rho\right)
\end{array}
\right)  +\left(
\begin{array}
[c]{c}%
3-2\left(  \sin2\theta+\sin2\rho\right) \\
-1+2\left(  \sin2\theta+\sin2\rho\right)
\end{array}
\right) \nonumber\\
&  \qquad\qquad+2\left(
\begin{array}
[c]{c}%
\cos2\left(  \theta-\rho\right)  +\sin2\left(  \theta-\rho\right) \\
\sin2\left(  \theta-\rho\right)  -\cos2\left(  \theta-\rho\right)
\end{array}
\right) \nonumber\\
&  \qquad=V^{-1}\left(
\begin{array}
[c]{c}%
0\\
\xi_{2}^{{}}%
\end{array}
\right)  =\left(
\begin{array}
[c]{c}%
\xi_{3}^{{}}\\
-\xi_{3}^{{}}%
\end{array}
\right)  ,\nonumber
\end{align}
where we used%
\begin{multline}
PQ=\frac{1}{4}\left(  \left(
\begin{array}
[c]{cc}%
1 & 0\\
0 & 1
\end{array}
\right)  +\left(
\begin{array}
[c]{cc}%
\cos2\theta+\cos2\rho & \sin2\theta+\sin2\rho\\
\sin2\theta+\sin2\rho & -\left(  \cos2\theta+\cos2\rho\right)
\end{array}
\right)  \right. \label{eqMom14bis}\\
\left.  +\left(
\begin{array}
[c]{cc}%
\cos2\left(  \theta-\rho\right)  & -\sin2\left(  \theta-\rho\right) \\
\sin2\left(  \theta-\rho\right)  & \cos2\left(  \theta-\rho\right)
\end{array}
\right)  \right)  .
\end{multline}
Noting again that the right-hand side of (\ref{eqMom.14}) takes the form
$\left( 
\begin{smallmatrix}
\xi\\
-\xi
\end{smallmatrix}
\right) $, we arrive at%
\begin{equation}
-2\left(  \sin2\theta+\sin2\rho\right)  +4\sin2\left(  \theta-\rho\right)  =0,
\label{eqMom.15}%
\end{equation}
which directly gives the second part of (\ref{eqMom.4}).

The specific solution is found by transforming (\ref{eqMom.15}) into%
\begin{equation}
\left(  1-2\cos2\rho\right)  \sin2\theta=-\left(  1+2\cos2\theta\right)
\sin2\rho. \label{eqMom.16}%
\end{equation}
Squaring and using (\ref{eqMom.3}) to eliminate $\rho$ yields%
\begin{equation}
2\cos^{2}2\theta+4\cos2\theta+\frac{3}{4}=0, \label{eqMom.17}%
\end{equation}
which gives
\begin{equation}
\cos2\theta=-1+\sqrt{\frac{5}{8}} \label{eqMom.18}%
\end{equation}
by the quadratic formula, using $\left|  \cos2\theta\right|  \leq1$ to choose
the positive radical. Working this back through (\ref{eqMom.3}) yields%
\begin{equation}
\cos2\rho=\frac{3}{2}-\sqrt{\frac{5}{8}}, \label{eqMom.19}%
\end{equation}
and substitution of (\ref{eqMom.18}) and (\ref{eqMom.19}) into (\ref{eqMom.16}%
) then shows that $\sin2\theta$ and $\sin2\rho$ must have the same sign, so
that the solutions stated in the proposition are the only ones possible, the
numerically exhibited pair being those for which $\sin2\theta$ and $\sin2\rho$
are both positive.
\end{proof}

\begin{figure}[tbp]
\setlength{\unitlength}{1bp}
\begin{picture}
(240,240) \put(0,0){\includegraphics[bb=0 1 360 359,height=240bp,width=240bp]
{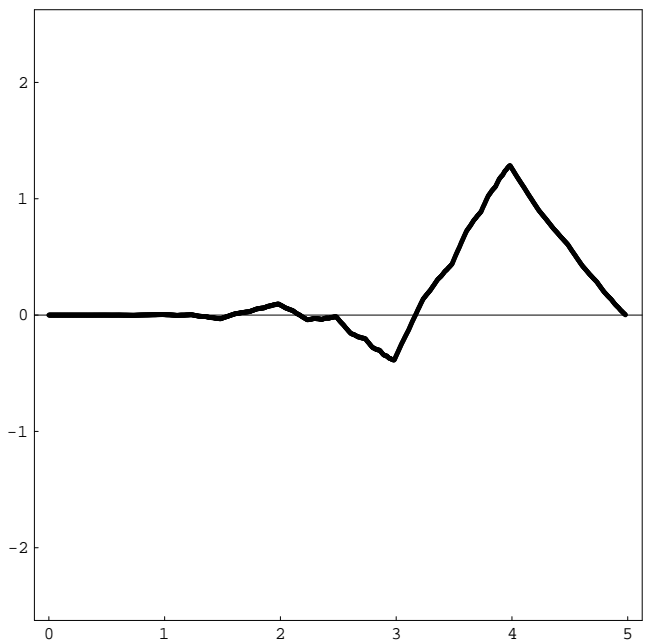}}
\end{picture}
\caption{The ultra-smooth wavelet scaling function: $\varphi_{\theta,\rho}%
^{{}}\left(  x\right)  $ for $\theta\approx0.284\,\pi$, $\rho\approx
0.124\,\pi$, with two vanishing moments of $\psi$ (see Proposition
\ref{ProMom.2}(b) and the discussion preceding it). The H\"{o}lder-Sobolev
exponent of this $\varphi$ is $\geq1.136$; see Remark \ref{RemExpNew.HoSo}.}%
\label{FigMomThree}%
\end{figure}

\begin{remark}
\label{RemExpNew.HoSo}The H\"{o}lder-Sobolev exponent is at least as good as
$0.84$ when we have a $(1+z)^{2}$ factor of $m_{0}(z)$, and at least as good
as $1.136$ if $(1+z)^{3}$ is a factor \cite{LaSu00} \cite{Vol95}; see Figure
\textup{\ref{FigMomThree}.}
\end{remark}

\section*{\label{App}Appendix (by Brian Treadway)}

\setcounter{equation}{0}\setcounter{theorem}{0}\renewcommand{\theequation
}{A.\arabic{equation}}\renewcommand{\thetheorem}{A.\arabic{theorem}}In this
appendix we show some cascade approximations of wavelet scaling functions for
the example with $g=3$ discussed in Sections \ref{Irr} and \ref{Exp} above.
The local method of cascade iteration works here just as it did in
\cite{BrJo99b}. It is just a matter of enumerating terms.

The direct method of iterating the relation (\ref{eqInt.9}), i.e.,%
\begin{equation}
\varphi\left(  x\right)  =\sum_{k=0}^{2g-1}a_{k}\varphi\left(  2x-k\right)  ,
\label{eqA.1}%
\end{equation}
proceeds by translating a distance $k$ to the right, multiplying by $a_{k}$,
summing over $k$, and scaling down by $2$. This takes an expression in which
every term has $n$ factors, each of which is an $a_{k}$ for some $k$, and
turns it into an expression in which every term has $n+1$ such factors. In
fact, every ordered product of $n$ $a_{k}$'s occurs exactly once, at some
dyadic point in the $n$'th stage. Multiplication of these factors is
commutative, but we forget that for the moment and take care to add the new
factors at the left.

For example: With $g=3$ there are $6$ coefficients, $a_{0},a_{1},\dots,a_{5}$.\medskip

$n=0$: $1$ point, $1=6^{0}$ term, no factors of $a_{k}$, so the vacuous term
is just $1$, and we have the Haar function:
\[%
\begin{array}
[c]{ll}%
x=0 & \varphi=1
\end{array}
\]

$n=1$: $6$ points, $6=6^{1}$ terms, as follows:
\[%
\begin{array}
[t]{ll}%
x=0 & \varphi=a_{0}\\
x=1/2 & \varphi=a_{1}\\
x=1 & \varphi=a_{2}%
\end{array}
\quad
\begin{array}
[t]{ll}%
x=3/2 & \varphi=a_{3}\\
x=2 & \varphi=a_{4}\\
x=5/2 & \varphi=a_{5}%
\end{array}
\]

$n=2$: $16$ points, $36=6^{2}$ terms, as follows:
\[%
\begin{array}
[t]{ll}%
x=0 & \varphi=a_{0}a_{0}\\
x=1/4 & \varphi=a_{0}a_{1}\\
x=1/2 & \varphi=a_{1}a_{0}+a_{0}a_{2}\\
x=3/4 & \varphi=a_{1}a_{1}+a_{0}a_{3}\\
x=1 & \varphi=a_{2}a_{0}+a_{1}a_{2}+a_{0}a_{4}\\
x=5/4 & \varphi=a_{2}a_{1}+a_{1}a_{3}+a_{0}a_{5}\\
x=3/2 & \varphi=a_{3}a_{0}+a_{2}a_{2}+a_{1}a_{4}\\
x=7/4 & \varphi=a_{3}a_{1}+a_{2}a_{3}+a_{1}a_{5}%
\end{array}
\quad
\begin{array}
[t]{ll}%
x=2 & \varphi=a_{4}a_{0}+a_{3}a_{2}+a_{2}a_{4}\\
x=9/4 & \varphi=a_{4}a_{1}+a_{3}a_{3}+a_{2}a_{5}\\
x=5/2 & \varphi=a_{5}a_{0}+a_{4}a_{2}+a_{3}a_{4}\\
x=11/4 & \varphi=a_{5}a_{1}+a_{4}a_{3}+a_{3}a_{5}\\
x=3 & \varphi=a_{5}a_{2}+a_{4}a_{4}\\
x=13/4 & \varphi=a_{5}a_{3}+a_{4}a_{5}\\
x=7/2 & \varphi=a_{5}a_{4}\\
x=15/4 & \varphi=a_{5}a_{5}%
\end{array}
\]

$n=3$: $36$ points, $216=6^{3}$ terms, grouping left as an exercise for the reader.\medskip

One could write down all $6^{n}$ terms at the outset of each stage, and then
just ask which ones go with which of the values of $x$. The answer is:
\begin{equation}%
\addtolength{\displayboxwidth}{\leftmargini}\begin{minipage}
[t]{\displayboxwidth}\raggedright The indices of the $a$'s, interpreted as digits
at the right of the fraction point in a non-unique positional binary number
system with six digits instead of two, give the $x$ to which a given term
should be assigned.%
\end{minipage}
\label{eqA.answer}%
\end{equation}
As an example from the list above, the term $a_{4}a_{1}$ goes with
$x=4\left(  2^{-1}\right)  +1\left(  2^{-2}\right)  =9/4$. All the others (at
any cascade stage) can be assigned in the same way.

In the ordering used above, if we go from one
stage to the next by the ``translate by $k$, multiply by $a_{k}$, sum
over $k$, and scale down by $2$'' method, terms will be built up by adding
factors from the left, while if we use the ``local linear combination'' method
from \cite[Appendix]{BrJo99b} and \cite[Section 6.5, pp.~204--206]{Dau92},
terms will be built up by adding factors from the right. In either method, the
full set of $6^{n}$ terms (or $4^{n}$ in \cite{BrJo99b}, where $g=2$) will be
obtained, without duplication, and each term will be assigned to the same $x$
regardless of which method is used.

At cascade stage $n$ there are values of $\varphi$ assigned to
values $x_{i}$ of $x$ that are consecutive integer
multiples of $2^{-n}$ ranging from $0$ to $\left(  q^{\left(  n\right)
}-1\right)  \cdot2^{-n}$. The iteration begins with the ordered list of values
of $\varphi$ at those points on the $x$-axis:%
\begin{equation}
\left(  \varphi^{\left(  n\right)  }\left(  0\right)  ,\varphi^{\left(
n\right)  }\left(  1\cdot2^{-n}\right)  ,\varphi^{\left(  n\right)  }\left(
2\cdot2^{-n}\right)  ,\dots,\varphi^{\left(  n\right)  }\left(  \left(
q^{\left(  n\right)  }-1\right)  \cdot2^{-n}\right)  \right)  . \label{eqA.2}%
\end{equation}
Sets of $g$ adjacent values in this list go into the computation of each point
at the next cascade stage, which is on a finer grid of $q^{\left(  n+1\right)
}$ consecutive integer multiples of $2^{-\left(  n+1\right)  }$. Each set of
$g$ adjacent values of $\varphi$ from the list (\ref{eqA.2}) yields $2$ values
on the $2^{-\left(  n+1\right)  }$ grid by a linear combination with two sets
of alternate $a_{k}$'s as coefficients. This can be expressed as a matrix
product (shown here in the case $g=3$):%
\begin{multline}
\left(
\begin{array}
[c]{ll}%
\varphi^{\left(  n+1\right)  }\left(  0\right)  & \varphi^{\left(  n+1\right)
}\left(  1\cdot2^{-\left(  n+1\right)  }\right) \\
\varphi^{\left(  n+1\right)  }\left(  2\cdot2^{-\left(  n+1\right)  }\right)
& \varphi^{\left(  n+1\right)  }\left(  3\cdot2^{-\left(  n+1\right)  }\right)
\\
\varphi^{\left(  n+1\right)  }\left(  4\cdot2^{-\left(  n+1\right)  }\right)
& \varphi^{\left(  n+1\right)  }\left(  5\cdot2^{-\left(  n+1\right)  }\right)
\\
\varphi^{\left(  n+1\right)  }\left(  6\cdot2^{-\left(  n+1\right)  }\right)
& \varphi^{\left(  n+1\right)  }\left(  7\cdot2^{-\left(  n+1\right)  }\right)
\\
\vdots & \vdots\\
\varphi^{\left(  n+1\right)  }\left(  \left(  q^{\left(  n+1\right)
}-6\right)  \cdot2^{-\left(  n+1\right)  }\right)  & \varphi^{\left(
n+1\right)  }\left(  \left(  q^{\left(  n+1\right)  }-5\right)  \cdot
2^{-\left(  n+1\right)  }\right) \\
\varphi^{\left(  n+1\right)  }\left(  \left(  q^{\left(  n+1\right)
}-4\right)  \cdot2^{-\left(  n+1\right)  }\right)  & \varphi^{\left(
n+1\right)  }\left(  \left(  q^{\left(  n+1\right)  }-3\right)  \cdot
2^{-\left(  n+1\right)  }\right) \\
\varphi^{\left(  n+1\right)  }\left(  \left(  q^{\left(  n+1\right)
}-2\right)  \cdot2^{-\left(  n+1\right)  }\right)  & \varphi^{\left(
n+1\right)  }\left(  \left(  q^{\left(  n+1\right)  }-1\right)  \cdot
2^{-\left(  n+1\right)  }\right)
\end{array}
\right) \label{eqA.3}\\
=\left(
\begin{array}
[c]{lll}%
0 & 0 & \varphi^{\left(  n\right)  }\left(  0\right) \\
0 & \varphi^{\left(  n\right)  }\left(  0\right)  & \varphi^{\left(  n\right)
}\left(  1\cdot2^{-n}\right) \\
\varphi^{\left(  n\right)  }\left(  0\right)  & \varphi^{\left(  n\right)
}\left(  1\cdot2^{-n}\right)  & \varphi^{\left(  n\right)  }\left(
2\cdot2^{-n}\right) \\
\varphi^{\left(  n\right)  }\left(  1\cdot2^{-n}\right)  & \varphi^{\left(
n\right)  }\left(  2\cdot2^{-n}\right)  & \varphi^{\left(  n\right)  }\left(
3\cdot2^{-n}\right) \\
\vdots & \vdots & \vdots\\
\varphi^{\left(  n\right)  }\left(  \left(  q^{\left(  n\right)  }-3\right)
\cdot2^{-n}\right)  & \varphi^{\left(  n\right)  }\left(  \left(  q^{\left(
n\right)  }-2\right)  \cdot2^{-n}\right)  & \varphi^{\left(  n\right)
}\left(  \left(  q^{\left(  n\right)  }-1\right)  \cdot2^{-n}\right) \\
\varphi^{\left(  n\right)  }\left(  \left(  q^{\left(  n\right)  }-2\right)
\cdot2^{-n}\right)  & \varphi^{\left(  n\right)  }\left(  \left(  q^{\left(
n\right)  }-1\right)  \cdot2^{-n}\right)  & 0\\
\varphi^{\left(  n\right)  }\left(  \left(  q^{\left(  n\right)  }-1\right)
\cdot2^{-n}\right)  & 0 & 0
\end{array}
\right) \\
\bullet\left(
\begin{array}
[c]{ll}%
a_{4} & a_{5}\\
a_{2} & a_{3}\\
a_{0} & a_{1}%
\end{array}
\right)  .
\end{multline}
Thus a $\left(  q^{\left(  n\right)  }+\left(  g-1\right)  \right)  \times g$
matrix (the first factor on the right in (\ref{eqA.3}) above), partitioned
from the list (\ref{eqA.2}) of $q^{\left(  n\right)  }$ points at stage $n$,
yields a $\left(  q^{\left(  n\right)  }+\left(  g-1\right)  \right)  \times2$
matrix (the left-hand side of (\ref{eqA.3})), which is then flattened out to
give the list of values for stage $n+1$. The number of points $q^{\left(
n+1\right)  }$ in the new list is the total number of entries in that $\left(
q^{\left(  n\right)  }+\left(  g-1\right)  \right)  \times2$ matrix,%
\begin{equation}
q^{\left(  n+1\right)  }=2\left(  q^{\left(  n\right)  }+\left(  g-1\right)
\right)  , \label{eqA.4}%
\end{equation}
and with $q^{\left(  0\right)  }=1$, this recursively yields the number
$q^{\left(  n\right)  }$ of points at each stage as%
\begin{equation}
q^{\left(  n\right)  }=\left(  2g-1\right)  \cdot2^{n}-2\left(  g-1\right)  .
\label{eqA.5}%
\end{equation}
In the case $g=3$, this is%
\begin{equation}
q^{\left(  n\right)  }=5\cdot2^{n}-4. \label{eqA.6}%
\end{equation}

The three steps in the local method, (1)~partitioning a list into rows of
$g=3$ points (first adding $g-1$ zeroes at each end), (2)~matrix
multiplication, and (3)~flattening the resulting matrix back into a single
list, are easily implemented in \textit{Mathematica} \cite{Wol96}. All that
then remains to compute a cascade approximation of a wavelet scaling function
$\varphi$ is to specify the numerical values of the coefficients $a_{k}$ and
repeat the procedure $n$ times, starting with the one-element list%
\begin{equation}
\left(
\begin{array}
[c]{c}%
1
\end{array}
\right)  . \label{eqA.7}%
\end{equation}

The same local cascade relation expressed in (\ref{eqA.3}) as giving two
values of $\varphi^{\left(  n+1\right)  }$ from three values of $\varphi
^{\left(  n\right)  }$ can also be set up to give four values of
$\varphi^{\left(  n+1\right)  }$ from four values of $\varphi^{\left(
n\right)  }$, or five values of $\varphi^{\left(  n+1\right)  }$ from five
values of $\varphi^{\left(  n\right)  }$, by combining overlapping ranges of
the initial and final lists. The terms involved on the successive $x$-grids
are shown in the diagrams below, and the corresponding matrices are given.%
\[%
\begin{array}
[c]{ll}%
\begin{array}
[c]{cc}%
\setlength{\unitlength}{1bp}%
\begin{picture}
(216,49)\put(0,0){\includegraphics[bb=0 0 432 98, height=49bp,width=216bp]
{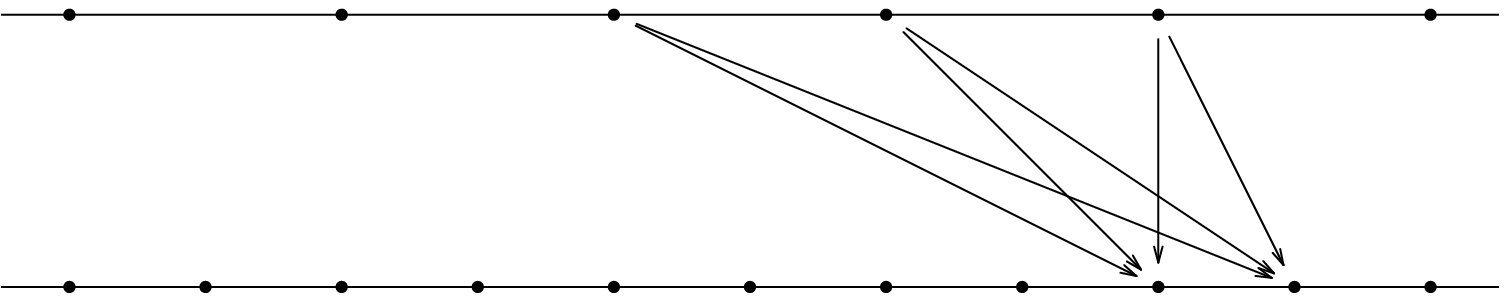}}%
\end{picture}
& \hspace{-16pt}%
\end{array}
& \left(
\begin{array}
[c]{cc}%
a_{4} & a_{5}\\
a_{2} & a_{3}\\
a_{0} & a_{1}%
\end{array}
\right)  \text{\quad(as in (\ref{eqA.3}))}\\%
\begin{array}
[c]{cc}%
\setlength{\unitlength}{1bp}%
\begin{picture}
(216,61)\put(0,0){\includegraphics[bb=0 0 432 98, height=49bp,width=216bp]
{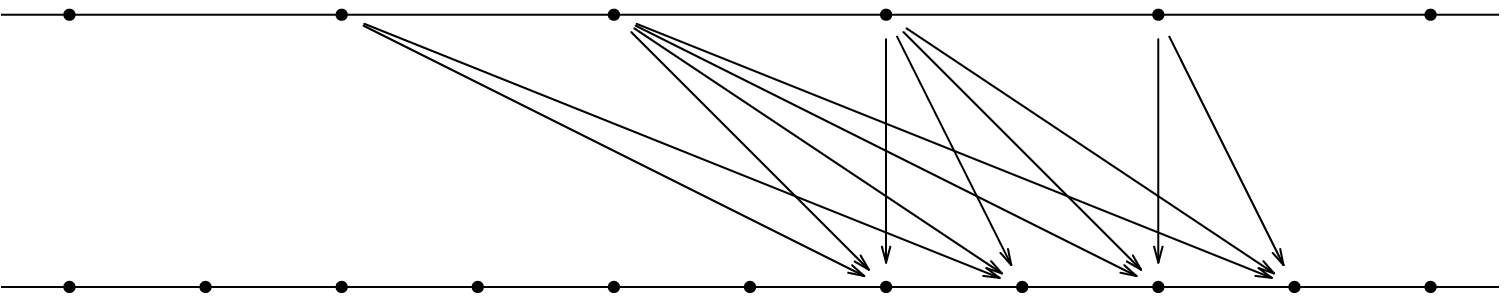}}%
\end{picture}
& \hspace{-16pt}%
\end{array}
& \left(
\begin{array}
[c]{cccc}%
a_{4} & a_{5} & 0 & 0\\
a_{2} & a_{3} & a_{4} & a_{5}\\
a_{0} & a_{1} & a_{2} & a_{3}\\
0 & 0 & a_{0} & a_{1}%
\end{array}
\right) \\%
\begin{array}
[c]{cc}%
\setlength{\unitlength}{1bp}%
\begin{picture}
(216,61)\put(0,0){\includegraphics[bb=0 0 432 98, height=49bp,width=216bp]
{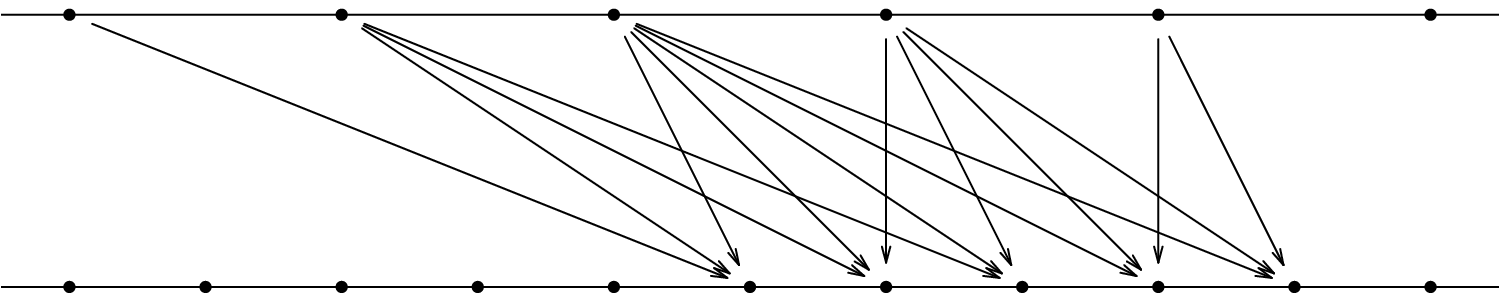}}%
\end{picture}
& \hspace{-16pt}%
\end{array}
& \left(
\begin{array}
[c]{ccccc}%
a_{5} & 0 & 0 & 0 & 0\\
a_{3} & a_{4} & a_{5} & 0 & 0\\
a_{1} & a_{2} & a_{3} & a_{4} & a_{5}\\
0 & a_{0} & a_{1} & a_{2} & a_{3}\\
0 & 0 & 0 & a_{0} & a_{1}%
\end{array}
\right) \\%
\begin{array}
[c]{cc}%
\setlength{\unitlength}{1bp}%
\begin{picture}
(216,61)\put(0,0){\includegraphics[bb=0 0 432 98, height=49bp,width=216bp]
{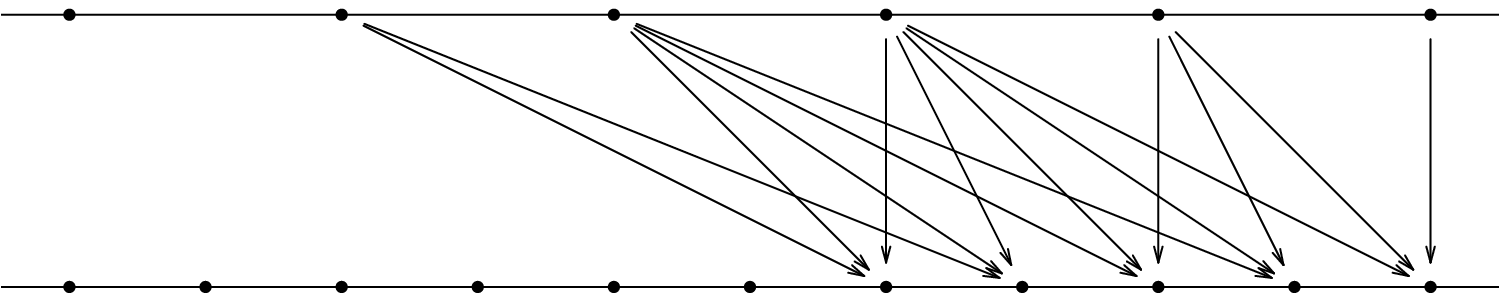}}%
\end{picture}
& \hspace{-16pt}%
\end{array}
& \left(
\begin{array}
[c]{ccccc}%
a_{4} & a_{5} & 0 & 0 & 0\\
a_{2} & a_{3} & a_{4} & a_{5} & 0\\
a_{0} & a_{1} & a_{2} & a_{3} & a_{4}\\
0 & 0 & a_{0} & a_{1} & a_{2}\\
0 & 0 & 0 & 0 & a_{0}%
\end{array}
\right)
\end{array}
\]
All of these represent the \emph{same} calculation of the $\left(  n+1\right)
$'st cascade stage fron the $n$'th stage: they merely collect different
locally related sets of points in the successive stages. The advantage of the
matrices that are square is that they allow successive stages to be expressed
as powers of the matrix.

An eigenvector decomposition yielding an explicit expression for the
$n\rightarrow\infty$ limit of the cascade stages, like that done in
\cite[Appendix]{BrJo99b}, could be done here using the $4\times4$ matrix above
(here the matrices are written to act on the left). The detailed calculation
is much more extensive than in the $2\times2$ case of \cite{BrJo99b}, and we
will not present it here, but note only that the two $5\times5$ matrices above
have simple eigenvalues and left eigenvectors (in addition to the row
$\left(
\begin{array}
[c]{cccc}%
1 & 1 & \dots & 1
\end{array}
\right)  $ that all the matrices have): for the first, eigenvalue $a_{5}$,
eigenvector $\left(
\begin{array}
[c]{ccccc}%
1 & 0 & 0 & 0 & 0
\end{array}
\right)  $, and for the second, eigenvalue $a_{0}$, eigenvector $\left(
\begin{array}
[c]{ccccc}%
0 & 0 & 0 & 0 & 1
\end{array}
\right)  $. Since the starting list for the cascade computation is $\left(
\begin{array}
[c]{c}%
1
\end{array}
\right)  $, to be ``padded'' with zeroes on the left and right, these two
eigenvectors occur explicitly at the first cascade stage and continue
thereafter. This shows that $\varphi^{\left(  n\right)  }$ diverges at one end
of its support interval or the other when one of these eigenvalues is greater
than $1$, growing like $a_{5}^{n}$ when $a_{5}>1$ or like $a_{0}^{n}$ when
$a_{0}>1$. The full eigenvector decomposition then has a term, generically
nonzero, that grows in the same way for other points $x_{i}$, and so the same
divergence occurs at points throughout the support interval. The regions where
$a_{0}$ or $a_{5}$ is greater than $1$, leading to this divergence in the
cascade iteration, are shown by shading in Figure \ref{FigSmooth}.

The coefficients here are not restricted, as they are in \cite{Wan00},
to $a_{i}\leq 1$, so we do not have cycles in the cascade iteration:
the terms grow indefinitely. As a result, the condition $a_{5}>1$ or
$a_{0}>1$ is sufficient, but not necessary, for divergence.

For the example with $g=3$ specified in (\ref{eqSpe.20and21})
above, the formulas for the coefficients $a_{0},a_{1},\dots,a_{5}$ of the
polynomial $m_{0}^{\left(  A\right)  }\left(  z\right)  $ in (\ref{eqExp.8})
may be derived as follows. From (\ref{eqSpe.9}) we have
\begin{equation}
A(z)=V(Q_{\theta}^{\perp}+zQ_{\theta}^{{}})(Q_{\rho}^{\perp}+zQ_{\rho}^{{}})
\label{eqA.8}%
\end{equation}
with%
\begin{equation}
V=\left(
\begin{array}
[c]{cc}%
1 & 1\\
1 & -1
\end{array}
\right)  , \label{eqA.8bis}%
\end{equation}%
\begin{equation}
Q_{\theta}=\left(
\begin{array}
[c]{cc}%
\cos^{2}\theta & \cos\theta\sin\theta\\
\cos\theta\sin\theta & \sin^{2}\theta
\end{array}
\right)  =\frac{1}{2}\left(  \left(
\begin{array}
[c]{cc}%
1 & 0\\
0 & 1
\end{array}
\right)  +\left(
\begin{array}
[c]{cc}%
\cos2\theta & \sin2\theta\\
\sin2\theta & -\cos2\theta
\end{array}
\right)  \right)  , \label{eqA.9}%
\end{equation}
and%
\begin{equation}
Q_{\theta}^{\perp}=Q_{\theta+\left(  \pi/2\right)  }^{{}} \label{eqA.10}%
\end{equation}
Then the coefficients $a_{0},a_{1},\dots,a_{5}$ are:%
\begin{equation}%
\begin{aligned}
a_{0}  &  =\frac{1}{4}(1-\cos2\theta-\sin2\theta-\cos2\rho-\sin2\rho
+\cos(2\theta-2\rho)+\sin(2\theta-2\rho)),\\ a_{1}  &  =\frac{1}{4}%
(1+\cos2\theta-\sin2\theta+\cos2\rho-\sin2\rho+\cos(2\theta-2\rho
)-\sin(2\theta-2\rho)),\\ a_{2}  &  =\frac{1}{2}(1-\cos(2\theta-2\rho
)-\sin(2\theta-2\rho)),\\ a_{3}  &  =\frac{1}{2}(1-\cos(2\theta-2\rho
)+\sin(2\theta-2\rho)),\\ a_{4}  &  =\frac{1}{4}(1+\cos2\theta+\sin
2\theta+\cos2\rho+\sin2\rho+\cos(2\theta-2\rho)+\sin(2\theta-2\rho)),\\ a_{5}
&  =\frac{1}{4}(1-\cos2\theta+\sin2\theta-\cos2\rho+\sin2\rho+\cos
(2\theta-2\rho)-\sin(2\theta-2\rho)).
\end{aligned}
\label{eqA.11}%
\end{equation}
These can be seen to meet the conditions (\ref{eqInt.7})--(\ref{eqInt.8}) for
the coefficients of a scaling function, The even- and odd-indexed coefficients
also sum to a constant separately:
\begin{equation}
\sum_{i=0}^{2}a_{2i}=\sum_{i=0}^{2}a_{2i+1}=1 \label{eqA.11bis}%
\end{equation}
(see \cite[eq.~(9.12)]{ReWe98}), which is what makes the constant vector
$\left(
\begin{array}
[c]{cccc}%
1 & 1 & \dots & 1
\end{array}
\right)  $ an eigenvector or the $a_{i}$-matrices above. Of course, the
coefficients $a_{i}$, as functions of $\theta$ and $\rho$, have the
periodicity (with period $\pi$ in both angles $\theta$ and $\rho$) of the
projections $Q_{\theta}$ and $Q_{\rho}$ they were derived from:%
\begin{equation}
a_{i}\left(  \theta,\rho\right)  =a_{i}\left(  \theta+m\pi,\rho+n\pi\right)
,\qquad m,n\in\mathbb{Z}. \label{eqA.11ter}%
\end{equation}
In addition, they are related in pairs by the reflection relation%
\begin{equation}
a_{i}\left(  \theta,\rho\right)  =a_{5-i}\left(  -\theta,-\rho\right)  ,\qquad
i=0,\dots,5. \label{eqA.12}%
\end{equation}
These relations carry through the successive stages of the local cascade
computation (\ref{eqA.3}) as the
periodicity
\begin{equation}
\varphi_{\theta,\rho}^{\left(  n\right)  }\left(  x\right)  =\varphi
_{\theta+m\pi,\rho+n\pi}^{\left(  n\right)  }\left(  x\right)  \label{eqA.13per}%
\end{equation}
and the reflection symmetry%
\begin{equation}
\varphi_{\theta,\rho}^{\left(  n\right)  }\left(  x\right)  =\varphi
_{\theta-\pi,\rho-\pi}^{\left(  n\right)  }\left(  x\right)  =\varphi
_{\pi-\theta,\pi-\rho}^{\left(  n\right)  }\left(  x_{f}^{\left(  n\right)
}-x\right)  , \label{eqA.13}%
\end{equation}
where $x_{f}$ is the last point to which a value is assigned by the $n$th
stage. By (\ref{eqA.6}) we have%
\begin{equation}
x_{f}^{\left(  n\right)  }=\left(  q^{\left(  n\right)  }-1\right)
\cdot2^{-n}=5-5\cdot2^{-n}. \label{eqA.14}%
\end{equation}
There are pairwise relations involving a translation in the $\left(
\theta,\rho\right)  $ plane by half the period:%
\begin{equation}
a_{2i+j}\left(  \theta,\rho\right)  =a_{2\left(  2-i\right)  +j}\left(
\theta-\frac{\pi}{2},\rho-\frac{\pi}{2}\right)  ,\qquad i=0,1,2,\qquad j=0,1.
\label{eqA.15}%
\end{equation}
There are also twofold and threefold affine symmetries, such as
the invariance of $a_{2}$ under the twofold transformation
\begin{equation}
\theta\longmapsto -\theta,\qquad 
\rho\longmapsto -\rho +\frac{\pi}{4},\label{eqA.twofold}%
\end{equation}
or the
invariance of $a_{0}$ under the threefold transformation
\begin{equation}
\theta\longmapsto -\rho +\frac{\pi}{4},\qquad 
\rho\longmapsto\theta -\rho +\frac{\pi}{2},\label{eqA.threefold}%
\end{equation}
which has $a_{0}$'s three local extrema as its fixed points
in the $\pi$-periodic context.

Contour plots of $a_{0}$ and $a_{2}$ in the $\theta,\rho$-plane are shown in
Figure \ref{FigAppContour} below; the other $a_{i}$'s can be derived from
these by the translation in (\ref{eqA.15}) or the rotation around the origin
in (\ref{eqA.12}), or around the point $\left(  \frac{\pi}{2},\frac{\pi}%
{2}\right)  $ when this is combined with translation by $\pi$ in both angles.

\begin{figure}[tbp]
\setlength{\unitlength}{1bp}
\begin{picture}
(360,338) \put(17,174){\includegraphics[bb=0 0 360
360,height=163bp,width=163bp] {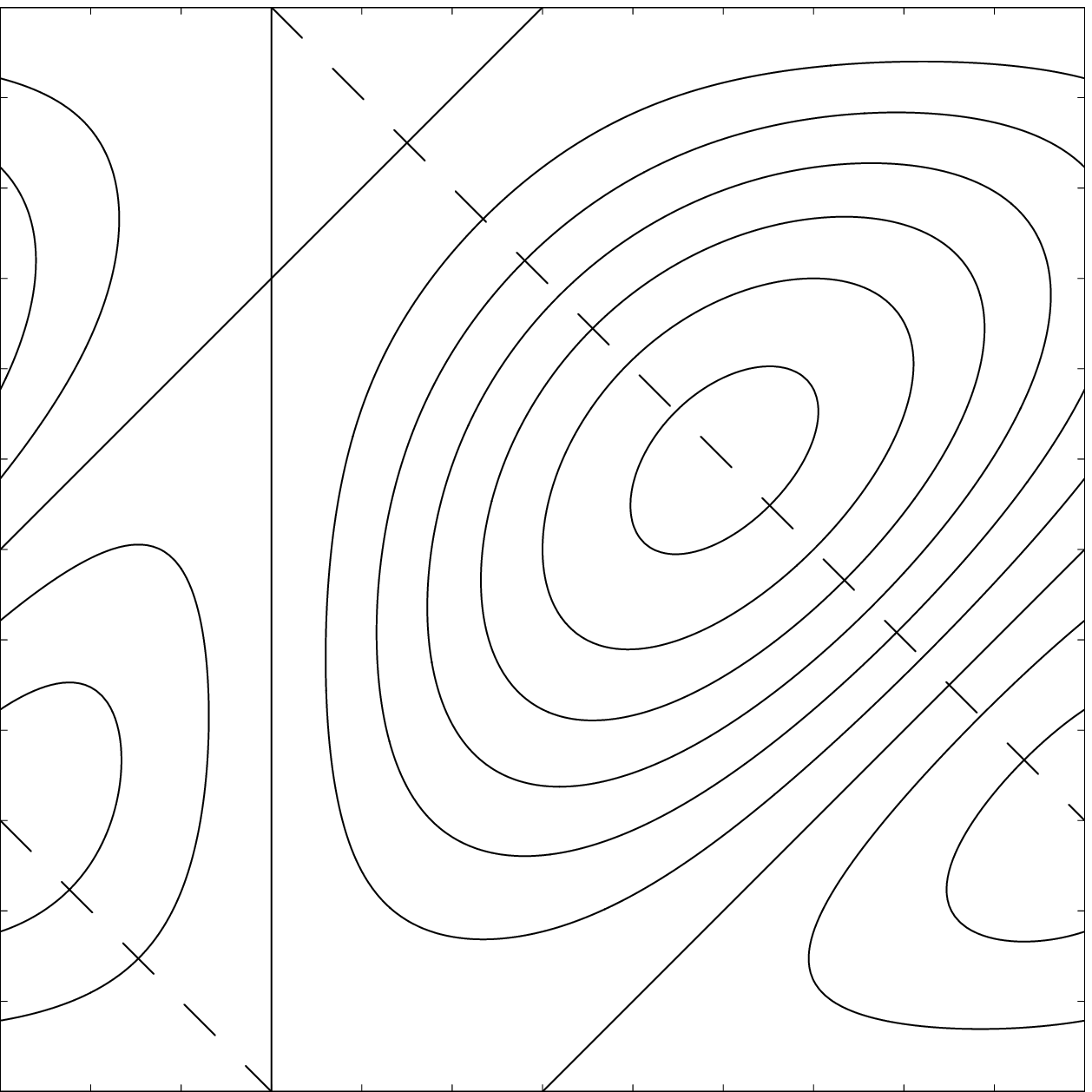}}  \put(197,174){\includegraphics
[bb=0 0 360 360,height=163bp,width=163bp] {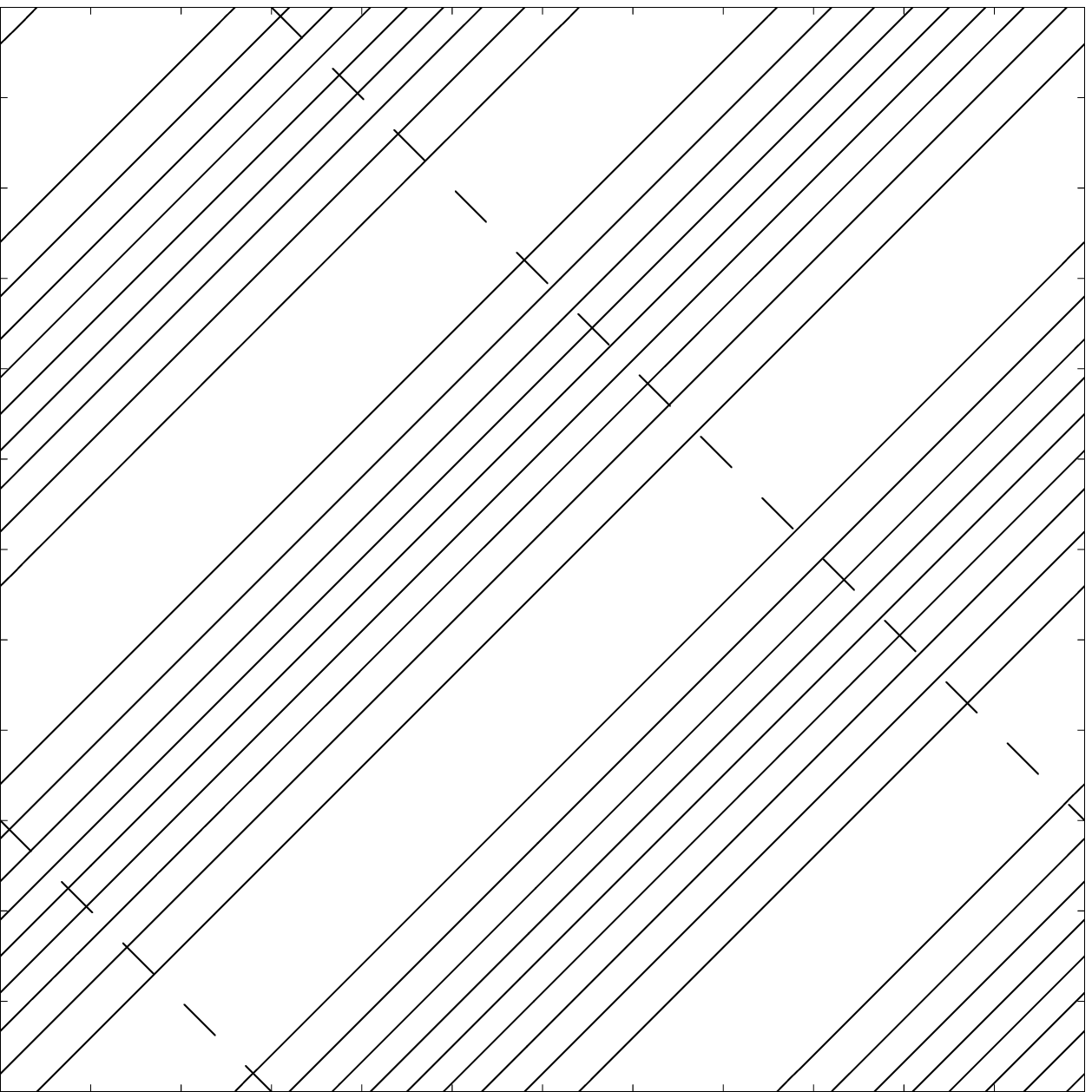}}  \put
(17,172){\makebox(0,0)[t] {$\scriptstyle0$}}  \put(58,172){\makebox(0,0)[t]
{$\scriptstyle\frac{\pi}{4}$}}  \put(98.67,172){\makebox(0,0)[t]
{$\scriptstyle\frac{\pi}{2}$}}  \put(139.33,172){\makebox(0,0)[t]
{$\scriptstyle\frac{3\pi}{4}$}}  \put(180,172){\makebox(0,0)[t] {$\scriptstyle
\pi$}}  \put(16,175){\makebox(0,0)[r] {$\scriptstyle0$}}  \put
(16,214.33){\makebox(0,0)[r] {$\scriptstyle\frac{\pi}{4}$}}  \put
(16,255){\makebox(0,0)[r] {$\scriptstyle\frac{\pi}{2}$}}  \put(3,236){\makebox
(0,0)[l] {$\rho$}}  \put(16,296.33){\makebox(0,0)[r] {$\scriptstyle\frac{3\pi
}{4}$}}  \put(16,337){\makebox(0,0)[tr] {$\scriptstyle\pi$}}  \put
(197,172){\makebox(0,0)[t] {$\scriptstyle0$}}  \put(238,172){\makebox(0,0)[t]
{$\scriptstyle\frac{\pi}{4}$}}  \put(278.67,172){\makebox(0,0)[t]
{$\scriptstyle\frac{\pi}{2}$}}  \put(319.33,172){\makebox(0,0)[t]
{$\scriptstyle\frac{3\pi}{4}$}}  \put(360,172){\makebox(0,0)[t] {$\scriptstyle
\pi$}}  \put(196,175){\makebox(0,0)[r] {$\scriptstyle0$}}  \put
(196,214.33){\makebox(0,0)[r] {$\scriptstyle\frac{\pi}{4}$}}  \put
(196,255){\makebox(0,0)[r] {$\scriptstyle\frac{\pi}{2}$}}  \put
(183,236){\makebox(0,0)[l] {$\rho$}}  \put(196,296.33){\makebox(0,0)[r]
{$\scriptstyle\frac{3\pi}{4}$}}  \put(196,337){\makebox(0,0)[tr]
{$\scriptstyle\pi$}}  \put(0,150){\makebox(180,14) {$\theta$}}  \put
(180,150){\makebox(180,14) {$\theta$}}  \put(16,136){\makebox(164,14) {Contour
plot of $a_{0}(\theta,\rho)$}}  \put(196,136){\makebox(164,12) {Contour plot
of $a_{2}(\theta,\rho)$}}  \put(0,14){\includegraphics[bb=0 9 360
213,height=102bp,width=180bp] {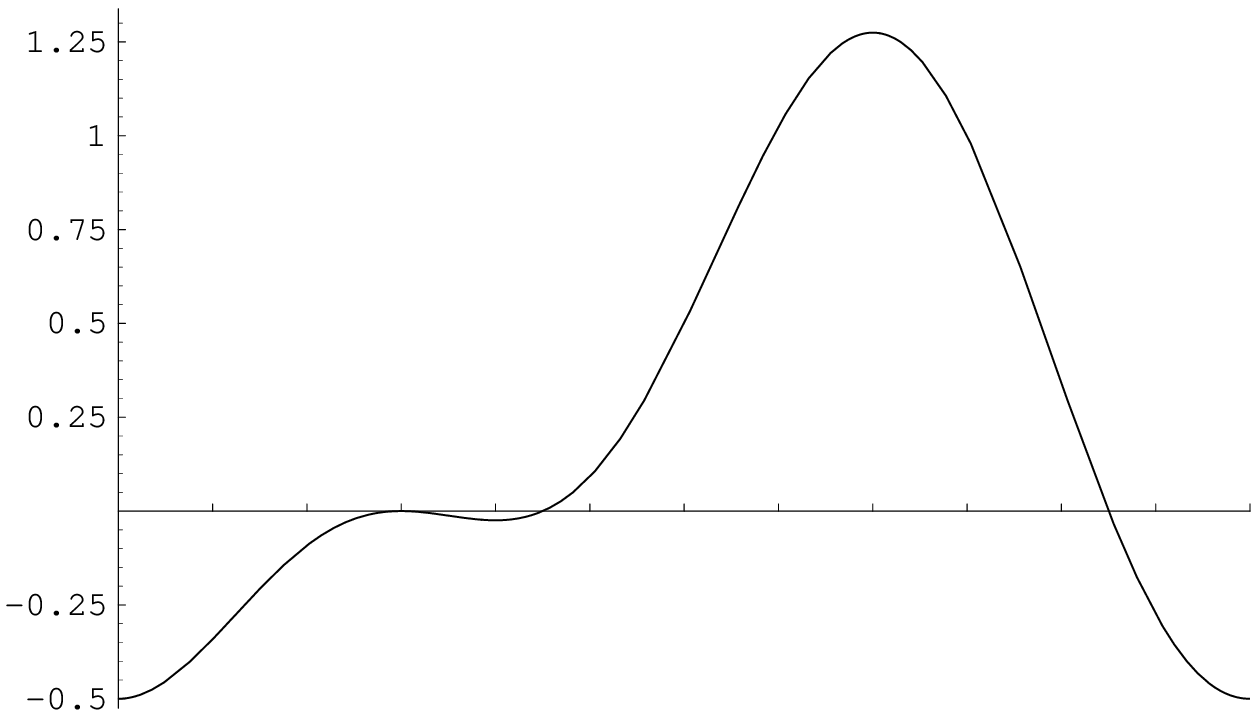}}  \put(107,86){\makebox
(0,0)[br]{$a_{0}(\theta,\frac{\pi}{4}-\theta)$}} \put(16,43){\makebox
(0,0)[r]{$\scriptscriptstyle0$}} \put(58,42){\makebox(0,0)[t]{$\scriptstyle
\frac{\pi}{4}$}} \put(98.67,42){\makebox(0,0)[t]{$\scriptstyle\frac{\pi}{2}$}}
\put(139.33,42){\makebox(0,0)[t]{$\scriptstyle\frac{3\pi}{4}$}} \put
(180,42){\makebox(0,0)[t]{$\scriptstyle\pi$}} \put(180,14){\includegraphics
[bb=0 9 360 213,height=102bp,width=180bp] {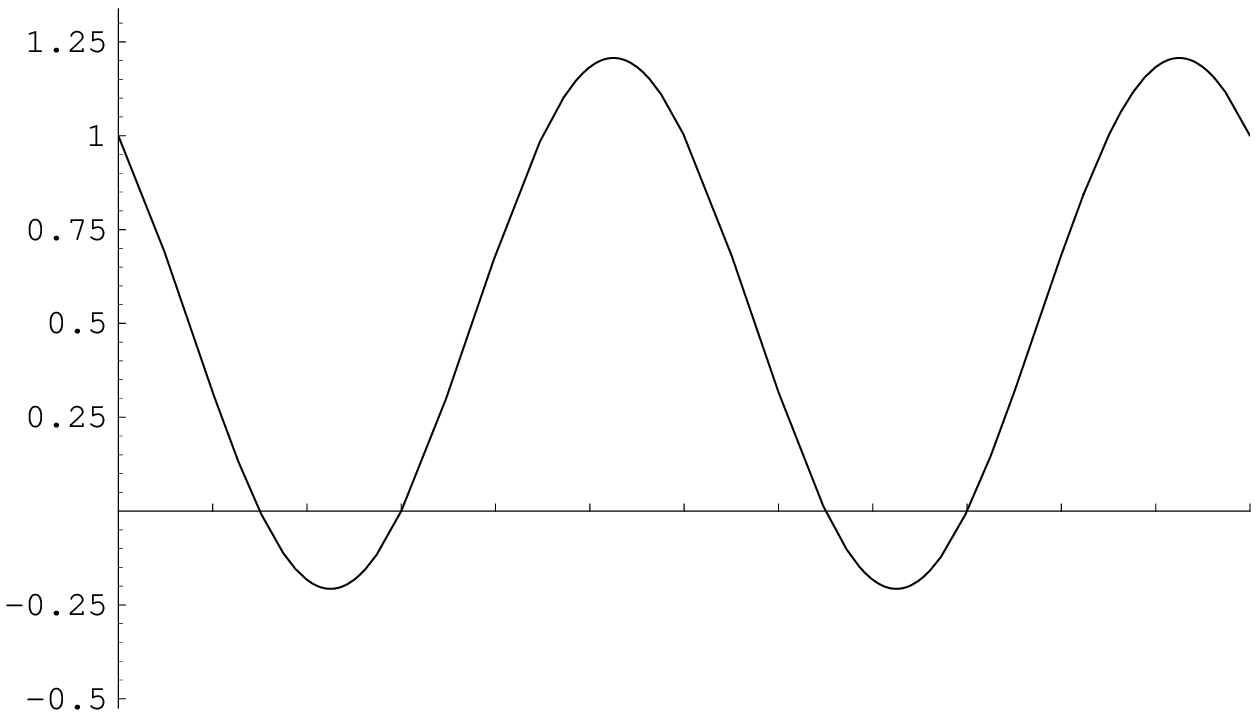}}  \put
(283,86){\makebox(0,0)[bl]{$a_{2}(\theta,\frac{\pi}{4}-\theta)$}}
\put(196,43){\makebox(0,0)[r]{$\scriptscriptstyle0$}} \put(238,42){\makebox
(0,0)[t]{$\scriptstyle\frac{\pi}{4}$}} \put(278.67,42){\makebox
(0,0)[t]{$\scriptstyle\frac{\pi}{2}$}} \put(319.33,42){\makebox
(0,0)[t]{$\scriptstyle\frac{3\pi}{4}$}} \put(360,42){\makebox
(0,0)[t]{$\scriptstyle\pi$}} \put(0,20){\makebox(180,14) {$\theta$}}
\put(180,20){\makebox(180,14) {$\theta$}}  \put(16,0){\makebox(164,14)
{Cross-section of $a_{0}$ on dashed line}}  \put(196,0){\makebox(164,14)
{Cross-section of $a_{2}$ on dashed line}}
\end{picture}
\caption[Contour plots of $a_{0}$ and $a_{2}$.]{Contour plots of $a_{0}$ and
$a_{2}$. For $a_{0}$, the three lines $\theta=\pi/4$, $\rho=0$ (or $\pi$), and
$\theta-\rho=\pi/2$ are contours of $a_{0}=0$, whose intersections are saddle
points; all the maxima and minima of $a_{0}$ are located on the dashed line
$\theta+\rho=\pi/4$ (or $5\pi/4$), displayed in the lower plot. For $a_{2}$,
contours of constant $a_{2}$ are diagonal lines of the form $\theta-\rho
=$~constant.}%
\label{FigAppContour}%
\end{figure}

Discussion of what variation $\varphi^{\left(  n\right)  }$ has between
$x_{i}$ and $x_{i+1}$ is somewhat metaphysical: it only matters that
$\varphi^{\left(  n\right)  }\left(  x_{i}\right)  $ is associated with the
point $x_{i}=i\cdot2^{-n}$. For example, if we said $\varphi^{\left(
n\right)  }\left(  x\right)  =\varphi^{\left(  n\right)  }\left(  i\cdot
2^{-n}\right)  $ for $i\cdot2^{-n}\leq x<\left(  i+\frac{1}{2}\right)
\cdot2^{-n}$ and $\varphi^{\left(  n\right)  }\left(  x\right)  =0$ for
$\left(  i+\frac{1}{2}\right)  \cdot2^{-n}\leq x<\left(  i+1\right)
\cdot2^{-n}$, then the same ``shape'' would occur in every bottom-level
interval of every cascade stage, preventing continuity from appearing in the
$n\rightarrow\infty$ limit even if it would otherwise have appeared. Other
variations within the shortest dyadic intervals at a given stage have a
similar arbitrariness; $\varphi^{\left(  n\right)  }$ is really only defined
pointwise for finite $n$. When $n$ goes to infinity, continuity arises in some
cases. On the other hand, for a given $x$ and finite $n$, $\varphi
_{\theta,\rho}^{\left(  n\right)  }\left(  x\right)  $ as a function of
$\theta$ and $\rho$ is continuous, since it is a polynomial (of order $n$) in
the $a_{i}$'s; when the order $n$ goes to infinity, singularities arise.

To decide when the corresponding scaling function $\varphi_{\theta,\rho
}\left(  x\right)  $ generates a wavelet in the \emph{strict sense} or merely
a \emph{tight frame,} as discussed in Section \ref{Int} above, we use the
method of Cohen \cite{Coh92b,CoRy95} (see also \cite{BEJ00}): We identify
cycles on $\mathbb{T}$ for the doubling map $z\mapsto z^{2}$, i.e., a finite
cyclic subset unequal to $\left\{  1\right\}  $ and invariant under $z\mapsto
z^{2}$. The result is that $\varphi\left(  x\right)  $ generates a ``strict''
wavelet if and only if%
\begin{equation}
\left\{  z\in\mathbb{T};m_{0}^{\left(  \varphi\right)  }\left(  -z\right)
=0\right\}  \label{eqAPJ.1}%
\end{equation}
does not contain a nontrivial cycle.

The cycles on $\mathbb{T}$ are not subgroups of $\mathbb{T}$ but rather
cyclic orbits on $\mathbb{T}$ under the $z\mapsto z^{2}$ action of one of the
cyclic groups $\mathbb{Z}_{k}$, $k=1,2,\dots$. Such a cyclic orbit $C_{k}$
with $k$ distinct points $z_{1},\dots,z_{k}$ must be of the form
$z_{1}\rightarrow z_{2}\rightarrow\dots\rightarrow z_{k}\rightarrow z_{1}$,
where $z_{i+1}^{{}}=z_{i}^{2}$ if $i=1,2,\dots,k-1$, and $z_{k}^{2}=z_{1}^{{}%
}$. Hence points $c$ in an orbit $C_{k}$ must satisfy $c^{2^{k}}=c$, and each
$c$ must be a $\left(  2^{k}-1\right)  $'th root of $1$. Different orbits must
be disjoint, and their union will be invariant under $z\mapsto z^{2}$ acting
on $\mathbb{T}$. The converse is not true. Note also that we can have
different $\left(  2^{k}-1\right)  $'th roots $c$ of $1$ defining different
cyclic orbits for the same $k$. If $k=1$ or $k=2$, then in each case there is
only one orbit, but if $k=3$, there are two choices. Since $m_{0}^{\left(
\theta,\rho\right)  }\left(  z\right)  $ is for each $\theta,\rho$ a
polynomial of degree at most $5$, the cardinality of a cycle contained in
(\ref{eqAPJ.1}) is at most $4.$ Thus, if $z$ is contained in such a cycle, we
must have one of the possibilities $z^{2}=z$, $z^{4}=z$, $z^{8}=z$. Hence the
cycles of length at most $3$ are $\left\{  1\right\}  $, $\left\{
\omega,\omega^{2}\right\}  $ where $\omega=e^{i2\pi/3}$, $\left\{  \zeta
,\zeta^{2},\zeta^{4}\right\}  $ where $\zeta:=e^{i2\pi/7}$, and $\left\{
\bar{\zeta},\bar{\zeta}^{2},\bar{\zeta}^{4}\right\}  =\left\{  \zeta^{6}%
,\zeta^{5},\zeta^{3}\right\}  $. But as $m_{0}\left(  -1\right)  =0$ always,
$\left(  z+1\right)  $ is always a factor of $m_{0}\left(  z\right)  $, and
since the cycle should be different from the trivial cycle $\left\{
1\right\}  $, we are reduced to the case $\left\{  \omega,\omega^{2}\right\}
$. The other cycles would make $m_{0}^{\left(  \theta,\rho\right)  }$
divisible by a polynomial of degree at least $4$.

Thus we have the following cases: $m_{0}^{\left(  \theta,\rho\right)  }\left(
z\right)  $ may be divisible by%
\begin{equation}
p_{3}\left(  z\right)  =\prod_{k=0}^{2}\left(  \omega^{k}+z\right)  =1+z^{3},
\label{eqAPJ.2}%
\end{equation}
by%
\begin{align}
p_{4}\left(  z\right)   &  =\left(  1+z\right)  \left(  \zeta+z\right)
\left(  \zeta^{2}+z\right)  \left(  \zeta^{4}+z\right) \label{eqAPJ.3}\\
&  =1+\bar{\beta}z-z^{2}+\beta z^{3}+z^{4},\nonumber
\end{align}
or by%
\begin{align}
p_{4}^{\left(  \#\right)  }\left(  z\right)  =\overline{p_{4}\left(  \bar
{z}\right)  }  &  =\left(  1+z\right)  \left(  \zeta^{3}+z\right)  \left(
\zeta^{5}+z\right)  \left(  \zeta^{6}+z\right) \label{eqAPJ.4}\\
&  =1+\beta z-z^{2}+\bar{\beta}z^{3}+z^{4},\nonumber
\end{align}
where $\zeta$ (as above) and $\beta$ are defined as%
\begin{equation}
\zeta:=e^{i\frac{2\pi}{7}},\qquad\beta:=1+\zeta+\zeta^{2}+\zeta^{4}=\frac
{1}{2}+i\frac{\sqrt{7}}{2}. \label{eqAPJ.4bis}%
\end{equation}

\begin{proposition}
\label{ProAPJ.1}There are only four cases of QMF-functions
\begin{equation}
a_{0}+a_{1}z+a_{2}z^{2}+a_{3}z^{3}+a_{4}z^{4}+a_{5}z^{5},\qquad a_{k}%
\in\mathbb{R}, \label{eqProAPJ.1.0}%
\end{equation}
which give tight frames that are not strict wavelets. In addition to the
$\mathcal{P}_{3}\left(  \mathbb{T},\mathrm{U}_{2}\left(  \mathbb{C}\right)
\right)  $-conditions, they satisfy%
\begin{equation}
\sum_{k=0}^{5}a_{k}=2. \label{eqProAPJ.1.1}%
\end{equation}
The four correspond to the three loops%
\begin{equation}
\frac{1}{\sqrt{2}}\left(
\begin{array}
[c]{cc}%
z & z^{2}\\
1 & -z
\end{array}
\right)  ,\qquad\frac{1}{\sqrt{2}}\left(
\begin{array}
[c]{cc}%
1 & z\\
z & -z^{2}%
\end{array}
\right)  ,\qquad\frac{1}{\sqrt{2}}\left(
\begin{array}
[c]{cc}%
z^{2} & 1\\
z^{2} & -1
\end{array}
\right)  . \label{eqProAPJ.1.2}%
\end{equation}
and the loop%
\begin{equation}
\frac{1}{\sqrt{2}}\left(
\begin{array}
[c]{cc}%
1 & z^{2}\\
1 & -z^{2}%
\end{array}
\right)  . \label{eqProAPJ.1.2bis}%
\end{equation}
The wavelet representation $T^{\left(  A\right)  }$ is irreducible for the
first two of the four, and reducible for the last two. The values of
$\lambda_{0}\left(  A\right)  $ are as follows: $1/2$, $1/2$, $0$, and $1$,
respectively. The first three have cycles of order $2$ and the last one a
cycle of order $4$. The corresponding system of coefficients is as follows:%
\begin{equation}%
\begin{array}
[c]{ccccccc}%
a_{0} & a_{1} & a_{2} & a_{3} & a_{4} & a_{5} & \\\cline{1-6}%
0 & 0 & 1 & 0 & 0 & 1 & \\
1 & 0 & 0 & 1 & 0 & 0 & \smash{\left.  \vphantom{%
\begin{matrix}
0\\0\\0%
\end{matrix}
}\right\rbrace }\text{ two-cycle}\\
0 & 1 & 0 & 0 & 1 & 0 & \\
1 & 0 & 0 & 0 & 0 & 1 & \text{\ four-cycle}%
\end{array}
\label{eqProAPJ.1.3}%
\end{equation}
and so all four cases are Haar wavelets. The scaling functions $\varphi\left(
x\right)  $ may be taken as in Table \textup{\ref{TableProAPJ.1.1}.}

\begin{table}[tbp]
\caption{$\varphi\left(  x\right)  $ in the four cases: Tight frames
corresponding to cycles of length two, and a four-cycle.}%
\label{TableProAPJ.1.1}
$%
\begin{array}
[c]{ll}%
\begin{array}
[c]{cc}%
\setlength{\unitlength}{1bp}%
\begin{picture}
(265,129)(-9,0) \put(0,0){\includegraphics[bb=8 0 337 164,
height=123bp,width=247bp] {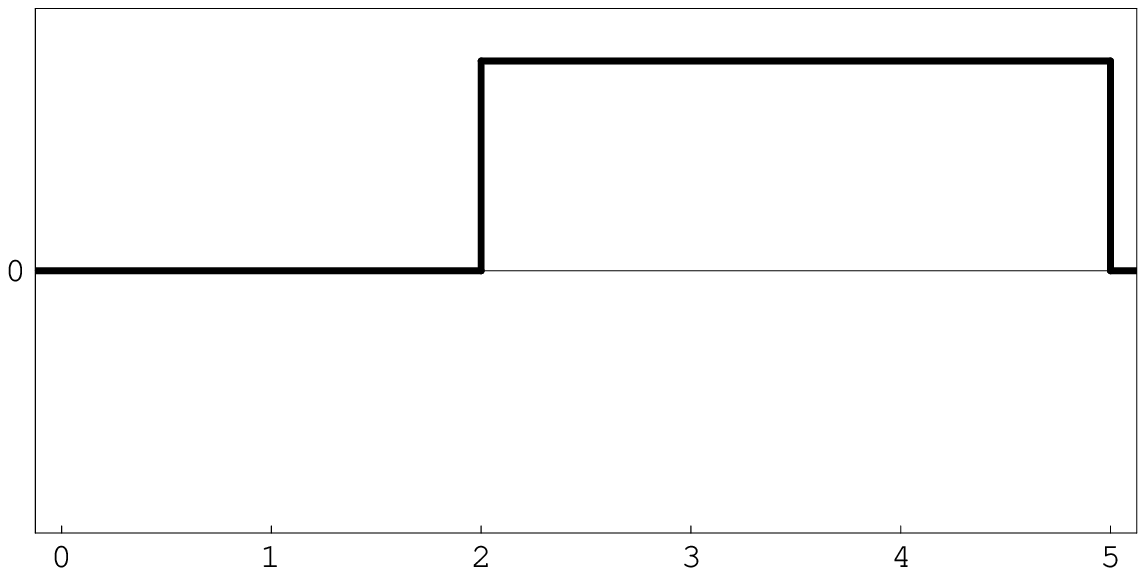}} \put(248,2){\makebox(0,12)[l]{$x$}}
\end{picture}
& \hspace{-16pt}%
\end{array}
&
\begin{array}
[c]{l}%
\theta=\pi/4,\\
\rho=\pi/2\mathpunct{;}%
\end{array}
\\%
\begin{array}
[c]{cc}%
\setlength{\unitlength}{1bp}%
\begin{picture}
(265,129)(-9,0) \put(0,0){\includegraphics[bb=8 0 337 164,
height=123bp,width=247bp] {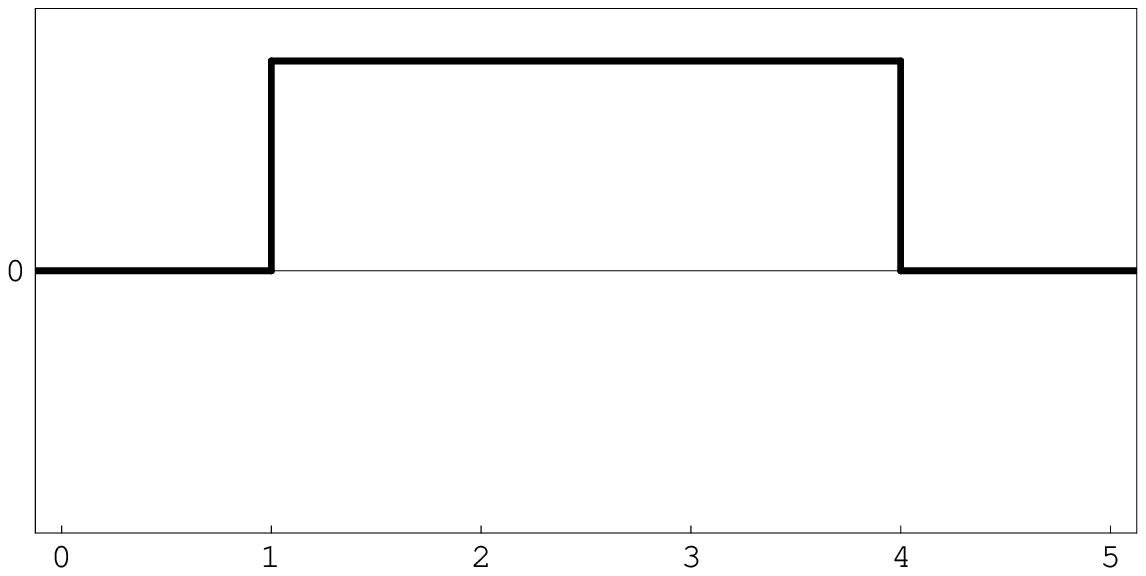}} \put(248,2){\makebox(0,12)[l]{$x$}}
\end{picture}
& \hspace{-16pt}%
\end{array}
&
\begin{array}
[c]{l}%
\theta=0\text{ \textup{(}or }\pi\text{\textup{),} }\\
\rho=0\text{ \textup{(}or }\pi\text{\textup{);}}%
\end{array}
\\%
\begin{array}
[c]{cc}%
\setlength{\unitlength}{1bp}%
\begin{picture}
(265,129)(-9,0) \put(0,0){\includegraphics[bb=8 0 337 164,
height=123bp,width=247bp] {scal03.eps}} \put(248,2){\makebox(0,12)[l]{$x$}}
\end{picture}
& \hspace{-16pt}%
\end{array}
&
\begin{array}
[c]{l}%
\theta=3\pi/4,\\
\rho=\pi/2\mathpunct{;}%
\end{array}
\\%
\begin{array}
[c]{cc}%
\setlength{\unitlength}{1bp}%
\begin{picture}
(265,129)(-9,0) \put(0,0){\includegraphics[bb=8 0 337 164,
height=123bp,width=247bp] {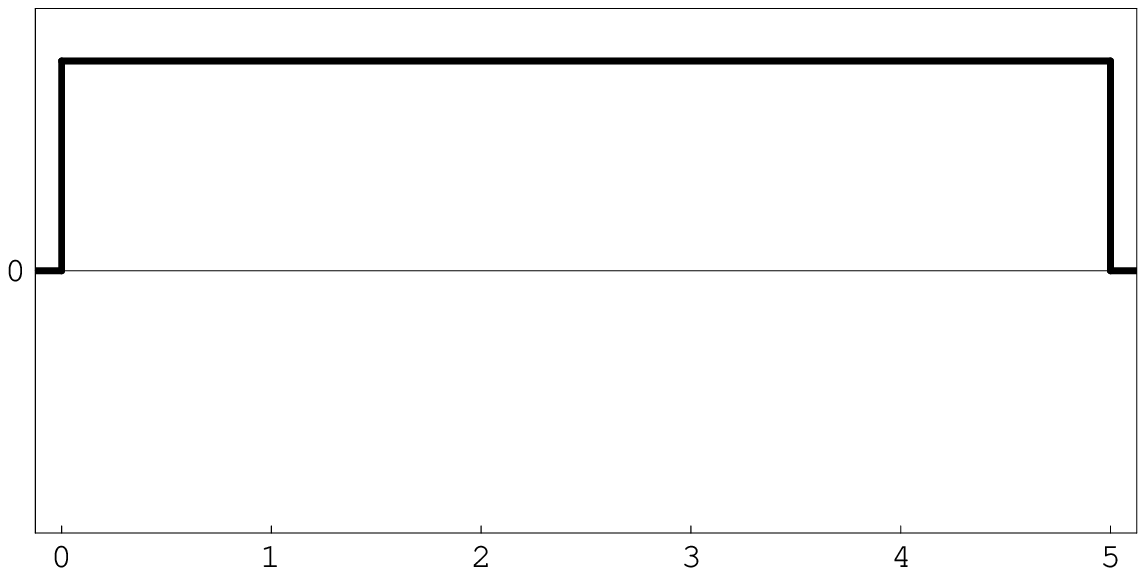}} \put(248,2){\makebox(0,12)[l]{$x$}}
\end{picture}
& \hspace{-16pt}%
\end{array}
&
\begin{array}
[c]{l}%
\\
\theta=\rho=\pi/2.
\end{array}
\end{array}
$\end{table}
\end{proposition}

\begin{proof}
If $m_{0}\left(  z\right)  $ is divisible by $1+z^{3}$, then its six
coefficients $a_{0},a_{1},\dots,a_{5}$ must be of the form $c_{0},c_{1}%
,c_{2},c_{0},c_{1},c_{2}$, and the associated loop $\mathbb{T}\rightarrow
\mathrm{U}_{2}\left(  \mathbb{C}\right)  $,%
\begin{equation}
\left(
\begin{array}
[c]{cc}%
c_{0}+c_{2}z+c_{1}z^{2} & c_{1}+c_{0}z+c_{2}z^{2}\\
\bar{c}_{2}+\bar{c}_{0}z+\bar{c}_{1}z^{2} & -\left(  \bar{c}_{1}+\bar{c}%
_{2}z+\bar{c}_{0}z^{2}\right)
\end{array}
\right)  . \label{eqProAPJ.1.4}%
\end{equation}
The corresponding $\mathrm{U}_{2}\left(  \mathbb{C}\right)  $-conditions then
yield:%
\begin{equation}
2\bar{c}_{0}c_{2}+\bar{c}_{1}c_{0}+\bar{c}_{2}c_{1}=0,\qquad\bar{c}_{0}%
c_{1}+\bar{c}_{1}c_{2}=0. \label{eqProAPJ.1.5}%
\end{equation}
Substitution of the second into the first yields $\bar{c}_{0}c_{2}=0$. Hence,
of the three numbers $c_{0},c_{1},c_{2}$, at most one, and therefore precisely
one, can be nonzero. But each of the three cases is determined up to scale,
and condition (\ref{eqProAPJ.1.1}) decides the scale. We are therefore led to
the three loops in (\ref{eqProAPJ.1.2}), and the rule (\ref{eqExp.8}) then
gives the three scaling functions $\varphi\left(  x\right)  $ which are listed
in (\ref{eqProAPJ.1.3}) and Table \ref{TableProAPJ.1.1}.

The cycle of the last line in Table \ref{TableProAPJ.1.1} is of order $4$. Let
$\mu:=e^{i2\pi/5}=\lambda^{3}$ ($\lambda:=e^{i2\pi/15}$). Then the cycle is
$\left\{  \mu,\mu^{2},\mu^{4},\mu^{3}\right\}  $, and $\prod_{k=0}^{4}\left(
\mu^{k}+z\right)  =z^{5}+1$, which is the $m_{0}\left(  z\right)  $ for the
last line of Table \ref{TableProAPJ.1.1}. (It is from a root of $1$ of order
$2^{l}-1$ ($=15$) for $l=4$.)

The other length-$3$ loops which would be possible are, as noted, $\left\{
\zeta,\zeta^{2},\zeta^{4}\right\}  $ and $\left\{  \bar{\zeta},\bar{\zeta}%
^{2},\bar{\zeta}^{4}\right\}  $, with $\zeta=e^{i2\pi/7}$. We claim that they
do \emph{not} in fact occur.

If one of them did occur, then the corresponding $m_{0}\left(  z\right)  $
would be divisible by either $p_{4}\left(  z\right)  $, or by $p_{4}^{\left(
\#\right)  }\left(  z\right)  $. But $p_{4}\left(  1\right)  =p_{4}^{\left(
\#\right)  }\left(  1\right)  =2$, so the factorization would be $m_{0}\left(
z\right)  =p_{4}\left(  z\right)  l\left(  z\right)  $ where $l\left(
z\right)  =a+\left(  1-a\right)  z$. (We have picked the normalization of
$m_{0}\left(  z\right)  $ given by $m_{0}\left(  1\right)  =2$ for
convenience.) From the formulas (\ref{eqA.8bis}) and (\ref{eqA.11}) we note
that the coefficients $a_{0},a_{1},\dots,a_{5}$ are real. Divisibility by
$p_{4}\left(  z\right)  $ means that $-1$, $-\zeta$, $-\zeta^{2}$, and
$-\zeta^{4}$ are roots of $m_{0}\left(  z\right)  $. So the complex conjugates
$-\bar{\zeta}$, $-\bar{\zeta}^{2}$, $-\bar{\zeta}^{4}$ are also roots. But
that would give us all seven points, $-1,-\zeta,-\zeta^{2},-\zeta^{3}%
,-\zeta^{4},-\zeta^{5},-\zeta^{6}$, as distinct roots of $m_{0}\left(
z\right)  $, which is impossible since $m_{0}$ is of degree at most $5$.
\end{proof}

In conclusion, when $g=3$, the variety of the wavelets which are only tight
frames sits on a finite subset of the full variety of all $\mathcal{P}%
_{3}\left(  \mathbb{T},\mathrm{U}_{2}\left(  \mathbb{C}\right)  \right)  $ examples.\medskip

\newcommand{\caplayout}{Layout of scaling function plots.
\emph{Round solid points}
(\raisebox{-1.75pt}{\huge$\bullet$}):
approximate locations of ``ultra-smooth'' wavelet scaling functions
in relation to scaling functions plotted.
\emph{Shading}: divergence due to $a_{0}>1$ or $a_{5}>1$.
\emph{Boxes}: marginal divergence due to $a_{0}=1$ or $a_{5}=1$.}
\begin{table}[tbp]
\caption{\protect\caplayout}%
\label{layout}
\newcommand{\laystrut}{\vphantom{abcdefghijkl}}
\begin{tabular}
[c]{ccccc}%
\phantom{$\rho=11\pi/12$} &
\begin{tabular}
[c]{|lll|}\hline
\llap{$\rho=11\pi/12\qquad$}al & bl & cl\laystrut\\
\llap{$5\pi/6\qquad$}ak & bk & ck\laystrut\\
\llap{$3\pi/4\qquad$}aj & bj & cj\laystrut\\
\llap{$2\pi/3\qquad$}ai & bi & ci\laystrut\\\hline
\end{tabular}
&
\begin{tabular}
[c]{|lll|}\hline
dl & el & fl\laystrut\\
dk & ek & fk\laystrut\\
dj & ej & fj\laystrut\\
di & ei & fi\laystrut\\\hline
\end{tabular}
&
\begin{tabular}
[c]{|lll|}\hline
gl & hl & il\laystrut\\
gk & hk & ik\rlap{\kern6pt\raisebox{6pt}[0pt][0pt]{\huge$\bullet$}}\laystrut\\
gj & hj & ij\laystrut\\
gi & \rlap{\kern-1.5pt\raisebox{-3pt}[0pt][0pt]{\includegraphics
[bb=0 0 15 13,height=11bp,width=13bp]{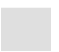}}}hi
 & \rlap{\kern-1.5pt\raisebox{-3pt}[0pt][0pt]{\includegraphics
[bb=0 0 15 13,height=11bp,width=13bp]{shade.eps}}}ii\laystrut\\\hline
\end{tabular}
&
\begin{tabular}
[c]{|lll|}\hline
jl & kl & ll\laystrut\\
jk & kk & lk\laystrut\\
\rlap{\kern-1.5pt\raisebox{0.25pt}[0pt][0pt]{\framebox
[13bp]{\rule{0bp}{4.5bp}}}}jj & kj & lj\laystrut\\
\rlap{\kern-1.5pt\raisebox{-3pt}[0pt][0pt]{\includegraphics
[bb=0 0 15 13,height=11bp,width=13bp]{shade.eps}}}ji
 & \rlap{\kern-1.5pt\raisebox{-3pt}[0pt][0pt]{\includegraphics
[bb=0 0 15 13,height=11bp,width=13bp]{shade.eps}}}ki
 & li\laystrut\\\hline
\end{tabular}
\\
& {\small p.~\pageref{P9}} & {\small p.~\pageref{P10}} 
& {\small p.~\pageref{P11}} & {\small p.~\pageref{P12}}\\
&
\begin{tabular}
[c]{|lll|}\hline
\llap{$7\pi/12\qquad$}ah & bh & ch\laystrut\\
\llap{$\pi/2\qquad$}ag & bg & cg\laystrut\\
\llap{$5\pi/12\qquad$}af & bf & cf\laystrut\\
\llap{$\pi/3\qquad$}ae & be & \rlap{\kern-1.5pt\raisebox{-3pt}[0pt][0pt]{\includegraphics
[bb=0 0 15 13,height=11bp,width=13bp]{shade.eps}}}ce\laystrut\\\hline
\end{tabular}
&
\begin{tabular}
[c]{|lll|}\hline
dh & eh & \rlap{\kern-1.5pt\raisebox{-3pt}[0pt][0pt]{\includegraphics
[bb=0 0 15 13,height=11bp,width=13bp]{shade.eps}}}fh\laystrut\\
\rlap{\kern-1.5pt\raisebox{0.25pt}[0pt][0pt]{\framebox
[13bp]{\rule{0bp}{4.5bp}}}}dg
 & \rlap{\kern-1.5pt\raisebox{-3pt}[0pt][0pt]{\includegraphics
[bb=0 0 15 13,height=11bp,width=13bp]{shade.eps}}}eg
 & \rlap{\kern-1.5pt\raisebox{-3pt}[0pt][0pt]{\includegraphics
[bb=0 0 15 13,height=11bp,width=13bp]{shade.eps}}}fg\laystrut\\
\rlap{\kern-1.5pt\raisebox{-3pt}[0pt][0pt]{\includegraphics
[bb=0 0 15 13,height=11bp,width=13bp]{shade.eps}}}df
 & \rlap{\kern-1.5pt\raisebox{-3pt}[0pt][0pt]{\includegraphics
[bb=0 0 15 13,height=11bp,width=13bp]{shade.eps}}}ef
 & \rlap{\kern-1.5pt\raisebox{-3pt}[0pt][0pt]{\includegraphics
[bb=0 0 15 13,height=11bp,width=13bp]{shade.eps}}}ff\laystrut\\
\rlap{\kern-1.5pt\raisebox{-3pt}[0pt][0pt]{\includegraphics
[bb=0 0 15 13,height=11bp,width=13bp]{shade.eps}}}de
 & \rlap{\kern-1.5pt\raisebox{-3pt}[0pt][0pt]{\includegraphics
[bb=0 0 15 13,height=11bp,width=13bp]{shade.eps}}}ee
 & \rlap{\kern-1.5pt\raisebox{-3pt}[0pt][0pt]{\includegraphics
[bb=0 0 15 13,height=11bp,width=13bp]{shade.eps}}}fe\laystrut\\\hline
\end{tabular}
&
\begin{tabular}
[c]{|lll|}\hline
gh & \rlap{\kern-1.5pt\raisebox{-3pt}[0pt][0pt]{\includegraphics
[bb=0 0 15 13,height=11bp,width=13bp]{shade.eps}}}hh
 & \rlap{\kern-1.5pt\raisebox{-3pt}[0pt][0pt]{\includegraphics
[bb=0 0 15 13,height=11bp,width=13bp]{shade.eps}}}ih\laystrut\\
\rlap{\kern-1.5pt\raisebox{0.25pt}[0pt][0pt]{\framebox
[13bp]{\rule{0bp}{4.5bp}}}}gg
 & \rlap{\kern-1.5pt\raisebox{-3pt}[0pt][0pt]{\includegraphics
[bb=0 0 15 13,height=11bp,width=13bp]{shade.eps}}}hg
 & \rlap{\kern-1.5pt\raisebox{-3pt}[0pt][0pt]{\includegraphics
[bb=0 0 15 13,height=11bp,width=13bp]{shade.eps}}}ig\laystrut\\
gf & \rlap{\kern-1.5pt\raisebox{-3pt}[0pt][0pt]{\includegraphics
[bb=0 0 15 13,height=11bp,width=13bp]{shade.eps}}}hf
 & if\laystrut\\
ge & he & ie\laystrut\\\hline
\end{tabular}
&
\begin{tabular}
[c]{|lll|}\hline
\rlap{\kern-1.5pt\raisebox{-3pt}[0pt][0pt]{\includegraphics
[bb=0 0 15 13,height=11bp,width=13bp]{shade.eps}}}jh
 & kh & lh\laystrut\\
\rlap{\kern-1.5pt\raisebox{0.25pt}[0pt][0pt]{\framebox
[13bp]{\rule{0bp}{4.5bp}}}}jg
 & kg & lg\laystrut\\
jf & kf & lf\laystrut\\
je & ke & le\laystrut\\\hline
\end{tabular}
\\
& {\small p.~\pageref{P5}} & {\small p.~\pageref{P6}} &
{\small p.~\pageref{P7}} & {\small p.~\pageref{P8}}\\
&
\begin{tabular}
[c]{|lll|}\hline
\llap{$\pi/4\qquad$}ad & bd & cd\laystrut\\
\llap{$\pi/6\qquad$}ac & bc & cc\laystrut\\
\llap{$\pi/12\qquad$}ab & bb & cb\laystrut\\
\llap{$\rho=0\qquad$}\raisebox{-25pt}[0pt][0pt]{\llap{$\theta=\;$}\rlap{$0$}%
}aa & \raisebox{-25pt}[0pt][0pt]{\rlap{$\frac{\pi}{12}$}%
}ba & \raisebox{-25pt}[0pt][0pt]{\rlap{$\frac{\pi}{6}$}}ca\laystrut\\\hline
\end{tabular}
&
\begin{tabular}
[c]{|lll|}\hline
\rlap{\kern-1.5pt\raisebox{0.25pt}[0pt][0pt]{\framebox
[13bp]{\rule{0bp}{4.5bp}}}}dd
 & ed & fd\laystrut\\
dc & ec & fc\laystrut\\
db\rlap{\kern0.5pt\raisebox{5pt}[0pt][0pt]{\huge$\bullet$}}
 & eb & ec\laystrut\\
\raisebox{-25pt}[0pt][0pt]{\rlap{$\frac{\pi}{4}$}}da
& \raisebox{-25pt}[0pt][0pt]{\rlap{$\frac{\pi}{3}$}}ea
 & \raisebox{-25pt}[0pt][0pt]{\rlap{$\frac{5\pi}{12}$}}fa\laystrut\\\hline
\end{tabular}
&
\begin{tabular}
[c]{|lll|}\hline
gd & hd & id\laystrut\\
gc & hc & ic\laystrut\\
gb & hb & ib\laystrut\\
\raisebox{-25pt}[0pt][0pt]{\rlap{$\frac{\pi}{2}$}}ga & \raisebox
{-25pt}[0pt][0pt]{\rlap{$\frac{7\pi}{12}$}}ha
& \raisebox{-25pt}[0pt][0pt]{\rlap{$\frac{2\pi}{3}$}}ia\laystrut\\\hline
\end{tabular}
&
\begin{tabular}
[c]{|lll|}\hline
jd & kd & ld\laystrut\\
jc & kc & lc\laystrut\\
jb & kb & lb\laystrut\\
\raisebox
{-25pt}[0pt][0pt]{\rlap{$\frac{3\pi}{4}$}}ja
& \raisebox{-25pt}[0pt][0pt]{\rlap{$\frac{5\pi}{6}$}}ka & \raisebox
{-25pt}[0pt][0pt]{\rlap{$\frac{11\pi}{12}$}}la\laystrut\\\hline
\end{tabular}
\\
& {\small p.~\pageref{P1}} & {\small p.~\pageref{P2}} &
{\small p.~\pageref{P3}} & {\small p.~\pageref{P4}}\\
& $\rule{0pt}{14pt}$ &  &  &
\end{tabular}
\end{table}

On the following pages are plots, for various values of the angles $\theta$
and $\rho$, of the wavelet scaling functions $\varphi_{\theta,\rho}\left(
x\right)  $ at the $8$th cascade level, computed by the local algorithm
described in (\ref{eqA.3}) above.
The layout of the plots is shown in the chart in Table \ref{layout}.
The plots for $\theta=\pi$ or $\rho=\pi$ beyond
the top and right of this
chart are the same as those for $\theta=0$ or $\rho=0$,
because of the periodicity (\ref{eqA.13per}).

The ``ultra-smooth'' scaling function with $m_{0}^{\left(  \theta,\rho\right)
}\left(  z\right)  $ divisible by $\left(  1+z\right)  ^{3}$, shown in Figure
\ref{FigMomThree} above, 
and its counterpart under the symmetry (\ref{eqA.13}),
lie at the positions shown by a round solid point
(\raisebox{-1.75pt}{\huge$\bullet$})
in both Table \ref{layout} and Figure \ref{FigSmooth}.

The scaling functions from the $g=2$ family in \cite{BrJo99b} appear as
subsets of the $g=3$ family here, supported on various subintervals of
$\left[  0,5\right]  $ of length $3$. The correspondence results, for
particular values of $\left(  \theta,\rho\right)  $, from the vanishing of two
of the $a_{i}$ coefficients, and the equality of the other four $a_{i}$'s to
the four coefficients of the $g=2$ family. The values of $\left(  \theta
,\rho\right)  $ corresponding to continuous scaling functions in the $g=2$
family \cite{CoHe92,CoHe94,Wan95,Wan96,DaLa92}
(see \cite[Remark 3.1]{BrJo99b}) are indicated in Table \ref{embed}.
The values of $\left(  \theta,\rho\right)  $ that give these known
continuous scaling functions are indicated graphically in Figure
\ref{FigSmooth}, along with the vanishing-moment points where the polynomial
$m_{0}^{\left(  \theta,\rho\right)  }\left(  z\right)  $ is divisible by
$\left(  1+z\right)  ^{2}$ and by $\left(  1+z\right)  ^{3}$ (see Proposition
\ref{ProMom.2}), and the tight-frame cases (see Proposition \ref{ProAPJ.1}).
Some regions of the $\left(  \theta,\rho\right)  $ plane where the cascade
approximants do not converge to a continuous scaling function
are also indicated in the same figure.

\begin{table}[tbp]
\caption{Embedding of the $g=2$ family in the $g=3$ family.}%
\label{embed}
\begin{tabular}
[c]{rll}%
\textbf{Support interval:} & $\left(  \theta,\rho\right)  $\textbf{ values:} &
\textbf{Continuous }$\varphi_{\theta,\rho}^{{}}\left(  x\right)  $\textbf{
at:}\\
$x\in\left[  0,3\right]  $ & $\left\{  \left(  \theta,\rho\right)
:\theta=3\pi/4\right\}  $ & $\theta=3\pi/4,\;\rho\in\left(  0,\pi/4\right)
\cup\left(  3\pi/4,\pi\right)  $\\
$x\in\left[  1,4\right]  $ & $\left\{  \left(  \theta,\rho\right)
:\rho=0\right\}  $ & $\rho=0,\;\theta\in\left(  \pi/4,\pi/2\right)
\cup\left(  \pi/2,3\pi/4\right)  $\\
$x\in\left[  2,5\right]  $ & $\left\{  \left(  \theta,\rho\right)  :\theta
=\pi/4\right\}  $ & $\theta=\pi/4,\;\rho\in\left(  0,\pi/4\right)  \cup\left(
3\pi/4,\pi\right)  $%
\end{tabular}
\end{table}

\newcommand{\capFigSmooth}
{\emph{Thin curved lines} (in the four corners): vanishing first
moment of $\psi$ (Proposition \ref{ProMom.2}(a)). \emph{Round solid points}
(\raisebox{-1.75pt}{\huge$\bullet$}):
vanishing second moment of $\psi$;
``ultra-smooth'' wavelet scaling function (Proposition \ref{ProMom.2}(b),
Figure \ref{FigMomThree}). \emph{Thick
straight lines}: embedding of (continuous portion of) $g=2$ family in $g=3$
family
(Table \ref{embed}).
\emph{Round open points}
(\raisebox{-1.75pt}{\huge$\circ$}):
translated Haar functions within $g=2$ family
(plots ``da'', ``dd'', ``dj'', ``ga'', ``ja'', ``jd'', ``jj'').
\emph{Square points}: tight frames (Proposition \ref{ProAPJ.1},
plots ``aa'', ``dg'', ``gg'', ``jg'').
\emph{Shading}: divergence due to $a_{0}>1$ or $a_{5}>1$.}

\begin{figure}[tbp]
\setlength{\unitlength}{1.47239bp}
\begin{picture}
(197,180)(0,158)  \put(21.25,178.25){\includegraphics[bb=0 0 360
360,height=227.75bp,width=227.75bp] {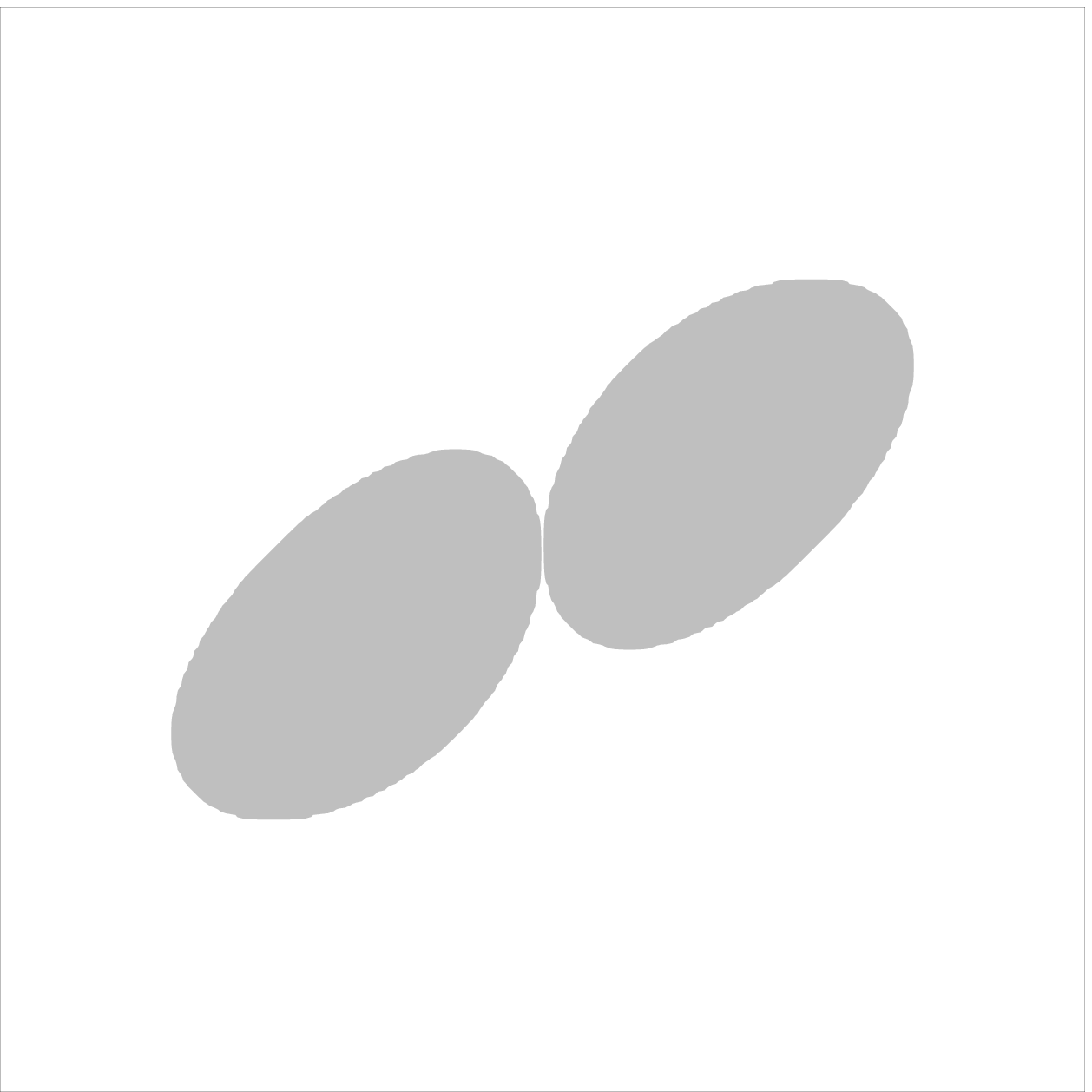}} \put
(17,174){\includegraphics[bb=0 0 360 360,height=240bp,width=240bp] {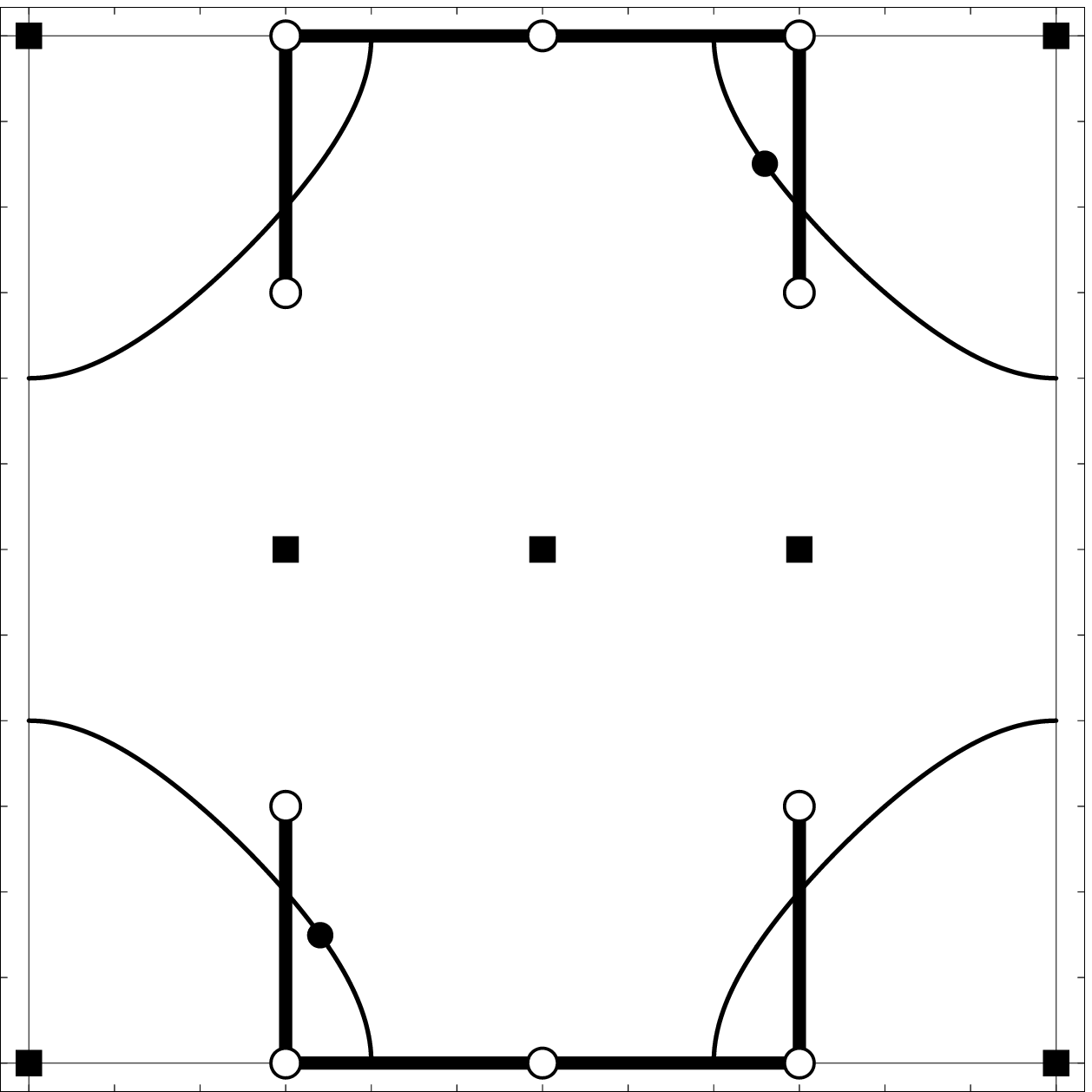}}
\put(21,172){\makebox(0,0)[t]{$0$ }}   \put(60,172){\makebox(0,0)[t]{$\frac
{\pi}{4}$}}   \put(98.67,172){\makebox(0,0)[t]{$\frac{\pi}{2}$}}
\put(137.33,172){\makebox(0,0)[t]{$\frac{3\pi}{4}$}}   \put(176,172){\makebox
(0,0)[t] {$\pi$}} \put(16,179){\makebox(0,0)[r] {$0$}}   \put
(16,216.33){\makebox(0,0)[r] {$\frac{\pi}{4}$}}   \put(16,255){\makebox
(0,0)[r] {$\frac{\pi}{2}$}} \put(3,234){\makebox(0,0)[l] {$\rho$}}
\put(16,294.33){\makebox(0,0)[r] {$\frac{3\pi}{4}$}}   \put(16,333){\makebox
(0,0)[tr] {$\pi$}}   \put(0,156){\makebox(180,10) {$\theta\qquad$}}
\end{picture}
\caption{\protect\capFigSmooth}%
\label{FigSmooth}%
\end{figure}

Putting all the $144$ pictures
together as illustrated in Table \ref{layout},
we get graphic support for
the observation that the two spin-vectors
in the factorization (\ref{eqA.8}) produce
more smoothness of $x\mapsto\varphi_{\theta ,\rho }\left( x\right) $ when they
are not aligned, i.e., off the diagonal
$\theta =\rho $. It also shows that,
close to one of the true Haar
wavelets, i.e., when $\varphi $ is the indicator
function of some $\left[ k,k+1\right) $, there
is a continuous $\varphi $, while close
to a mock Haar wavelet (i.e.,
one that is only a tight frame)
it appears that the graph of the
scaling functions have Hausdorff
dimension $>1$, hence the
``fractal'' appearance.

\begin{acknowledgements}
We are very grateful to Ola Bratteli, Ken Davidson, and David Kribs for
enlightening discussions, and to Brian Treadway and Cymie Wehr for expert
typesetting, and graphics artwork.
\end{acknowledgements}

\newpage

\begin{figure}[tbp]
\setlength{\unitlength}{1bp}
\begin{picture}(360,524)
\put(0,405){\includegraphics[%
height=119bp,width=120bp]{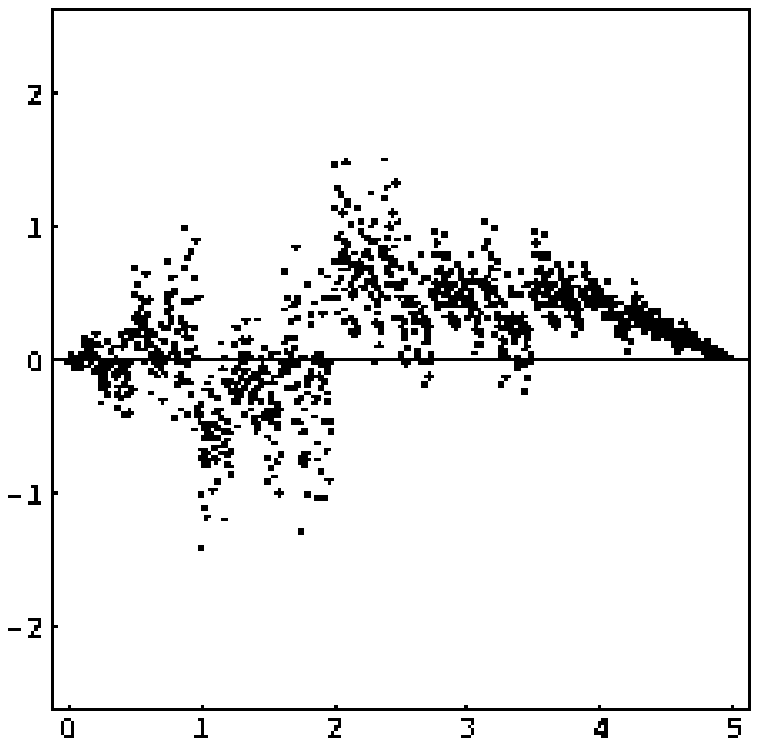}}
\put(120,405){\includegraphics[%
height=119bp,width=120bp]{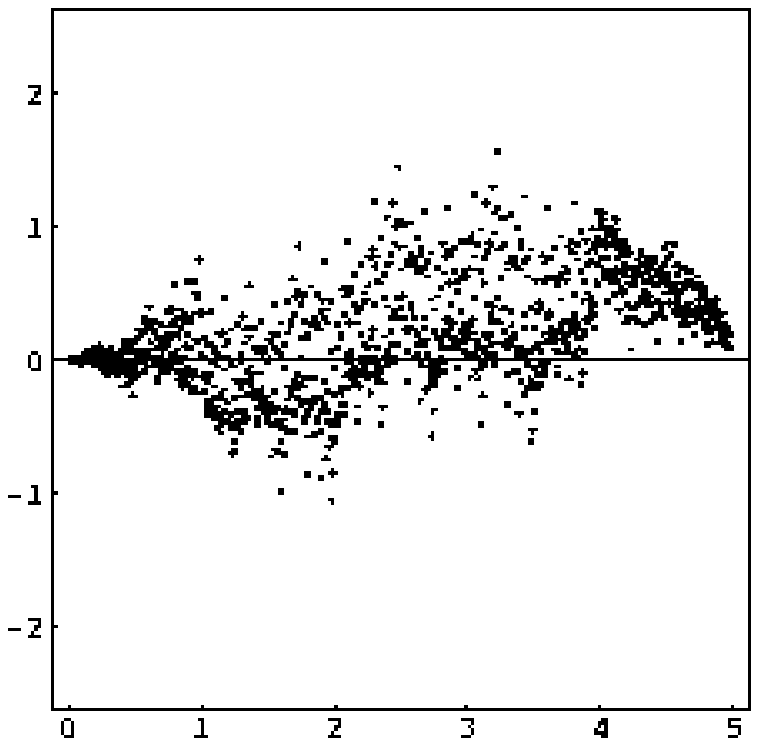}}
\put(240,405){\includegraphics[%
height=119bp,width=120bp]{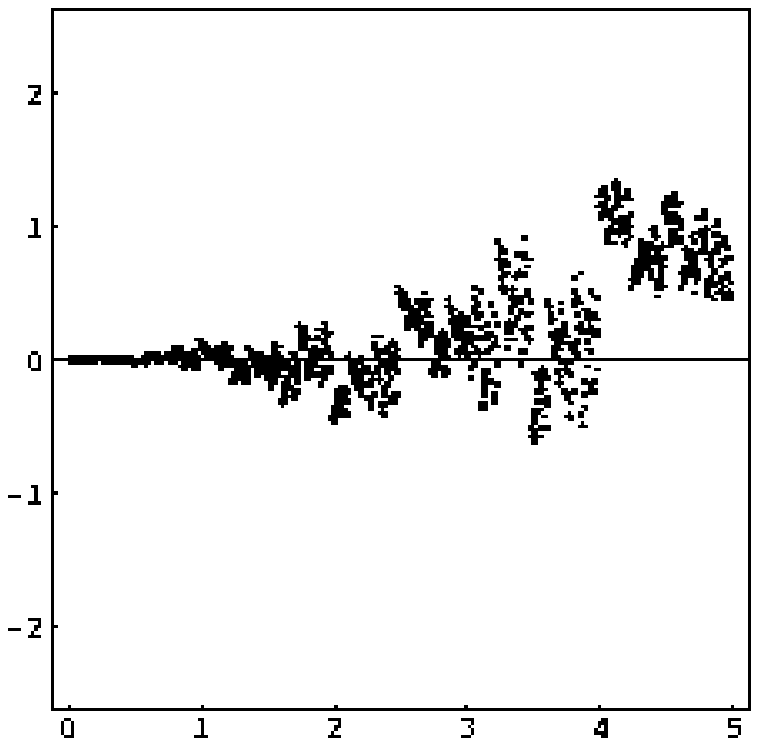}}
\put(0,393){\makebox(120,12){ad: $\theta=0,\;\rho=\pi/4$}}
\put(120,393){\makebox(120,12){bd: $\theta=\pi/12,\;\rho=\pi/4$}}
\put(240,393){\makebox(120,12){cd: $\theta=\pi/6,\;\rho=\pi/4$}}
\put(0,274){\includegraphics[%
height=119bp,width=120bp]{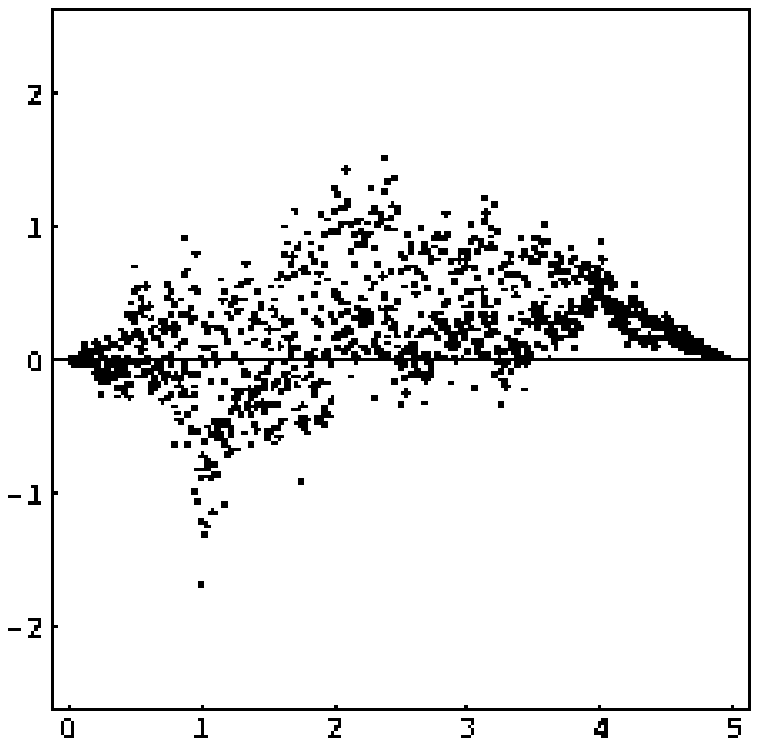}}
\put(120,274){\includegraphics[%
height=119bp,width=120bp]{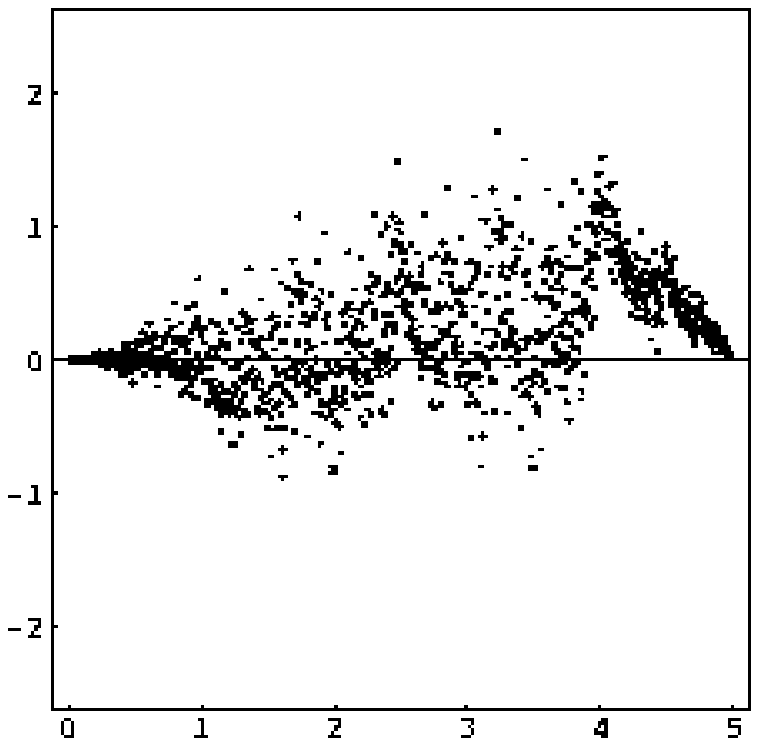}}
\put(240,274){\includegraphics[%
height=119bp,width=120bp]{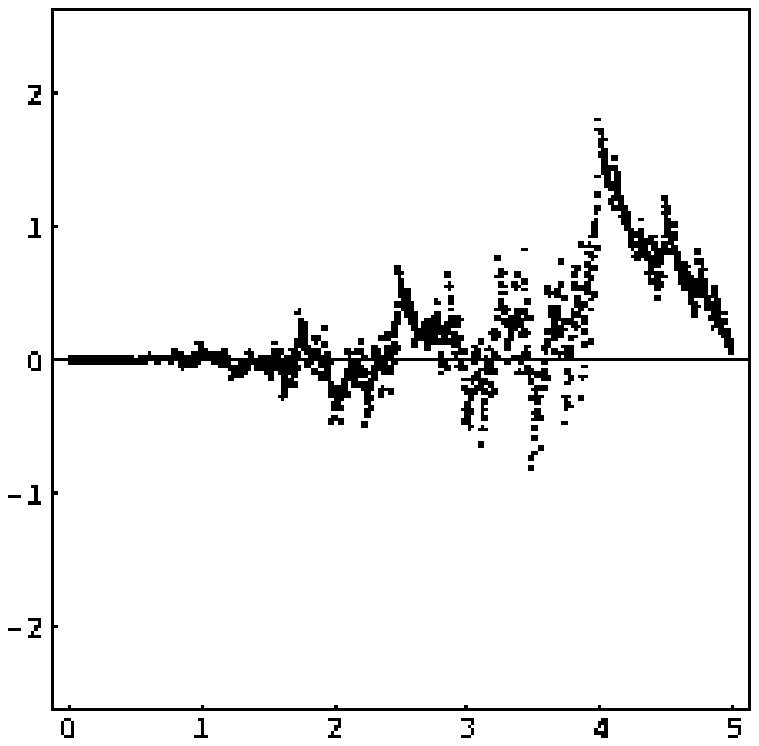}}
\put(0,262){\makebox(120,12){ac: $\theta=0,\;\rho=\pi/6$}}
\put(120,262){\makebox(120,12){bc: $\theta=\pi/12,\;\rho=\pi/6$}}
\put(240,262){\makebox(120,12){cc: $\theta=\pi/6,\;\rho=\pi/6$}}
\put(0,143){\includegraphics[%
height=119bp,width=120bp]{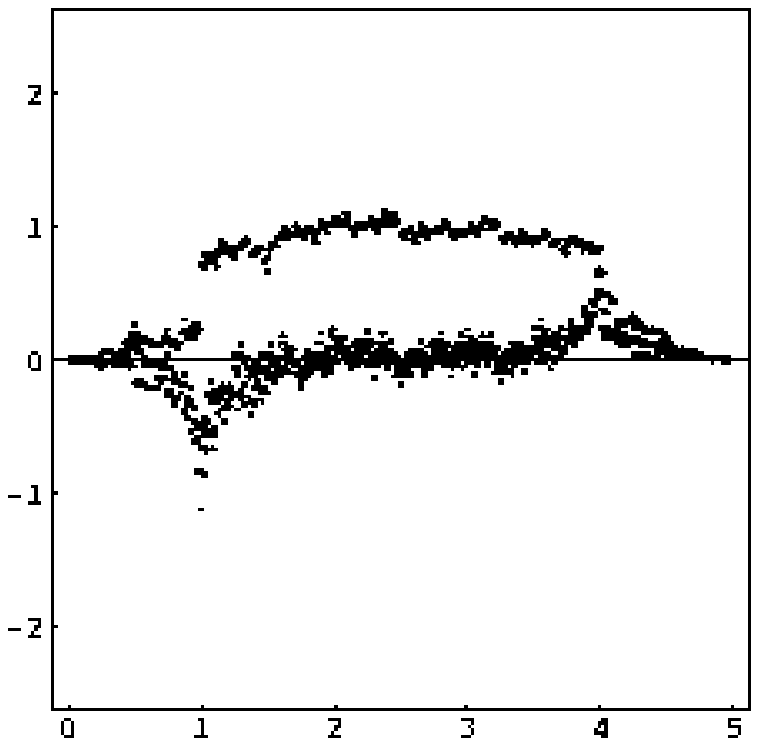}}
\put(120,143){\includegraphics[%
height=119bp,width=120bp]{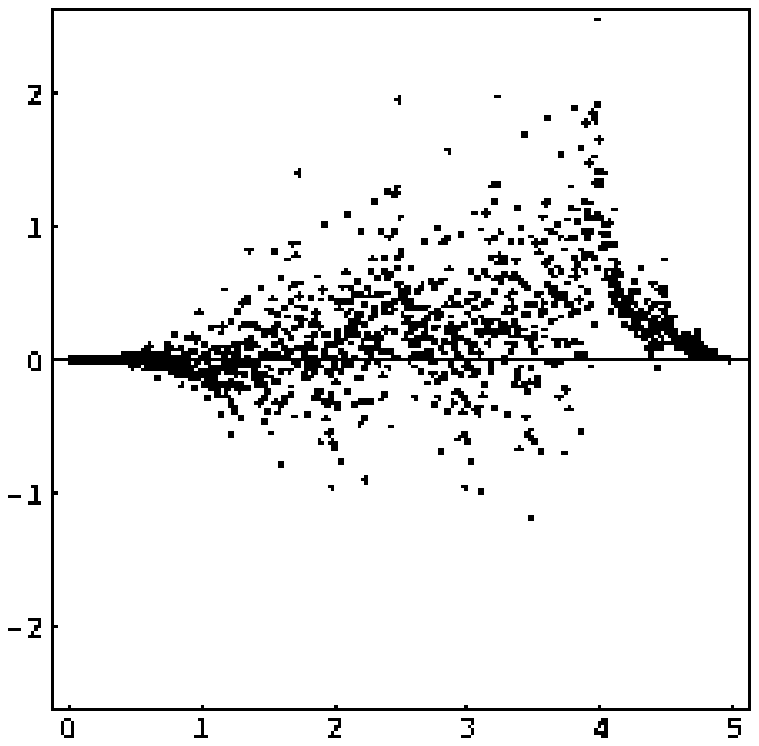}}
\put(240,143){\includegraphics[%
height=119bp,width=120bp]{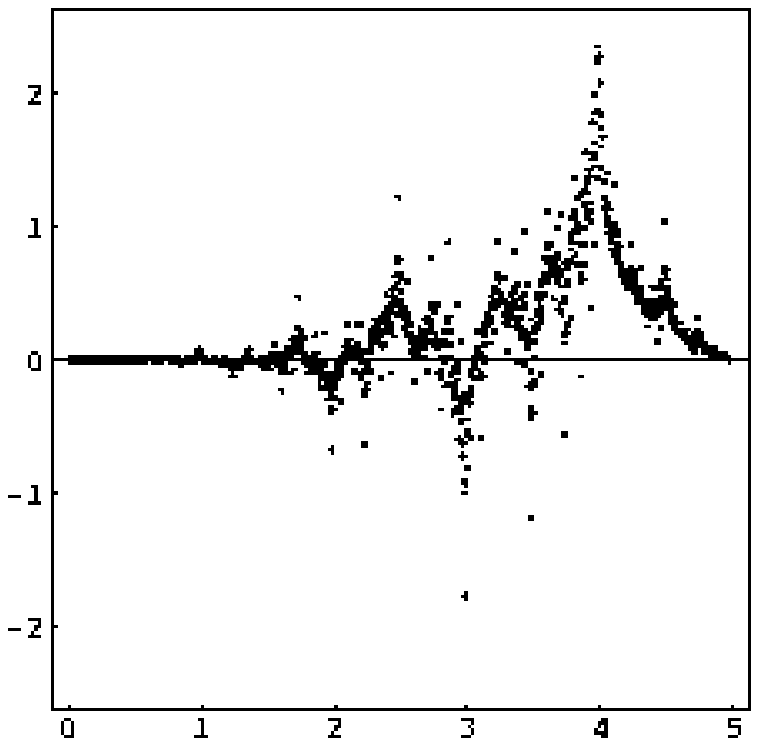}}
\put(0,131){\makebox(120,12){ab: $\theta=0,\;\rho=\pi/12$}}
\put(120,131){\makebox(120,12){bb: $\theta=\pi/12,\;\rho=\pi/12$}}
\put(240,131){\makebox(120,12){cb: $\theta=\pi/6,\;\rho=\pi/12$}}
\put(0,12){\includegraphics[%
height=119bp,width=120bp]{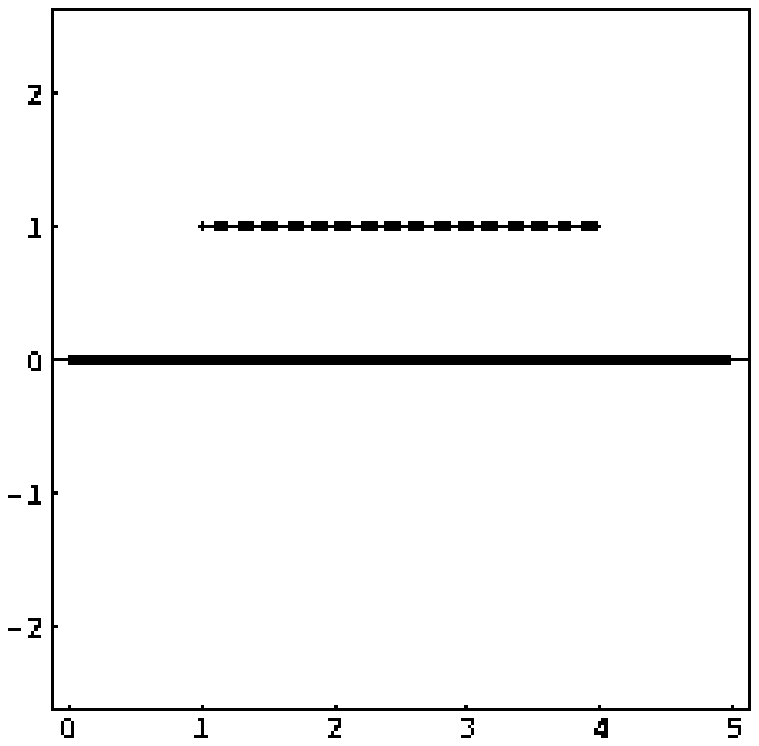}}
\put(120,12){\includegraphics[%
height=119bp,width=120bp]{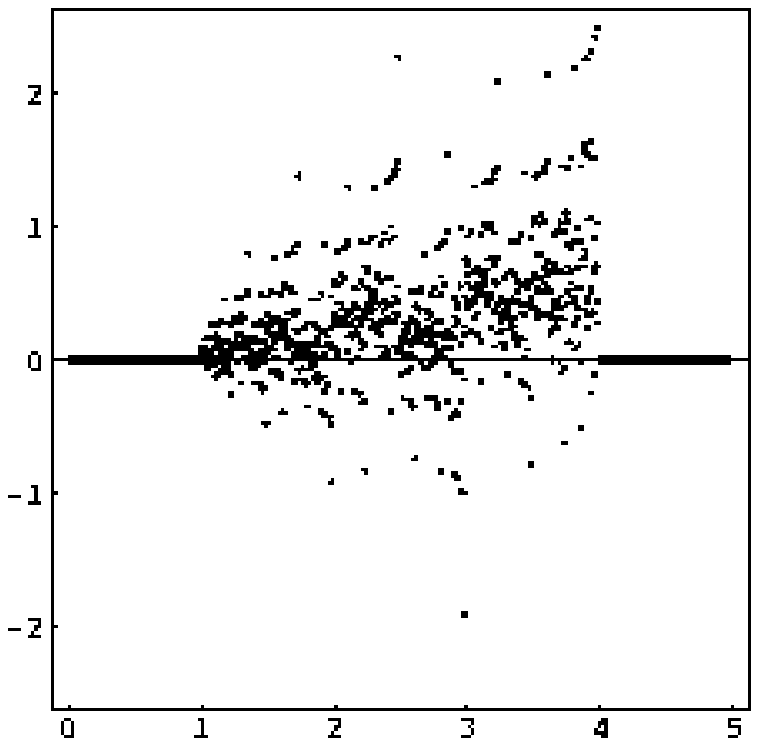}}
\put(240,12){\includegraphics[%
height=119bp,width=120bp]{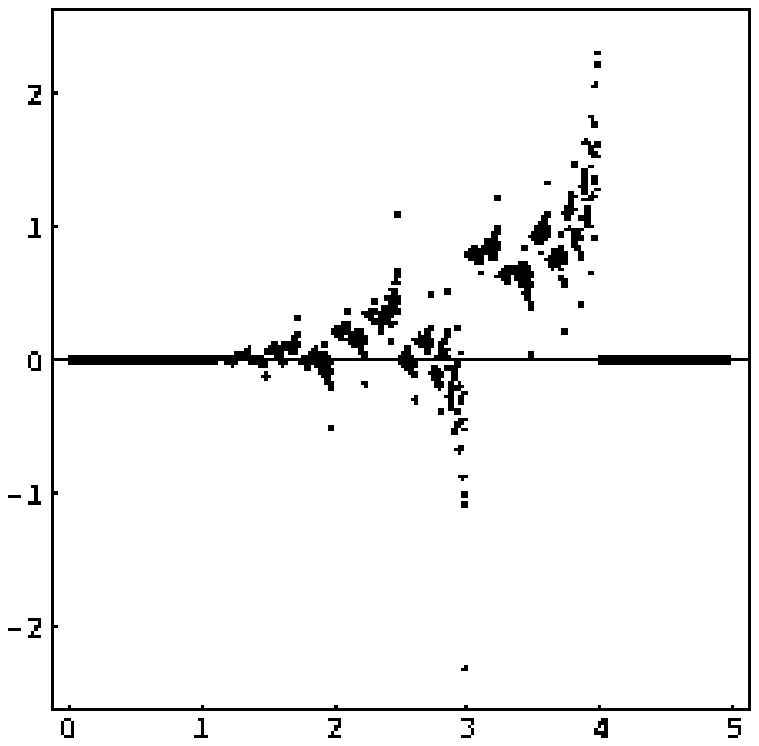}}
\put(0,0){\makebox(120,12){aa: $\theta=0,\;\rho=0$}}
\put(120,0){\makebox(120,12){ba: $\theta=\pi/12,\;\rho=0$}}
\put(240,0){\makebox(120,12){ca: $\theta=\pi/6,\;\rho=0$}}
\end{picture}
\label{P1}\end{figure}

\begin{figure}[tbp]
\setlength{\unitlength}{1bp}
\begin{picture}(360,524)
\put(0,405){\includegraphics[%
height=119bp,width=120bp]{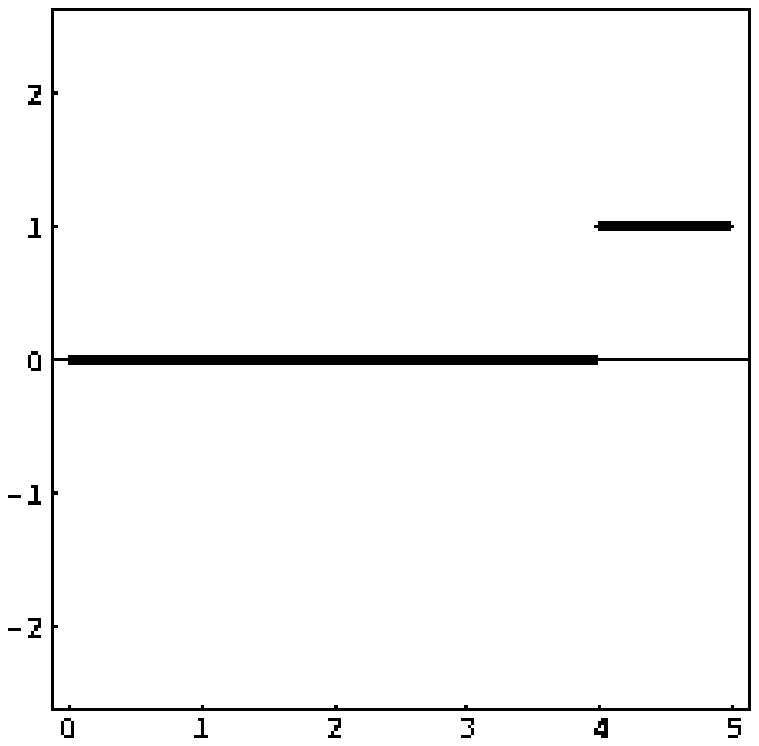}}
\put(120,405){\includegraphics[%
height=119bp,width=120bp]{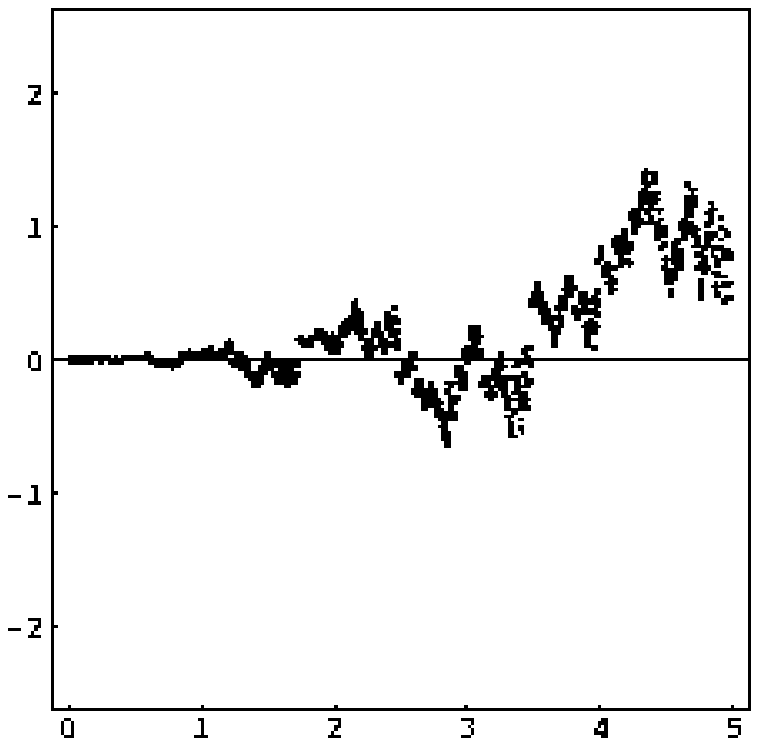}}
\put(240,405){\includegraphics[%
height=119bp,width=120bp]{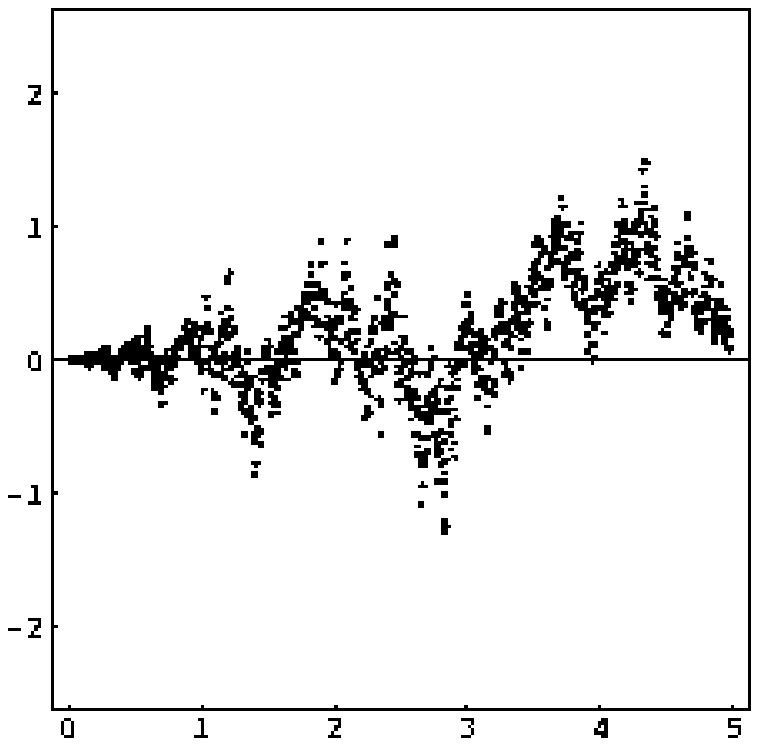}}
\put(0,393){\makebox(120,12){dd: $\theta=\pi/4,\;\rho=\pi/4$}}
\put(120,393){\makebox(120,12){ed: $\theta=\pi/3,\;\rho=\pi/4$}}
\put(240,393){\makebox(120,12){fd: $\theta=5\pi/12,\;\rho=\pi/4$}}
\put(0,274){\includegraphics[%
height=119bp,width=120bp]{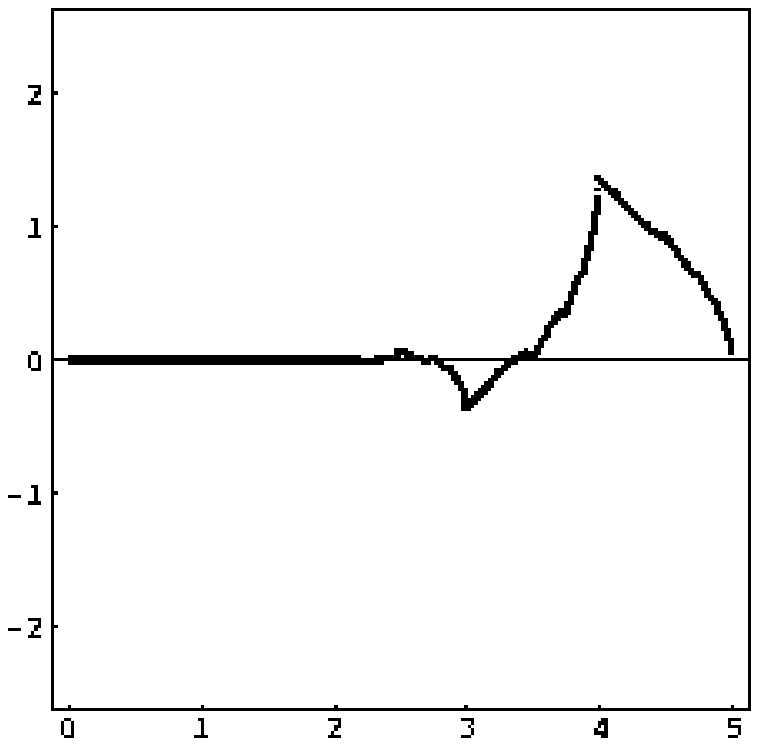}}
\put(120,274){\includegraphics[%
height=119bp,width=120bp]{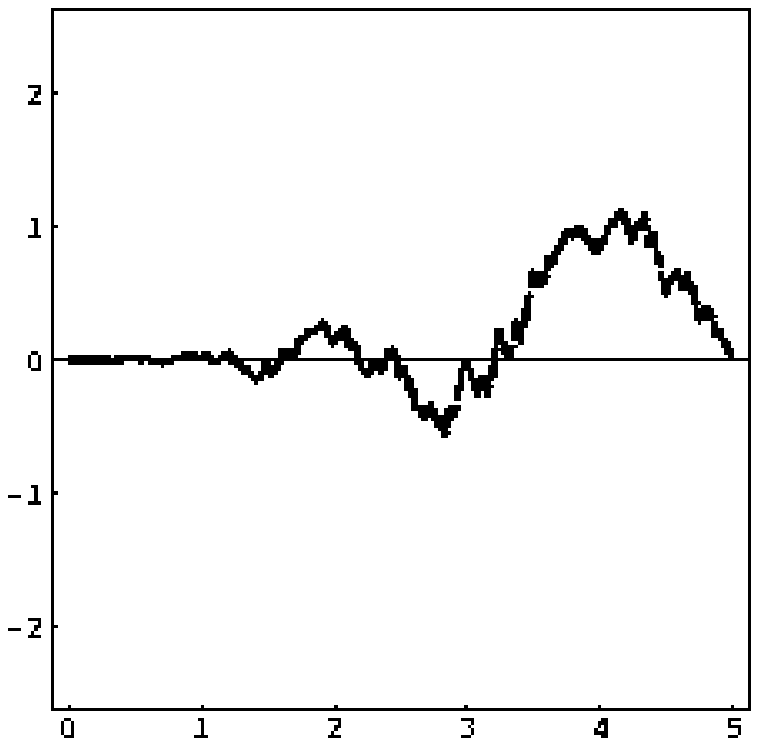}}
\put(240,274){\includegraphics[%
height=119bp,width=120bp]{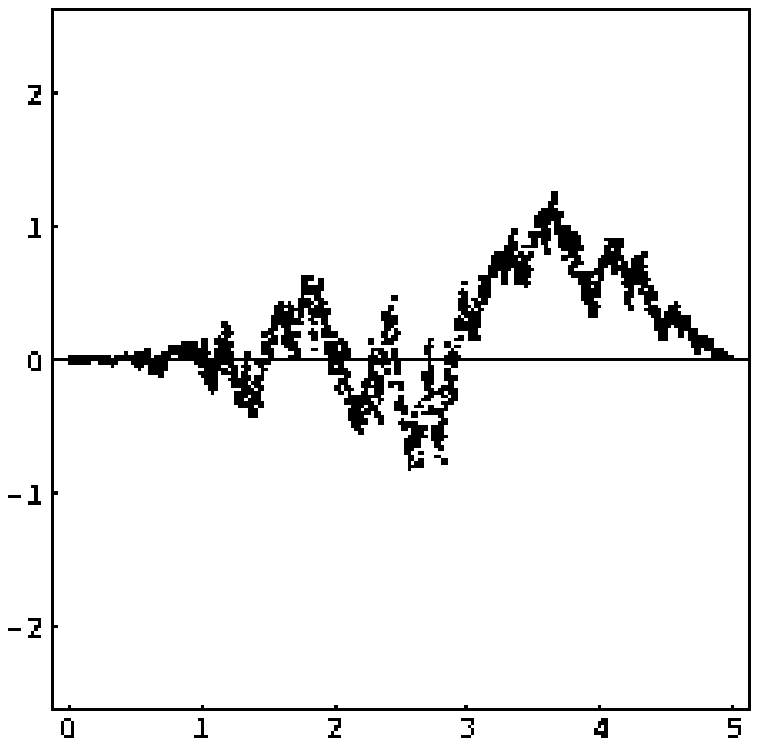}}
\put(0,262){\makebox(120,12){dc: $\theta=\pi/4,\;\rho=\pi/6$}}
\put(120,262){\makebox(120,12){ec: $\theta=\pi/3,\;\rho=\pi/6$}}
\put(240,262){\makebox(120,12){fc: $\theta=5\pi/12,\;\rho=\pi/6$}}
\put(0,143){\includegraphics[%
height=119bp,width=120bp]{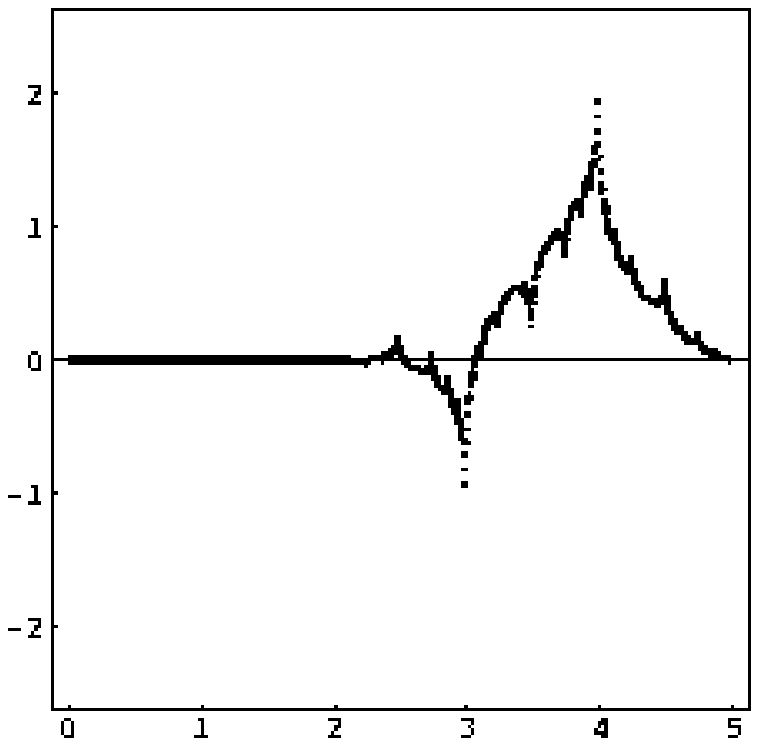}}
\put(120,143){\includegraphics[%
height=119bp,width=120bp]{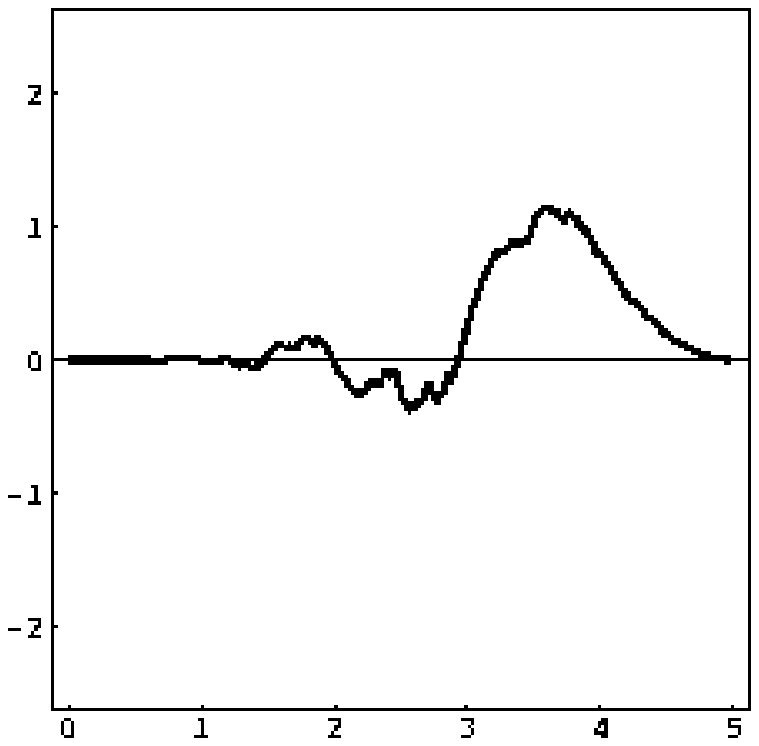}}
\put(240,143){\includegraphics[%
height=119bp,width=120bp]{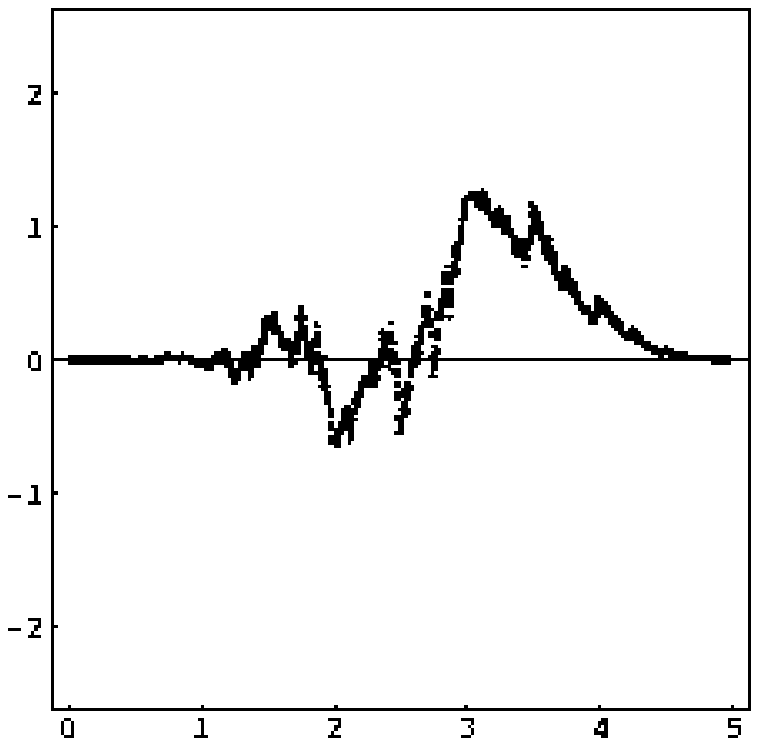}}
\put(0,131){\makebox(120,12){db: $\theta=\pi/4,\;\rho=\pi/12$}}
\put(120,131){\makebox(120,12){eb: $\theta=\pi/3,\;\rho=\pi/12$}}
\put(240,131){\makebox(120,12){fb: $\theta=5\pi/12,\;\rho=\pi/12$}}
\put(0,12){\includegraphics[%
height=119bp,width=120bp]{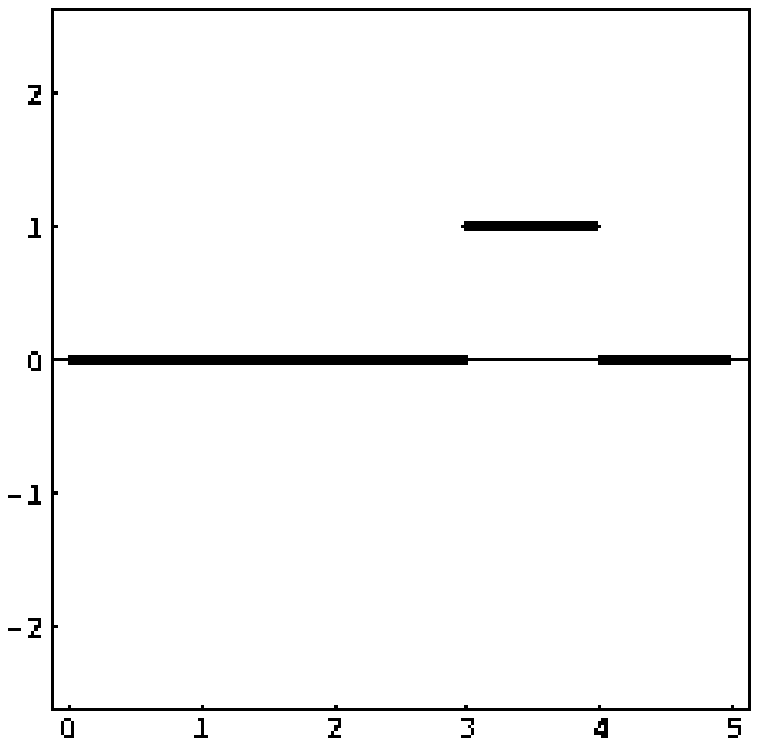}}
\put(120,12){\includegraphics[%
height=119bp,width=120bp]{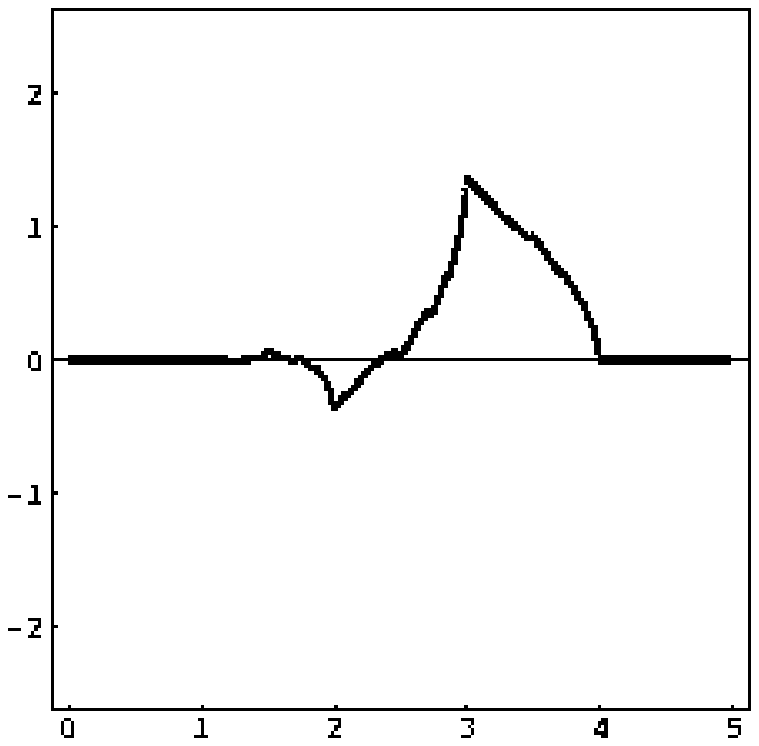}}
\put(240,12){\includegraphics[%
height=119bp,width=120bp]{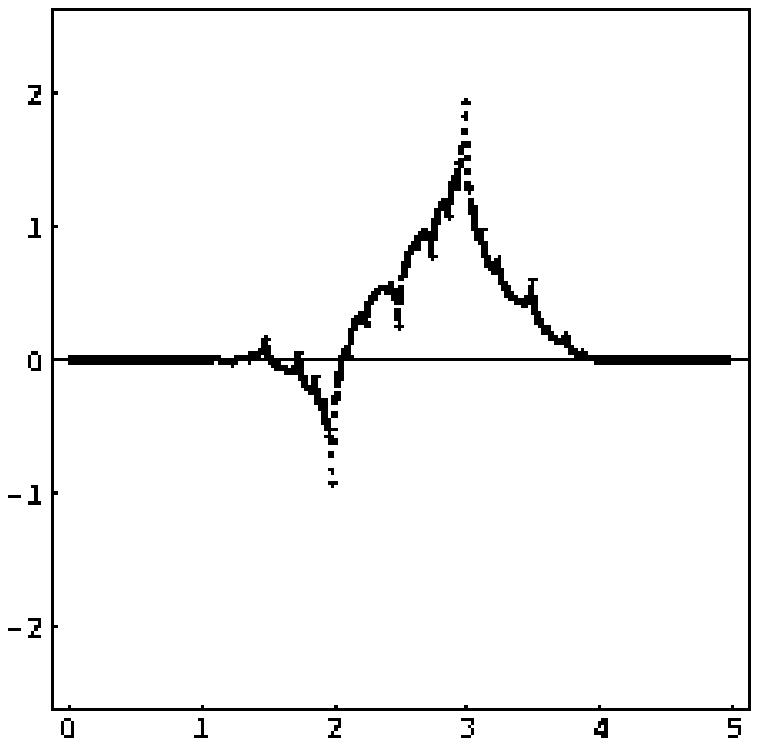}}
\put(0,0){\makebox(120,12){da: $\theta=\pi/4,\;\rho=0$}}
\put(120,0){\makebox(120,12){ea: $\theta=\pi/3,\;\rho=0$}}
\put(240,0){\makebox(120,12){fa: $\theta=5\pi/12,\;\rho=0$}}
\end{picture}
\label{P2}\end{figure}

\begin{figure}[tbp]
\setlength{\unitlength}{1bp}
\begin{picture}(360,524)
\put(0,405){\includegraphics[%
height=119bp,width=120bp]{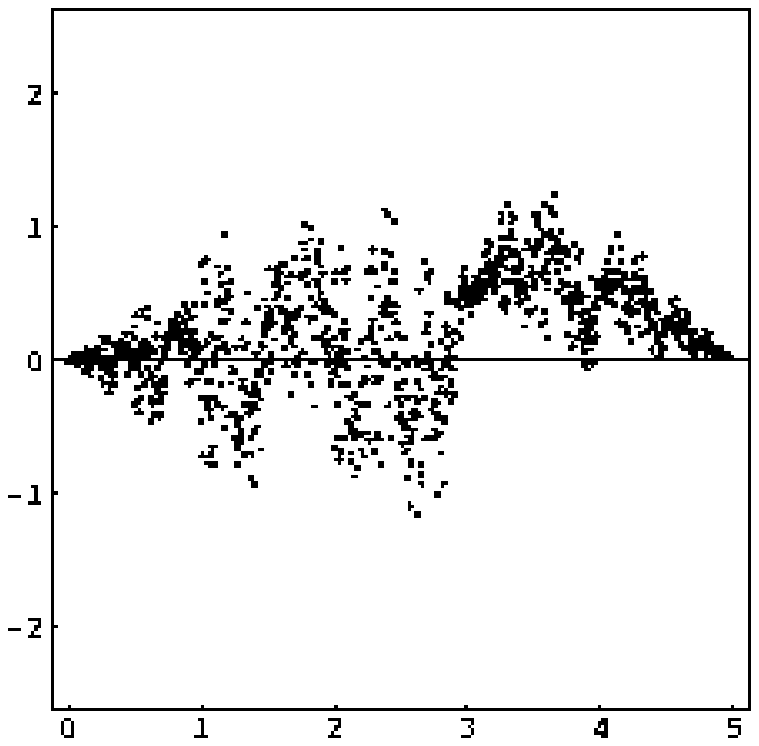}}
\put(120,405){\includegraphics[%
height=119bp,width=120bp]{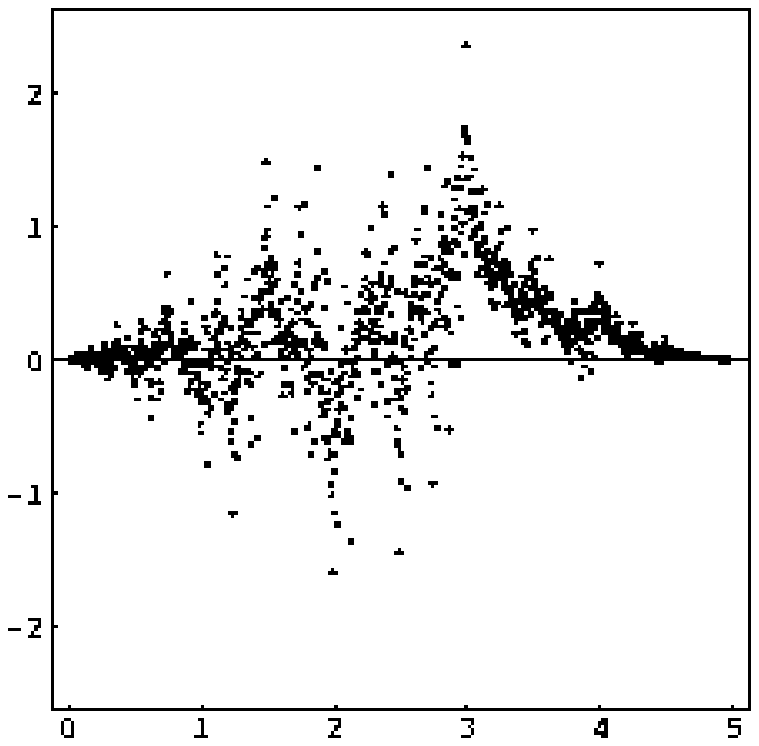}}
\put(240,405){\includegraphics[%
height=119bp,width=120bp]{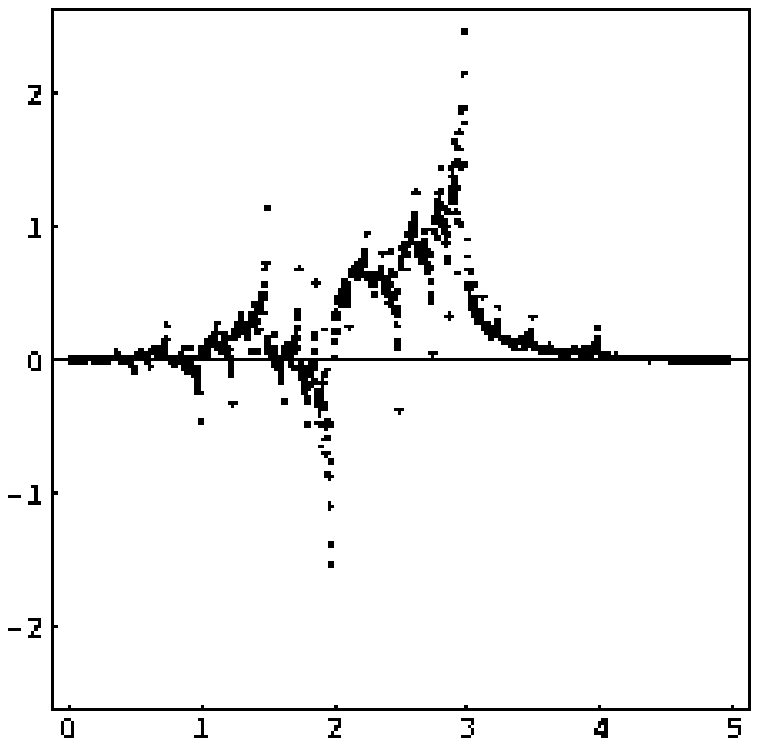}}
\put(0,393){\makebox(120,12){gd: $\theta=\pi/2,\;\rho=\pi/4$}}
\put(120,393){\makebox(120,12){hd: $\theta=7\pi/12,\;\rho=\pi/4$}}
\put(240,393){\makebox(120,12){id: $\theta=2\pi/3,\;\rho=\pi/4$}}
\put(0,274){\includegraphics[%
height=119bp,width=120bp]{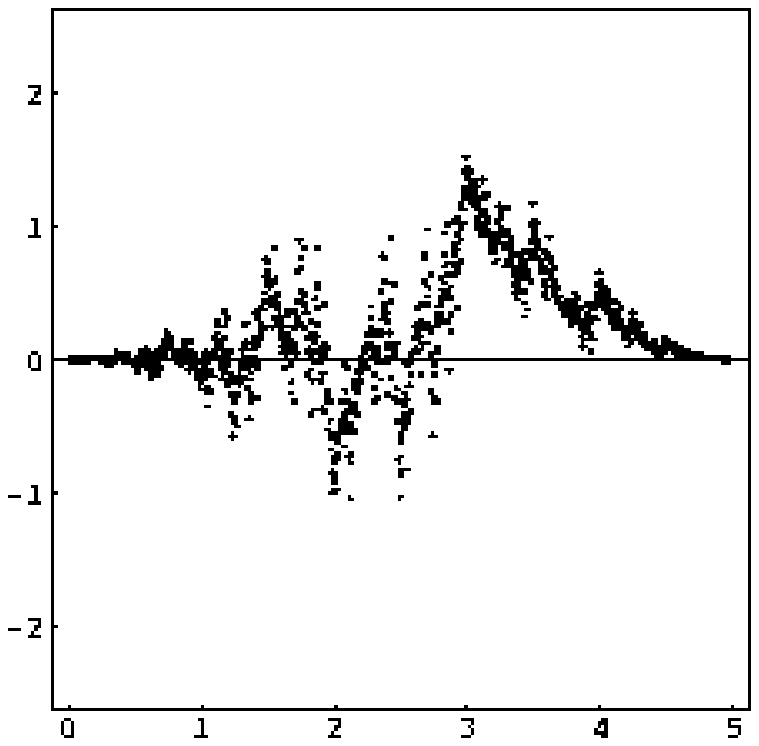}}
\put(120,274){\includegraphics[%
height=119bp,width=120bp]{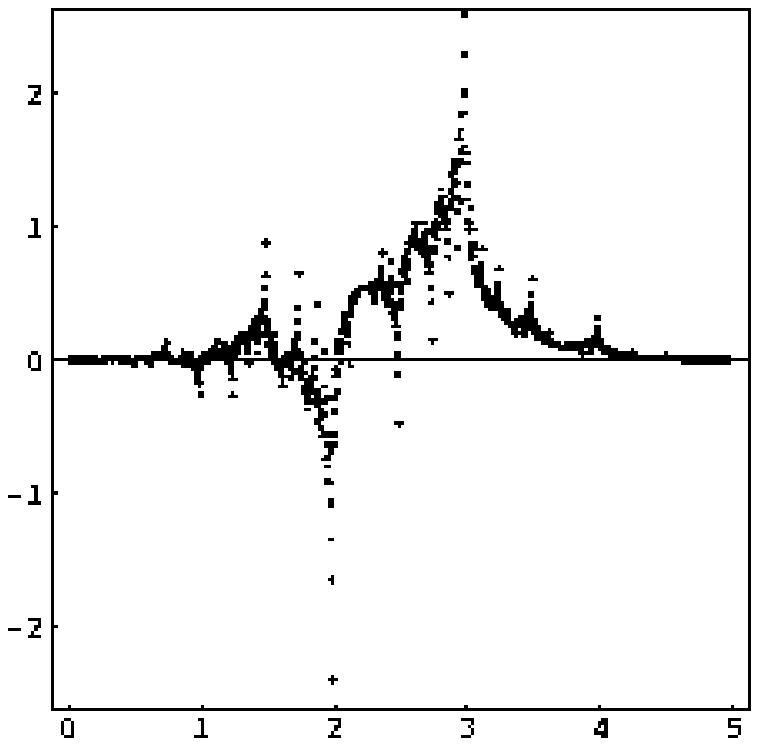}}
\put(240,274){\includegraphics[%
height=119bp,width=120bp]{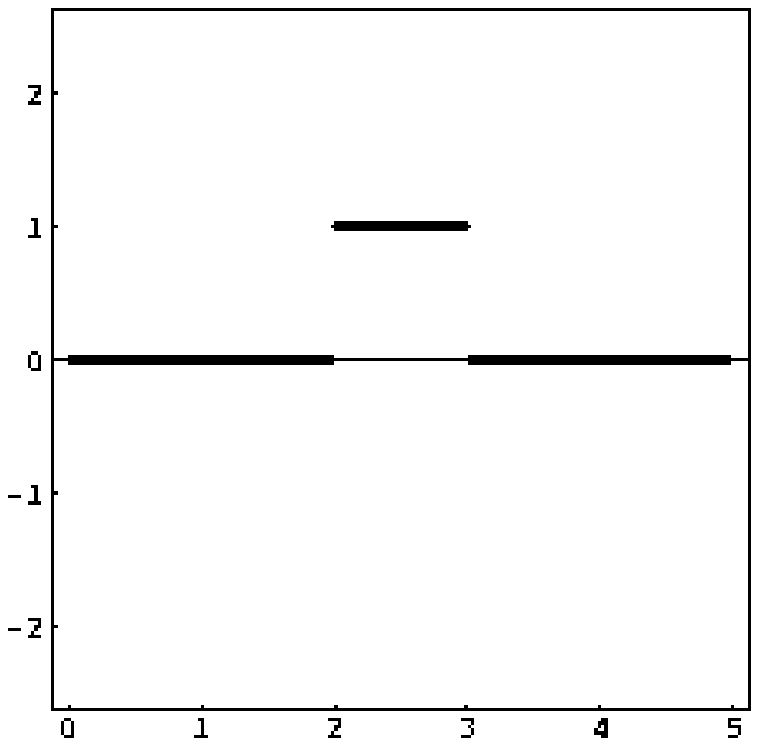}}
\put(0,262){\makebox(120,12){gc: $\theta=\pi/2,\;\rho=\pi/6$}}
\put(120,262){\makebox(120,12){hc: $\theta=7\pi/12,\;\rho=\pi/6$}}
\put(240,262){\makebox(120,12){ic: $\theta=2\pi/3,\;\rho=\pi/6$}}
\put(0,143){\includegraphics[%
height=119bp,width=120bp]{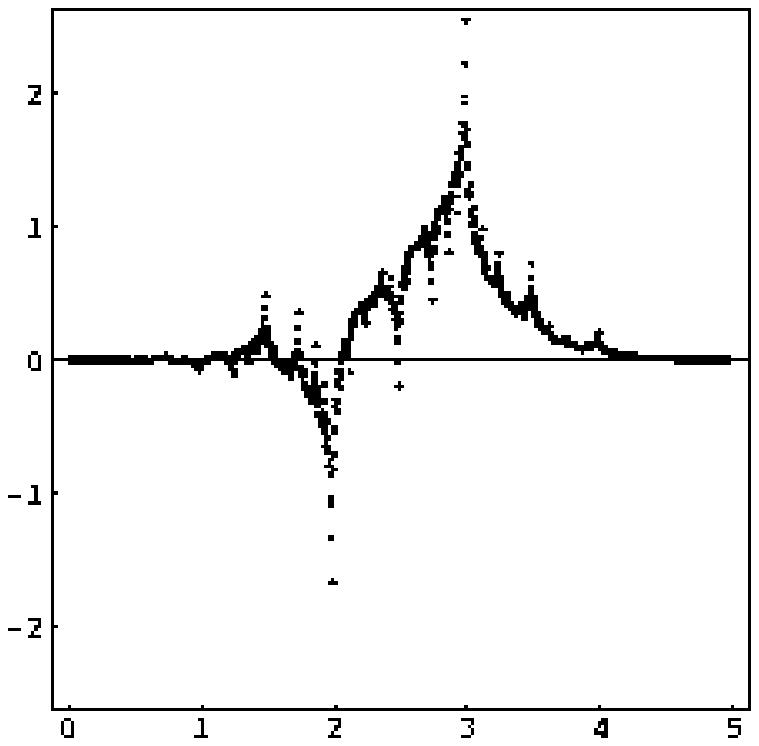}}
\put(120,143){\includegraphics[%
height=119bp,width=120bp]{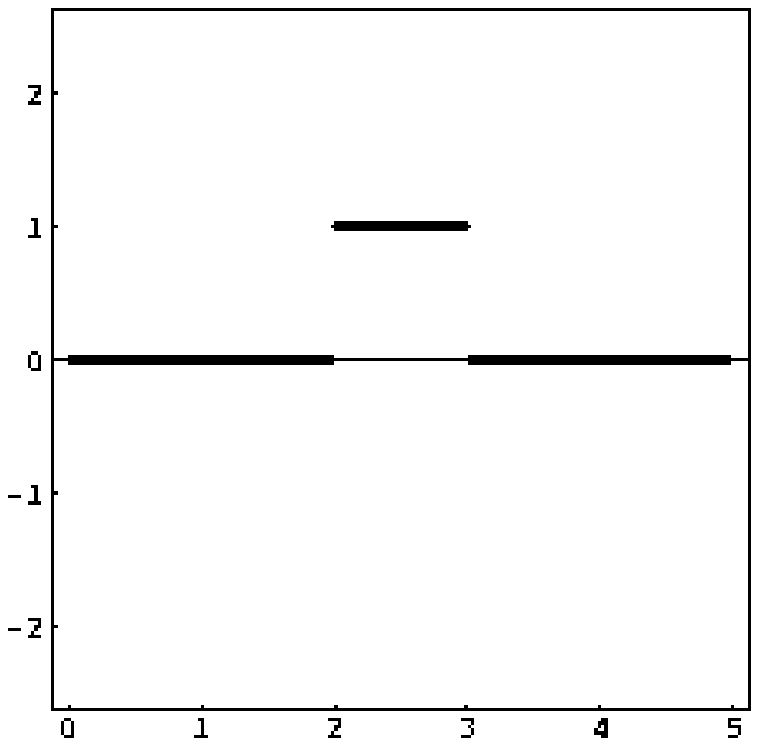}}
\put(240,143){\includegraphics[%
height=119bp,width=120bp]{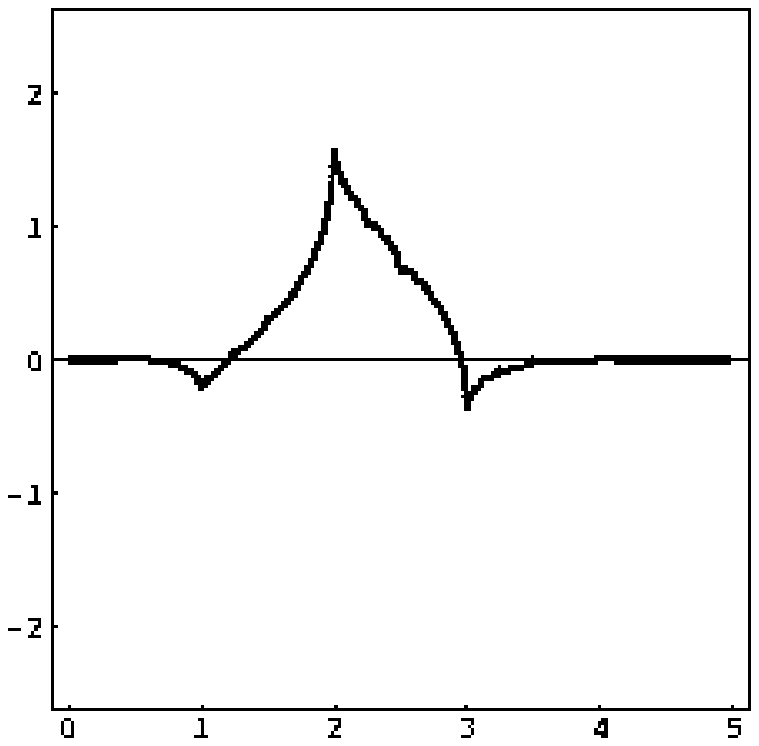}}
\put(0,131){\makebox(120,12){gb: $\theta=\pi/2,\;\rho=\pi/12$}}
\put(120,131){\makebox(120,12){hb: $\theta=7\pi/12,\;\rho=\pi/12$}}
\put(240,131){\makebox(120,12){ib: $\theta=2\pi/3,\;\rho=\pi/12$}}
\put(0,12){\includegraphics[%
height=119bp,width=120bp]{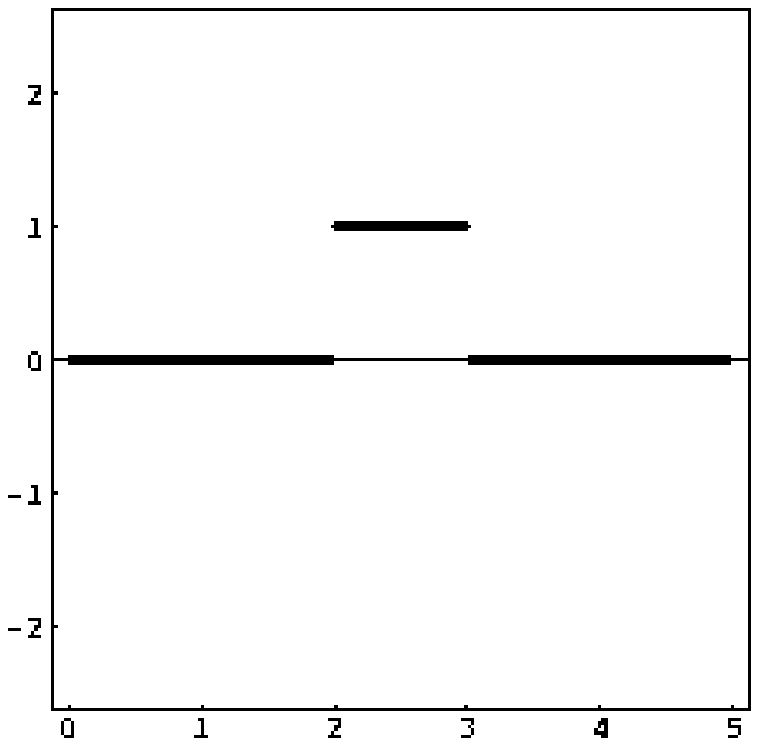}}
\put(120,12){\includegraphics[%
height=119bp,width=120bp]{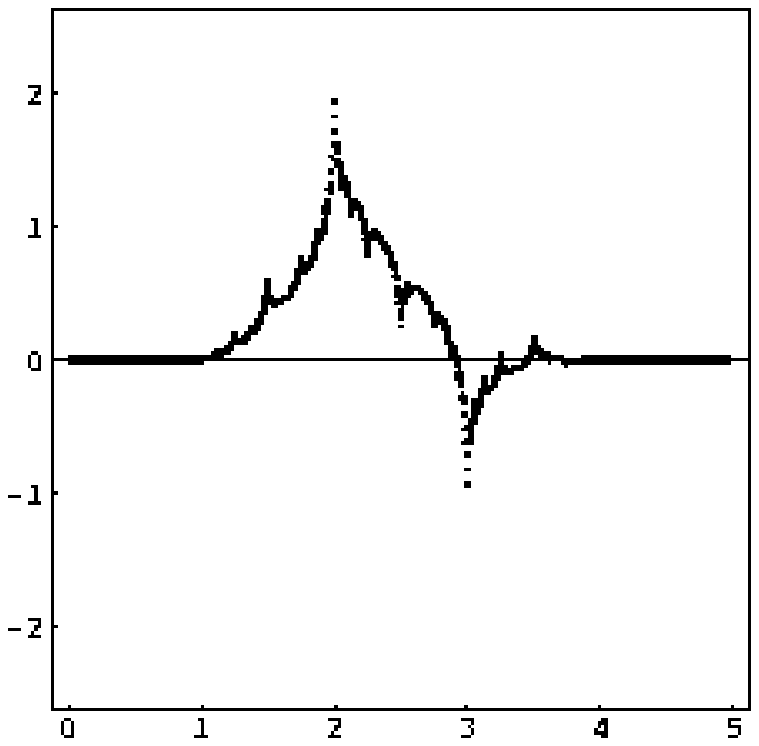}}
\put(240,12){\includegraphics[%
height=119bp,width=120bp]{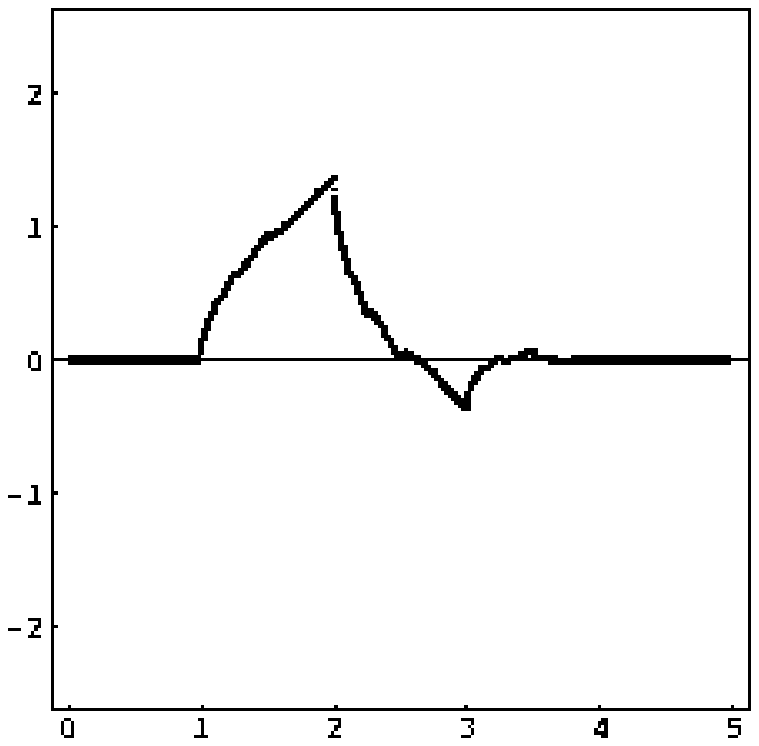}}
\put(0,0){\makebox(120,12){ga: $\theta=\pi/2,\;\rho=0$}}
\put(120,0){\makebox(120,12){ha: $\theta=7\pi/12,\;\rho=0$}}
\put(240,0){\makebox(120,12){ia: $\theta=2\pi/3,\;\rho=0$}}
\end{picture}
\label{P3}\end{figure}

\begin{figure}[tbp]
\setlength{\unitlength}{1bp}
\begin{picture}(360,524)
\put(0,405){\includegraphics[%
height=119bp,width=120bp]{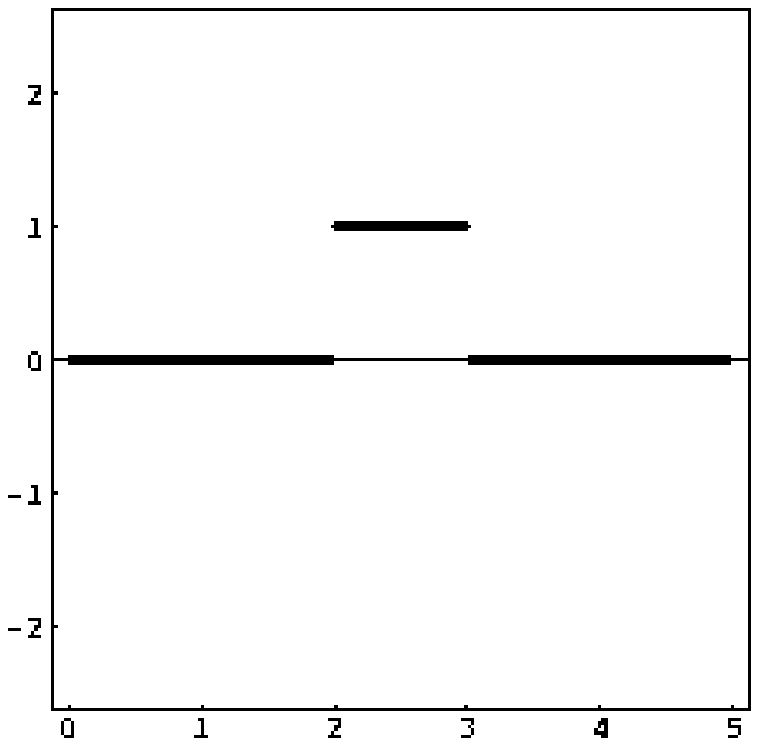}}
\put(120,405){\includegraphics[%
height=119bp,width=120bp]{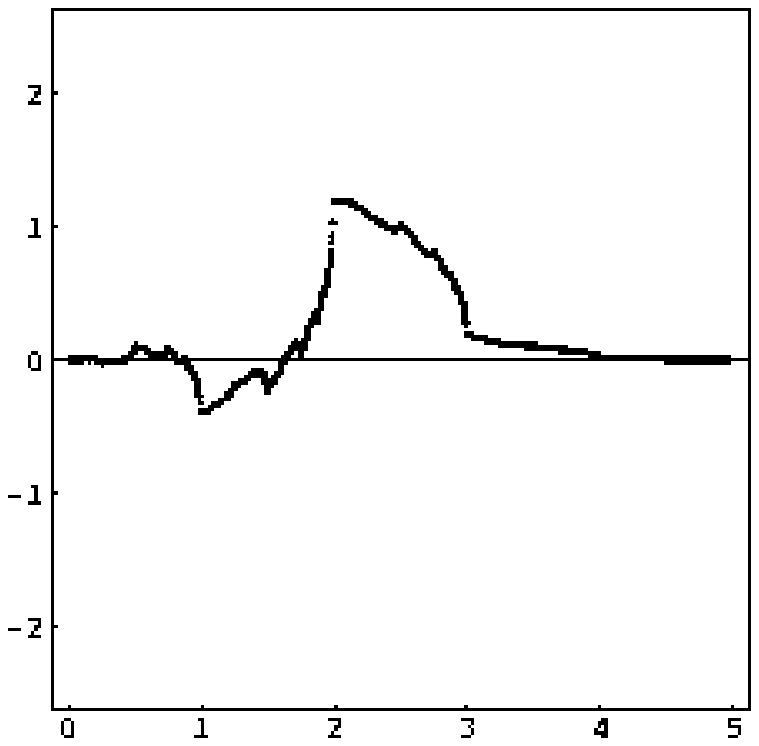}}
\put(240,405){\includegraphics[%
height=119bp,width=120bp]{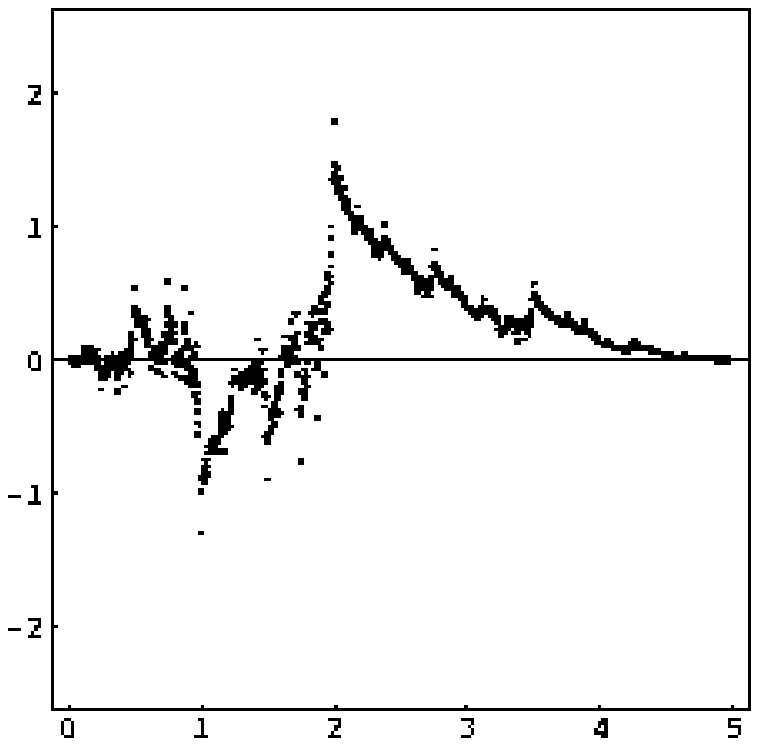}}
\put(0,393){\makebox(120,12){jd: $\theta=3\pi/4,\;\rho=\pi/4$}}
\put(120,393){\makebox(120,12){kd: $\theta=5\pi/6,\;\rho=\pi/4$}}
\put(240,393){\makebox(120,12){ld: $\theta=11\pi/12,\;\rho=\pi/4$}}
\put(0,274){\includegraphics[%
height=119bp,width=120bp]{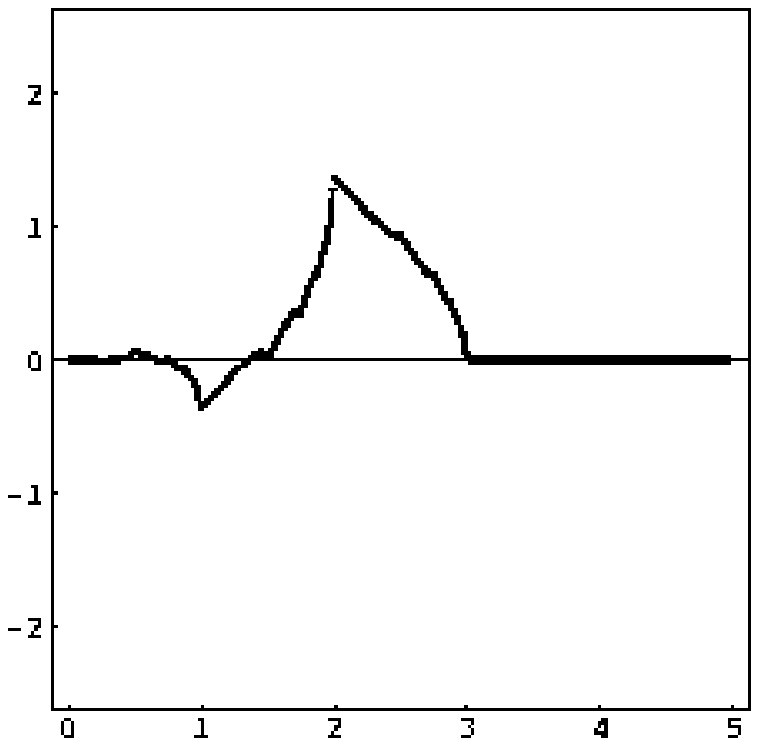}}
\put(120,274){\includegraphics[%
height=119bp,width=120bp]{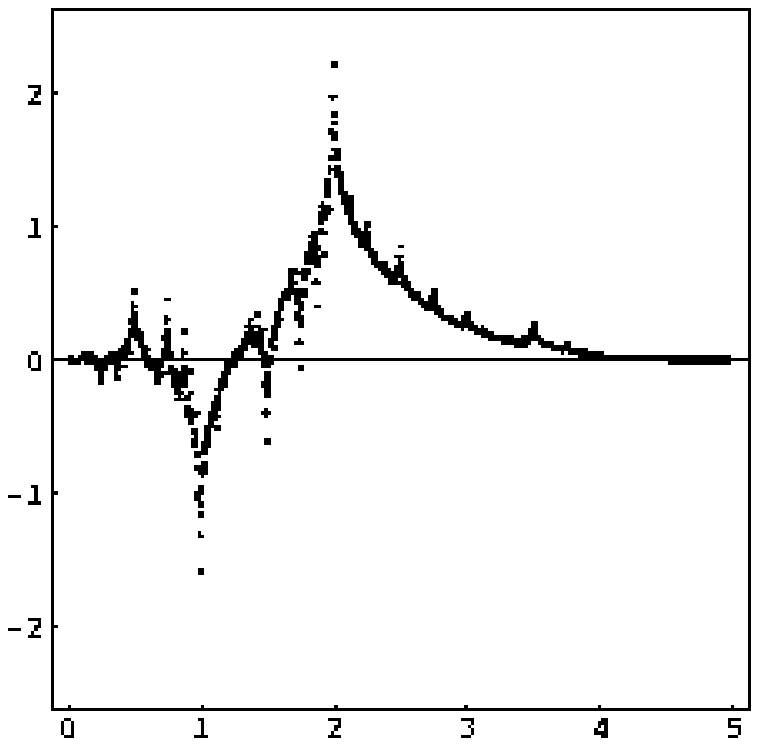}}
\put(240,274){\includegraphics[%
height=119bp,width=120bp]{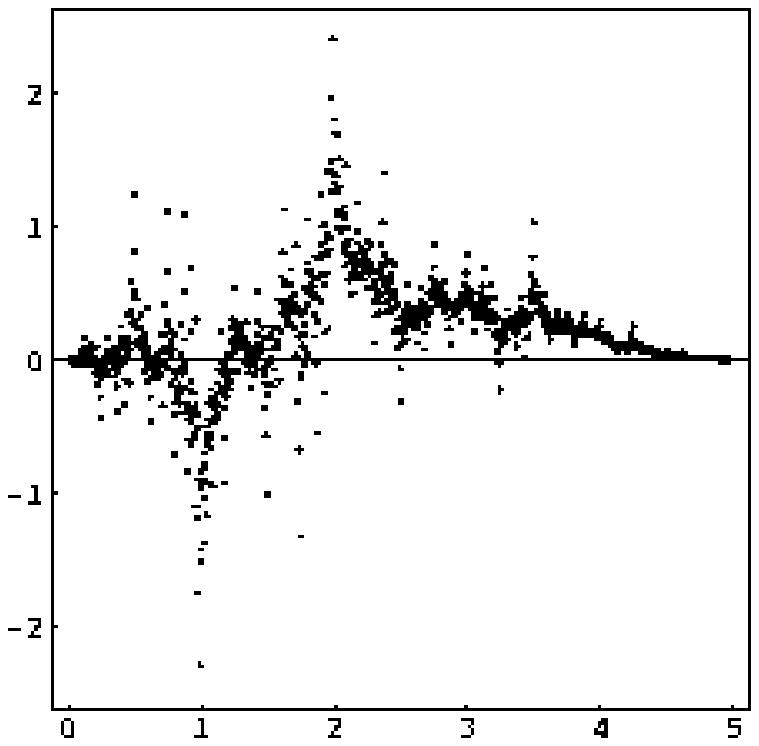}}
\put(0,262){\makebox(120,12){jc: $\theta=3\pi/4,\;\rho=\pi/6$}}
\put(120,262){\makebox(120,12){kc: $\theta=5\pi/6,\;\rho=\pi/6$}}
\put(240,262){\makebox(120,12){lc: $\theta=11\pi/12,\;\rho=\pi/6$}}
\put(0,143){\includegraphics[%
height=119bp,width=120bp]{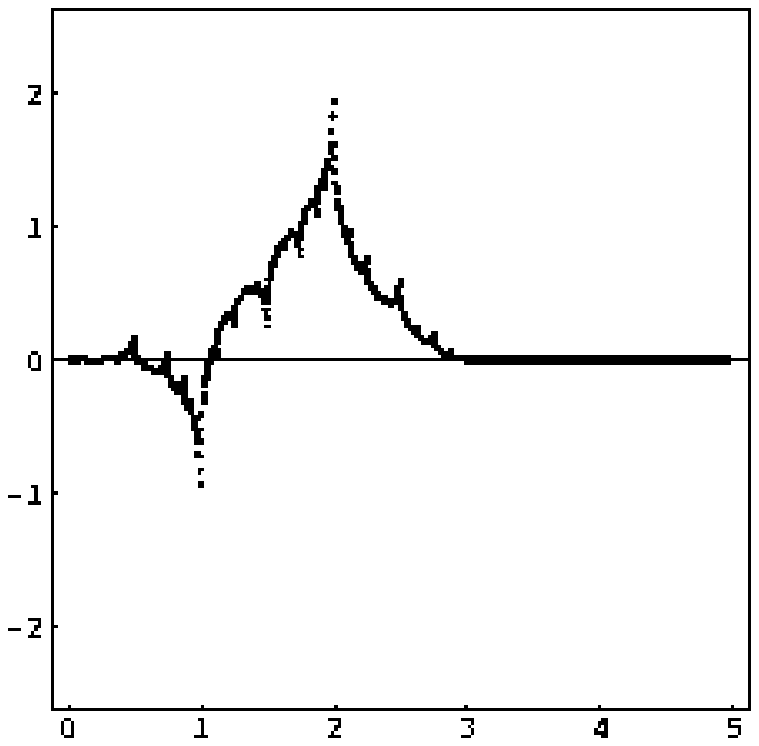}}
\put(120,143){\includegraphics[%
height=119bp,width=120bp]{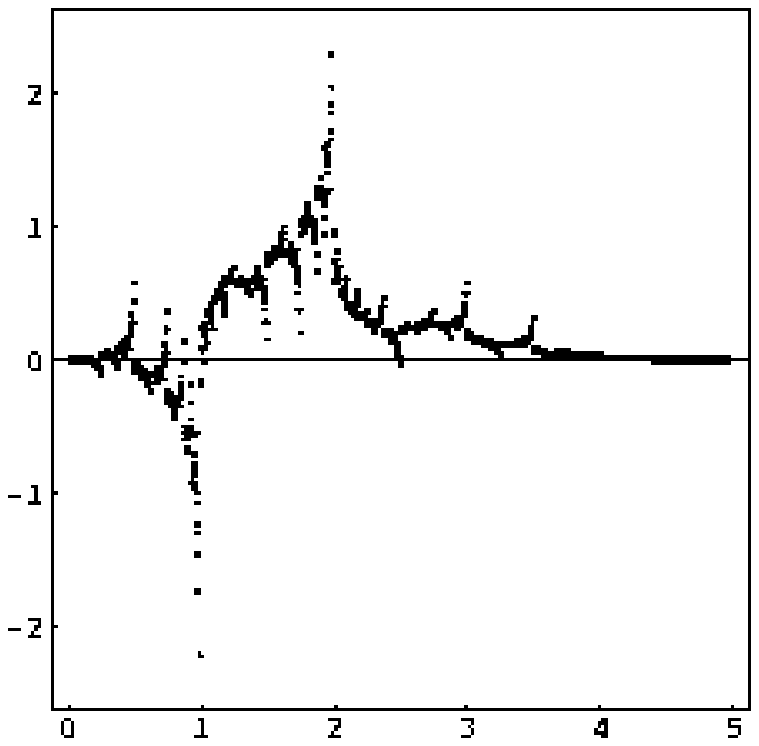}}
\put(240,143){\includegraphics[%
height=119bp,width=120bp]{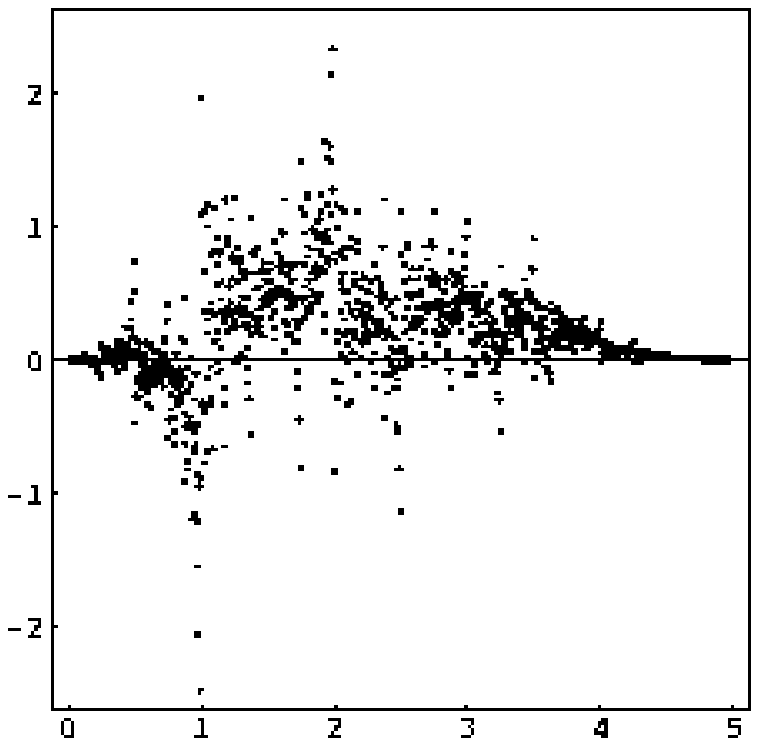}}
\put(0,131){\makebox(120,12){jb: $\theta=3\pi/4,\;\rho=\pi/12$}}
\put(120,131){\makebox(120,12){kb: $\theta=5\pi/6,\;\rho=\pi/12$}}
\put(240,131){\makebox(120,12){lb: $\theta=11\pi/12,\;\rho=\pi/12$}}
\put(0,12){\includegraphics[%
height=119bp,width=120bp]{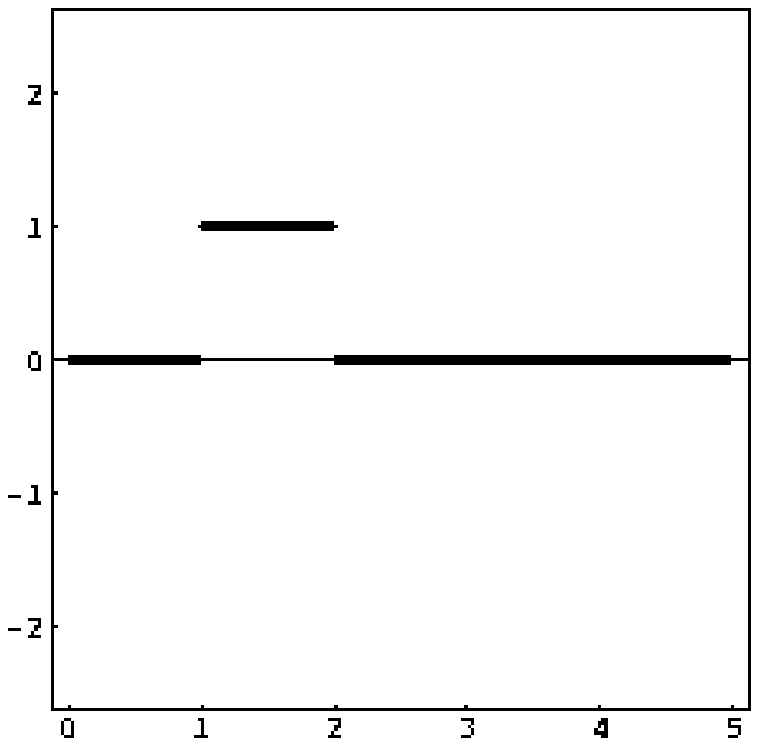}}
\put(120,12){\includegraphics[%
height=119bp,width=120bp]{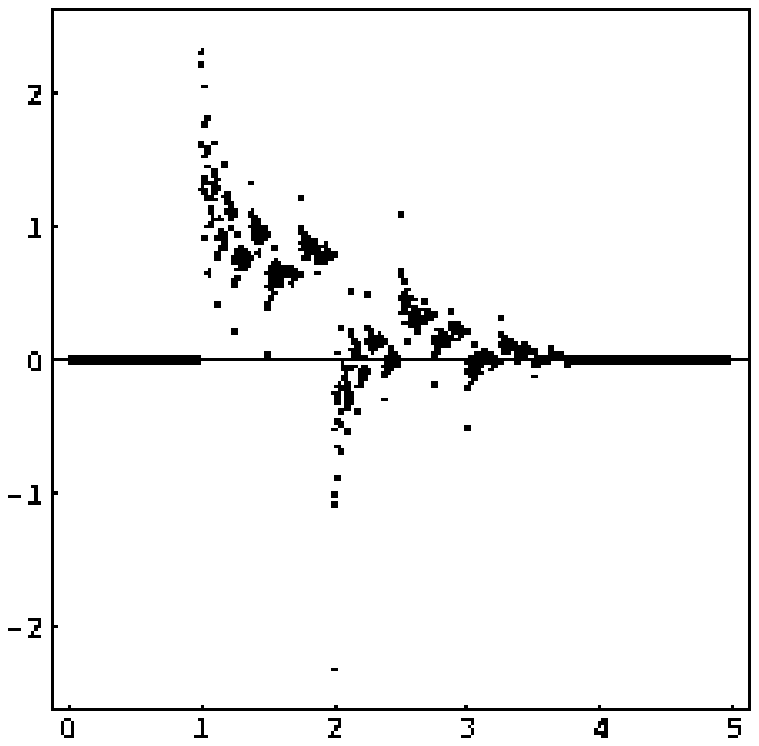}}
\put(240,12){\includegraphics[%
height=119bp,width=120bp]{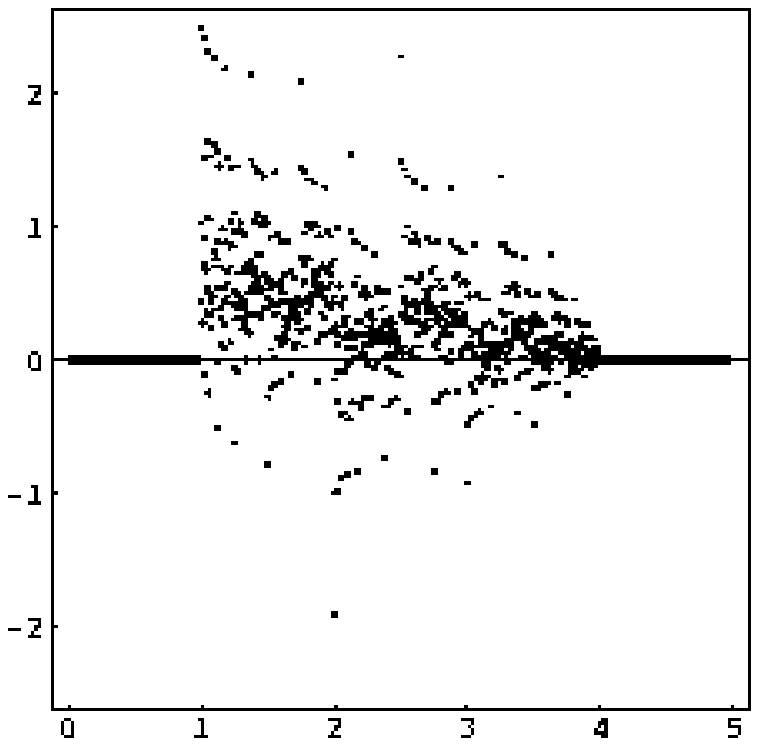}}
\put(0,0){\makebox(120,12){ja: $\theta=3\pi/4,\;\rho=0$}}
\put(120,0){\makebox(120,12){ka: $\theta=5\pi/6,\;\rho=0$}}
\put(240,0){\makebox(120,12){la: $\theta=11\pi/12,\;\rho=0$}}
\end{picture}
\label{P4}\end{figure}

\begin{figure}[tbp]
\setlength{\unitlength}{1bp}
\begin{picture}(360,524)
\put(0,405){\includegraphics[%
height=119bp,width=120bp]{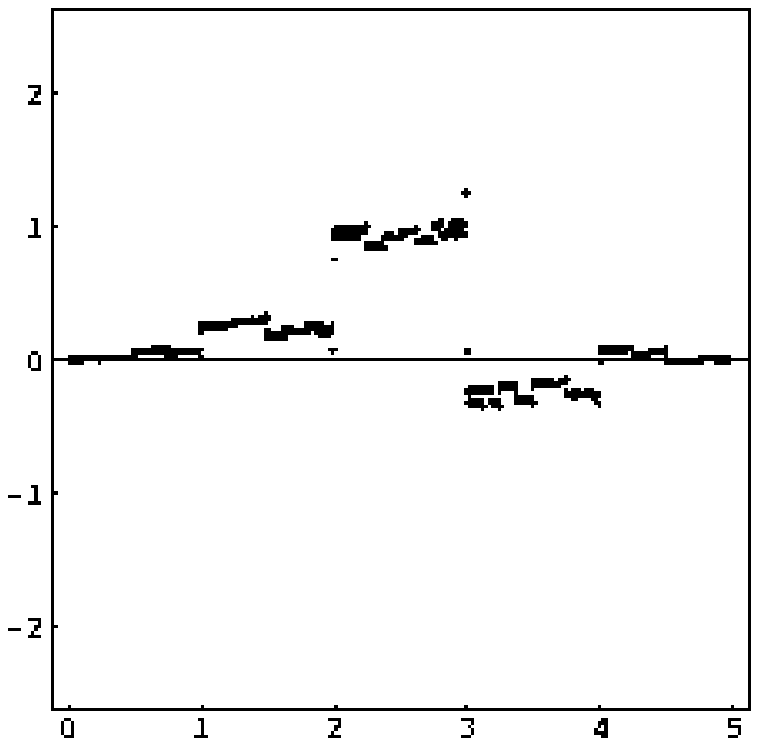}}
\put(120,405){\includegraphics[%
height=119bp,width=120bp]{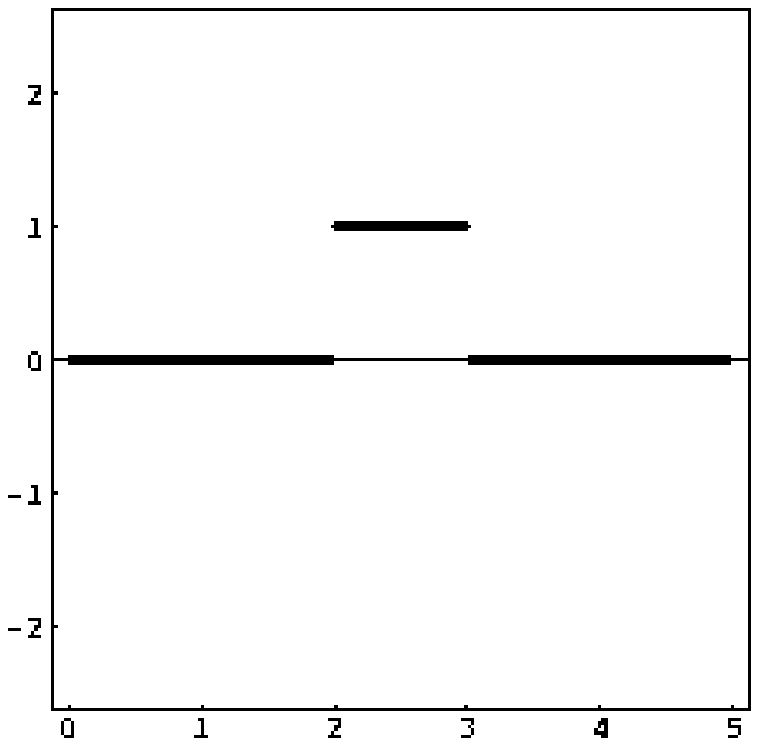}}
\put(240,405){\includegraphics[%
height=119bp,width=120bp]{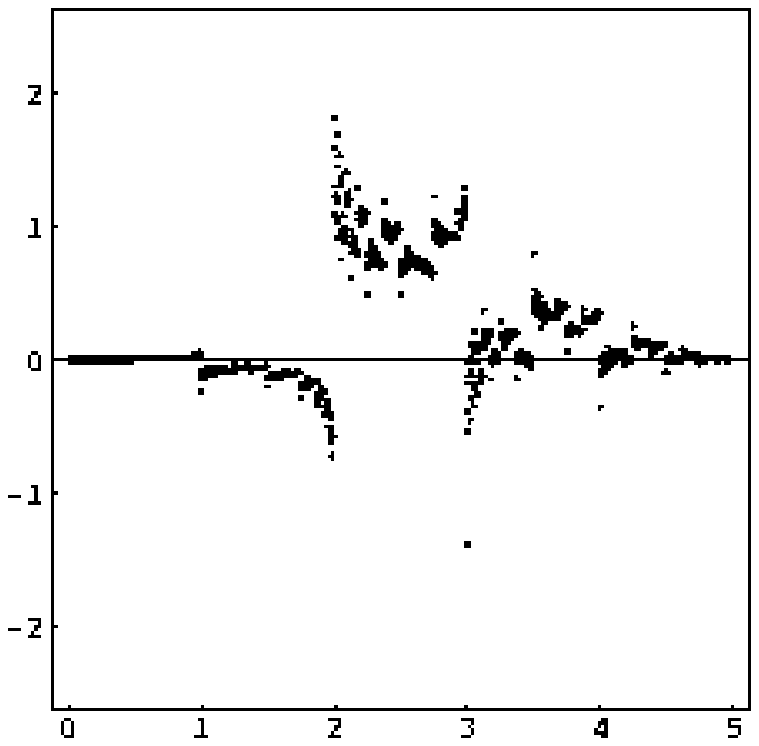}}
\put(0,393){\makebox(120,12){ah: $\theta=0,\;\rho=7\pi/12$}}
\put(120,393){\makebox(120,12){bh: $\theta=\pi/12,\;\rho=7\pi/12$}}
\put(240,393){\makebox(120,12){ch: $\theta=\pi/6,\;\rho=7\pi/12$}}
\put(0,274){\includegraphics[%
height=119bp,width=120bp]{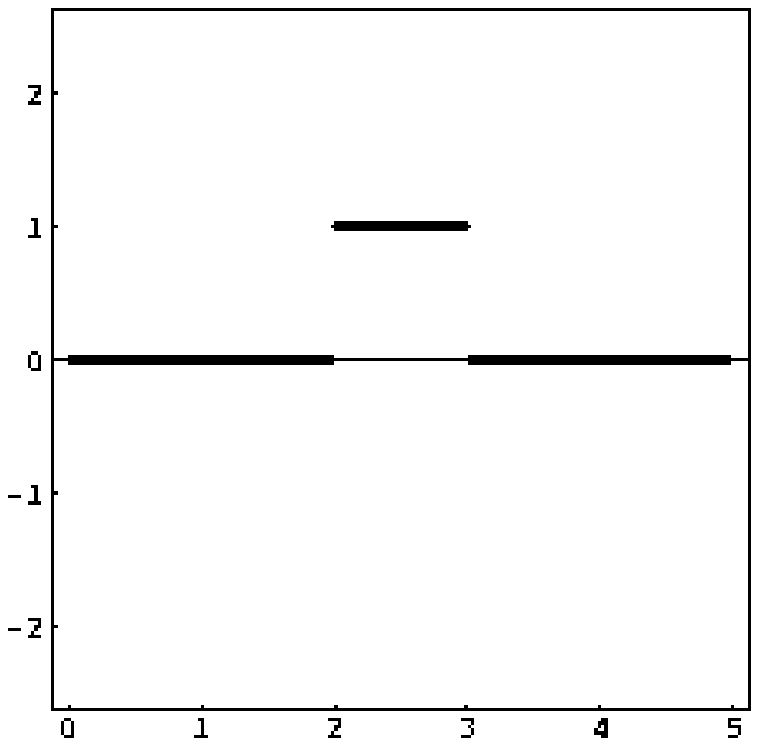}}
\put(120,274){\includegraphics[%
height=119bp,width=120bp]{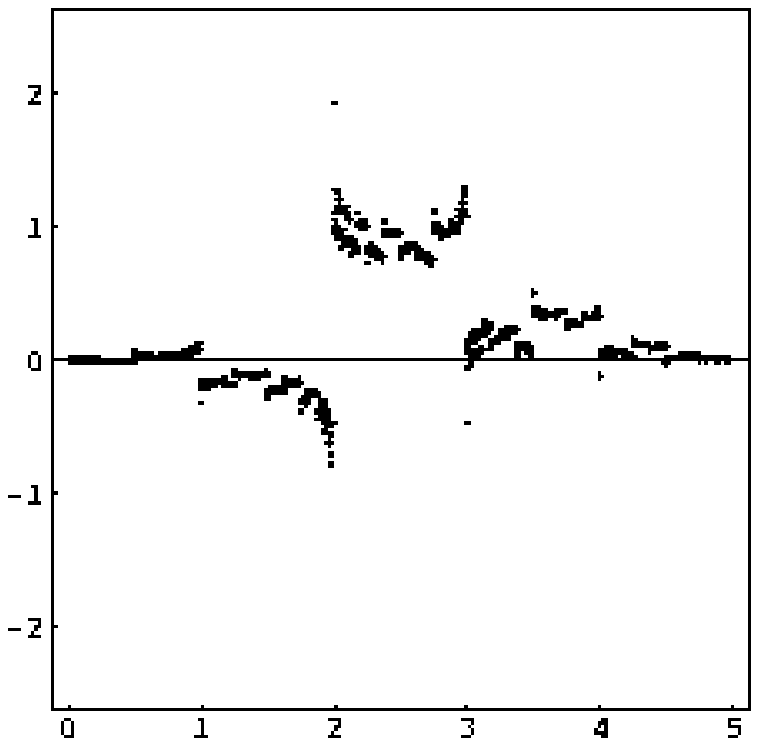}}
\put(240,274){\includegraphics[%
height=119bp,width=120bp]{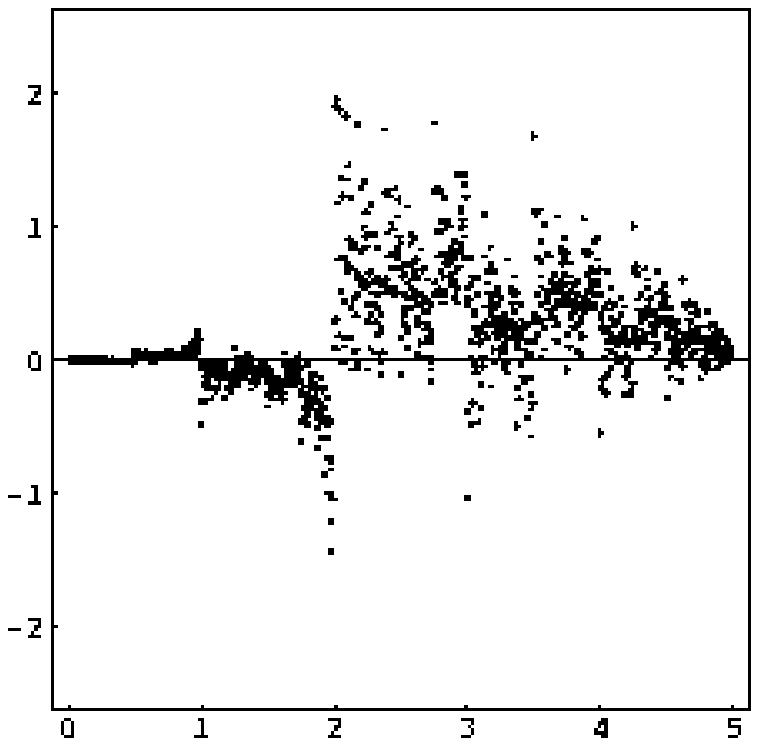}}
\put(0,262){\makebox(120,12){ag: $\theta=0,\;\rho=\pi/2$}}
\put(120,262){\makebox(120,12){bg: $\theta=\pi/12,\;\rho=\pi/2$}}
\put(240,262){\makebox(120,12){cg: $\theta=\pi/6,\;\rho=\pi/2$}}
\put(0,143){\includegraphics[%
height=119bp,width=120bp]{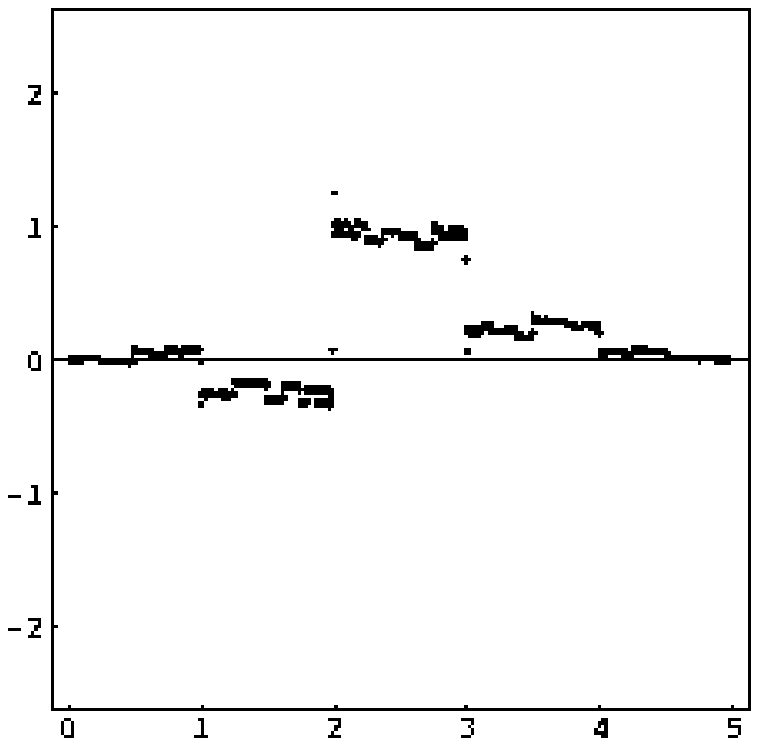}}
\put(120,143){\includegraphics[%
height=119bp,width=120bp]{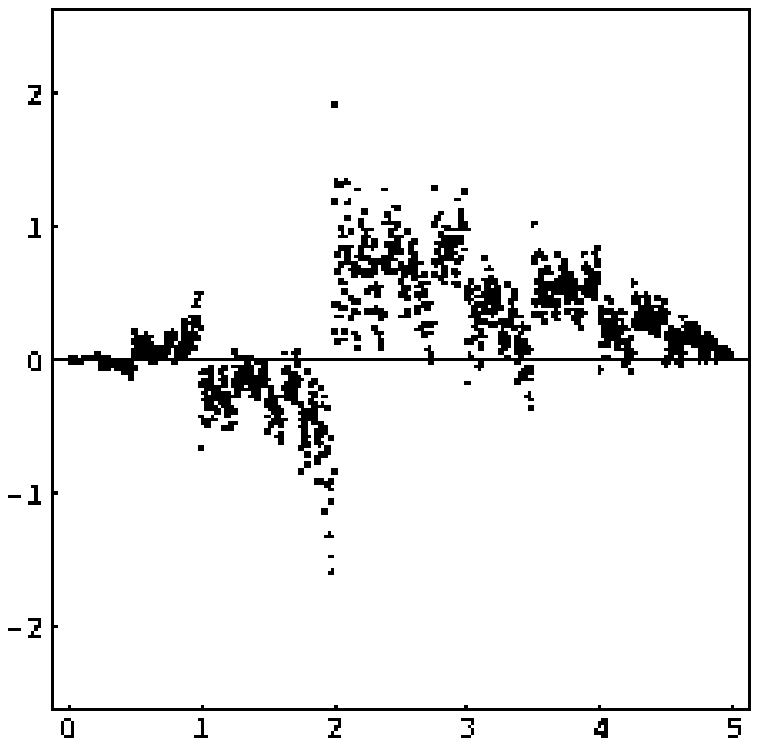}}
\put(240,143){\includegraphics[%
height=119bp,width=120bp]{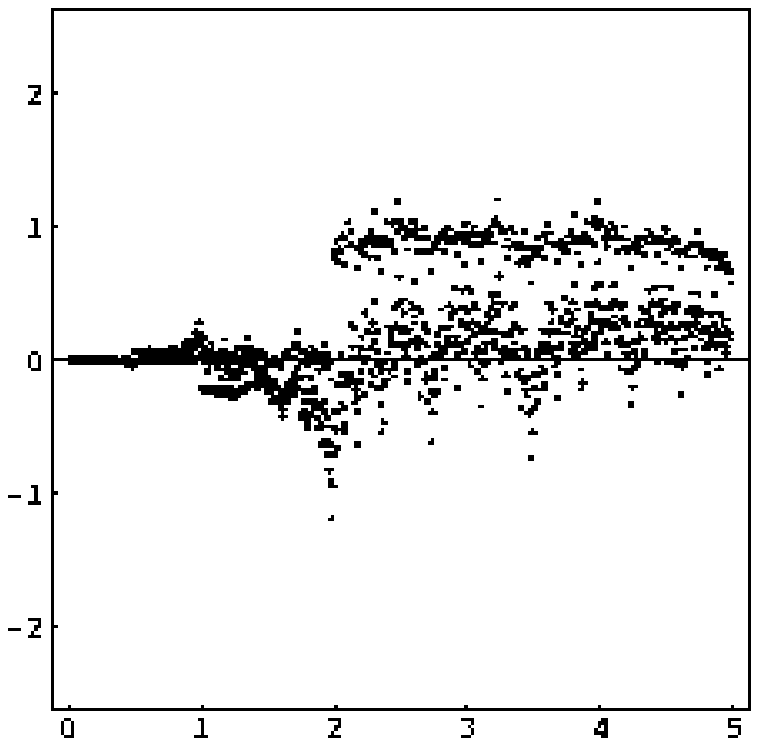}}
\put(0,131){\makebox(120,12){af: $\theta=0,\;\rho=5\pi/12$}}
\put(120,131){\makebox(120,12){bf: $\theta=\pi/12,\;\rho=5\pi/12$}}
\put(240,131){\makebox(120,12){cf: $\theta=\pi/6,\;\rho=5\pi/12$}}
\put(0,12){\includegraphics[%
height=119bp,width=120bp]{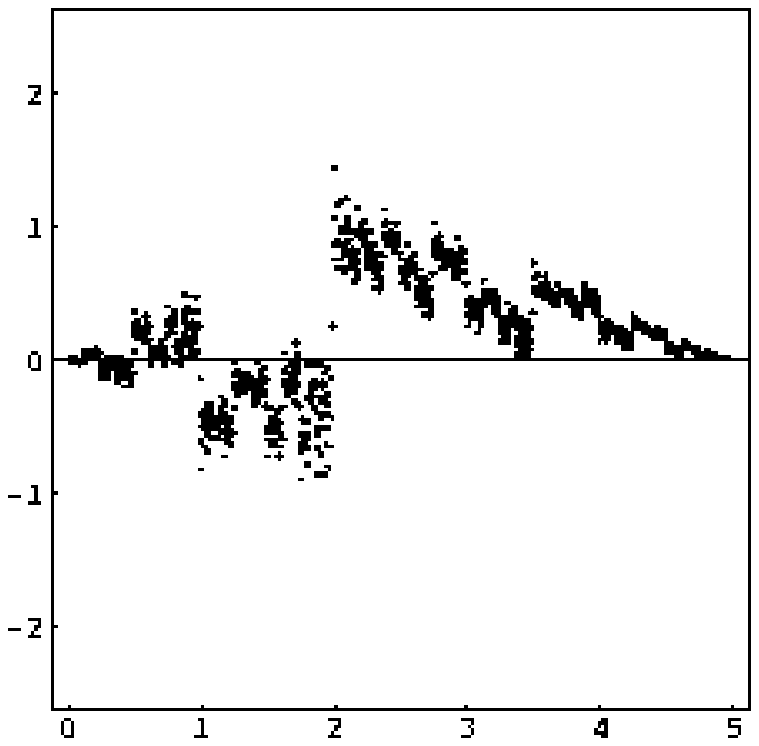}}
\put(120,12){\includegraphics[%
height=119bp,width=120bp]{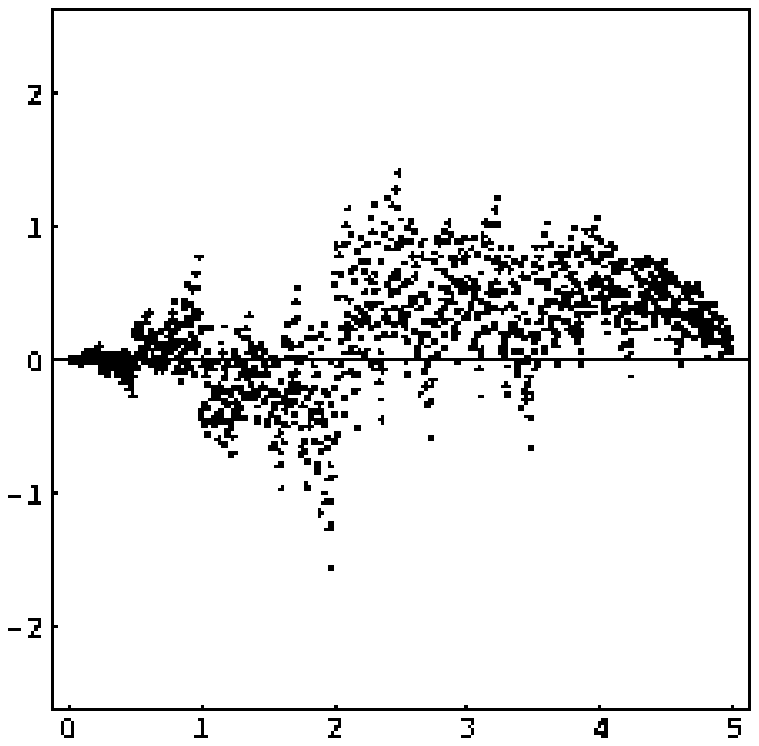}}
\put(240,12){\includegraphics[%
height=119bp,width=120bp]{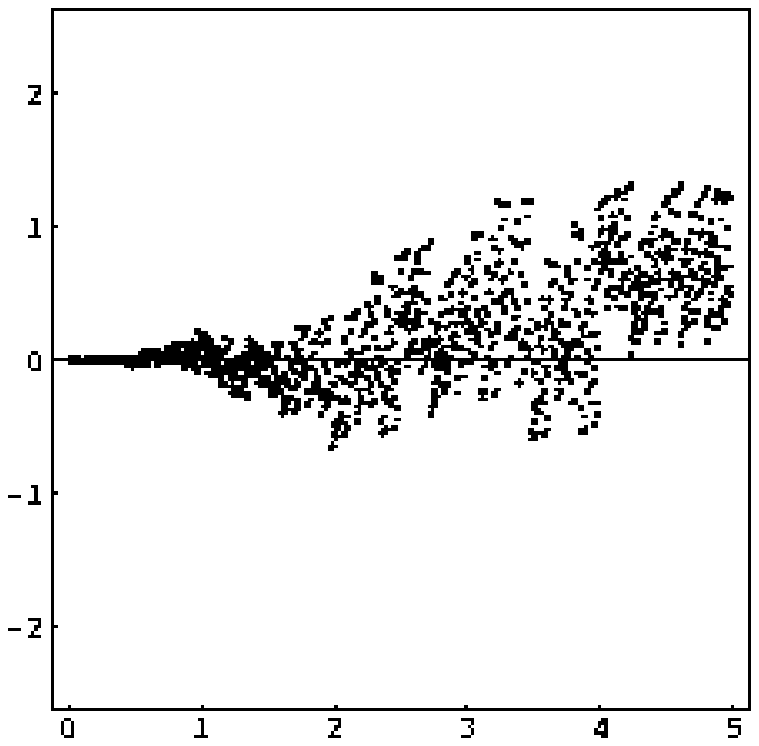}}
\put(0,0){\makebox(120,12){ae: $\theta=0,\;\rho=\pi/3$}}
\put(120,0){\makebox(120,12){be: $\theta=\pi/12,\;\rho=\pi/3$}}
\put(240,0){\makebox(120,12){ce: $\theta=\pi/6,\;\rho=\pi/3$}}
\end{picture}
\label{P5}\end{figure}

\begin{figure}[tbp]
\setlength{\unitlength}{1bp}
\begin{picture}(360,524)
\put(0,405){\includegraphics[%
height=119bp,width=120bp]{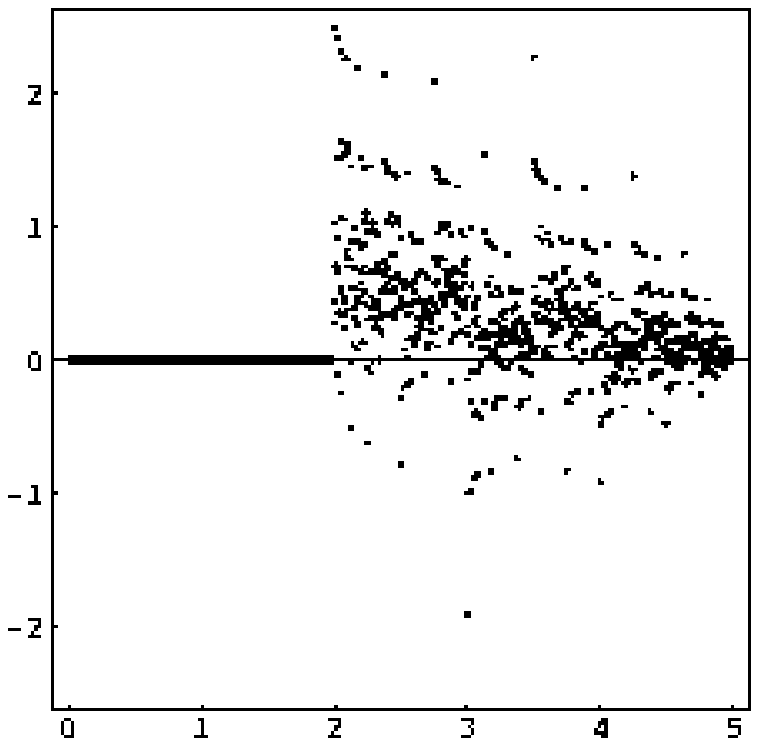}}
\put(120,405){\includegraphics[%
height=119bp,width=120bp]{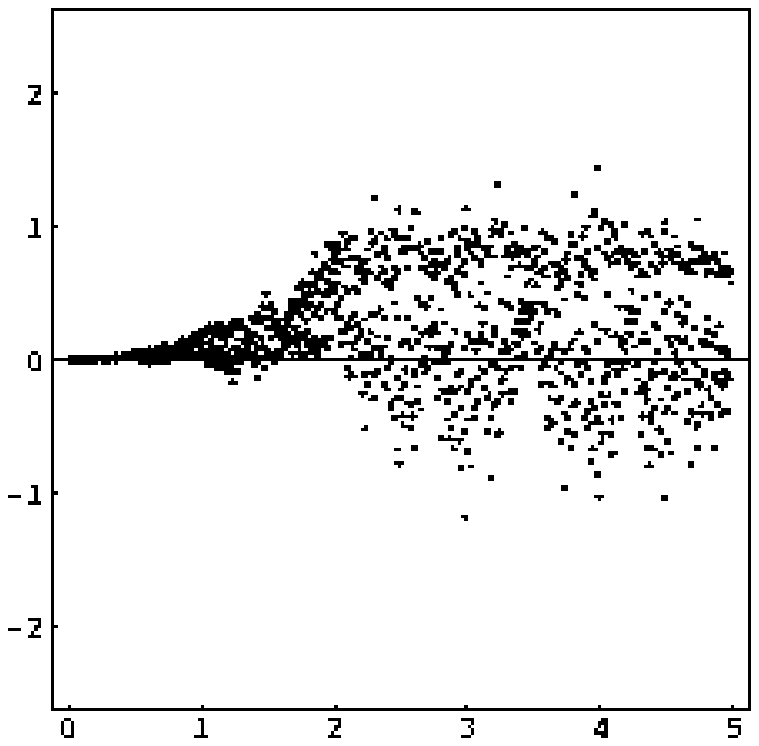}}
\put(240,405){\includegraphics[%
height=119bp,width=120bp]{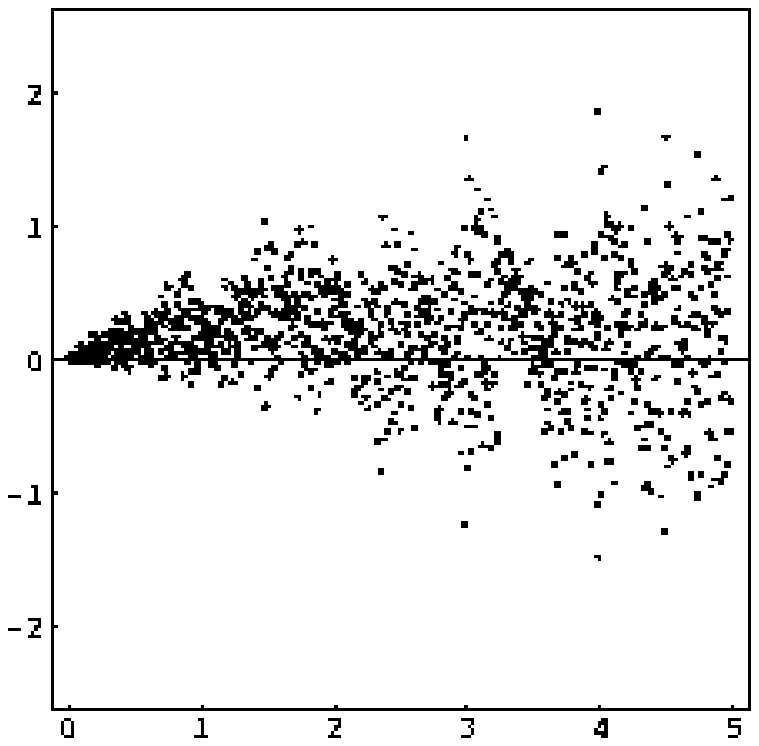}}
\put(0,393){\makebox(120,12){dh: $\theta=\pi/4,\;\rho=7\pi/12$}}
\put(120,393){\makebox(120,12){eh: $\theta=\pi/3,\;\rho=7\pi/12$}}
\put(240,393){\makebox(120,12){fh: $\theta=5\pi/12,\;\rho=7\pi/12$}}
\put(0,274){\includegraphics[%
height=119bp,width=120bp]{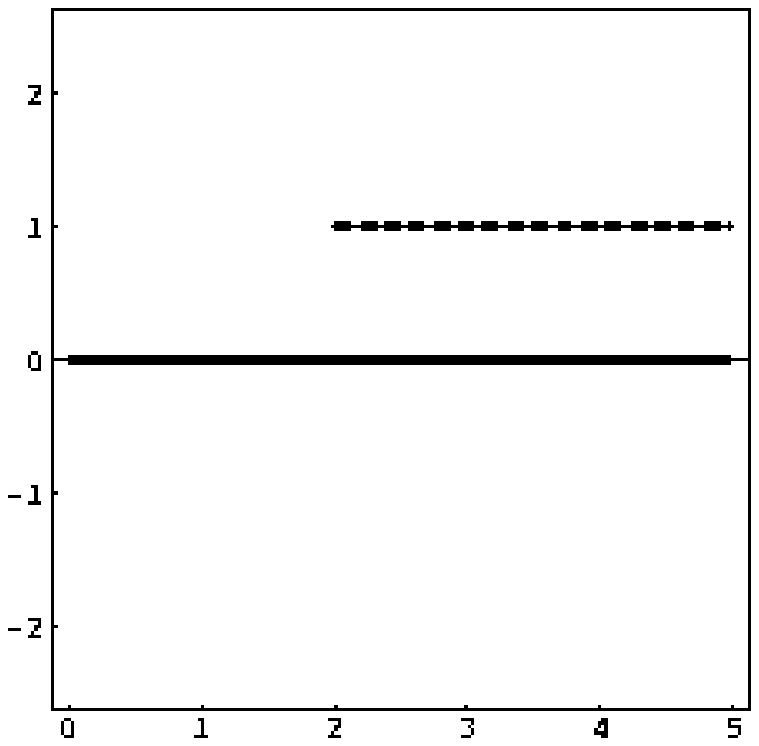}}
\put(120,274){\includegraphics[%
height=119bp,width=120bp]{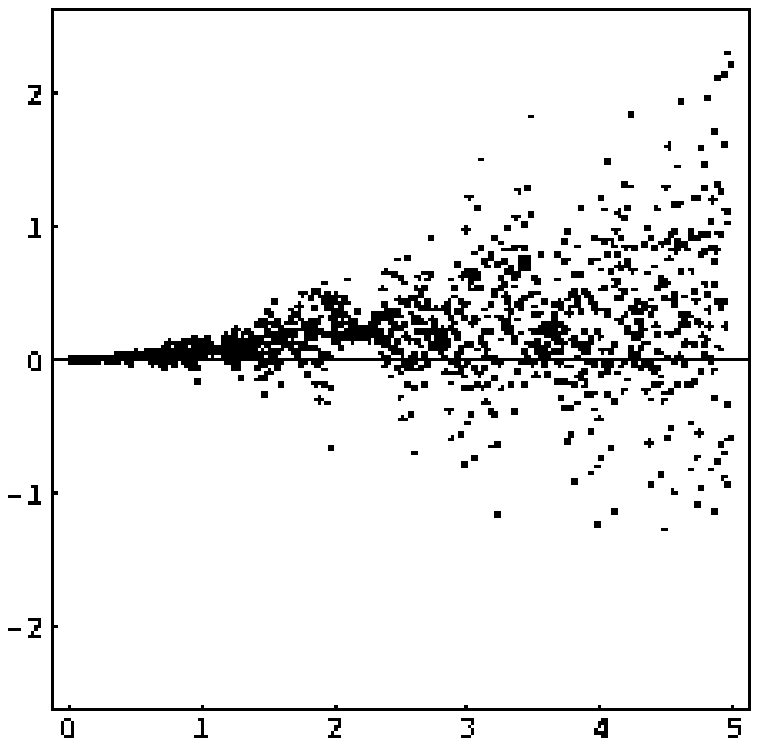}}
\put(240,274){\includegraphics[%
height=119bp,width=120bp]{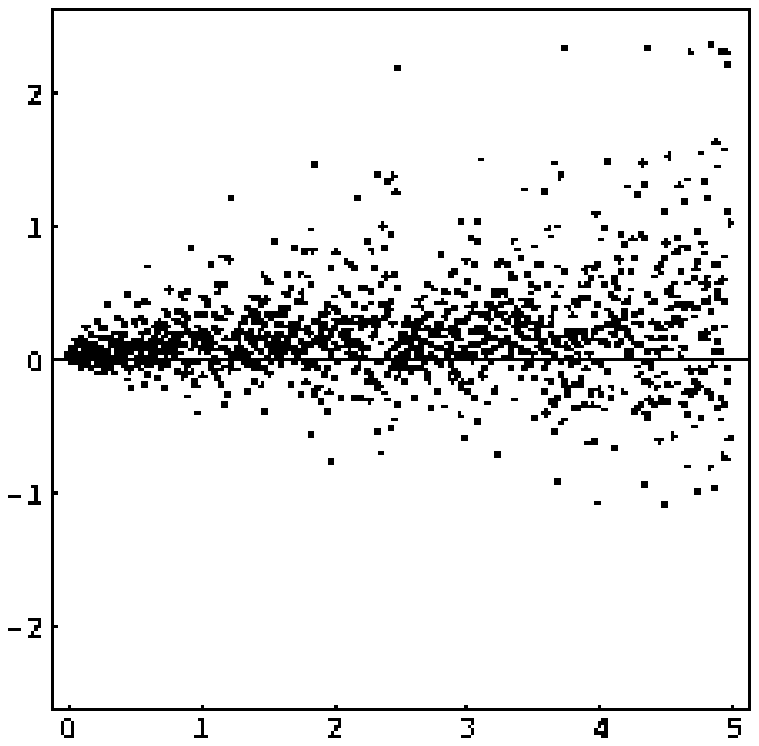}}
\put(0,262){\makebox(120,12){dg: $\theta=\pi/4,\;\rho=\pi/2$}}
\put(120,262){\makebox(120,12){eg: $\theta=\pi/3,\;\rho=\pi/2$}}
\put(240,262){\makebox(120,12){fg: $\theta=5\pi/12,\;\rho=\pi/2$}}
\put(0,143){\includegraphics[%
height=119bp,width=120bp]{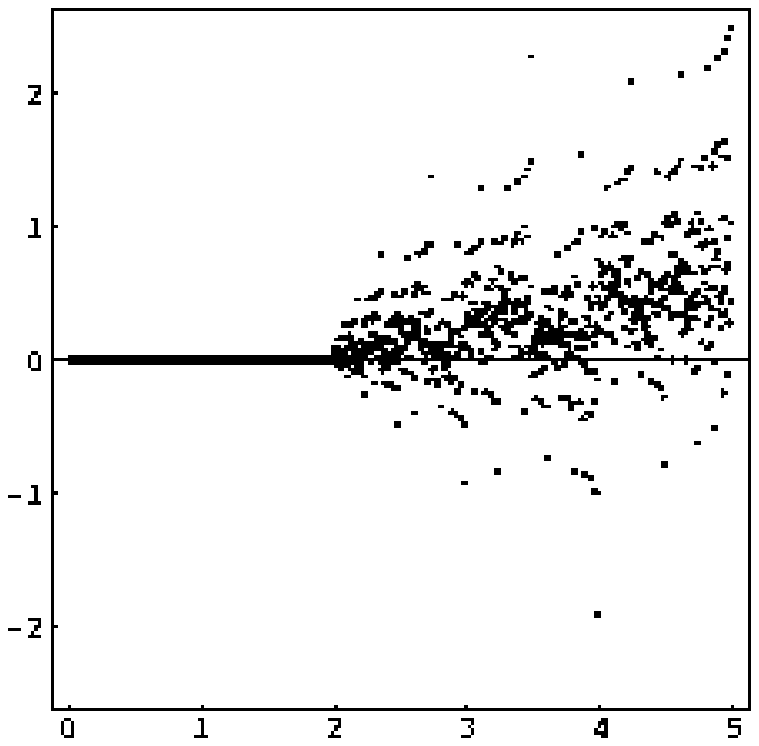}}
\put(120,143){\includegraphics[%
height=119bp,width=120bp]{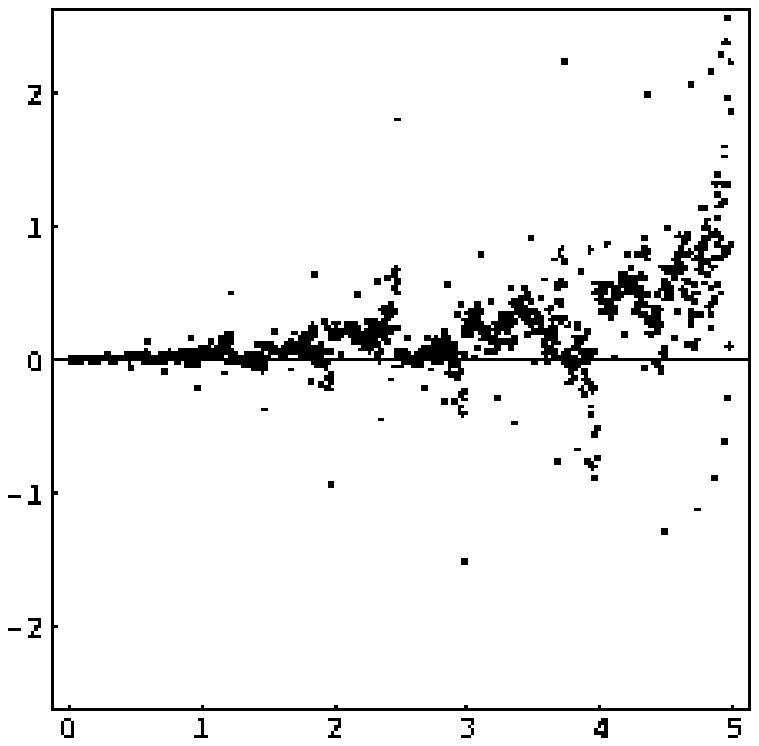}}
\put(240,143){\includegraphics[%
height=119bp,width=120bp]{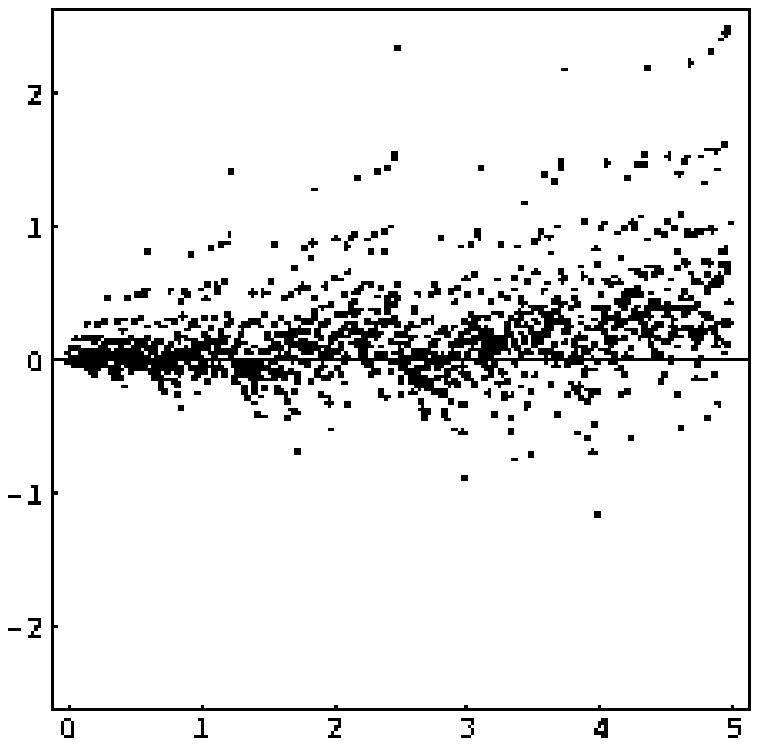}}
\put(0,131){\makebox(120,12){df: $\theta=\pi/4,\;\rho=5\pi/12$}}
\put(120,131){\makebox(120,12){ef: $\theta=\pi/3,\;\rho=5\pi/12$}}
\put(240,131){\makebox(120,12){ff: $\theta=5\pi/12,\;\rho=5\pi/12$}}
\put(0,12){\includegraphics[%
height=119bp,width=120bp]{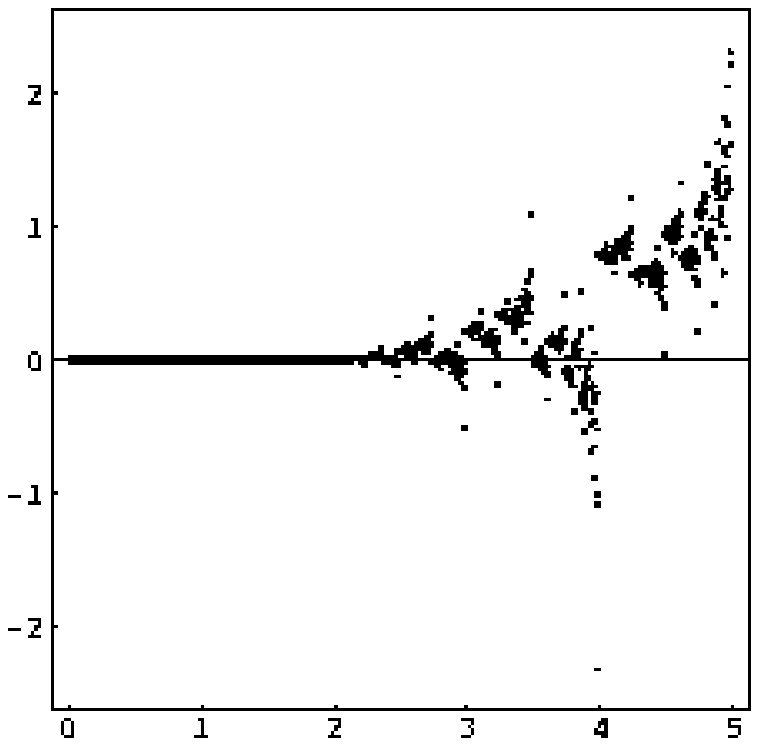}}
\put(120,12){\includegraphics[%
height=119bp,width=120bp]{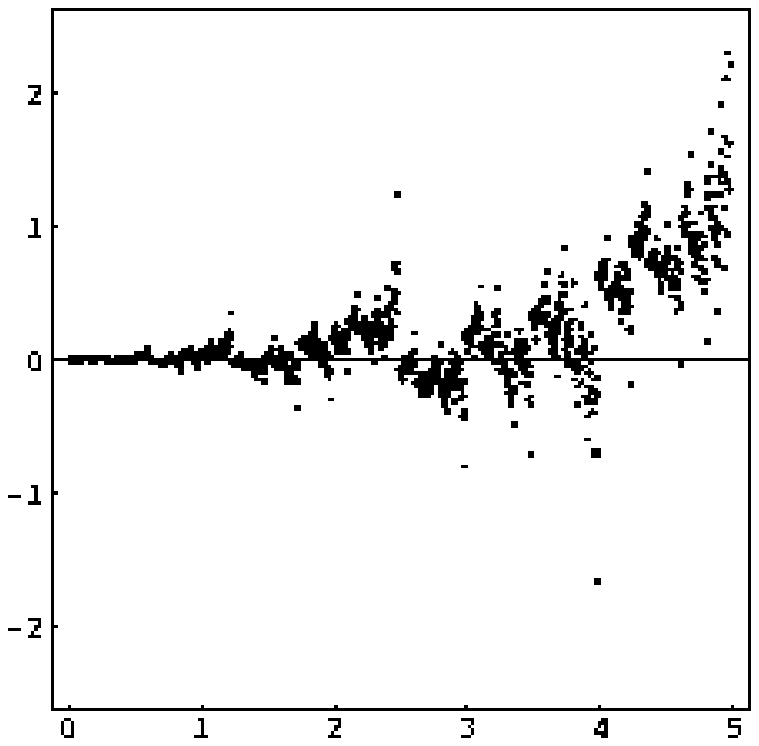}}
\put(240,12){\includegraphics[%
height=119bp,width=120bp]{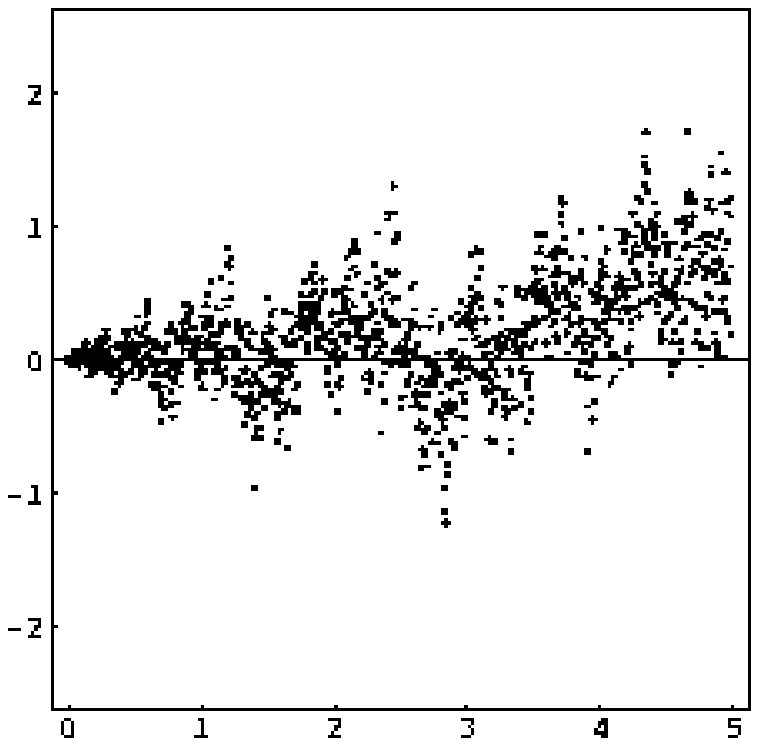}}
\put(0,0){\makebox(120,12){de: $\theta=\pi/4,\;\rho=\pi/3$}}
\put(120,0){\makebox(120,12){ee: $\theta=\pi/3,\;\rho=\pi/3$}}
\put(240,0){\makebox(120,12){fe: $\theta=5\pi/12,\;\rho=\pi/3$}}
\end{picture}
\label{P6}\end{figure}

\begin{figure}[tbp]
\setlength{\unitlength}{1bp}
\begin{picture}(360,524)
\put(0,405){\includegraphics[%
height=119bp,width=120bp]{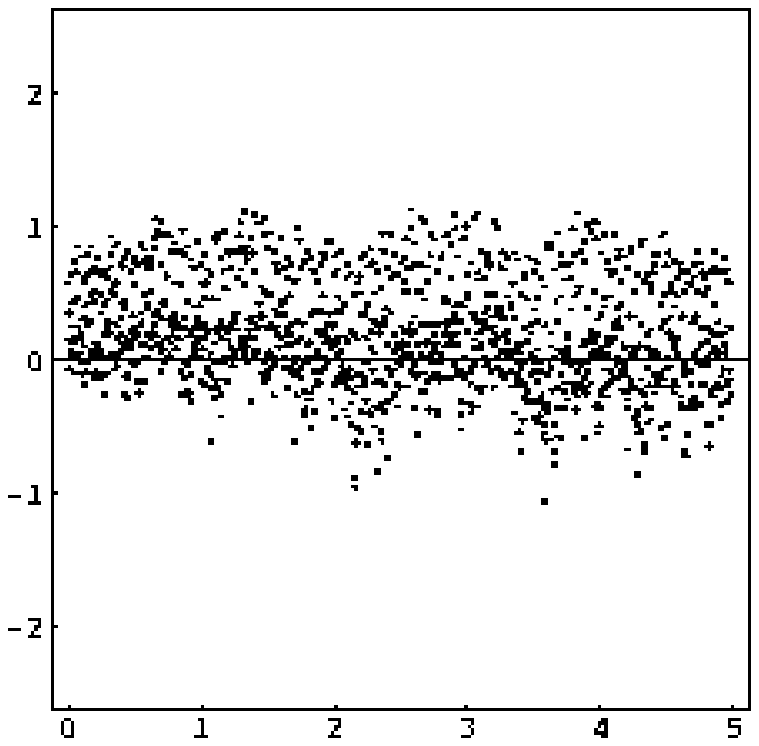}}
\put(120,405){\includegraphics[%
height=119bp,width=120bp]{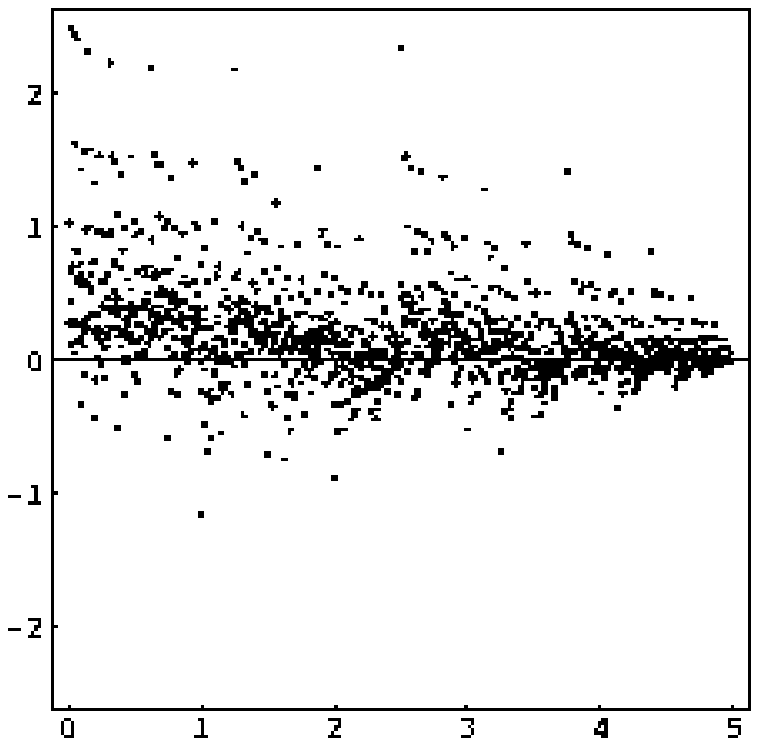}}
\put(240,405){\includegraphics[%
height=119bp,width=120bp]{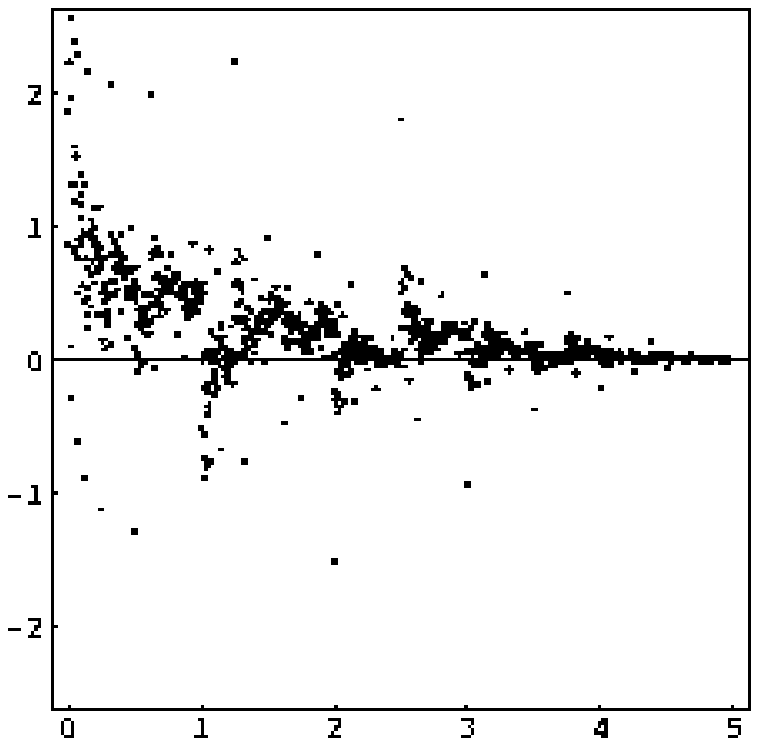}}
\put(0,393){\makebox(120,12){gh: $\theta=\pi/2,\;\rho=7\pi/12$}}
\put(120,393){\makebox(120,12){hh: $\theta=7\pi/12,\;\rho=7\pi/12$}}
\put(240,393){\makebox(120,12){ih: $\theta=2\pi/3,\;\rho=7\pi/12$}}
\put(0,274){\includegraphics[%
height=119bp,width=120bp]{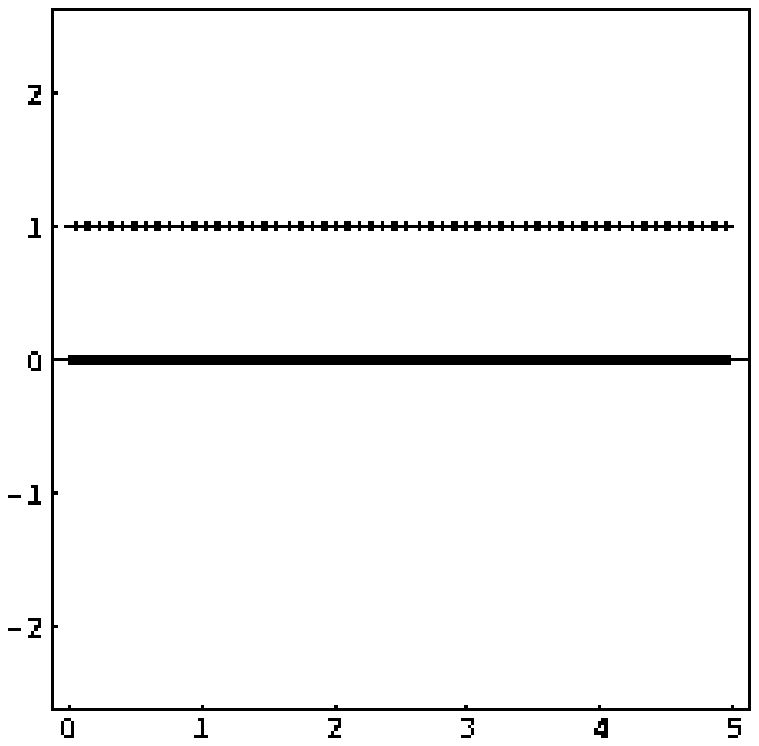}}
\put(120,274){\includegraphics[%
height=119bp,width=120bp]{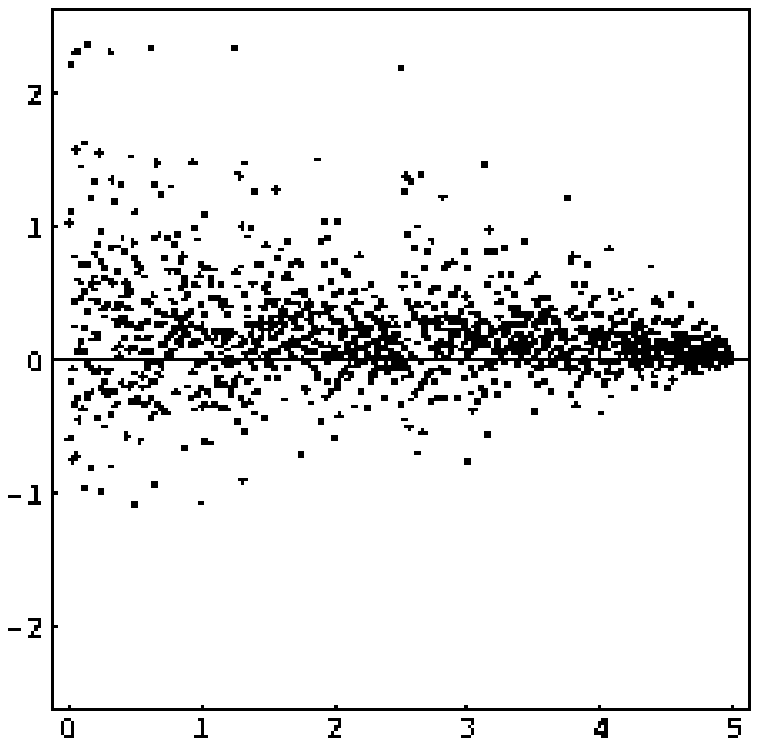}}
\put(240,274){\includegraphics[%
height=119bp,width=120bp]{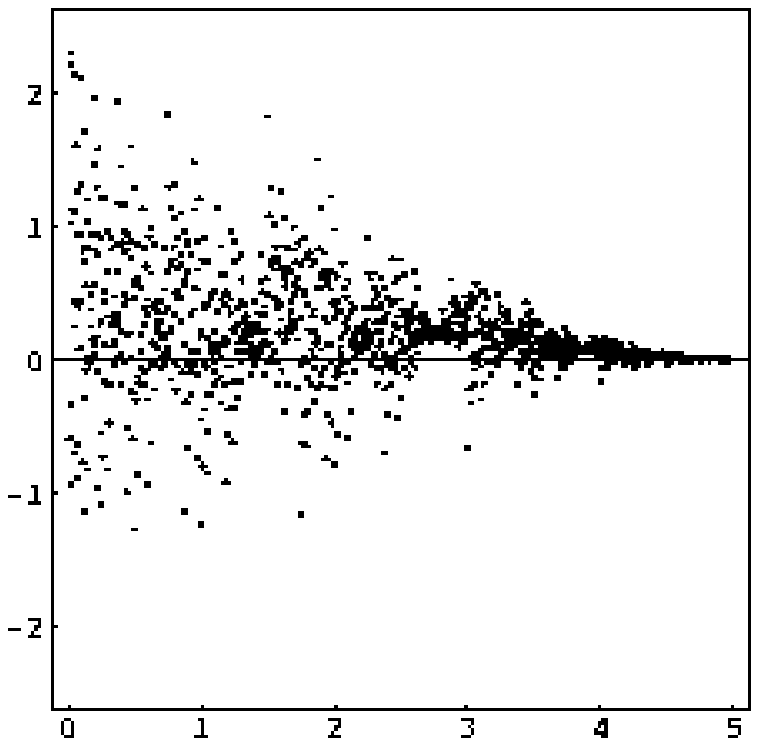}}
\put(0,262){\makebox(120,12){gg: $\theta=\pi/2,\;\rho=\pi/2$}}
\put(120,262){\makebox(120,12){hg: $\theta=7\pi/12,\;\rho=\pi/2$}}
\put(240,262){\makebox(120,12){ig: $\theta=2\pi/3,\;\rho=\pi/2$}}
\put(0,143){\includegraphics[%
height=119bp,width=120bp]{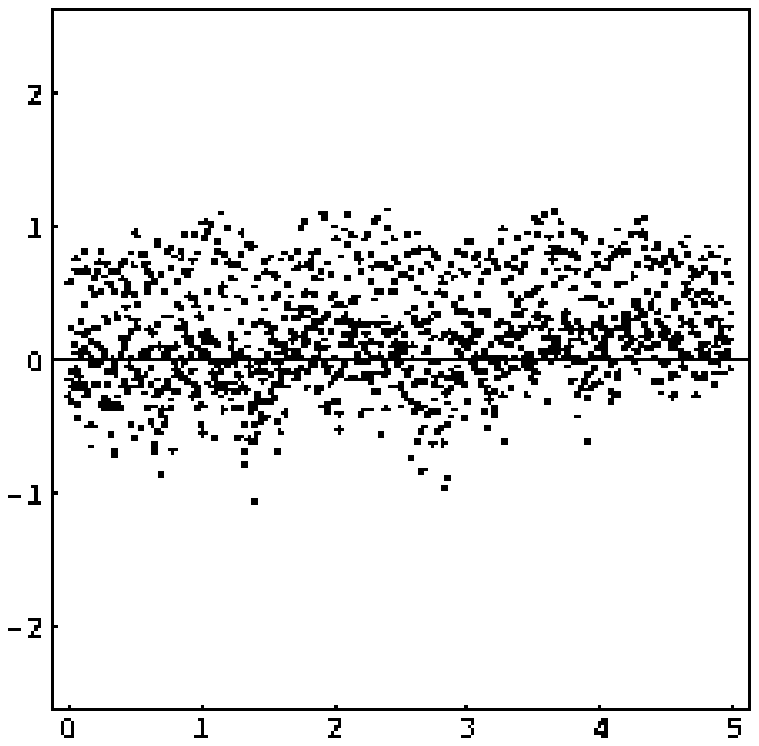}}
\put(120,143){\includegraphics[%
height=119bp,width=120bp]{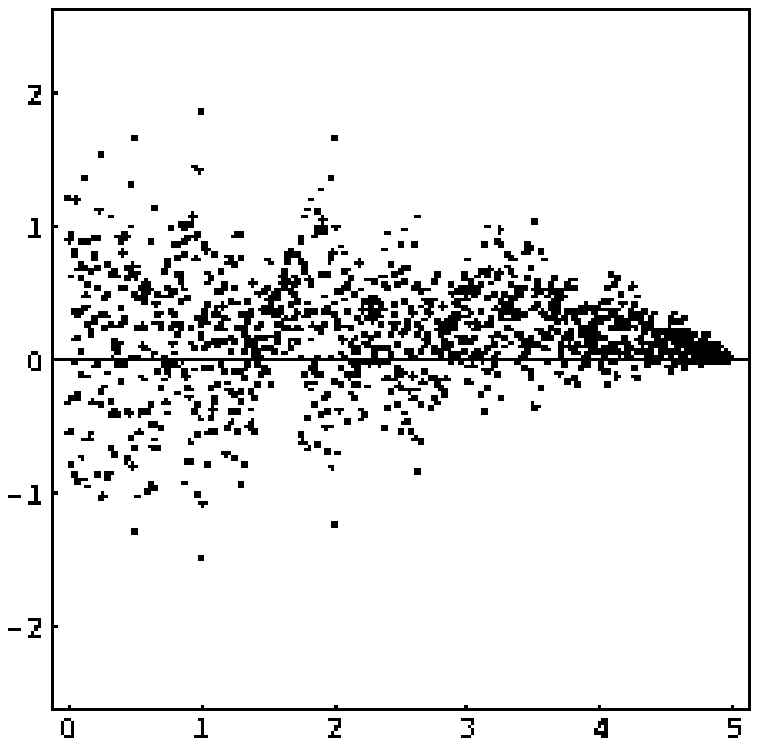}}
\put(240,143){\includegraphics[%
height=119bp,width=120bp]{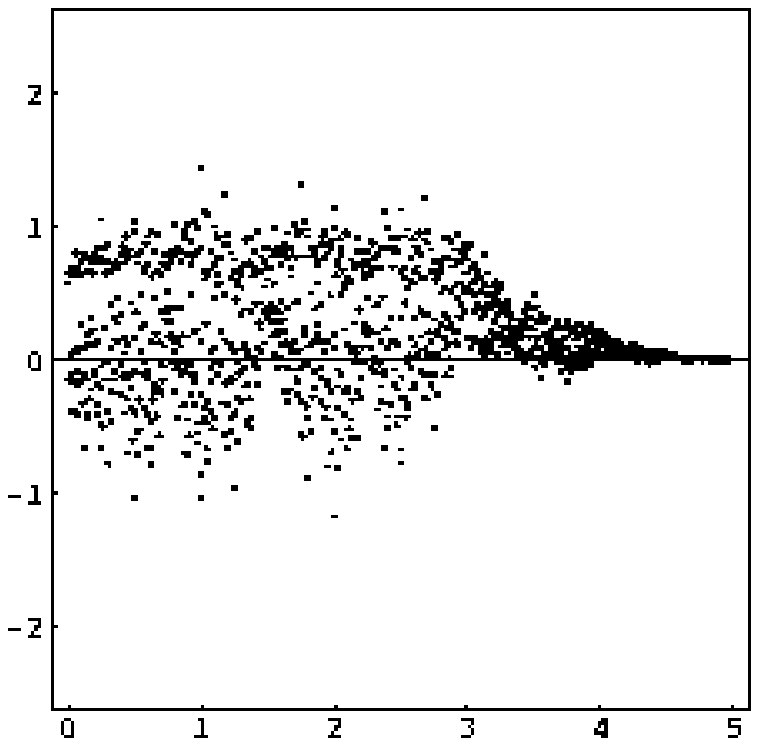}}
\put(0,131){\makebox(120,12){gf: $\theta=\pi/2,\;\rho=5\pi/12$}}
\put(120,131){\makebox(120,12){hf: $\theta=7\pi/12,\;\rho=5\pi/12$}}
\put(240,131){\makebox(120,12){if: $\theta=2\pi/3,\;\rho=5\pi/12$}}
\put(0,12){\includegraphics[%
height=119bp,width=120bp]{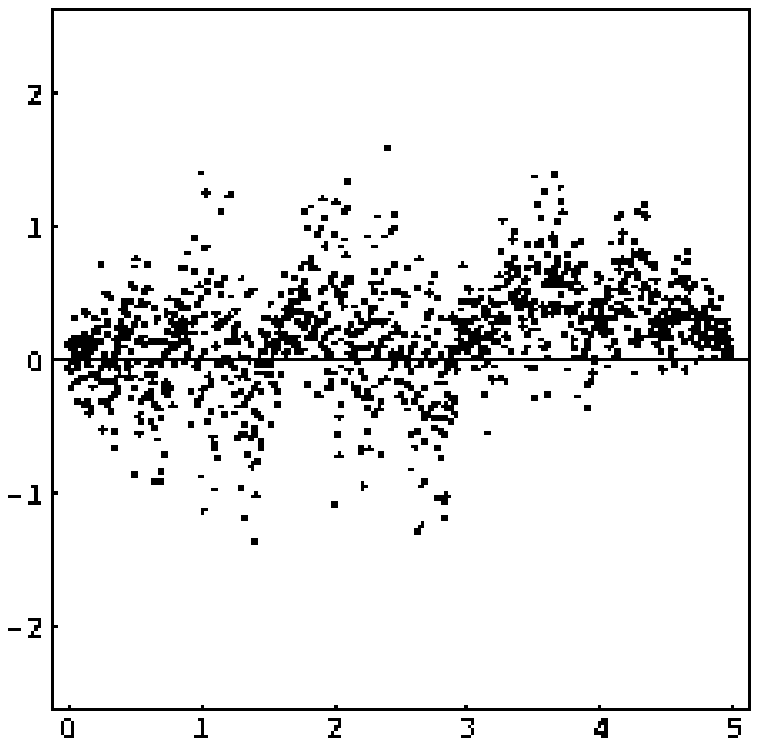}}
\put(120,12){\includegraphics[%
height=119bp,width=120bp]{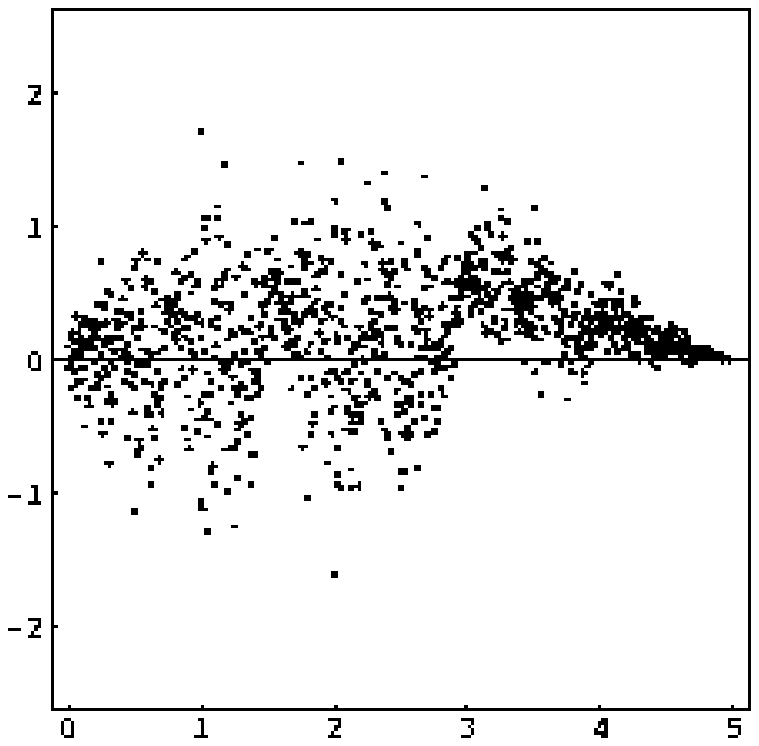}}
\put(240,12){\includegraphics[%
height=119bp,width=120bp]{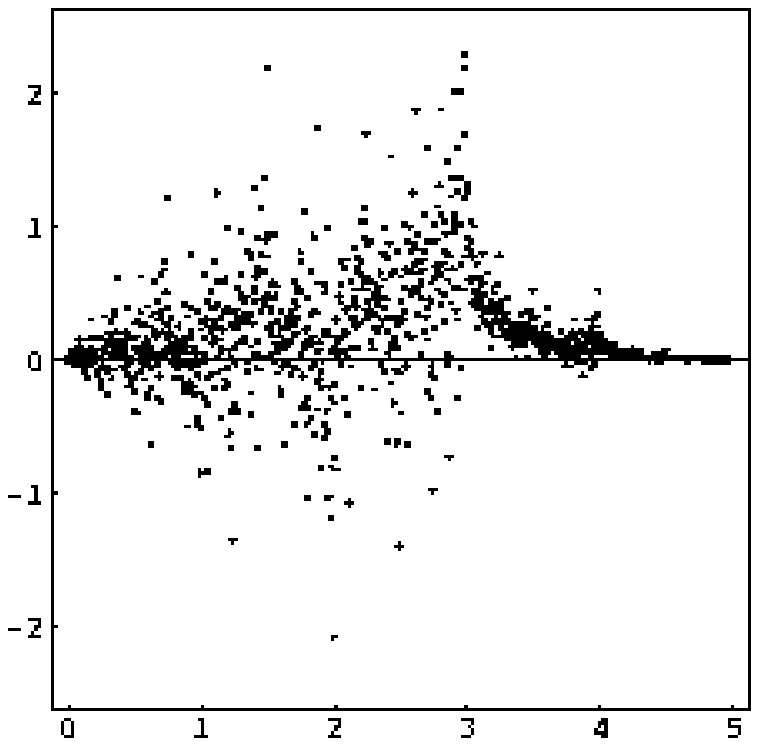}}
\put(0,0){\makebox(120,12){ge: $\theta=\pi/2,\;\rho=\pi/3$}}
\put(120,0){\makebox(120,12){he: $\theta=7\pi/12,\;\rho=\pi/3$}}
\put(240,0){\makebox(120,12){ie: $\theta=2\pi/3,\;\rho=\pi/3$}}
\end{picture}
\label{P7}\end{figure}

\begin{figure}[tbp]
\setlength{\unitlength}{1bp}
\begin{picture}(360,524)
\put(0,405){\includegraphics[%
height=119bp,width=120bp]{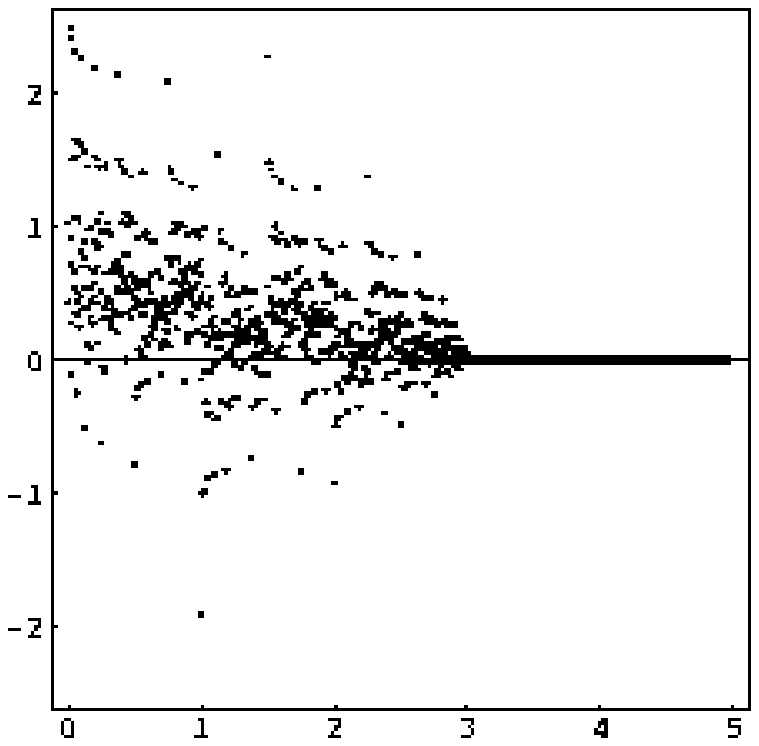}}
\put(120,405){\includegraphics[%
height=119bp,width=120bp]{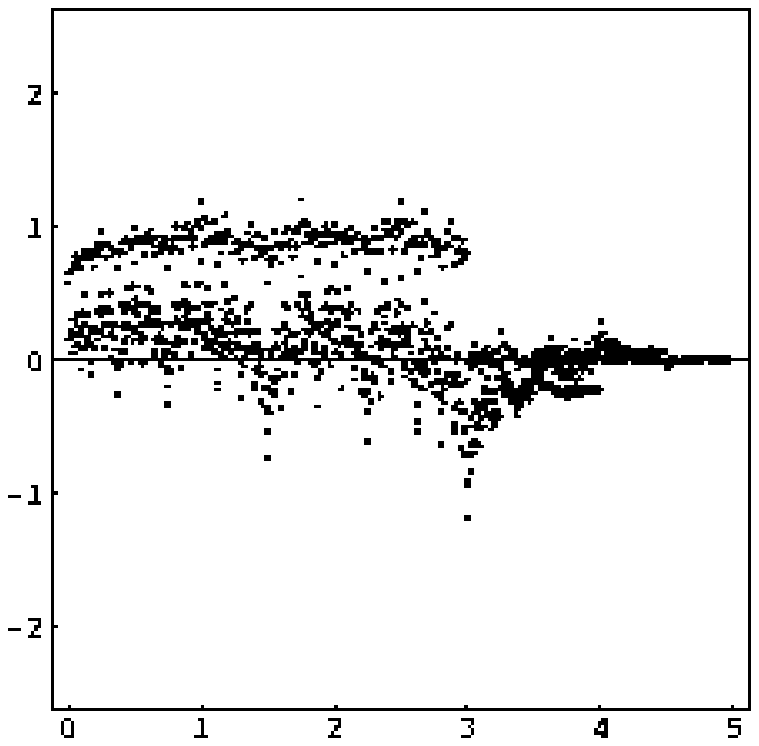}}
\put(240,405){\includegraphics[%
height=119bp,width=120bp]{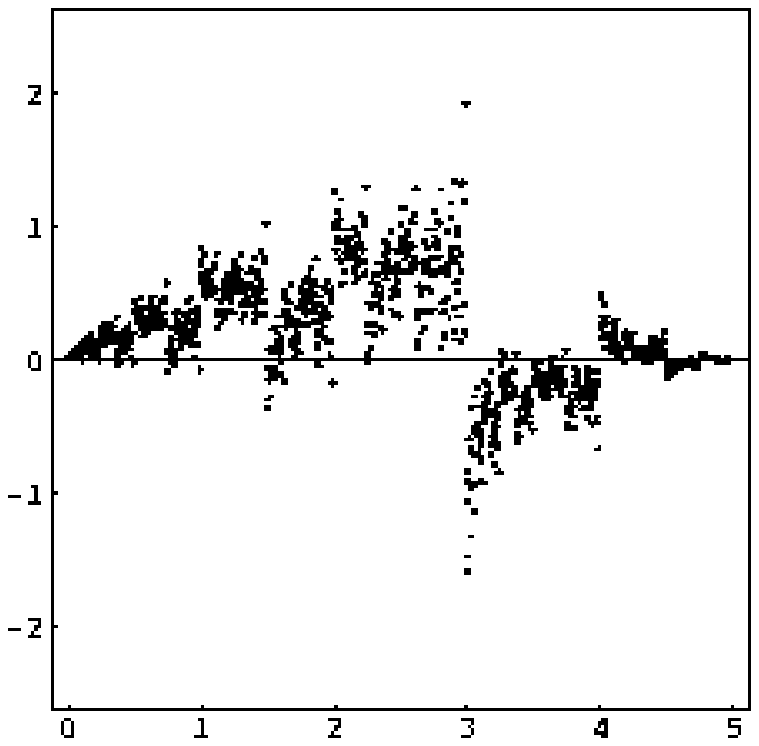}}
\put(0,393){\makebox(120,12){jh: $\theta=3\pi/4,\;\rho=7\pi/12$}}
\put(120,393){\makebox(120,12){kh: $\theta=5\pi/6,\;\rho=7\pi/12$}}
\put(240,393){\makebox(120,12){lh: $\theta=11\pi/12,\;\rho=7\pi/12$}}
\put(0,274){\includegraphics[%
height=119bp,width=120bp]{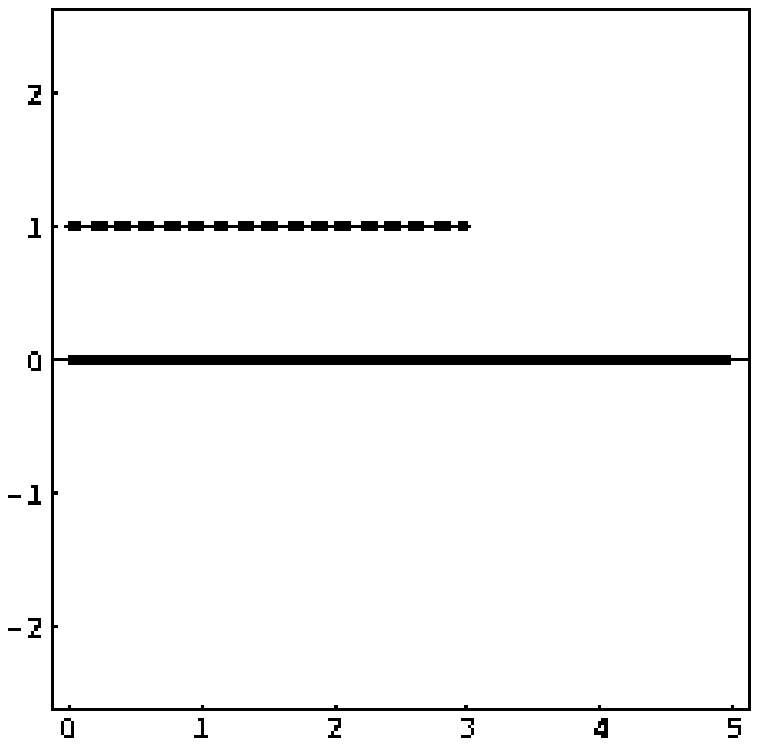}}
\put(120,274){\includegraphics[%
height=119bp,width=120bp]{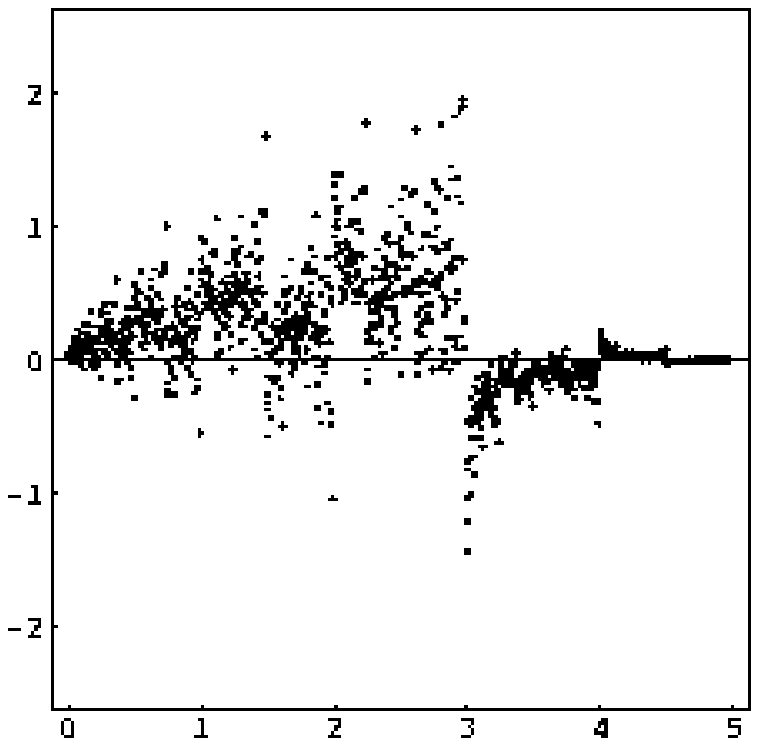}}
\put(240,274){\includegraphics[%
height=119bp,width=120bp]{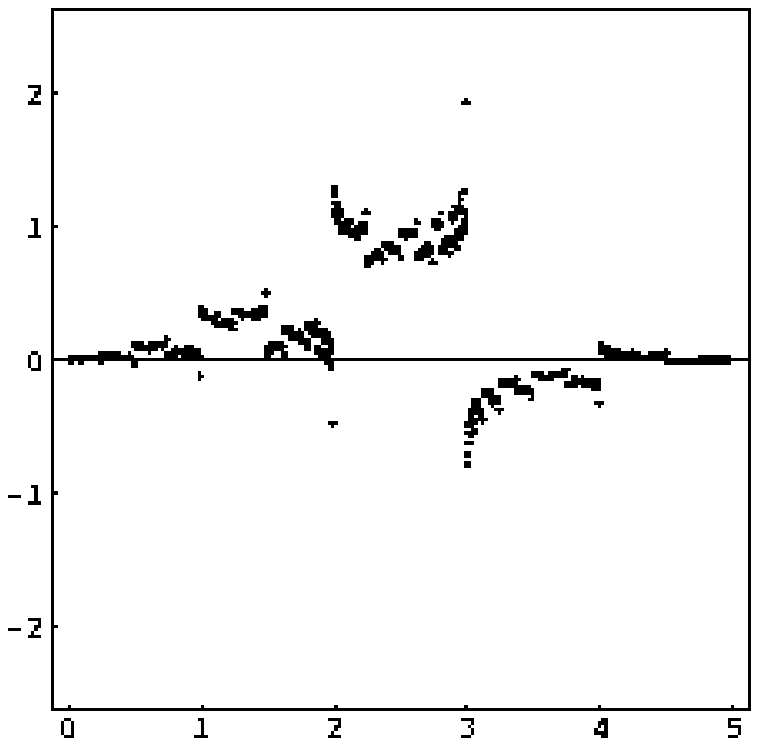}}
\put(0,262){\makebox(120,12){jg: $\theta=3\pi/4,\;\rho=\pi/2$}}
\put(120,262){\makebox(120,12){kg: $\theta=5\pi/6,\;\rho=\pi/2$}}
\put(240,262){\makebox(120,12){lg: $\theta=11\pi/12,\;\rho=\pi/2$}}
\put(0,143){\includegraphics[%
height=119bp,width=120bp]{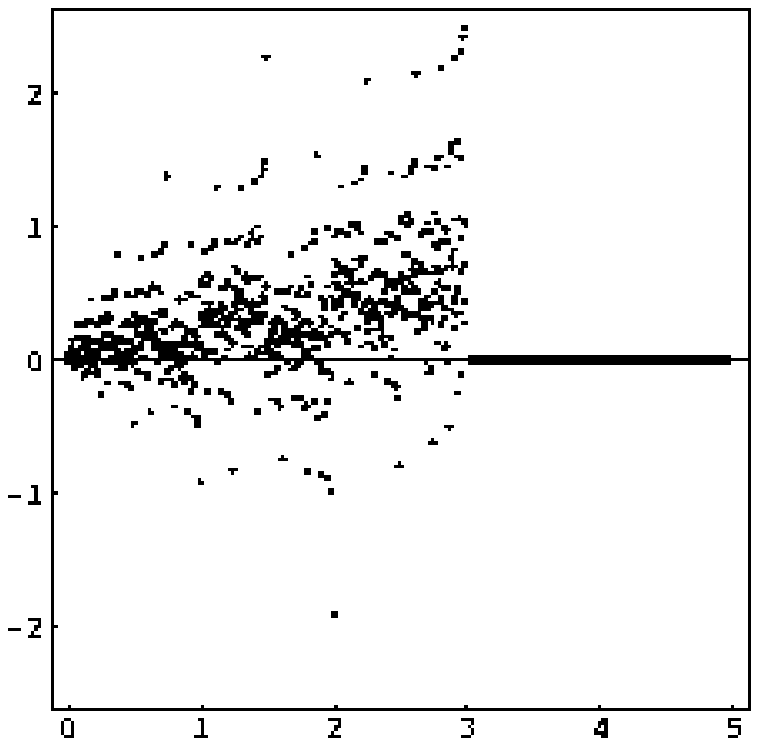}}
\put(120,143){\includegraphics[%
height=119bp,width=120bp]{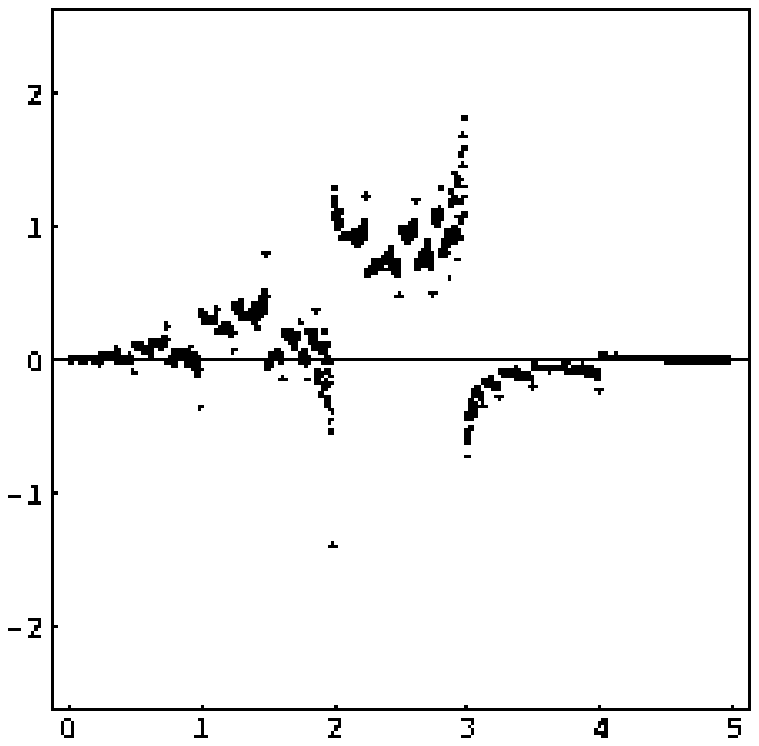}}
\put(240,143){\includegraphics[%
height=119bp,width=120bp]{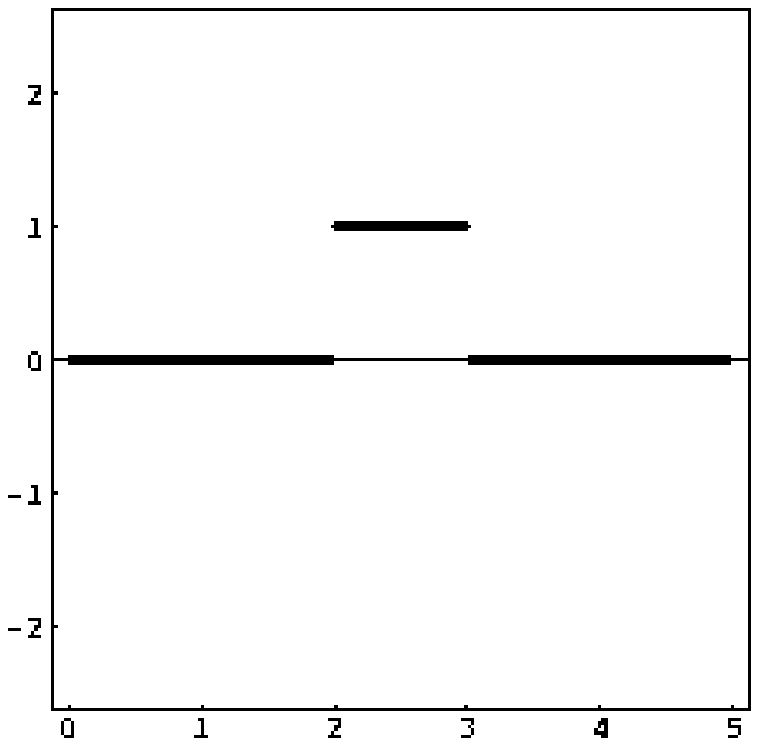}}
\put(0,131){\makebox(120,12){jf: $\theta=3\pi/4,\;\rho=5\pi/12$}}
\put(120,131){\makebox(120,12){kf: $\theta=5\pi/6,\;\rho=5\pi/12$}}
\put(240,131){\makebox(120,12){lf: $\theta=11\pi/12,\;\rho=5\pi/12$}}
\put(0,12){\includegraphics[%
height=119bp,width=120bp]{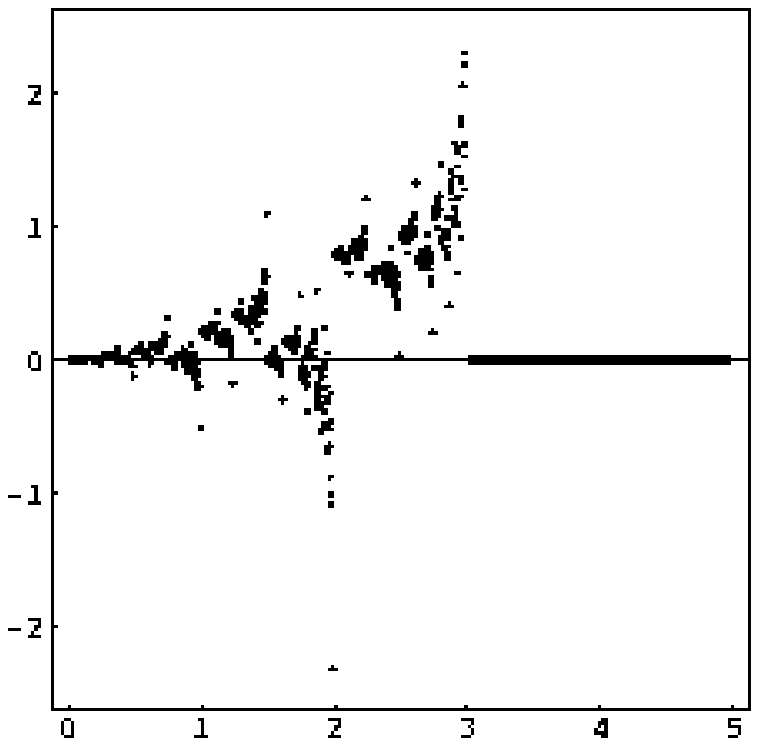}}
\put(120,12){\includegraphics[%
height=119bp,width=120bp]{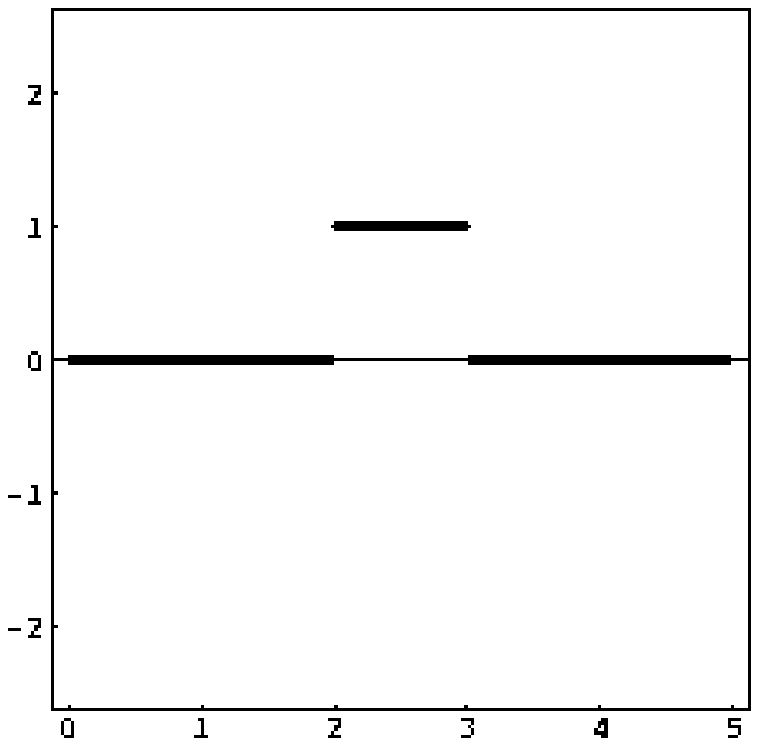}}
\put(240,12){\includegraphics[%
height=119bp,width=120bp]{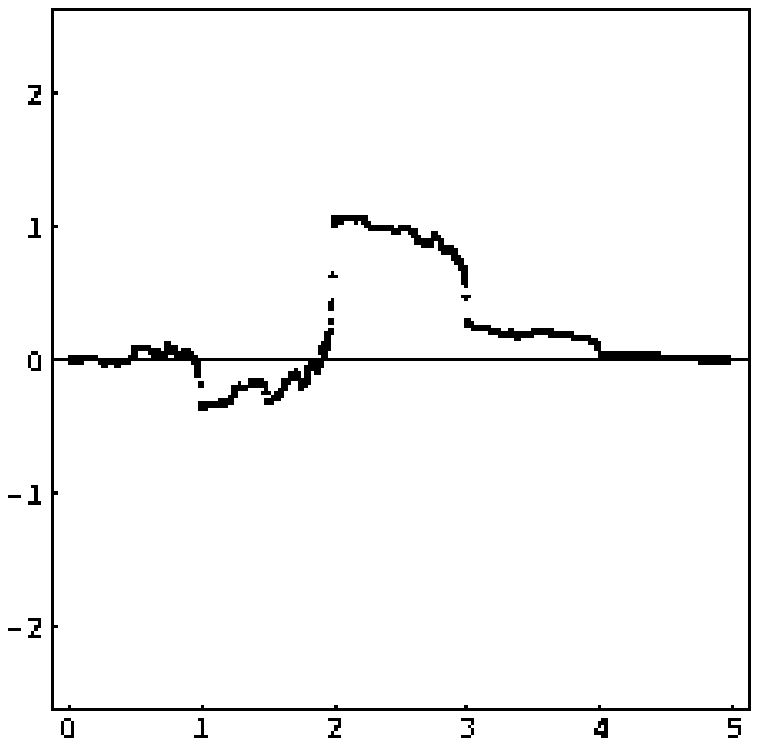}}
\put(0,0){\makebox(120,12){je: $\theta=3\pi/4,\;\rho=\pi/3$}}
\put(120,0){\makebox(120,12){ke: $\theta=5\pi/6,\;\rho=\pi/3$}}
\put(240,0){\makebox(120,12){le: $\theta=11\pi/12,\;\rho=\pi/3$}}
\end{picture}
\label{P8}\end{figure}

\begin{figure}[tbp]
\setlength{\unitlength}{1bp}
\begin{picture}(360,524)
\put(0,405){\includegraphics[%
height=119bp,width=120bp]{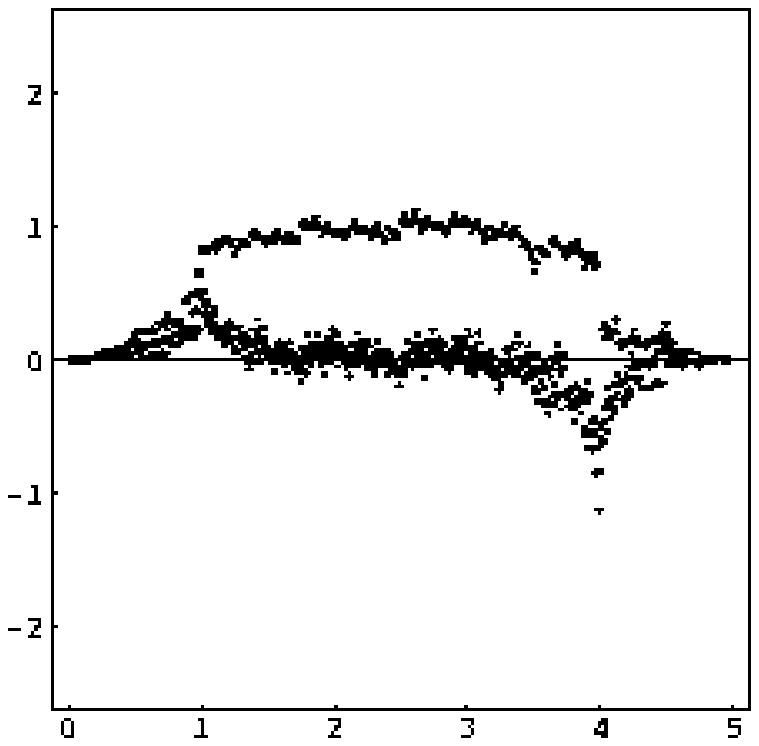}}
\put(120,405){\includegraphics[%
height=119bp,width=120bp]{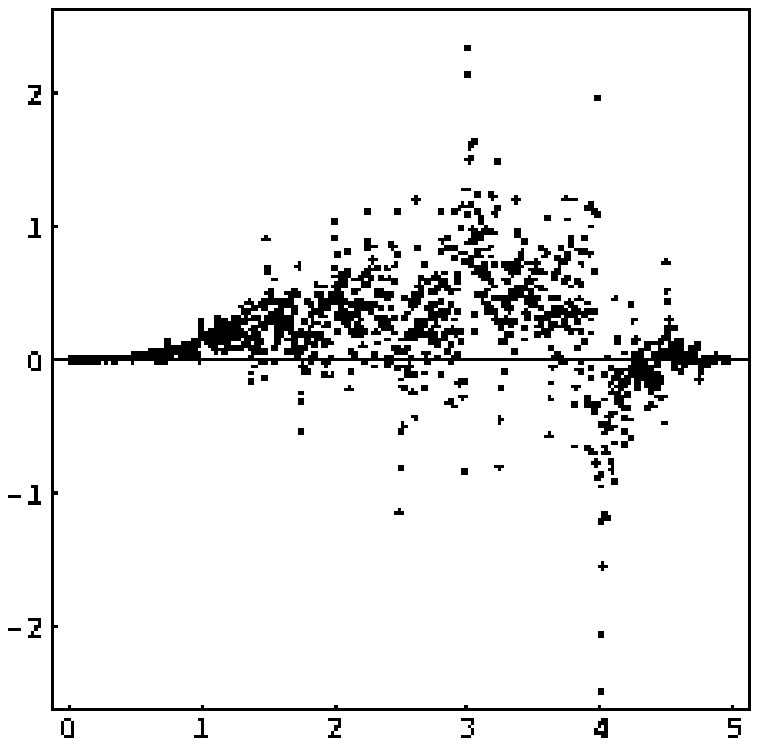}}
\put(240,405){\includegraphics[%
height=119bp,width=120bp]{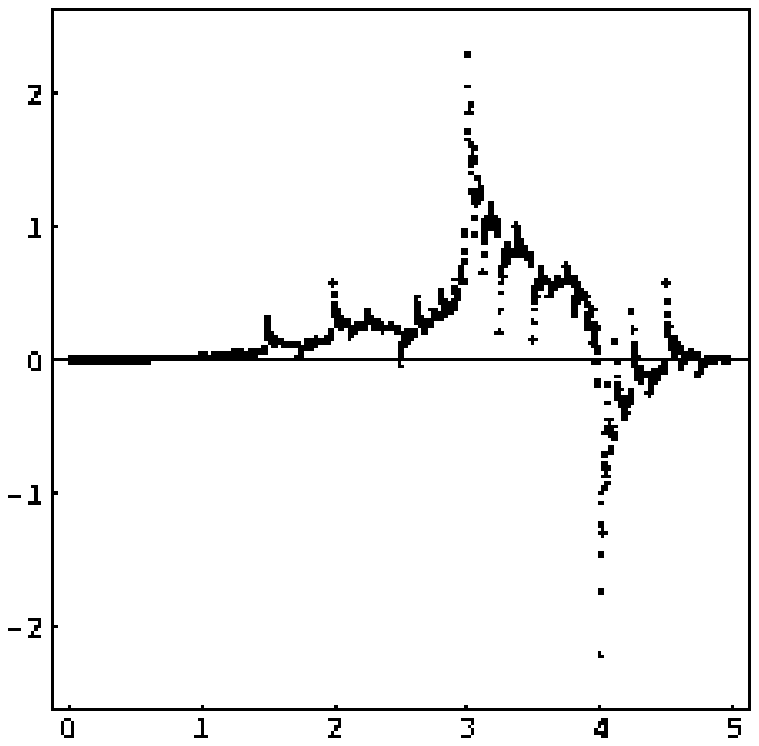}}
\put(0,393){\makebox(120,12){al: $\theta=0,\;\rho=11\pi/12$}}
\put(120,393){\makebox(120,12){bl: $\theta=\pi/12,\;\rho=11\pi/12$}}
\put(240,393){\makebox(120,12){cl: $\theta=\pi/6,\;\rho=11\pi/12$}}
\put(0,274){\includegraphics[%
height=119bp,width=120bp]{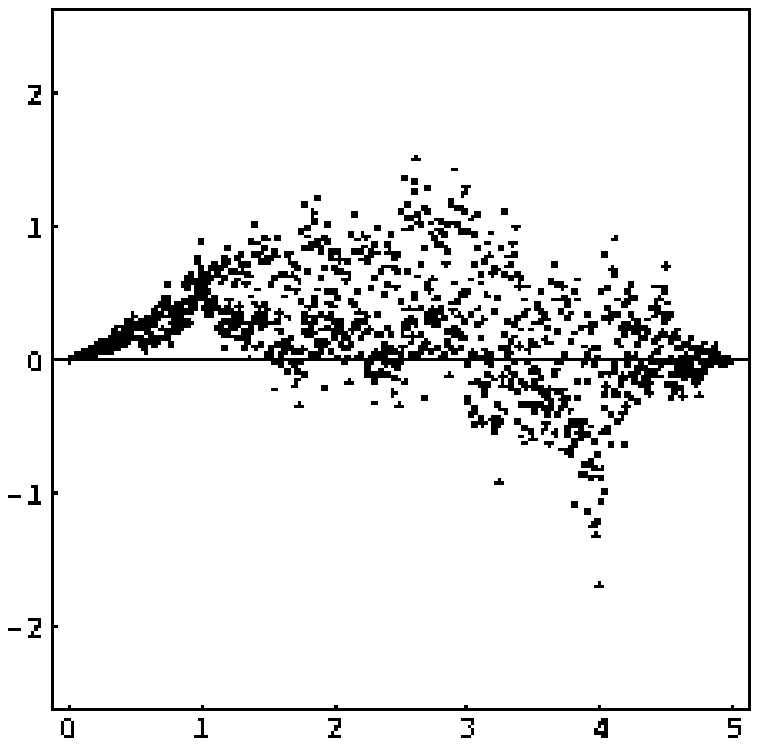}}
\put(120,274){\includegraphics[%
height=119bp,width=120bp]{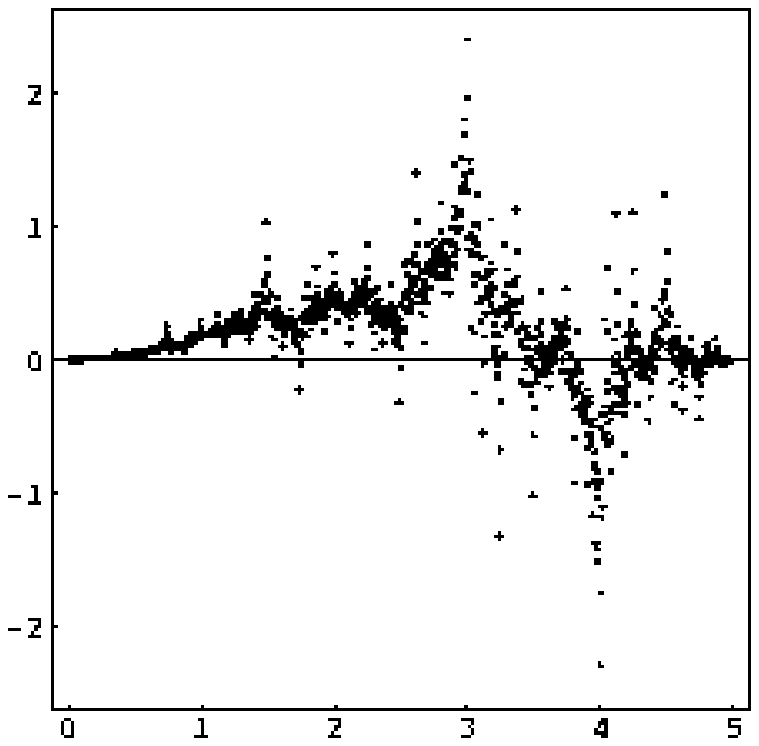}}
\put(240,274){\includegraphics[%
height=119bp,width=120bp]{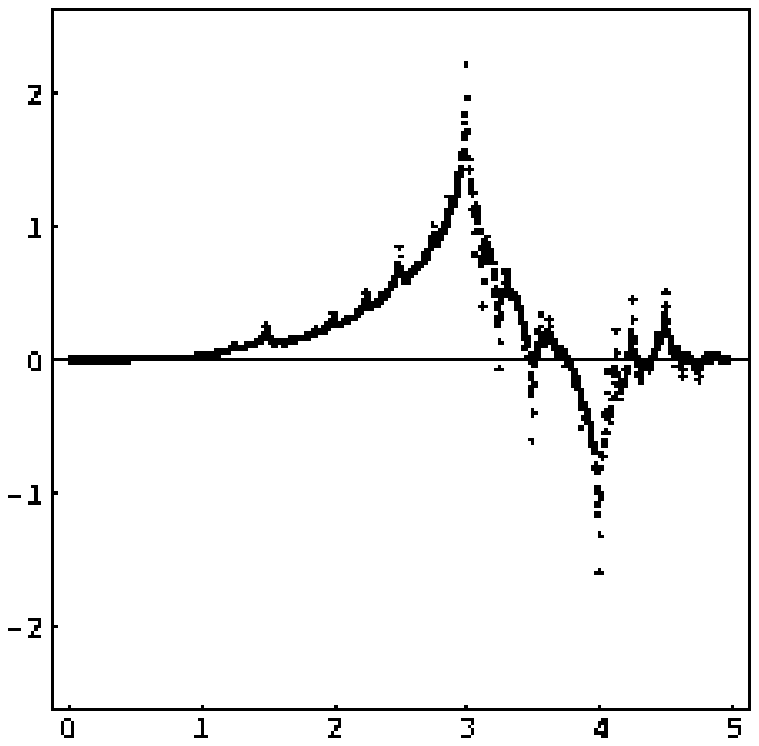}}
\put(0,262){\makebox(120,12){ak: $\theta=0,\;\rho=5\pi/6$}}
\put(120,262){\makebox(120,12){bk: $\theta=\pi/12,\;\rho=5\pi/6$}}
\put(240,262){\makebox(120,12){ck: $\theta=\pi/6,\;\rho=5\pi/6$}}
\put(0,143){\includegraphics[%
height=119bp,width=120bp]{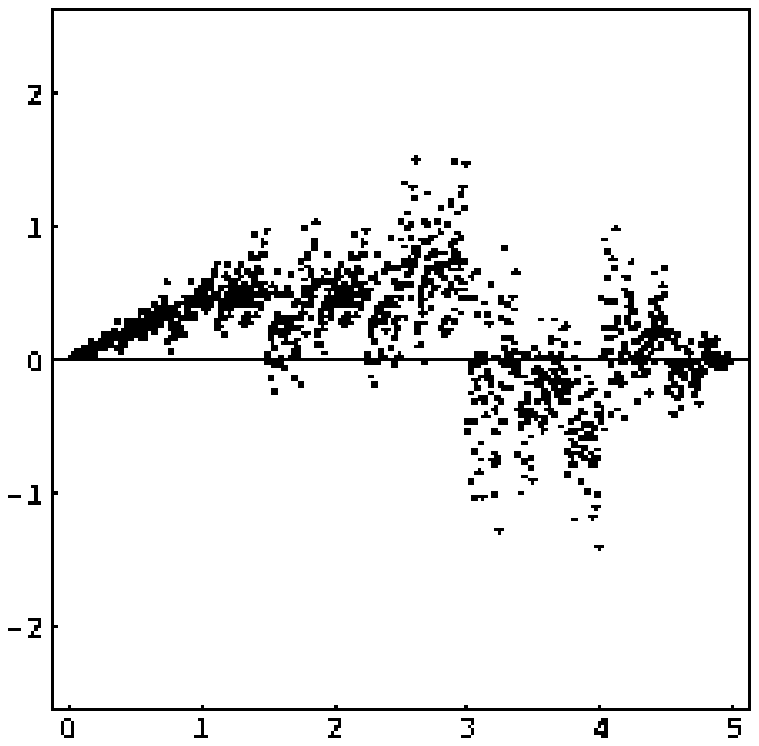}}
\put(120,143){\includegraphics[%
height=119bp,width=120bp]{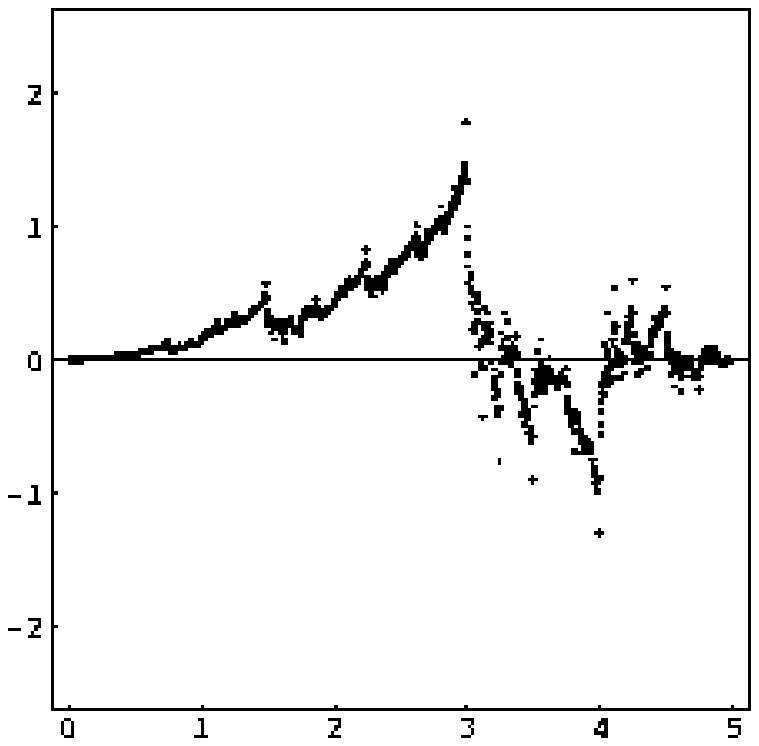}}
\put(240,143){\includegraphics[%
height=119bp,width=120bp]{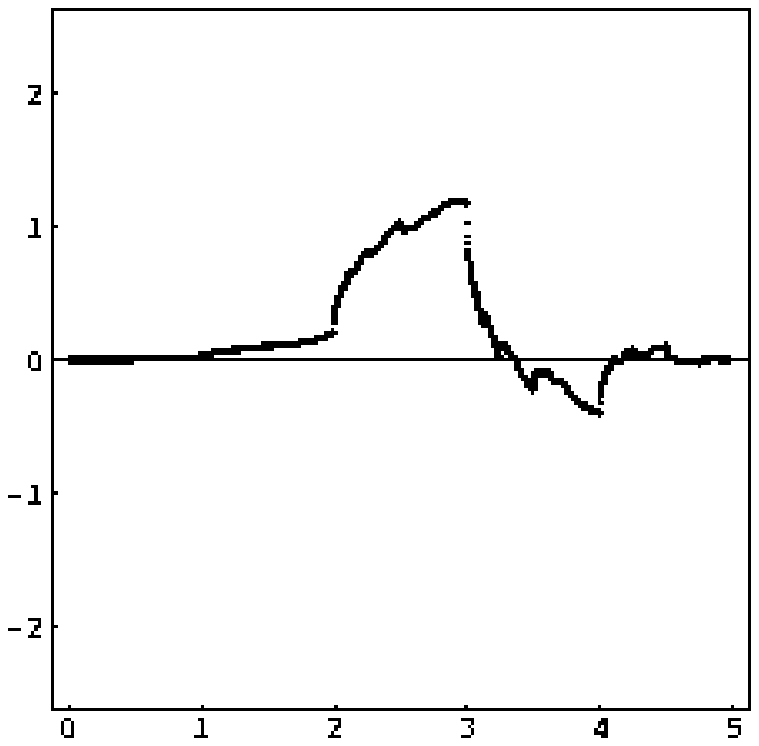}}
\put(0,131){\makebox(120,12){aj: $\theta=0,\;\rho=3\pi/4$}}
\put(120,131){\makebox(120,12){bj: $\theta=\pi/12,\;\rho=3\pi/4$}}
\put(240,131){\makebox(120,12){cj: $\theta=\pi/6,\;\rho=3\pi/4$}}
\put(0,12){\includegraphics[%
height=119bp,width=120bp]{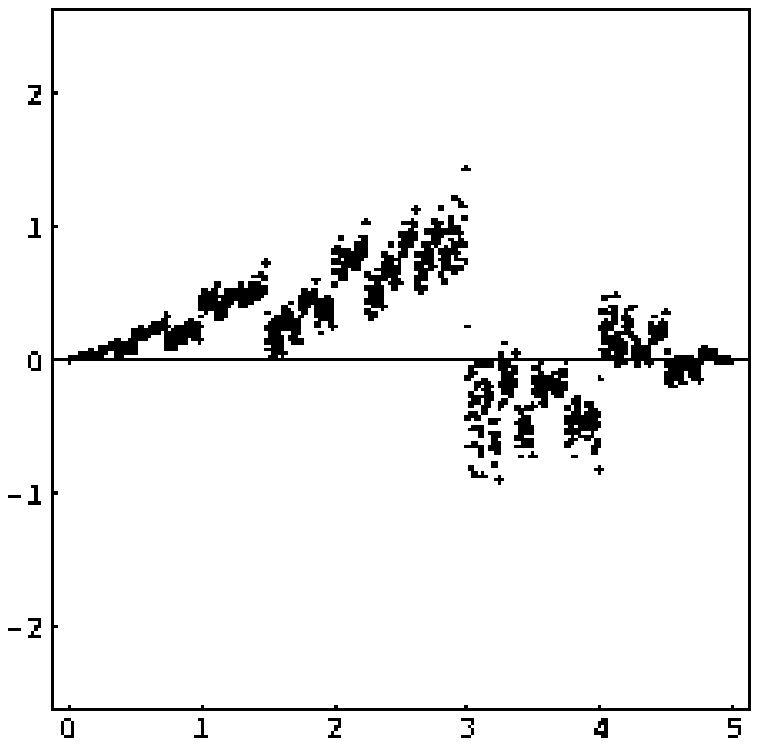}}
\put(120,12){\includegraphics[%
height=119bp,width=120bp]{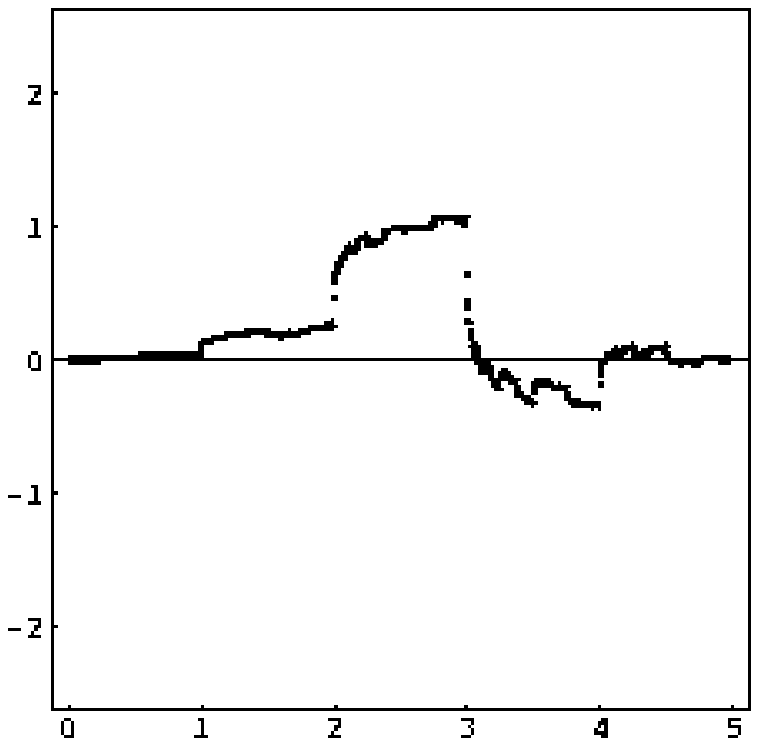}}
\put(240,12){\includegraphics[%
height=119bp,width=120bp]{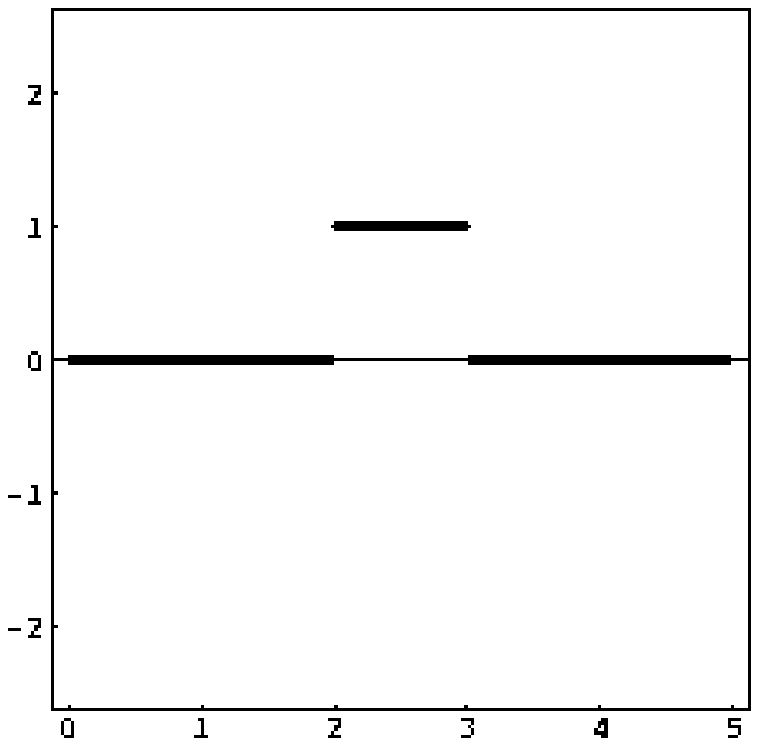}}
\put(0,0){\makebox(120,12){ai: $\theta=0,\;\rho=2\pi/3$}}
\put(120,0){\makebox(120,12){bi: $\theta=\pi/12,\;\rho=2\pi/3$}}
\put(240,0){\makebox(120,12){ci: $\theta=\pi/6,\;\rho=2\pi/3$}}
\end{picture}
\label{P9}\end{figure}

\begin{figure}[tbp]
\setlength{\unitlength}{1bp}
\begin{picture}(360,524)
\put(0,405){\includegraphics[%
height=119bp,width=120bp]{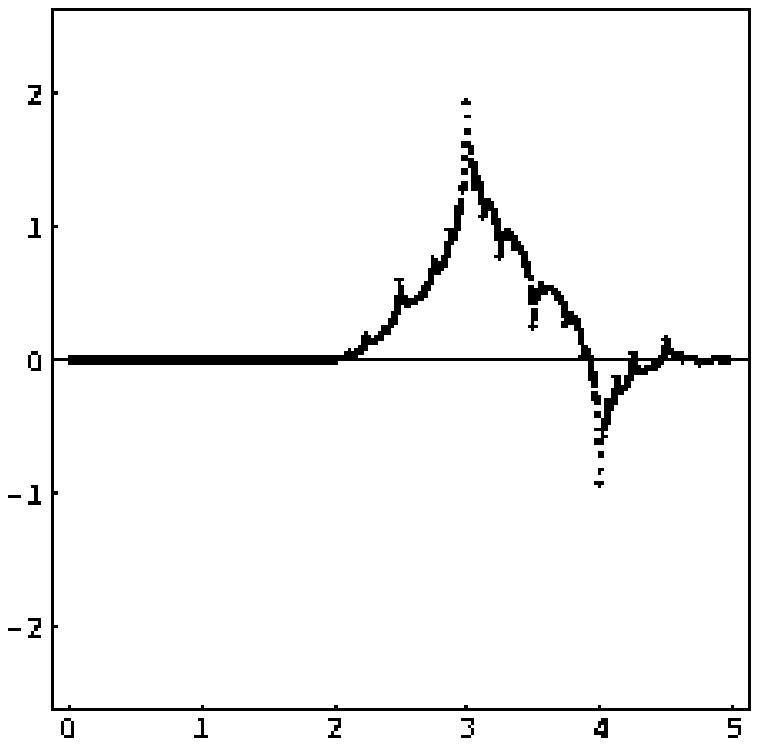}}
\put(120,405){\includegraphics[%
height=119bp,width=120bp]{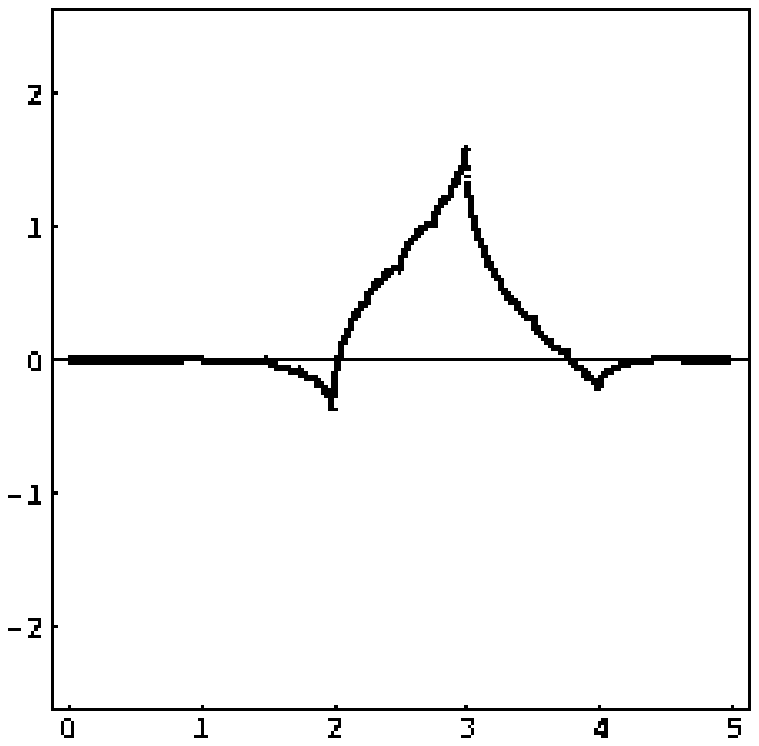}}
\put(240,405){\includegraphics[%
height=119bp,width=120bp]{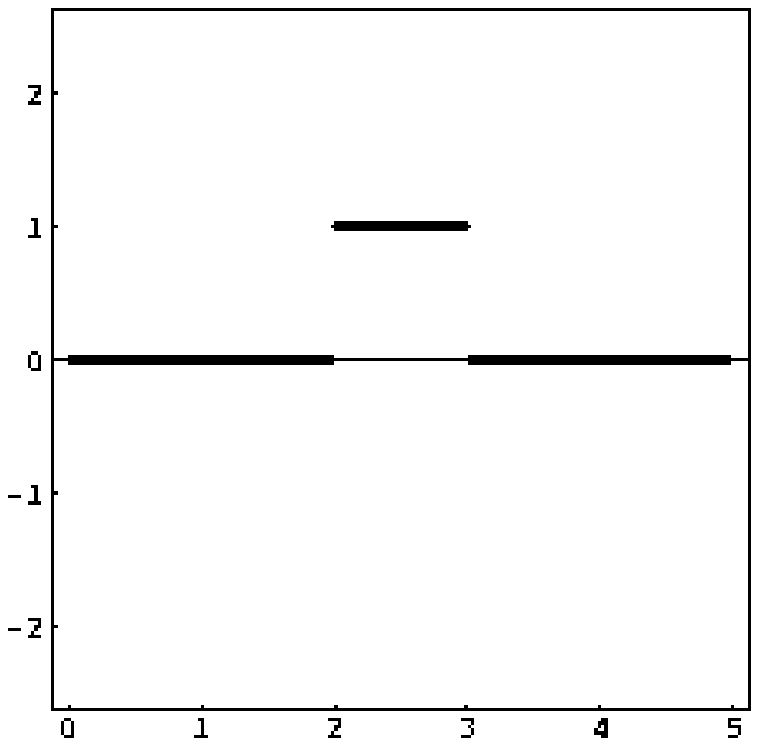}}
\put(0,393){\makebox(120,12){dl: $\theta=\pi/4,\;\rho=11\pi/12$}}
\put(120,393){\makebox(120,12){el: $\theta=\pi/3,\;\rho=11\pi/12$}}
\put(240,393){\makebox(120,12){fl: $\theta=5\pi/12,\;\rho=11\pi/12$}}
\put(0,274){\includegraphics[%
height=119bp,width=120bp]{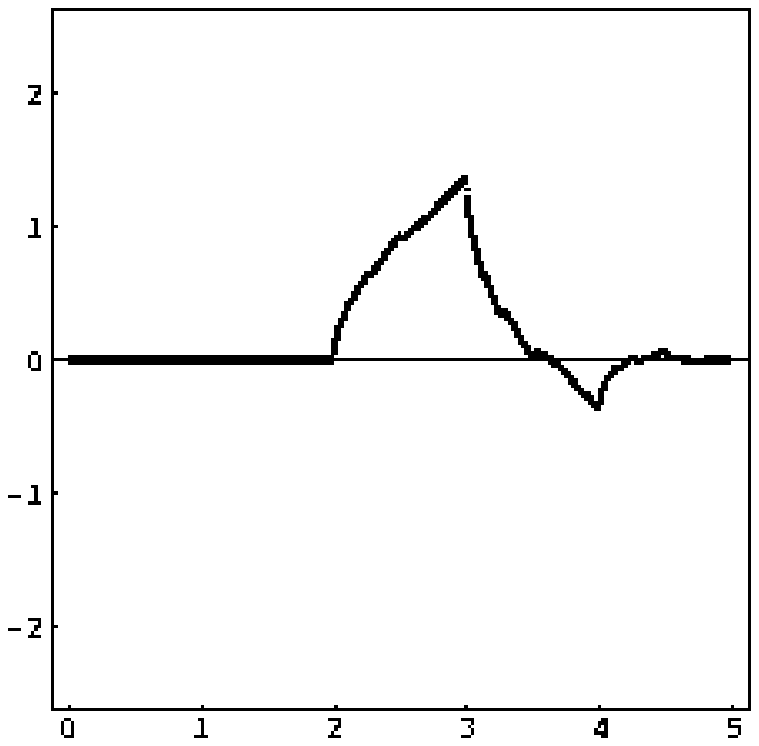}}
\put(120,274){\includegraphics[%
height=119bp,width=120bp]{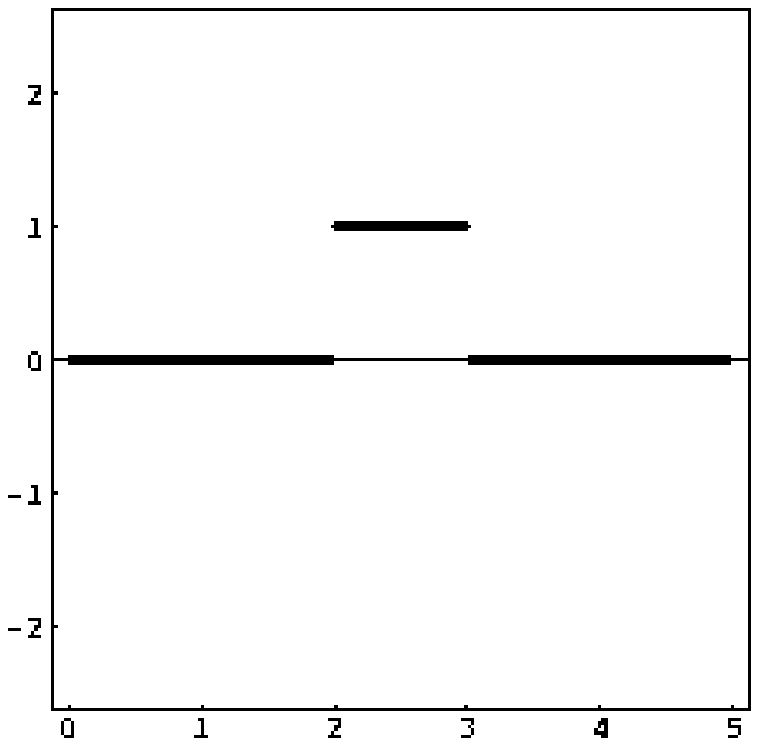}}
\put(240,274){\includegraphics[%
height=119bp,width=120bp]{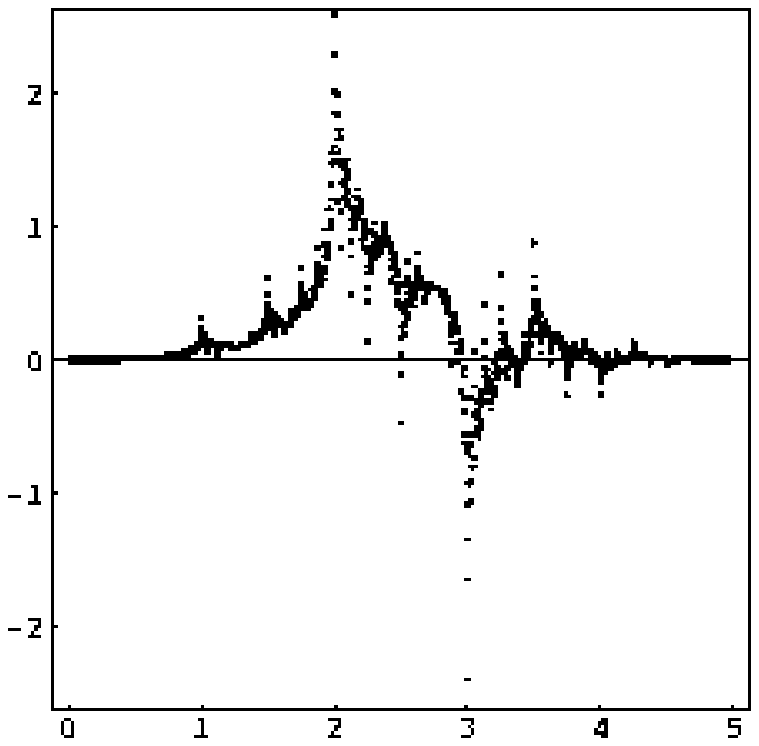}}
\put(0,262){\makebox(120,12){dk: $\theta=\pi/4,\;\rho=5\pi/6$}}
\put(120,262){\makebox(120,12){ek: $\theta=\pi/3,\;\rho=5\pi/6$}}
\put(240,262){\makebox(120,12){fk: $\theta=5\pi/12,\;\rho=5\pi/6$}}
\put(0,143){\includegraphics[%
height=119bp,width=120bp]{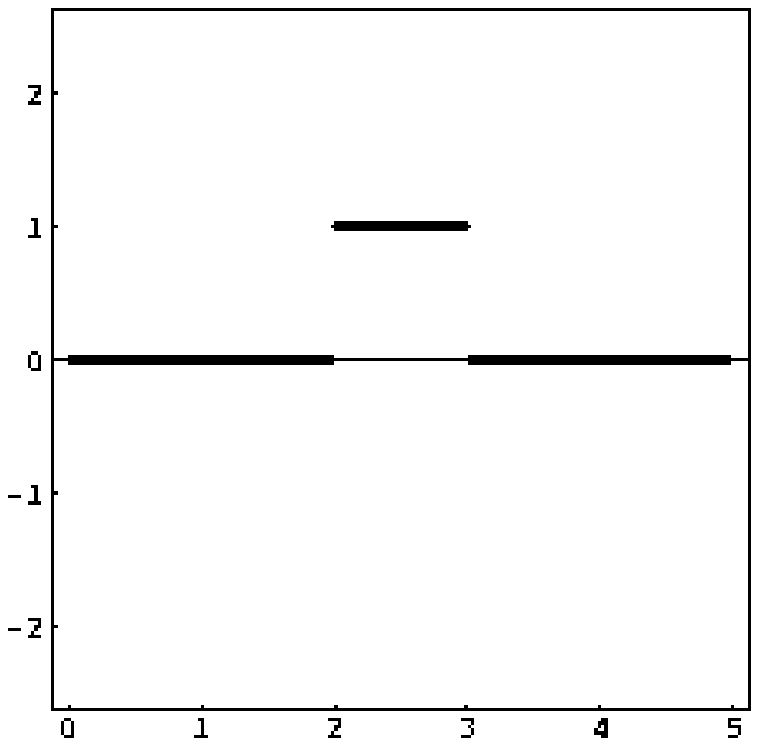}}
\put(120,143){\includegraphics[%
height=119bp,width=120bp]{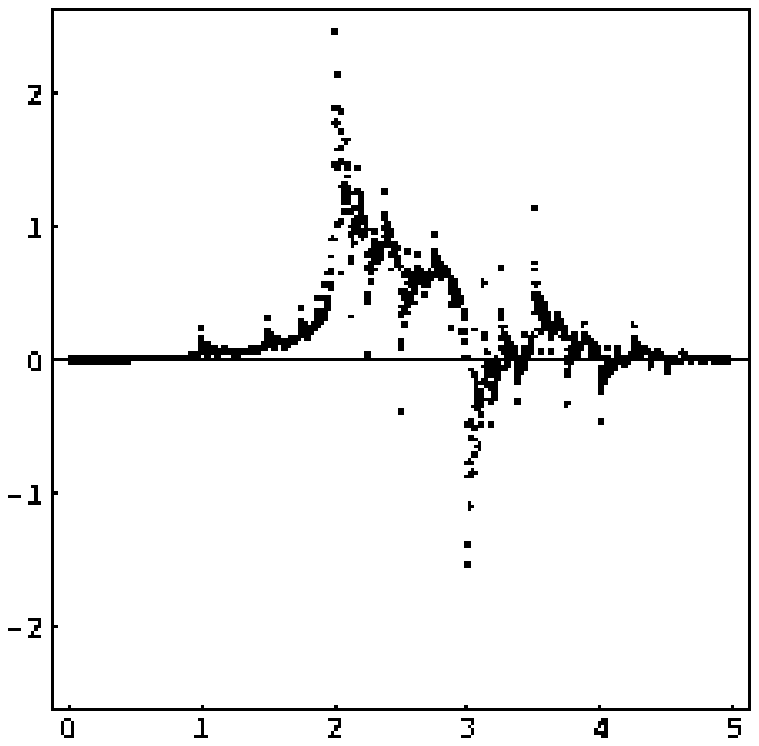}}
\put(240,143){\includegraphics[%
height=119bp,width=120bp]{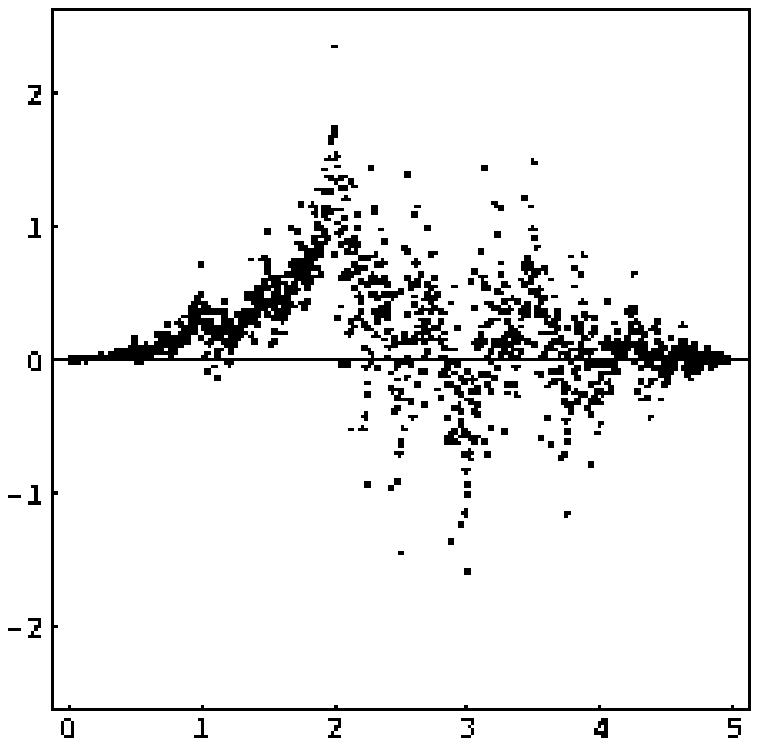}}
\put(0,131){\makebox(120,12){dj: $\theta=\pi/4,\;\rho=3\pi/4$}}
\put(120,131){\makebox(120,12){ej: $\theta=\pi/3,\;\rho=3\pi/4$}}
\put(240,131){\makebox(120,12){fj: $\theta=5\pi/12,\;\rho=3\pi/4$}}
\put(0,12){\includegraphics[%
height=119bp,width=120bp]{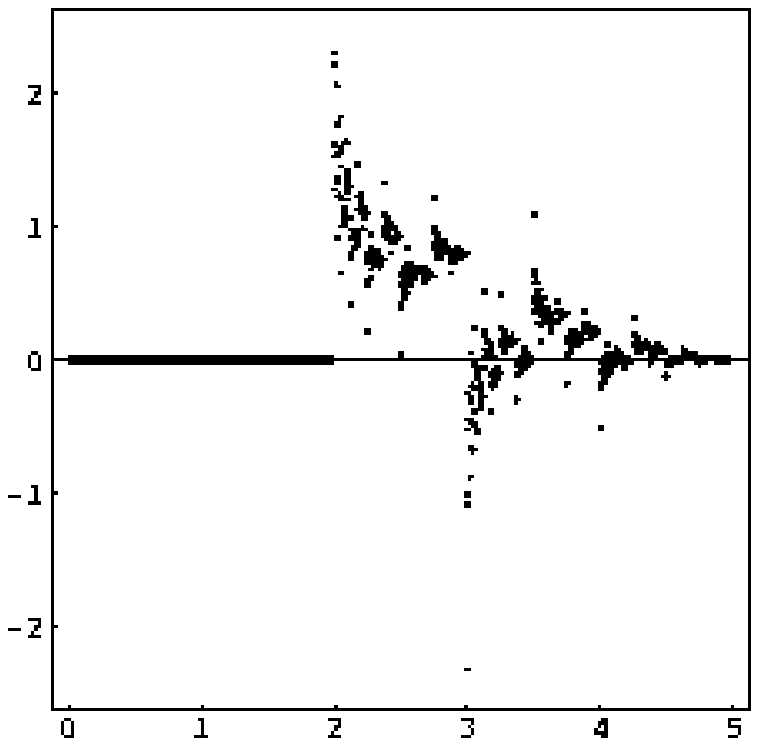}}
\put(120,12){\includegraphics[%
height=119bp,width=120bp]{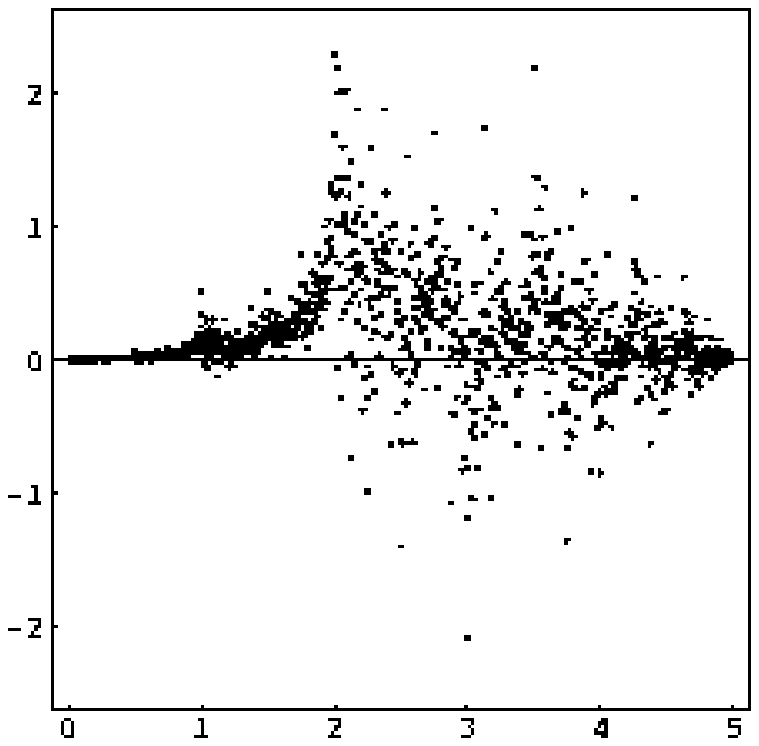}}
\put(240,12){\includegraphics[%
height=119bp,width=120bp]{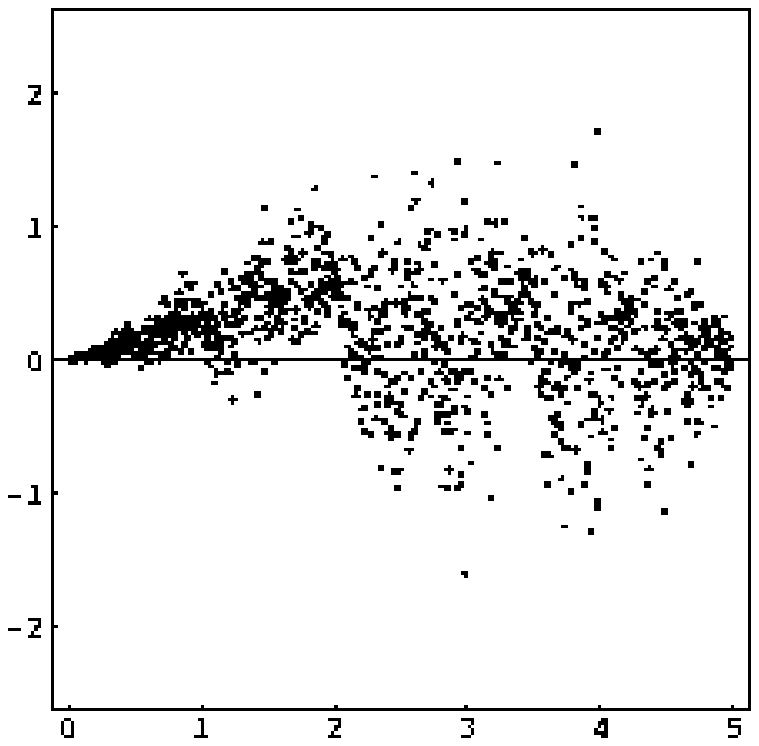}}
\put(0,0){\makebox(120,12){di: $\theta=\pi/4,\;\rho=2\pi/3$}}
\put(120,0){\makebox(120,12){ei: $\theta=\pi/3,\;\rho=2\pi/3$}}
\put(240,0){\makebox(120,12){fi: $\theta=5\pi/12,\;\rho=2\pi/3$}}
\end{picture}
\label{P10}\end{figure}

\begin{figure}[tbp]
\setlength{\unitlength}{1bp}
\begin{picture}(360,524)
\put(0,405){\includegraphics[%
height=119bp,width=120bp]{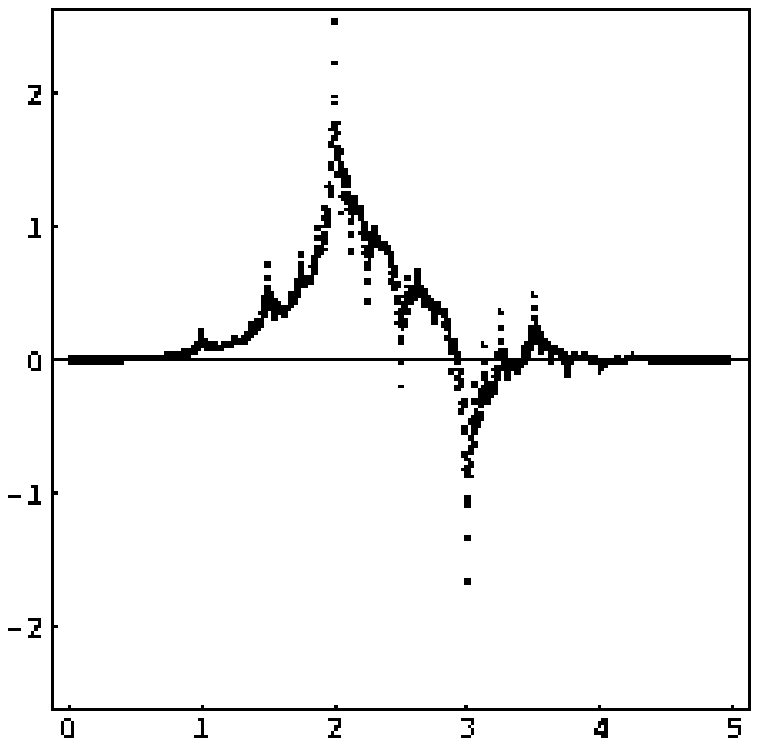}}
\put(120,405){\includegraphics[%
height=119bp,width=120bp]{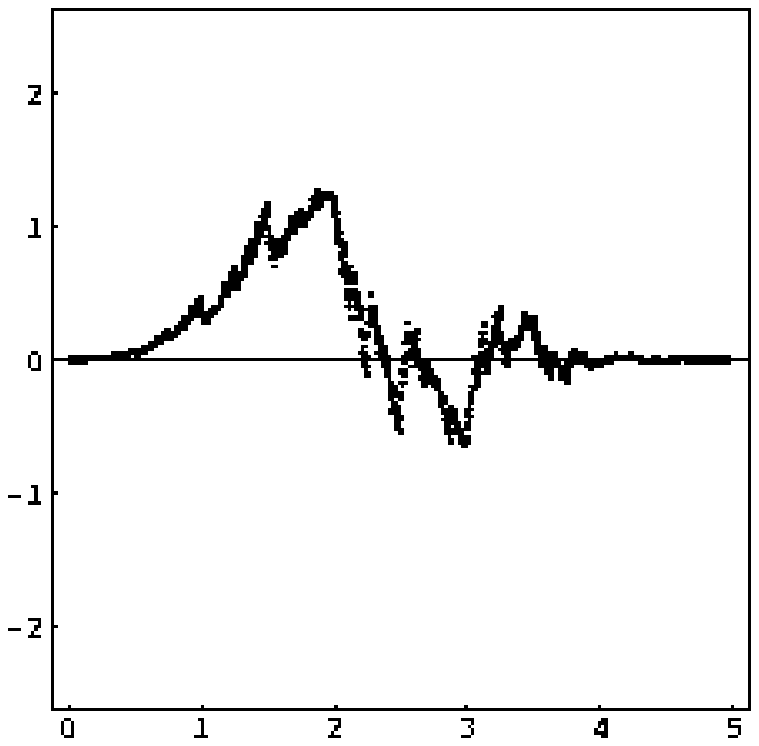}}
\put(240,405){\includegraphics[%
height=119bp,width=120bp]{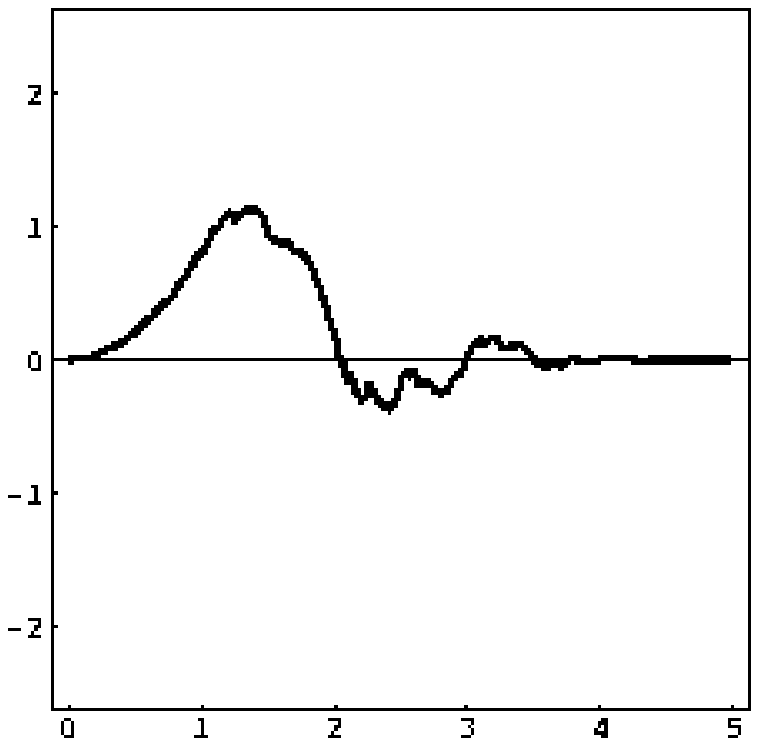}}
\put(0,393){\makebox(120,12){gl: $\theta=\pi/2,\;\rho=11\pi/12$}}
\put(120,393){\makebox(120,12){hl: $\theta=7\pi/12,\;\rho=11\pi/12$}}
\put(240,393){\makebox(120,12){il: $\theta=2\pi/3,\;\rho=11\pi/12$}}
\put(0,274){\includegraphics[%
height=119bp,width=120bp]{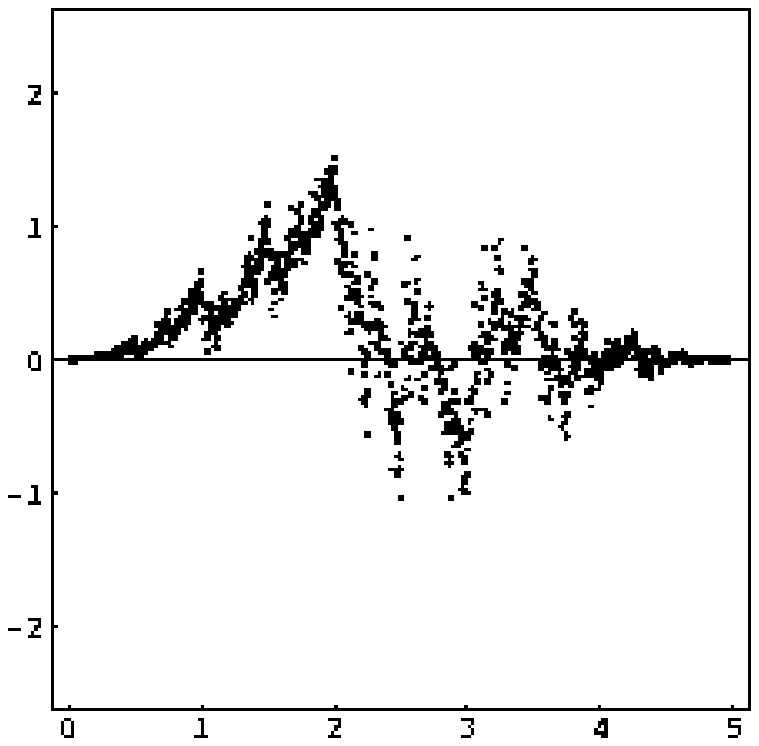}}
\put(120,274){\includegraphics[%
height=119bp,width=120bp]{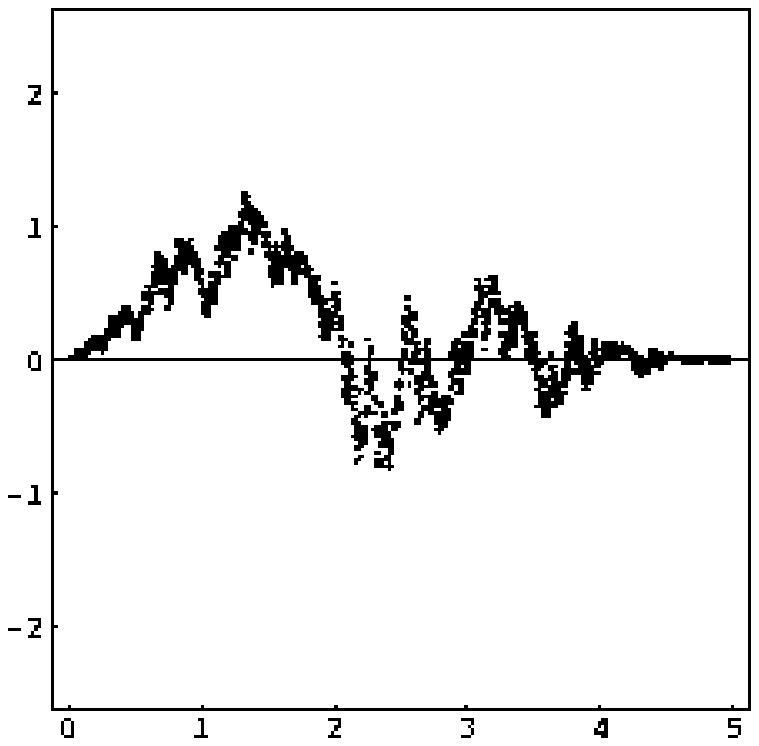}}
\put(240,274){\includegraphics[%
height=119bp,width=120bp]{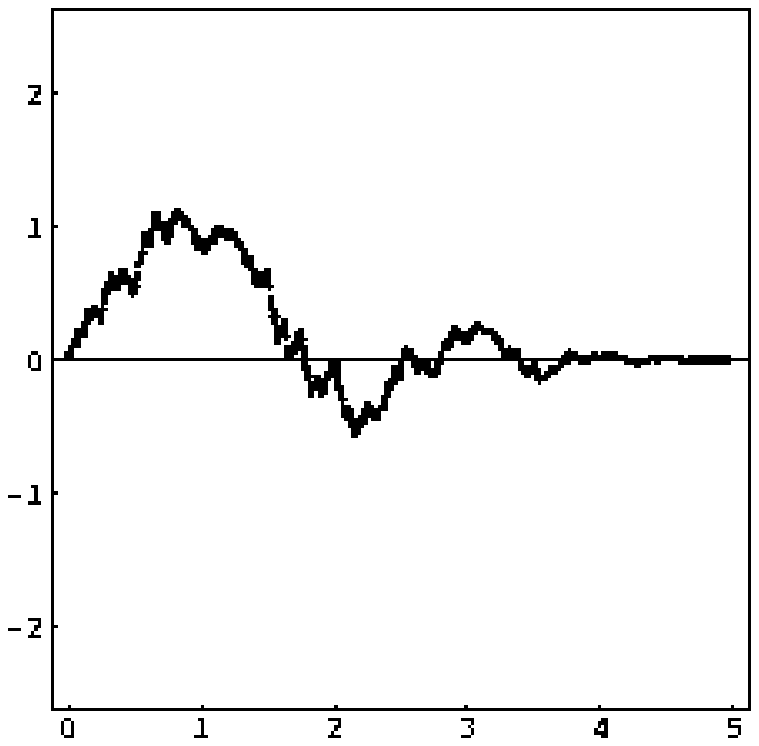}}
\put(0,262){\makebox(120,12){gk: $\theta=\pi/2,\;\rho=5\pi/6$}}
\put(120,262){\makebox(120,12){hk: $\theta=7\pi/12,\;\rho=5\pi/6$}}
\put(240,262){\makebox(120,12){ik: $\theta=2\pi/3,\;\rho=5\pi/6$}}
\put(0,143){\includegraphics[%
height=119bp,width=120bp]{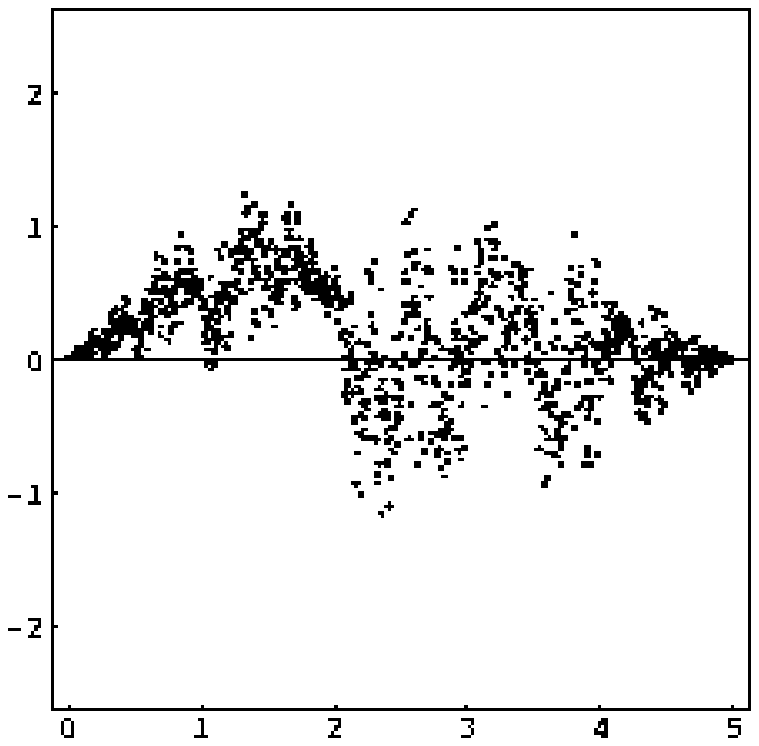}}
\put(120,143){\includegraphics[%
height=119bp,width=120bp]{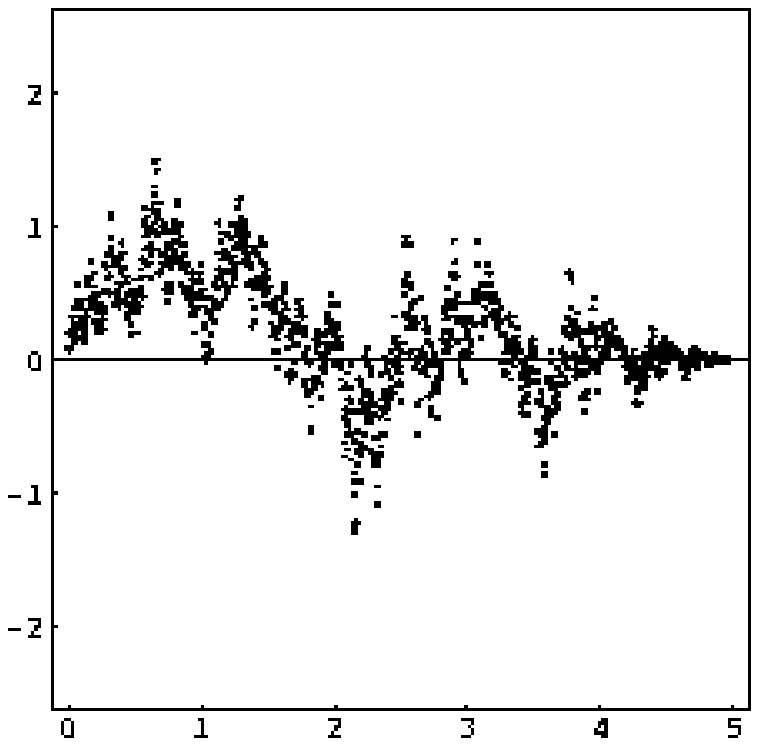}}
\put(240,143){\includegraphics[%
height=119bp,width=120bp]{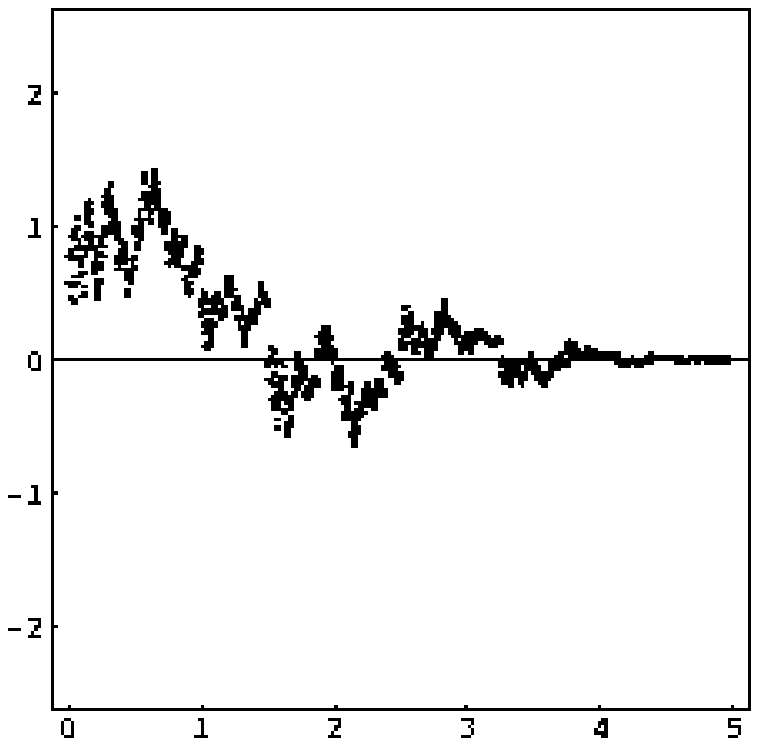}}
\put(0,131){\makebox(120,12){gj: $\theta=\pi/2,\;\rho=3\pi/4$}}
\put(120,131){\makebox(120,12){hj: $\theta=7\pi/12,\;\rho=3\pi/4$}}
\put(240,131){\makebox(120,12){ij: $\theta=2\pi/3,\;\rho=3\pi/4$}}
\put(0,12){\includegraphics[%
height=119bp,width=120bp]{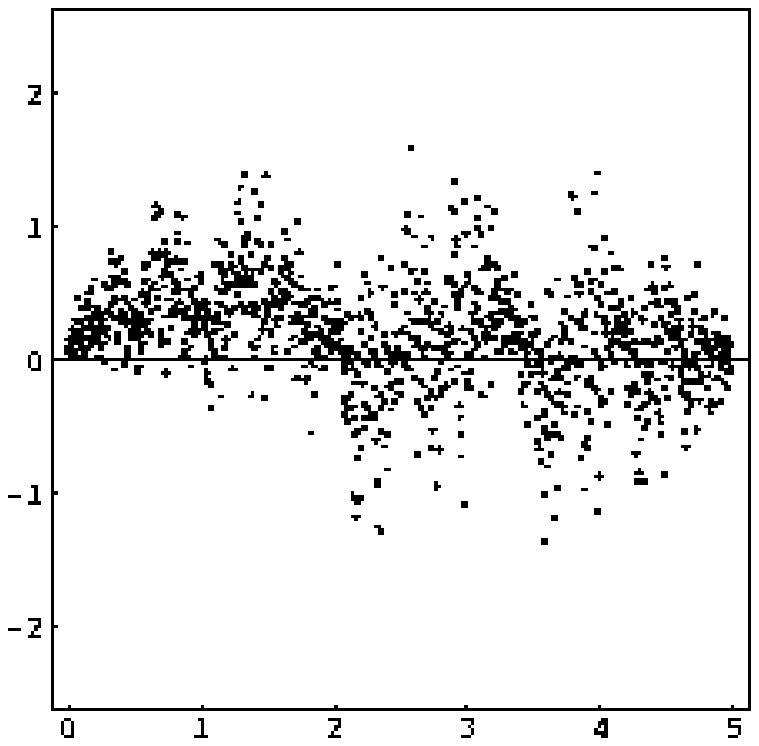}}
\put(120,12){\includegraphics[%
height=119bp,width=120bp]{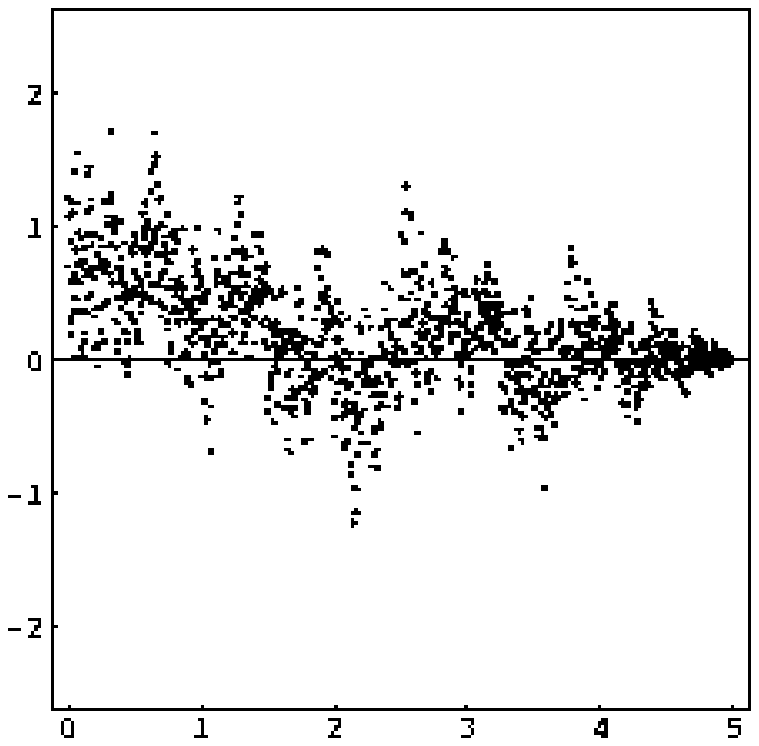}}
\put(240,12){\includegraphics[%
height=119bp,width=120bp]{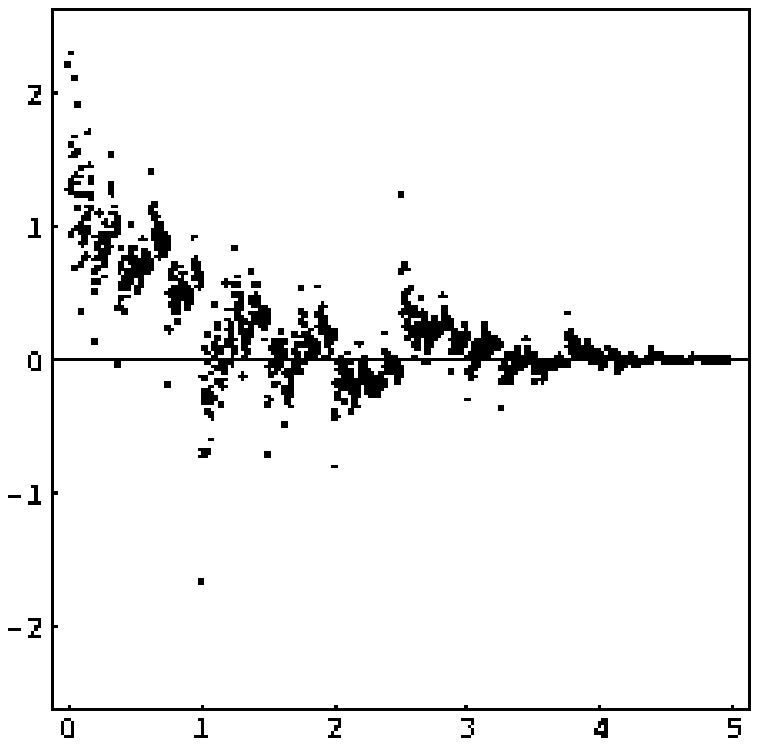}}
\put(0,0){\makebox(120,12){gi: $\theta=\pi/2,\;\rho=2\pi/3$}}
\put(120,0){\makebox(120,12){hi: $\theta=7\pi/12,\;\rho=2\pi/3$}}
\put(240,0){\makebox(120,12){ii: $\theta=2\pi/3,\;\rho=2\pi/3$}}
\end{picture}
\label{P11}\end{figure}

\begin{figure}[tbp]
\setlength{\unitlength}{1bp}
\begin{picture}(360,524)
\put(0,405){\includegraphics[%
height=119bp,width=120bp]{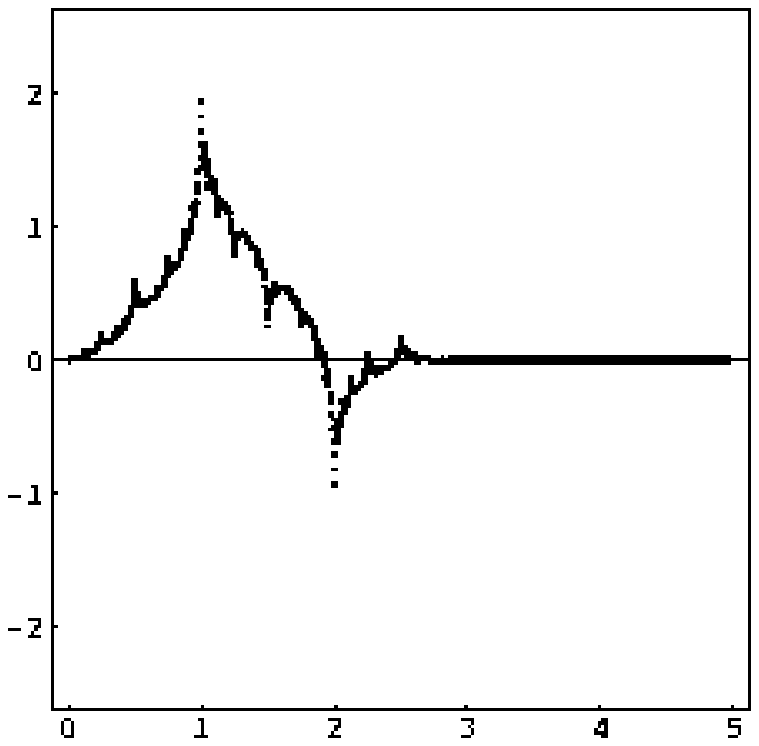}}
\put(120,405){\includegraphics[%
height=119bp,width=120bp]{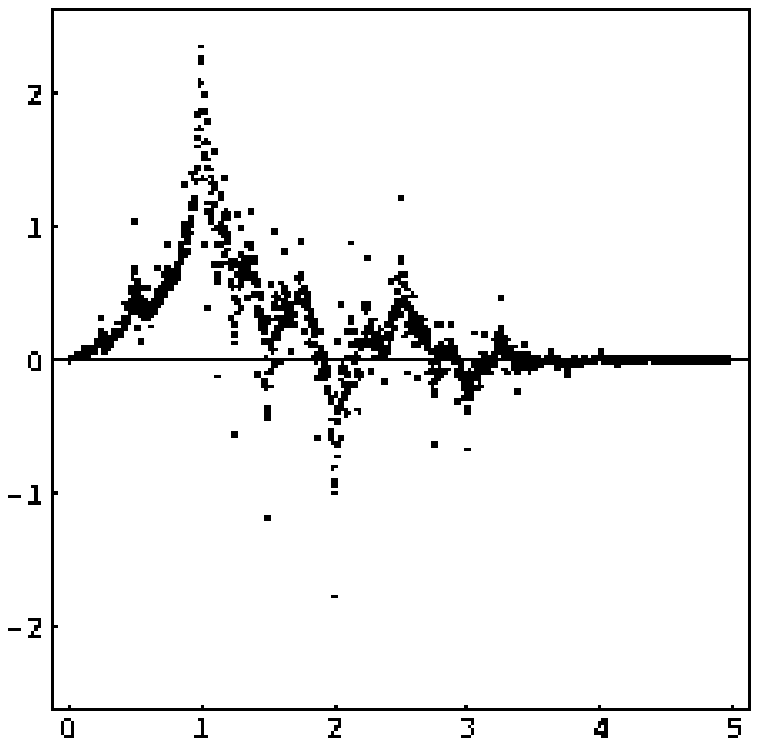}}
\put(240,405){\includegraphics[%
height=119bp,width=120bp]{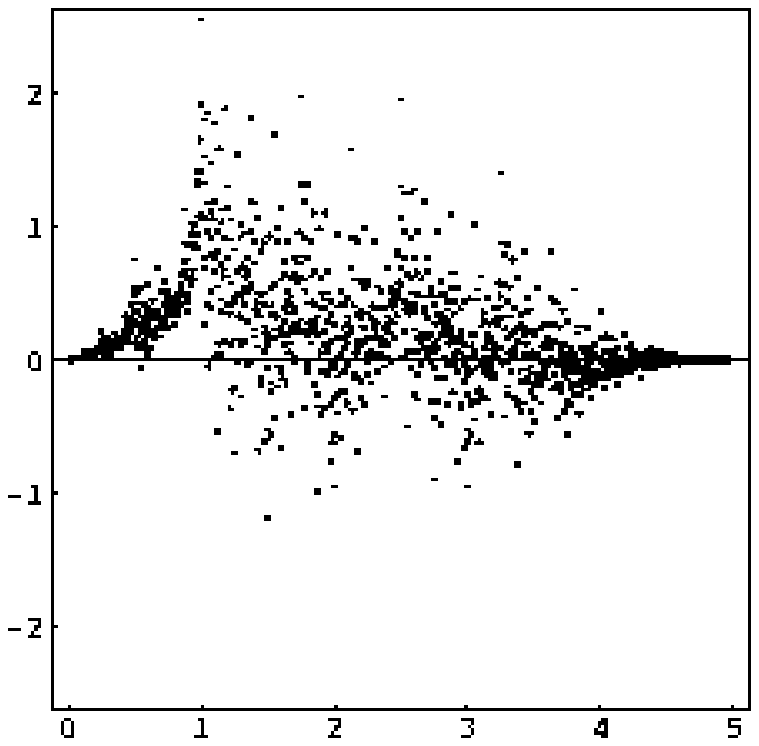}}
\put(0,393){\makebox(120,12){jl: $\theta=3\pi/4,\;\rho=11\pi/12$}}
\put(120,393){\makebox(120,12){kl: $\theta=5\pi/6,\;\rho=11\pi/12$}}
\put(240,393){\makebox(120,12){ll: $\theta=11\pi/12,\;\rho=11\pi/12$}}
\put(0,274){\includegraphics[%
height=119bp,width=120bp]{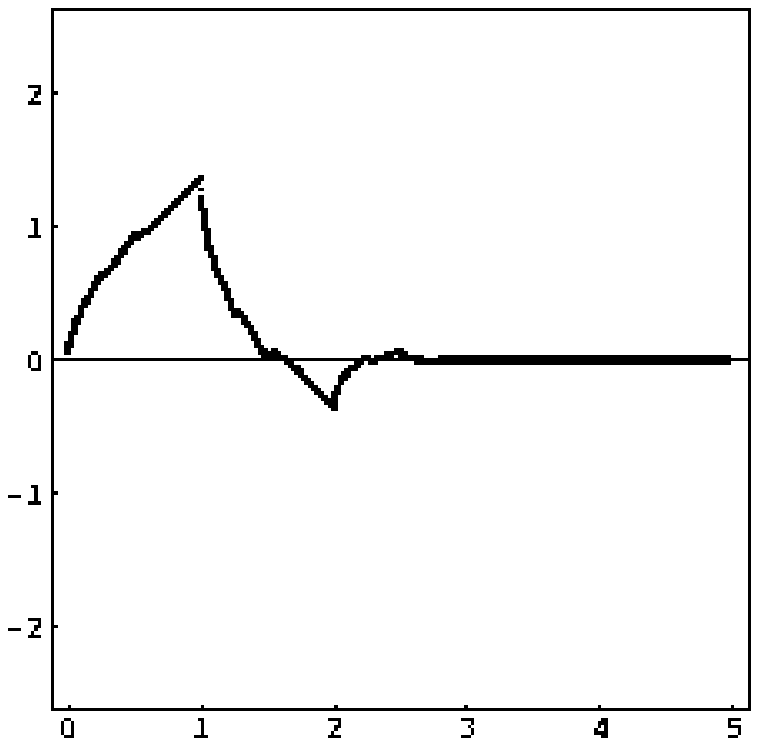}}
\put(120,274){\includegraphics[%
height=119bp,width=120bp]{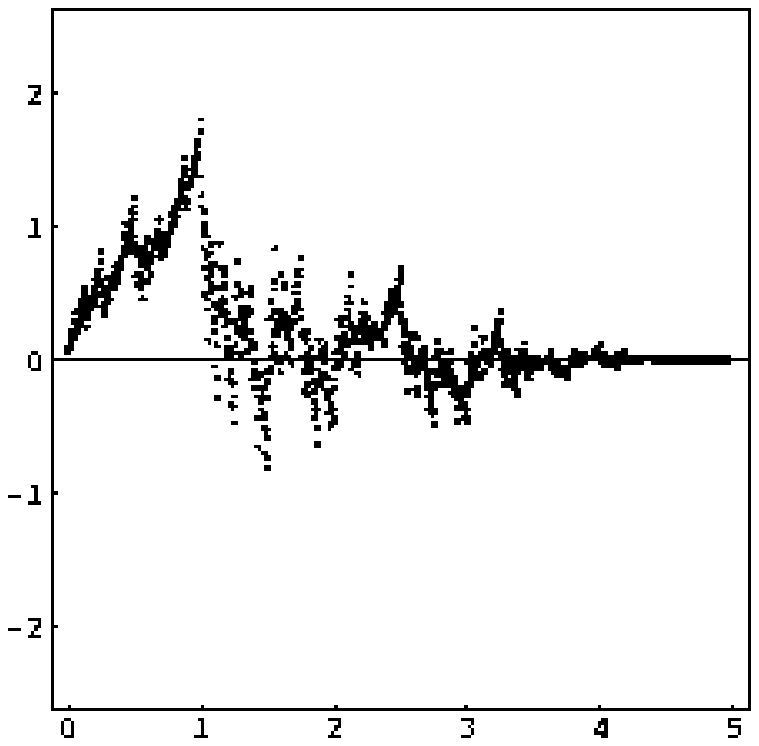}}
\put(240,274){\includegraphics[%
height=119bp,width=120bp]{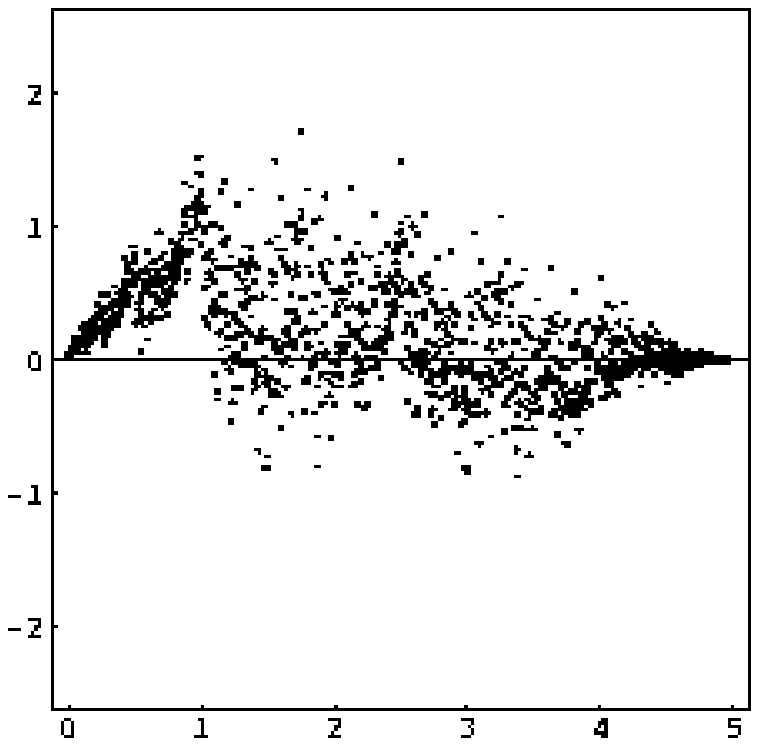}}
\put(0,262){\makebox(120,12){jk: $\theta=3\pi/4,\;\rho=5\pi/6$}}
\put(120,262){\makebox(120,12){kk: $\theta=5\pi/6,\;\rho=5\pi/6$}}
\put(240,262){\makebox(120,12){lk: $\theta=11\pi/12,\;\rho=5\pi/6$}}
\put(0,143){\includegraphics[%
height=119bp,width=120bp]{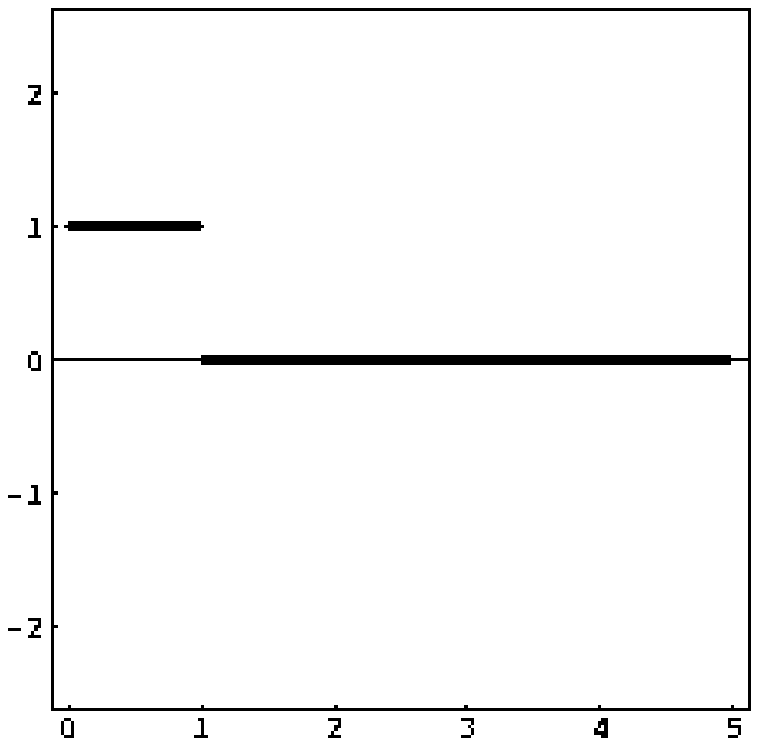}}
\put(120,143){\includegraphics[%
height=119bp,width=120bp]{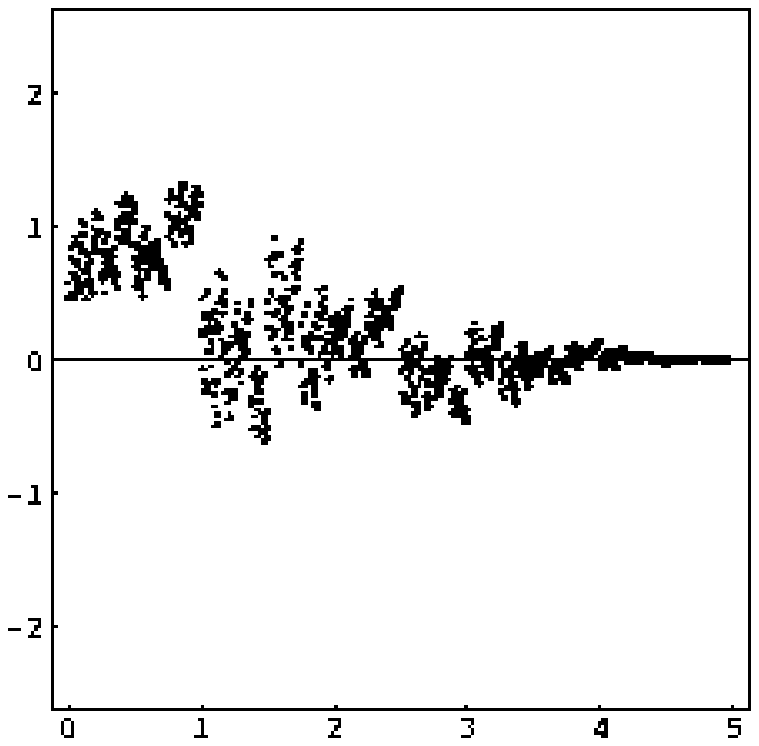}}
\put(240,143){\includegraphics[%
height=119bp,width=120bp]{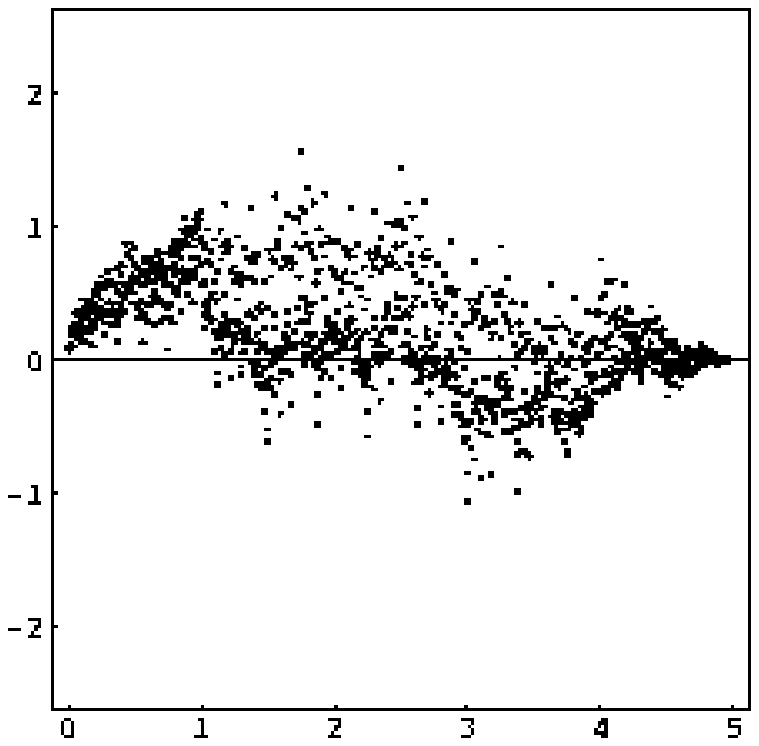}}
\put(0,131){\makebox(120,12){jj: $\theta=3\pi/4,\;\rho=3\pi/4$}}
\put(120,131){\makebox(120,12){kj: $\theta=5\pi/6,\;\rho=3\pi/4$}}
\put(240,131){\makebox(120,12){lj: $\theta=11\pi/12,\;\rho=3\pi/4$}}
\put(0,12){\includegraphics[%
height=119bp,width=120bp]{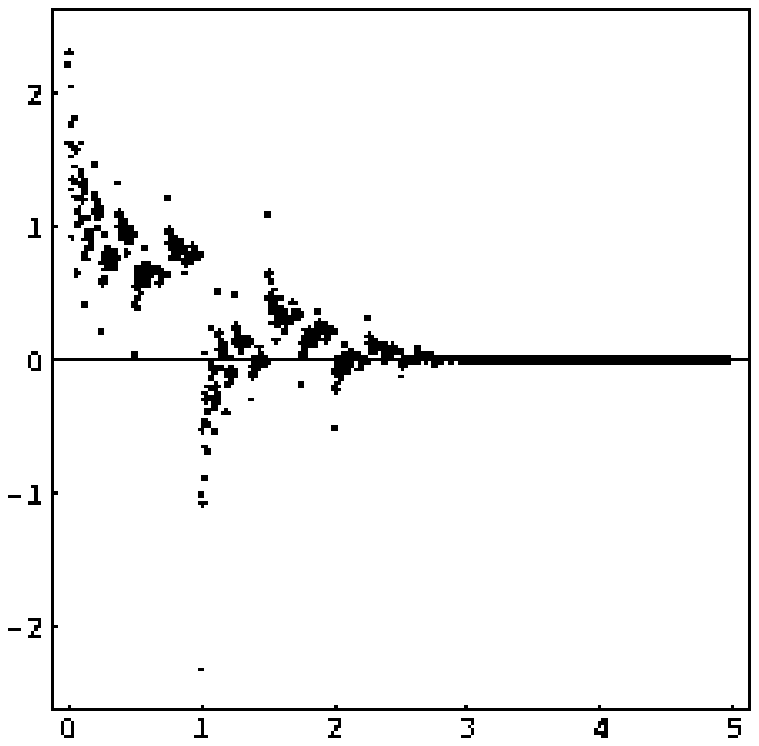}}
\put(120,12){\includegraphics[%
height=119bp,width=120bp]{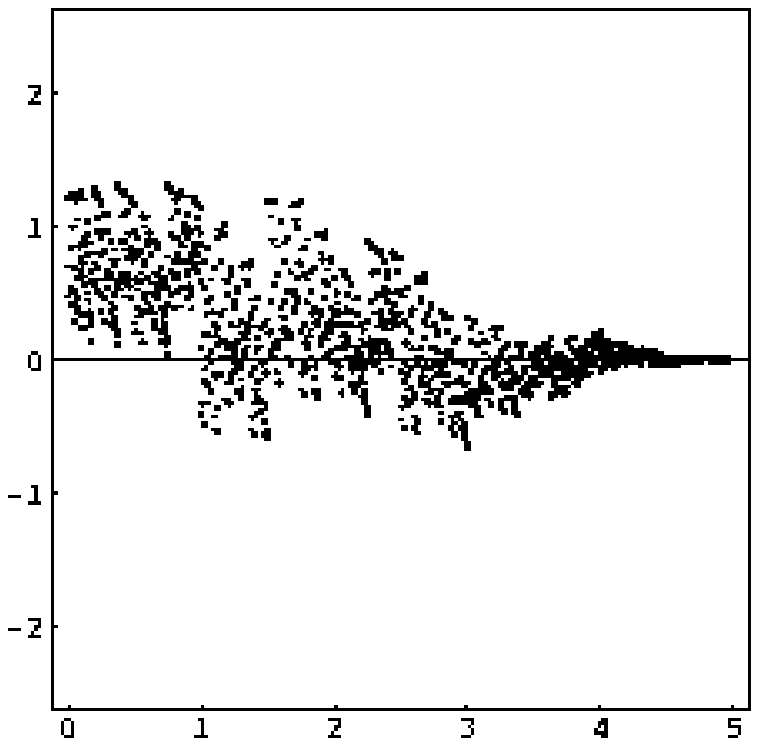}}
\put(240,12){\includegraphics[%
height=119bp,width=120bp]{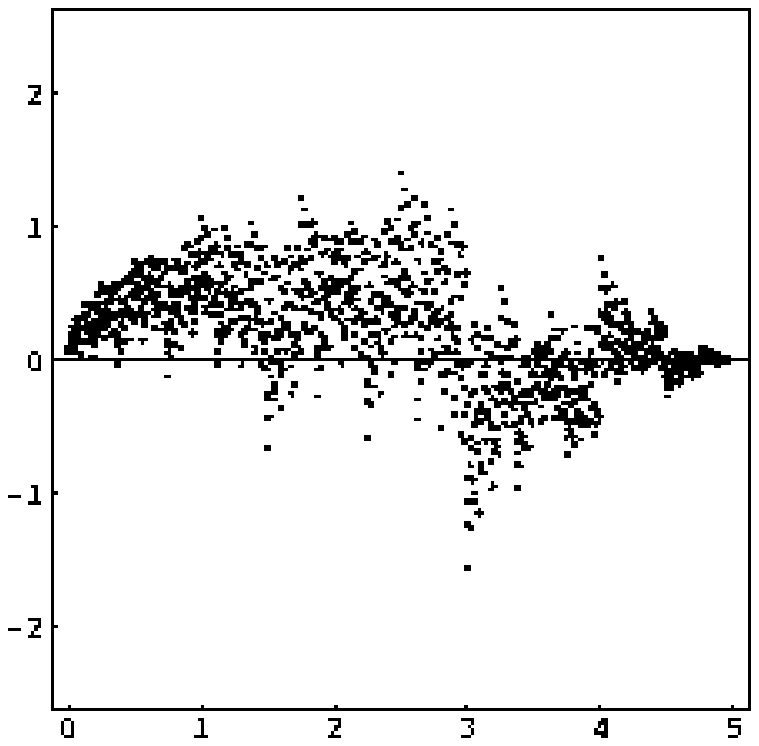}}
\put(0,0){\makebox(120,12){ji: $\theta=3\pi/4,\;\rho=2\pi/3$}}
\put(120,0){\makebox(120,12){ki: $\theta=5\pi/6,\;\rho=2\pi/3$}}
\put(240,0){\makebox(120,12){li: $\theta=11\pi/12,\;\rho=2\pi/3$}}
\end{picture}
\label{P12}\end{figure}

\clearpage

\bibliographystyle{bftalpha}
\bibliography{jorgen}
\end{document}